\documentclass[11pt,leqno]{article}
\usepackage[usenames,dvipsnames,svgnames]{xcolor}

\usepackage{hyperref}

\usepackage[toc]{appendix}

\newcommand{\addappendix}{%
  \section*{\appendixname}
  \addcontentsline{toc}{section}{\appendixname}
  \stepcounter{section}
  \renewcommand{\thesection}{A}
 \stepcounter{equation}
  }

\usepackage{amssymb,amsthm,amsmath,amscd,amsfonts,bbm,mathrsfs,fullpage,comment}
\usepackage{moreverb}
\usepackage[all]{xy}

 \def\noproof{{\unskip\nobreak\hfill\penalty50\hskip2em\hbox{}%
      \nobreak\hfill$\Box$\parfillskip=0pt%
     \finalhyphendemerits=0\par}}
\def\enddemo{\ifmmode\eqno\Box\else\noproof\vskip0.8truecm\fi}

\newtheorem{theo}{Theorem}[section]
\newtheorem{theorem}[theo]{Theorem}
\newtheorem{df}[theo]{Definition}

\newtheorem{lemma/def}[theo]{Lemma/Definition}

\newtheorem{prop}[theo]{Proposition}
\newtheorem{corollary}[theo]{Corollary}

\newtheorem{assum}[theo]{Assumption}
\newtheorem{claim}[theo]{Claim}
\newtheorem{remark}[theo]{Remark}
\newtheorem{remarks}[theo]{Remarks}

\newtheorem{lemma}[theo]{Lemma}

\newcommand{\lra}{\longrightarrow}
\newcommand{\lla}{\longleftarrow}
\newcommand{\hra}{\hookrightarrow}

\DeclareMathOperator{\Ban}{Ban}
\DeclareMathOperator{\BS}{BS}
\DeclareMathOperator{\Spec}{Spec}

\DeclareMathOperator{\Res}{Res}
\DeclareMathOperator{\Infl}{Infl}

\DeclareMathOperator{\Ind}{Ind}
\DeclareMathOperator{\cInd}{c-Ind}
\DeclareMathOperator{\ab}{ab}

\DeclareMathOperator{\aug}{aug}
\DeclareMathOperator{\fgaug}{fg aug}

\DeclareMathOperator{\res}{res}

\DeclareMathOperator{\rec}{rec}
\DeclareMathOperator{\charac}{char}

\DeclareMathOperator{\ord}{ord}
\DeclareMathOperator{\sord}{-ord}

\DeclareMathOperator{\id}{id}

\DeclareMathOperator{\ev}{ev}

\DeclareMathOperator{\Hom}{Hom}
\DeclareMathOperator{\Aut}{Aut}
\DeclareMathOperator{\Vol}{Vol}

\DeclareMathOperator{\uHom}{{\underline{\Hom}}}

\DeclareMathOperator{\Ext}{Ext}

\DeclareMathOperator{\End}{End}

\DeclareMathOperator{\fl}{fl}
\DeclareMathOperator{\pr}{pr}

\DeclareMathOperator{\rank}{rank}
\DeclareMathOperator{\sm}{sm}
\DeclareMathOperator{\open}{open}

\DeclareMathOperator{\supp}{supp}
\DeclareMathOperator{\tor}{tor}
\DeclareMathOperator{\wcdot}{\, \cdot\, }

\DeclareMathOperator{\image}{im}

\DeclareMathOperator{\coker}{coker}

\DeclareMathOperator{\Maps}{Maps}

\DeclareMathOperator{\Frob}{Frob}
\DeclareMathOperator{\Sh}{Sh}
\DeclareMathOperator{\opp}{opp}
\DeclareMathOperator{\Ord}{Ord}
\DeclareMathOperator{\wOrd}{\widetilde{\Ord}}

\DeclareMathOperator{\Tr}{Tr}

\DeclareMathOperator{\Gal}{Gal}
\DeclareMathOperator{\Norm}{N}
\DeclareMathOperator{\Nrd}{Nrd}
\DeclareMathOperator{\Mod}{Mod}

\DeclareMathOperator{\sets}{\underline{sets}}

\DeclareMathOperator{\spaces}{\underline{Top}}

\DeclareMathOperator{\Ram}{Ram}

\DeclareMathOperator{\GL}{GL}

\DeclareMathOperator{\PGL}{PGL}

\DeclareMathOperator{\ad}{ad}
\DeclareMathOperator{\bg}{big}
\DeclareMathOperator{\ext}{ext}

\DeclareMathOperator{\Orth}{O}

\DeclareMathOperator{\Ann}{Ann}

\DeclareMathOperator{\incl}{incl}
\DeclareMathOperator{\an}{an}
\DeclareMathOperator{\Stab}{Stab}

\DeclareMathOperator{\et}{et}

\DeclareMathOperator{\tr}{tr}
\DeclareMathOperator{\fin}{fin}
\DeclareMathOperator{\lfin}{lfin}

\DeclareMathOperator{\AR}{AR-}

\newcommand{\ARp}{\AR \fP\mbox{-}}

\newcommand{\pad}{\fP\mbox{-}\ad}

\DeclareMathOperator{\BC}{BC}
\DeclareMathOperator{\can}{can}

\DeclareMathOperator{\Man}{Man}
\DeclareMathOperator{\Sm}{Sm}

\DeclareMathOperator{\ladm}{ladm}
\DeclareMathOperator{\adm}{adm}

\DeclareMathOperator{\cont}{cont}

\DeclareMathOperator{\gal}{gal}

\DeclareMathOperator{\cycl}{cycl}

\DeclareMathOperator{\prolim}{\underset{\lla}{\lim}}
\DeclareMathOperator{\prolimk}{\underset{\underset{k}{\lla}}{\lim}}
\DeclareMathOperator{\prolimm}{\underset{\underset{m}{\lla}}{\lim}}
\DeclareMathOperator{\prolimn}{\underset{\underset{n}{\lla}}{\lim}}
\DeclareMathOperator{\dlim}{\underset{\lra}{\lim}}
\DeclareMathOperator{\colim}{colim}

\DeclareMathOperator{\barF}{\overline{F}}
\DeclareMathOperator{\barG}{\overline{G}}

\DeclareMathOperator{\einhalb}{\mbox{\small{$\frac{1}{2}$}}}

\newcommand{\io}{{\iota}}
\newcommand{\La}{{\Lambda}}
\newcommand{\la}{{\lambda}}

\newcommand{\fa}{{\mathfrak a}}
\newcommand{\fb}{{\mathfrak b}}

\newcommand{\fG}{{\mathfrak G}}
\newcommand{\fH}{{\mathfrak H}}
\newcommand{\fh}{{\mathfrak h}}

\newcommand{\fd}{{\mathfrak d}}
\newcommand{\ff}{{\mathfrak f}}

\newcommand{\fm}{{\mathfrak m}}
\newcommand{\fM}{{\mathfrak M}}
\newcommand{\fn}{{\mathfrak n}}
\newcommand{\fp}{{\mathfrak p}}
\newcommand{\fP}{{\mathfrak P}}
\newcommand{\fU}{{\mathfrak U}}
\newcommand{\fq}{{\mathfrak q}}

\newcommand{\bC}{{\mathbb C}}
\newcommand{\bD}{{\mathbb D}}

\newcommand{\bG}{{\mathbb G}}
\newcommand{\bH}{{\mathbb H}}

\newcommand{\bN}{{\mathbb N}}

\newcommand{\bQ}{{\mathbb Q}}
\newcommand{\bR}{{\mathbb R}}
\newcommand{\bT}{{\mathbb T}}
\newcommand{\bZ}{{\mathbb Z}}
\newcommand{\barQ}{{\overline{\mathbb Q}}}

\newcommand{\bA}{{\mathbb A}}

\newcommand{\hatW}{{\widehat{W}}}

\newcommand{\barB}{{\overline{B}}}

\newcommand{\barT}{{\overline{T}}}

\newcommand{\barchi}{{\overline{\chi}}}

\newcommand{\barnu}{{\overline{\nu}}}

\newcommand{\cA}{{\mathcal A}}
\newcommand{\cB}{{\mathcal B}}
\newcommand{\cC}{{\mathcal C}}
\newcommand{\cD}{{\mathcal D}}

\newcommand{\cF}{{\mathcal F}}
\newcommand{\cG}{{\mathcal G}}

\newcommand{\cI}{{\mathcal I}}
\newcommand{\cJ}{{\mathcal J}}

\newcommand{\cL}{{\mathcal L}}
\newcommand{\cM}{{\mathcal M}}
\newcommand{\cN}{{\mathcal N}}
\newcommand{\cO}{{\mathcal O}}

\newcommand{\cR}{{\mathcal R}}
\newcommand{\cS}{{\mathcal S}}
\newcommand{\cT}{{\mathcal T}}
\newcommand{\cU}{{\mathcal U}}
\newcommand{\cV}{{\mathcal V}}

\newcommand{\cY}{{\mathcal Y}}

\newcommand{\sA}{{\mathscr A}}
\newcommand{\sB}{{\mathscr B}}

\newcommand{\sD}{{\mathscr D}}

\newcommand{\sF}{{\mathscr F}}
\newcommand{\sG}{{\mathscr G}}
\newcommand{\sH}{{\mathscr H}}
\newcommand{\sI}{{\mathscr I}}
\newcommand{\sJ}{{\mathscr J}}
\newcommand{\sK}{{\mathscr K}}
\newcommand{\sL}{{\mathscr L}}

\newcommand{\sS}{{\mathscr S}}

\newcommand{\sX}{{\mathscr X}}
\newcommand{\sY}{{\mathscr Y}}

\newcommand{\wsF}{{\widetilde{\mathscr F}}}
\newcommand{\wsS}{{\widetilde{\mathscr S}}}

\newcommand{\wsK}{{\widetilde{\mathscr K}}}

\newcommand{\wE}{\widetilde{E}}
\newcommand{\wG}{\widetilde{G}}
\newcommand{\wH}{\widetilde{H}}
\newcommand{\wY}{\widetilde{Y}}
\newcommand{\wX}{\widetilde{X}}

\newcommand{\wphi}{\widetilde{\varphi}}

\newcommand{\wcO}{{\widetilde{\mathcal O}}}

\newcommand{\wK}{\widetilde{K}}

\newcommand{\wPhi}{\widetilde{\Phi}}
\newcommand{\oE}{{\vec{\mathcal E}}}

\newcommand{\whH}{\widehat{H}}

\newcommand{\wu}{\widetilde{u}}

\newcommand{\ep}{{\epsilon}}
\newcommand{\vep}{{\varepsilon}}

\newcommand{\barrho}{{\overline{\rho}}}

\newcommand{\PSh}{{\mathcal PSh}} 
\newcommand{\wvphi}{\widetilde{\varphi}}
\newcommand{\wio}{\widetilde{\iota}}
\DeclareMathOperator{\pt}{pt}
\newcommand{\wsX}{{\widetilde{\mathscr X}}}
\newcommand{\varep}{{\varepsilon}}
\newcommand{\wf}{\widetilde{f}}

\newcommand{\caInd}{{\mathcal Ind}}
\newcommand{\barsK}{{\overline{\sK}}}
\newcommand{\barsS}{{\overline{\sS}}}
\newcommand{\barS}{{\overline{S}}}
\newcommand{\barK}{{\overline{K}}}
\newcommand{\barL}{{\overline{L}}}
\DeclareMathOperator{\fOrd}{{\mathbb{O}rd}}
\DeclareMathOperator{\fcoh}{fc}
\newcommand{\barfh}{{\overline{\mathfrak h}}}
\newcommand{\cZ}{{\mathcal Z}}

\newcommand{\barp}{{\overline{p}}}
\DeclareMathOperator{\eis}{eis}
\DeclareMathOperator{\cl}{cl}

\newcommand{\bu}{{\bullet}}

\newcommand{\noi}{\noindent}

\begin{document}

\title{On certain cohomology groups attached to $\fp^{\infty}$-towers of quaternionic Hilbert modular varieties}

\author{By Michael Spie{\ss}} 

\maketitle
\date{}
\begin{abstract}

For a totally real number field $F$ and a nonarchimedean prime $\fp$ of $F$ lying above a prime number $p$ we introduce certain sheaf cohomology groups that intertwine the $\fp^{\infty}$-tower of a quaternionic Hilbert modular variety associated to a quaternion algebra $D$ over $F$ that is split at $\fp$ and a $p$-adically admissible representation of $\PGL_2(F_{\fp})$. Applied to infinitesimal $p$-adic deformations of the local factor at $\fp$ of a cuspidal automorphic representation $\pi$ of $D^*(\bA)$ this yields a natural construction of infinitesimal deformations of the Galois representation attached to $\pi$. 
 \end{abstract}

\tableofcontents

\section*{Introduction}
\addcontentsline{toc}{section}{Introduction}

In this paper we introduce certain categories of ({\'e}tale) sheaves associated to a ''locally pro-{\'e}tale tower'' of algebraic varieties equipped with an action of a locally profinite group $G$. Typically, the sheaves we consider are sheaves on a site without a final object. By ''intertwining'' the tower with a representation of $G$ we construct new cohomology theories and thereby interesting Galois representations that are usually not motivic, i.e.\ that do not occur in the {\'e}tale cohomology of a single algebraic variety. Our main motivation to develop the theory presented here is to produce custom fit infinitesimal deformations of Galois representations associated to automorphic forms. Our method can be applied to construct ''big'' Galois representations associated to $p$-adic families of automorphic forms as well. Of course the construction of such  Galois representations has a long history starting with the work of Hida and Wiles (see e.g.\ \cite{hida1, hida2}, \cite{wiles2}). Our approach differs from earlier methods insofar that we work completely within the framework of topos theory as developed in \cite{sga4}. 

The results of this paper have been applied in the companion paper \cite{spiess2} where we prove the equality of automorphic and arithmetic $\sL$-invariants associated to Hilbert modular forms (of parallel weight $(2,\ldots, 2)$). We like to point out that together with (Thm.\ 6.7, \cite{spiess1}) one obtains a proof of the exceptional zero conjecture for modular elliptic curves over totally real fields (see \cite{hida3}) that avoids the use of $p$-adic families of Hilbert modular forms completely.\footnote{Furthermore we like to point out that the results of the main technical section 4.2 of \cite{spiess1} have been clarified and simplified in (\cite{dasspi}, \S 3).} Further applications of the method introduced here will be given in a sequel to this paper. 

In this work we apply our construction only to the case of a tower of quaternionic Hilbert modular varieties over a single place $\fp$ (i.e.\ we vary the level only at one place). Thus in order to describe our method and results in some detail we focus for the rest of the introduction on this situation. Specifically we consider the following set-up. We fix a totally real field $F$ and a quaternion algebra $D$ over $F$. We assume $D$ is not totally definite and let $\Sigma=\{v_1, \ldots, v_d\}$, $d\ge 1$ be the set of archimedian places of $F$ that split $D$ (i.e.\ $D\otimes_F F_v$ is isomorphic to the matrix algebra $M_2(\bR)$ for $v\in \Sigma$). Let $G=D^*/F^*$ viewed as an algebraic group over $F$ and let $K_{\infty, +}$ be a maximal compact and connected subgroup of $G(F\otimes \bR)$. We denote the adele ring of $F$ by $\bA=\prod_{v}' F_v$ and let $\bA_f=\prod_{v\nmid \infty}' F_v$ be the finite adeles. For a small enough compact open subgroup $K_f$ of $G(\bA_f)$ the group $G(F)$ acts properly discontinuously on the symmetric space $G(\bA)/(K_f\times K_{\infty, +})$. The quaternionic Hilbert modular  variety $X_{K_f}$ of level $K_f$ is a $d$-dimensional smooth quasi-projective variety defined over a specific algebraic number field $L$ (= the reflex field of the corresponding Shimura datum). The associated complex analytic manifold $X_{K_f}^{\an}$ is the space of left $G(F)$-orbits in $G(\bA)/K_f\times K_{\infty, +}$. If $d=1$ (resp.\ $d=[F:\bQ]$) then $L=F$ (resp.\ $L=\bQ$). 

We fix a nonarchimedean place $\fp$ of $F$ that is unramified in $D$ so that we can (and will) identify $G_{\fp}:=G(F_{\fp})$ with the group $\PGL_2(F_{\fp})$. We assume that $K_f$ is of the form $K_f= K \times K_f^{\fp}$ where $K=K_{\fp}$ (resp.\ $K_f^{\fp}$) is a compact open subgroup of $G_{\fp}$ (resp.\ of $G(\bA_f^{\fp})$ where $\bA_f^{\fp}=\prod_{v\nmid \infty\fp}' F_v$). In the following we keep $K_f^{\fp}$ fixed and vary $K$. So we consider the $\fp^{\infty}$-tower $X:=\{X_K\}_{K}$ where we have abbreviated $X_K  = X_{K \times K_f^{\fp}}$. The tower is equipped with an action of $G_{\fp}$. 

We introduce a topos $\Sh(X_{\et})$, i.e.\ a category of {\'e}tale sheaves on a site associated to $X$ and study certain cohomology groups defined in terms of $\Ext$-groups of sheaves in $\Sh(X_{\et})$. The category $\Sh(X_{\et})$ may be viewed as the category of $G_{\fp}$-equivariant {\'e}tale sheaves on 
\[
X_{\infty}:= \, \prolim_{K} X_{K} \, =\,  \prolim_{K} X_{K \times K_f^{\fp}}
\] 
or also as sheaves on the ''quotient'' $G_{\fp}\backslash X_{\infty}$ (note that the latter does not exists as an algebraic variety). Given a discrete (i.e.\ smooth) representation $M$ of $G_{\fp}$ -- say over an Artinian local ring $A$ with finite residue field of characteristic $p$ -- we can associate a sheaf of $A$-modules on $X_{\et}$, i.e.\ an object of $\Sh(X_{\et}, A)$ that will be denoted by $M_{X, G_{\fp}}$. We then consider the cohomology groups 
\begin{equation}
\label{sschemecohom}
H_A^{\bu}(X_{\barQ}; M, \sF) : =\, \Ext_{\Sh((X_{\barQ})_{\et}, A)}^{\bu}(M_{X, G_{\fp}}, \sF)
\end{equation}
with coefficients in an arbitrary sheaf $\sF$. These cohomology groups are equipped with an action of the absolute Galois group $\fG=\Gal(\barQ/L)$ of $L$. They are also equipped with an action of a Hecke algebra. For the latter let us assume that $K_f^{\fp} = K_0(\fn)^{\fp}$ where $\fn$ is an ideal of $\cO_F$ that is relatively prime to $\fp$ and the discriminant $\fd$ of $D$ (for the definition of $K_0(\fn)^{\fp}$ see section \ref{subsection:qhmsv}). Then the groups \eqref{sschemecohom} are equipped with an action of a spherical Hecke algebra $\bT=\bZ[T_{\fq}\mid \fq \nmid \fp\fn\fd]$. 

We also consider a variant of these cohomology groups for the $\fp^{\infty}$-tower $X^{\an}:=\{X_{K}^{\an}\}_{K}$ of complex analytic manifolds defined similarly by
\begin{equation}
\label{sschemecohom2}
H_{\cO}^{\bu}(X^{\an}; M, \sF) : =\, \Ext_{\Sh(X^{\an}, \cO)}^{\bu}(M_{X^{\an}, G_{\fp}}, \sF).
\end{equation}
Here it is useful to allow more generally a complete noetherian local ring $\cO$ with finite residue field of characteristic $p$ as coefficient ring. 
\medskip

\noi We now list some basic properties of these cohomology theories.
\medskip

\noi {\bf Comparison with the {\'e}tale cohomology of a single variety $X_K$.} Let $K$ be a (sufficiently small) compact open subgroup of $G_{\fp}$ and let $M= \cInd_{K}^{G_{\fp}} A$ (compact induction). For $\sF\in \Sh(X_{\et}, A)$ let $\sF_K$ be the ''restriction'' of $\sF\in \Sh((X_K)_{\et}, A)$ to $X_K$.  Then we have
\begin{equation*}
\label{sschemecohom3}
H_A^{\bu}(X_{\barQ}; \cInd_{K}^{G_{\fp}} A, \sF) \, \cong \, H^{\bu}((X_K)_{\barQ}, \sF_K).
\end{equation*}

\noi {\bf Covering spectral sequence.} There exists a spectral sequence 
\begin{equation*}
\label{sschemecohom4}
E_2^{rs}= \Ext_{G_{\fp}, A}^r(M, H^s((X_{\barQ})_{\infty}, \sF)) \Rightarrow   H_A^{r+s}(X_{\barQ}; M, \sF)
\end{equation*}
where $H^{\bu}((X_{\barQ})_{\infty}, \sF)= \dlim_{K} H^{\bu}((X_K)_{\barQ}, \sF_K)$. 
Here by $\Ext_{G_{\fp}, A}^{\bu}(\wcdot, \wcdot)$ we denote the $\Ext$-groups in the category $\Mod_A^{\sm}(G_{\fp})$  of discrete $A[G_{\fp}]$-modules. If $M= \cInd_{K}^{G_{\fp}} A$ for $K\subseteq G_{\fp}$ sufficiently small the spectral sequence can be identified with the Hochschild-Serre spectral sequence (see \cite{artin}, 3.4.7) for the ''covering'' $X_{\infty}\to X_K$ with ''group of deck transformations'' $K$. 
\medskip

\noi {\bf Comparison between {\'e}tale and analytic cohomology.} Assume that $\sF\in \Sh(X_{\et}, A)$ is constructible, i.e.\ $\sF_K$ is constructible for one (hence every) sufficiently small compact open subgroup $K$ of $G_{\fp}$. As in the case of usual {\'e}tale cohomology one can associated to $\sF$ a sheaf $\sF^{\an}\in \Sh(X^{\an}, A)$ and there exists a comparison isomorphism
\begin{equation}
\label{sschemecohom5}
H_A^{\bu}(X_{\barQ}; M, \sF) \, \cong \, H_A^{\bu}(X^{\an}; M, \sF^{\an}).
\end{equation}
\medskip

\noi {\bf Comparison with group cohomology.} Let $(D^*)_+$ be the subgroup of elements of $D^*$ with totally positive reduced norm and put $G(F)_+ = (D^*)_+/F^*$. If $M$ is a discrete $\cO[G_{\fp}]$-module that is projective as an $\cO$-module then we have 
\begin{equation}
\label{sschemecohom6}
H_{\cO}^{\bu}(X^{\an}; M, \cO) \, \cong \, H^{\bu}(G(F)_+, \cA_{\cO}(M, K_f^{\fp}; \cO)).
\end{equation}
Here the coefficients on the left hand side are $\cO=\cO_{X^{\an}, G_{\fp}}$, the sheaf on $X^{\an}$ associated to the trivial $G_{\fp}$-module $\cO$. The coefficients $\cA_{\cO}(M, K_f^{\fp}; \cO)$ on the right hand side is 
 the $\cO[G(F)_+]$-module of maps $G(\bA_f^{\fp})/K_f^{\fp}\times  M\to \cO$ that are $\cO$-linear in the second component. The cohomology groups on the right hand side of \eqref{sschemecohom6} have been introduced and studied in \cite{spiess1}. The comparison isomorphisms \eqref{sschemecohom5} and \eqref{sschemecohom6} imply in particular that they are equipped with a canonical Galois (i.e.\ $\fG$-) action in case where $\cO=A$ is an Artinian local ring.\medskip

\noi {\bf Relation to $\fp$-ordinary cohomology} Let $B_{\fp}$ be the standard Borel subgroup of $G_{\fp}= \PGL_2(F_{\fp})$ and let $T_{\fp}\cong \bG_m(F_{\fp})= F_{\fp}^*$ be the maximal torus in $B_{\fp}$. In this paper we focus mainly on the case $M=\Ind_{B_{\fp}}^{G_{\fp}} W$ (parabolic induction) where $W$ is a discrete $A[T_{\fp}]$-module (or more generally a discrete $\cO[T_{\fp}]$-module). We abbreviate
\begin{eqnarray*}
\label{sschemecohom7}
\bH_A^{\bu}(X_{\barQ}; W, \sF) & := & H_A^{\bu}(X_{\barQ}; \Ind_{B_{\fp}}^{G_{\fp}} W, \sF), \\
 \bH_{\cO}^{\bu}(X^{\an}; W, \sF) & := & H_{\cO}^{\bu}(X^{\an}; \Ind_{B_{\fp}}^{G_{\fp}} W, \sF).\nonumber
\end{eqnarray*}
We also define the $n$-th $\fp$-{\it ordinary cohomology group} 
$\fOrd_A^n(X_{\barQ}, \sF)$ of $X$ with coefficient in $\sF \in \Sh((X_K)_{\et}, A)$ as the $n$-th right derived functor of $\sF \mapsto \Ord(H^0((X_{\barQ})_{\infty}, \sF))$ where $\Ord=\Ord_{B_{\fp}}: \Mod_A^{\sm}(G_{\fp}) \to \Mod_A^{\ladm}(T_{\fp})$ denotes Emerton's functor of ordinary parts \cite{emerton2}. These cohomology groups are linked by a spectral sequence   
\begin{equation}
\label{sschemecohom8}
E_2^{rs} = \Ext_{A, T_{\fp}}^r(W^{\io}, \fOrd_A^s(X_{\barQ}, \sF)) \, \Longrightarrow\, \bH_A^{r+s}(X_{\barQ}; W, \sF)
\end{equation}
(here $W^{\io}=W$ with reversed $T_{\fp}$-action given by $t\cdot w = t^{-1} w$).  If $\sF$ is constructible then we show
\begin{equation}
\label{sschemecohom9}
\fOrd_A^n(X_{\barQ}, \sF) \, \cong \, \dlim_m H^n((X_{K_1(m)})_{\barQ}, \sF_{K_1(m)})^{\ord}.
\end{equation}
Here $K_1(m):=K_1(\fp^m)$ is defined as usual as the subgroup of $A\in \PGL_2(\cO_{\fp})$ with $A \equiv \left(\begin{smallmatrix} 1 & *\\ 0 & 1\end{smallmatrix}\right)\bmod \fp^n$ (and modulo the center) and the superscript ''$\ord$'' denotes the ordinary part in the sense of Hida. If $\sF$ is constructible then the groups \eqref{sschemecohom9} are admissible $A[T_{\fp}]$-modules. Moreover if $W$ is admissible then the cohomology groups $\bH_A^{\bu}(X_{\barQ}; W, \sF)$ are equipped with a canonical profinite topology and the Galois action is continuous. If $W$ is finite then the groups $\bH_A^{\bu}(X_{\barQ}; W, \sF)$ are finite as well. 
\medskip

\noi {\bf $\varpi$-adic cohomology} We will also introduce a version of $\ell$-adic cohomology adapted to our context. Let $\cO$ be the valuation ring of a $p$-adic field $E$ and put $\cO_m = \cO/(\varpi^m)$ where $\varpi$ denotes a uniformizer. Similarly to the classical case one can define the notion of an {\it {\'e}tale constructible $\varpi$-adic sheaf} $\sF= (\sF_n)_n$ on $X$. If $W$ is a $\varpi$-adically admissible $\cO[T_{\fp}]$-module without $\varpi$-torsion (see \cite{emerton2}, 2.4.7)  then we set 
\begin{equation}
\label{sschemecohom10}
\bH_{\varpi-\ad}^n(X_{\barQ}; W, \sF) : =\, \prolimm \bH^n(X_{\barQ}; W\otimes_{\cO} \cO_m, \sF_m).
\end{equation}
This defines a contravariant $\delta$-functor (in the component $\sF$). The groups \eqref{sschemecohom10} are again equipped with a canonical profinite topology and the Galois action is continuous. 

By inverting $\varpi$ we can replace $W$ in \eqref{sschemecohom10} by an admissible $E$-Banach space representation $V$ of $T_{\fp}$. Concretely, we define 
\begin{equation}
\label{sschemecohom11}
\bH_{\varpi-\ad}^{\bu}(X_{\barQ}; V, \sF\otimes E) : =\, \bH_{\varpi-\ad}^{\bu}(X_{\barQ}; W, \sF)\otimes_{\cO} E
\end{equation}
where $W$ is a $T_{\fp}$-stable lattice in $V$. The cohomology groups \eqref{sschemecohom11} define a $\delta$-functor in both variables (contravariant in $V$).  If $\sF= \cO$ is constant and if $V$ is a finite-dimensional and discrete $T_{\fp}$-representation admitting a $T_{\fp}$-stable lattice then we have 
(Prop.\ \ref{prop:compwwhat} and Lemma \ref{lemma:cohhmsdlim})
\begin{equation}
\label{sschemecohom12}
\bH_{\varpi-\ad}^{\bu}(X_{\barQ}; V, E) \, =\, \bH_{\varpi-\ad}^{\bu}(X^{\an}; V, E)\, =\, \bH_{E}^{\bu}(X^{\an}; V, E).
\end{equation}

 \noi {\bf Relation to automorphic forms} Let $\chi: T_{\fp}\cong F_{\fp}^* \to E^*$ be a unitary quasicharacter (i.e.\ $\ker(\chi)$ is open and the image of $\chi$ is contained in $\cO^*$). We consider the cohomology groups \eqref{sschemecohom12} for the one-dimensional discrete $T_{\fp}$-representation $V=E(\chi)$. We show that a cuspidal automorphic representation $\pi=\bigotimes_v' \pi_v$ of $G(\bA)$ of parallel weight $2, 2 , \ldots, 2$ occurs in the cohomology group $\bH_{E}^n(X^{\an}; E(\chi), E)$ (i.e.\ the $\pi$-isotypical component $\bH_{E}^n(X^{\an}; E(\chi), E)_{\pi}$ is non-trivial) if and only $n=d$ or $n=d+1$, the local component of $\pi$ at $\fp$ is equal to $\Ind_{B_{\fp}}^{G_{\fp}} \chi$ and if the conductor of $\pi^{\infty, \fp} = \bigotimes_{v\nmid\fp\infty}' \pi_v$ divides $\fn$ (see Prop.\ \ref{prop:automorphpart} for the precise statement). 
 
If we consider instead e.g.\ the case $V=\wE[\Theta]$ where $\wE=E[\vep]$ is the ring dual numbers over $E$ and $\Theta: F_{\fp}^* \to \wE^*$ is an infinitesimal deformation of $\chi$ of the form $\Theta(x) = \chi(x)\cdot (1+\psi(x) \vep)$, $x\in F_{\fp}^*$ with $\psi\in \Hom_{\cont}(F_{\fp}^*, \bZ_p)$, then for $n= d, d+1$ and a suitable choice of $\psi$ the $\wE[\fG]$-module
\begin{equation}
\label{sschemecohom13}
\bH_{\varpi-\ad}^n(X_{\barQ}; \wE[\Theta], E)_{\pi}
\end{equation}
is a non-trivial deformation of $\bH_{E}^n(X^{\an}; E(\chi), E)_{\pi}$.\footnote{For a precise statement in the case where $\chi=1$ see \cite{spiess2}; the techniques employed there can be generalized to prove that for an arbitrary quasicharacter $\chi$ there exists an -- up to scalar -- unique $\psi$ with given kernel $\ker(\psi|_{\cO_{\fp}*})$ such that \eqref{sschemecohom13} is a free $\wE$-module of rank $\dim_E \bH_{E}^n(X^{\an}; E(\chi), E)_{\pi}$.} Here the subscript $\pi$ indicates that we localize with respect to $\ker(\la_{\pi}: \bT\to E)$.
  
\noi {\bf Action of the decomposition group at $\fp$} Assume that $d=1$ so that each $X_K$ is either a Shimura curve (if $D$ is a division algebra) or a modular curve (if $F=\bQ$ and if $D=M_2(\bQ)$ is the $2\times 2$-matrix algebra). For the action of the decomposition group $\fG_{\fp}=\Gal(\barQ_p/F_{\fp})$ on $\bH_{\pi}:=\bH_{\varpi-\ad}^1(X_{\barQ}; \wE[\Theta], E)_{\pi}$ we show that there exists a short exact sequence of $\wE[\fG_{\fp}]$-modules
\[
0\lra \bH_{\pi}^0\lra \bH_{\pi}\lra \bH_{\pi}^{\et}\to 0 
\]
so that the $\fG_{\fp}$-action on $\bH_{\pi}^0$ factors through the maximal abelian quotient of $\fG_{\fp}$ and can be described explicitly in terms of the local reciprocity map and of $\Theta$ (for a more general and precise statement see Theorems \ref{theorem:reciprocity} and \ref{theorem:reciprocity2}).

Now we describe briefly the content of each section. In Chapters \ref{section:general} and \ref{section:borelind} we deal with certain aspects of the representation theory of $G=\PGL_2(F)$ and its maximal torus $T$ over a complete local ring $\cO$. Here $F$ is now a fixed $p$-adic field. Specifically in section \ref{subsection:ladm} we establish some basic properties of the category of locally admissible representations of certain locally profinite groups and monoids. Sections \ref{subsection:lfin} and \ref{subsection:gammaord} consists of preparatory material from commutative algebra that will be used in sections \ref{subsection:parind} and \ref{subsection:derord} where we review and study the functor $\Ord: \Mod_{\cO}^{\sm}(G) \to \Mod_{\cO}^{\ladm}(T)$ and its right derived functors. In section \ref{subsection:ladmext} we investigate the $\Ext$-groups in the category of locally admissible representations of $T$ and in section \ref{subsection:admbanach}
we recall the notions of $\varpi$-adically admissible and of admissible Banach space representations of $T$. 

The key technical result of the first two chapters is Lemma \ref{lemma:injacyclic}. It is a purely representation-theoretic result contributing towards the proof of \eqref{sschemecohom9}. It also implies (see Cor.\ \ref{corollary:adordn}) that the $n$-th right derived functor of $\Ord$ when applied to a locally admissible representation $V$ of $G$ over an Artinian local ring $A$ agrees with Emerton's $n$-th higher functor of ordinary parts introduced in \cite{emerton3}. We note though that rather than in (\cite{emerton3}, Conj.\ 3.7.2) we consider $\Ord$ here as a functor on the whole category of discrete representations $G$. However our methods can also be used to give a proof of Emerton's conjecture for $G=\PGL_2(F)$ (see Remark \ref{remarks:derivedord} (b); however the proof is simpler than that of Cor.\ \ref{corollary:adordn} as it does not rely on Lemma \ref{lemma:injacyclic}).

In chapter \ref{section:sspaces} we develop a general framework for constructing the topoi $\Sh(X_{\et})$ and $\Sh(X^{\an})$ and the cohomology groups \eqref{sschemecohom} and \eqref{sschemecohom2}. For that we consider again an arbitrary locally profinite group $G$ together with a set $\sK$ of ''small enough'' open subgroups of $G$ satisfying certain conditions and a site $\cC$ that has fibre products, equalizers and coproducts. In section \ref{subsection:gsets} we discuss some basic properties of the site $\sS=\sS_{G, \sK}$ (equipped with the canonical topology) of discrete left $G$-sets $S$ such that the stabilizer of each element in $S$ is contained in $\sK$. In section \ref{subsection:ssschemes} we introduce the notion of an $\sS$-{\it object in $\cC$}. It is a continuous functor $X:\sS\to \cC$ that commutes with coproducts. We show that there exists a natural construction of a site $\cC/\colim X$ so that $X$ ''lifts'' to a functor $\wX:  \sS\to \cC/\colim X$. For an arbitrary coefficient ring $R$ the topos $\Sh(X, R)$ is defined as the category of sheaves of $R$-modules on $\cC/\colim X$. We show that this construction is independent of the choice of $\sK$. Moreover if $\sS$ has a final object $\pt$, i.e.\ if $G\in \sK$ then $\Sh(X, R)$ can be identified with category of sheaves $R$-modules on the site $\cC/X_G$ (where $X_G:=X(\pt)$).
The category of sheaves $R$-modules on $\sS$ can be identified with the category $\Mod_R^{\sm}(G)$ of discrete $R[G]$-modules. Under this identification we denote the pair of adjoint functors $(\wX_s, \wX^s)$ by
\begin{equation*}
\label{stopos1}
\Mod_R^{\sm}(G) \lra \Sh(X, R), \,\,\, M\mapsto M_{X, G}, \qquad  H^0(X_{\infty}, \wcdot): \Sh(X, R) \lra \Mod_R^{\sm}(G).
\end{equation*}
An $\sS$-object $X$ will be called {\it exact} if it commutes with equalizers. In this case we show that  the functor $M \mapsto M_{X, G}$ is exact. In section \ref{subsection:ssmorphism} we introduce the notion 
of a morphism $f: X \to Y$ of $\sS$-objects and establish basic properties of the induced morphism of topoi $(f^*, f_*): \Sh(X, R) \to \Sh(Y,R)$. In section \ref{subsection:functinG} we consider more generally functorial properties 
of $\Sh(X, R)$ with respect to both changing $X$ and the underlying group $G$. Here the advantage of working with the category $\sS$ rather than with just the collection of open subgroups $\sK$ of $G$ comes fully to bear. In section \ref{subsection:cohomology} the cohomology groups \eqref{sschemecohom} are introduced and their basic properties discussed in the context of arbitrary $\sS$-objects. In the subsequent sections \ref{subsection:ssspacescheme} and \ref{subsection:galcohomology} we introduce the notion of an $\sS$-space and $\sS$-scheme by specifying the underlying site $\cC$ and we discuss the Galois action on the groups \eqref{sschemecohom}. In the final section \ref{subsection:uniformizable} of the chapter the isomorphism \eqref{sschemecohom6} will be established within a rather general framework. 

In Chapter \ref{section:ordinary} the general theory developed in the previous chapter is specified to the case $G=\PGL_2(F)$ (where $F$ is a $p$-adic field) and $M$ is parabolically induced. Furthermore we consider only the case where $R=\cO$ is a complete noetherian local ring. In section \ref {subsection:principal} we introduce ordinary cohomology $\fOrd_{\cO}^{\bu}(X, \sF)$ of an $\sS$-space (or scheme) $X$. If $X$ is exact then we prove that a spectral sequence of the type \eqref{sschemecohom8} exists. Moreover
we show that \eqref{sschemecohom9} holds if either $X$ is an $\sS$-scheme, $\cO=A$ is Artinian and $\sF$ is constructible or if $X$ is an $\sS$-space, $\dim(\cO)\le 1$ and the sheaf $\sF$ satisfies certain finiteness conditions. In section \ref{subsection:varpiadic} we introduce the $\varpi$-adic cohomology groups \eqref{sschemecohom10} and \eqref{sschemecohom11}.

Finally, in Chapter \ref{section:cohomqhmsv} we consider the case of a quaternionic Hilbert modular $\sS$-variety $X$. By applying $X$ to $S= G_{\fp}/K$ with $K$ sufficiently small the variety $X(S)$ is equal to the Hilbert modular variety $X_K$ considered at the beginning of the introduction. In section \ref{subsection:qhmsv} we prove the existence of $X$ and show that it is an exact $\sS$-scheme. Occasionally, we also consider the case $d=0$, i.e.\ when $D$ is a totally definite quaternion algebra. In this case there exists only an associated ($0$-dimensional) $\sS$-space $X=X^{\an}$ and it is not exact anymore. If $d\ge 1$ then the varieties $X_K$ are proper except if $D=M_2(F)$. In this case one can view the collection of Borel-Serre compactifications of each $X(S)^{\an}$ for $S\in \sS$ also as an $\sS$-space. We denote it by $X^{\BS}$. In section \ref{subsection:borelserre} we study the cohomology of the boundary $\partial X=X^{\BS}\setminus X$ and in two subsequent we establish a variant of \eqref{sschemecohom6} for $\varpi$-adic cohomology, 
prove \eqref{sschemecohom12} and investigate which automorphic representations $\pi$ ''occur'' in the cohomology $\bH_{\varpi-\ad}^{\bu}(X_{\barQ}; E(\chi), E)=\bH_{E}^{\bu}(X^{\an}; E(\chi), E)$ when $\chi$ is a unitary quasicharacter. In section \eqref{subsection:shimuracurve} we consider the case $d=0,1$. We study in particular the groups 
\begin{equation}
\label{sschemecohom14}
\fOrd_{\varpi-\ad}^d(X^{\an}, \cO): =\, \prolimm \fOrd_{\cO_m}^d(X^{\an}, \cO_m) 
\end{equation}
where $\cO$ denotes again the valuation ring of a $p$-adic field $E$. We prove that \eqref{sschemecohom14} is a $\varpi$-adically admissible $\cO[T_{\fp}]$-module and that the dual of the associated $E$-Banach space $\fOrd_{\varpi-\ad}^d(X^{\an}, \cO)\otimes_{\cO} E$ is a module under a certain Iwasawa algebra $\La\otimes_{\cO} E$ and as such is finitely generated and projective (see Prop.\ \ref{prop:freetotdeg} for the precise statement). 

In the final two sections we consider the case $d=1$. We consider in particular the action of the decomposition group $\fG_{\fp}$ on $\bH_{\varpi-\ad}^1(X_{\barQ}; A(\chi), E)$ where $A$ is an Artinian local $E$-algebra and $\chi: T_{\fp}\to A^*$ is a bounded character. Using \eqref{sschemecohom8} and \eqref{sschemecohom9} it will be shown that for any $E$-Banach space representation $V$ of $T_{\fp}$ we have
\[
\bH_{\varpi-\ad}^1(X_{\barQ}; V, E) = \Hom_{E[T_{\fp}]}(V^{\io},\fOrd_{\varpi-\ad}^1(X_{\barQ}, \cO)\otimes_{\cO} E)
\]
and that $\fOrd_{\varpi-\ad}^1(X_{\barQ}, \cO)$ is the $\varpi$-adic completion of $\dlim_n \prolimm H^1((X_{K_1(n)})_{\barQ}, \cO_m)^{\ord}$.
This together with a trick of Hida allows us to deduce Theorem \ref{theorem:reciprocity} from a result of Wiles (\cite{wiles1}, Thm.\ 2.2).

\paragraph{{\it Acknowledgement}} I thank Lennart Gehrmann, Eike Lau, Vytautas Paskunas and Thomas Zink for helpful communications. 

\section{Representations of locally profinite groups}
\label{section:general}

\subsection{Notation and preliminary remarks}
\label{subsection:preremrep} 

Unless otherwise stated all rings are commutative with $1\ne 0$. If $\varphi: R\to R'$ is a ringhomomorphism and $N$ an $R$-module then we put $N_{R'} = N\otimes_R R'$. If 
$M$ is an $R'$-module then we can view it as an $R$-module via $\varphi$ and as such denote it by $M_{\varphi}$.

Let $R$ be noetherian ring and let $\cG$ be a locally profinite group or a locally profinite monoid. By the latter we mean that $\cG$ is a topological monoid so that its maximal subgroup $\cG^0$ is open and locally profinite. If $\cG$ is only a monoid then we will assume that $g\cG^0 = \cG^0g$ holds for every $g\in \cG$. 

The category of left $R[\cG]$-modules will be denoted by $\Mod_R(\cG)$. 
A left $R[\cG]$-module $V$ is called discrete (or smooth) if $V$ is discrete as a $\cG$-module, i.e.\ if the stabilizer $\Stab_{\cG}(v)=\{g\in \cG^0\mid gv=g\}$ is open in $\cG$ for every $v\in V$. A discrete $R[\cG]$-module $V$ is called admissible if $V^U$ is a finitely generated $R$-module for every open subgroup $U$ of $\cG^0$. 

The full subcategory of $\Mod_R(\cG)$ of discrete $R[\cG]$-modules will be denoted by $\Mod_R^{\sm}(\cG)$. Furthermore we denote by $\Mod_{R, f}^{\sm}(\cG)$ respectively $\Mod_R^{\adm}(\cG)$ the subcategory of discrete $R[\cG]$-modules that are finitely generated as $R$-modules resp.\ admissible. The category $\Mod_R^{\sm}(\cG)$ is abelian and $\Mod_{R, f}^{\sm}(\cG)$ is a Serre subcategory of $\Mod_R^{\sm}(\cG)$.

Then the category $\Mod_R^{\sm}(\cG)$ has enough injectives. Indeed, the category $\Mod_R(\cG)$ has enough injectives and the embedding $\Mod_R^{\sm}(\cG)\hra \Mod_R(\cG)$ is exact and left adjoint to the functor 
\[
\Mod_R(\cG)\lra \Mod_R^{\sm}(\cG), \quad V \mapsto V_{\sm}\, =\, \bigcup_{U\le \cG^0, U \open} V^U.
\]
So if $V\in \Mod_R^{\sm}(\cG)$ and $V\hra I$ is a monomorphism in $\Mod_R(\cG)$ with $I$ injective then it factors as $V\to I_{\sm} \to I$ and $I_{\sm}$ is an injective object of $\Mod_R^{\sm}(\cG)$. 

We note the following useful  

\begin{lemma}
\label{lemma:torsionfree}
Let $G$ be a second countable locally profinite group and let $R$ be a principal ideal domain. Let $W\in \Mod_{\cO}^{\adm}(G)$ and assume that $W$ is torsionfree as an $R$-module. Then $W$ is a free $R$-module. 
\end{lemma}

\begin{proof} Since $G$ is second countable there exists a sequence of open subgroups of $K_0\supseteq K_1 \supseteq K_2 \supseteq \ldots \supseteq K_n \supseteq \ldots$ f $G$ with $\bigcap_{n\ge 0} K_n = 1$. The $R$-submodule $W^{K_n}$ of $W$ is finitely generated and torsionfree hence free. Also the fact that $W$ is torsionfree implies that $P_n:= W^{K_{n+1}}/W^{K_n}$ is a torsionfree, hence free $R$-module as well. Therefore we have $W^{K_{n+1}}\cong W^{K_n}\oplus P_n$ as $R$-modules. It follows $W=\bigcup_{n\ge 0} W^{K_n}\cong W^{K_0}\oplus \bigoplus_{n\ge 0} P_n$,
so $W$ is free as an $R$-module.
\end{proof}

\subsection{Locally admissible $R[\cG]-$modules}
\label{subsection:ladm}

Let $\cG$ be a locally profinite group or a locally profinite monoid. A discrete $R[\cG]$-module $V$ is called {\it locally admissible} if $V$ is the union of its admissible submodules. Equivalently, $V$ is locally admissible if the submodule $R[\cG] v$ is admissible for every $v\in V$. The full subcategory of $\Mod_R^{\sm}(\cG)$ of locally admissible $R[\cG]$-modules will be denoted by $\Mod_R^{\ladm}(\cG)$. We will show that the category $\Mod_R^{\ladm}(\cG)$ is abelian under the following assumption.

\begin{assum}
\label{assum:normalsubgroup}
Every open subgroup $U$ of $\cG$ contains an open normal subgroup of $\cG$.\footnote{A subgroup $H$ of a monoid $\cG$ is normal if $gH = Hg$ for every $g\in G$.}
Moreover for any open normal subgroup $N$ of $\cG$ the quotient group (resp.\ monoid) $\cG/N$ is finitely generated.\end{assum}

Note that \ref{assum:normalsubgroup} holds if $\cG$ is a profinite group but not e.g.\ for $\cG=\PGL_2(\bQ_p)$. 

\begin{lemma}
\label{lemma:tadm}
Assume that $\cG$ satisfies condition \ref{assum:normalsubgroup} and let $V\in \Mod_R^{\sm}(\cG)$. The following conditions are equivalent

\noi (i) $V$ is a finitely generated and admissible $R[\cG]$-module.

\noi (ii) $V$ is finitely generated as an $R$-module.
\end{lemma}

\begin{proof} Assume that (i) holds and let $v_1, \ldots, v_m\in V$ with $V=\sum_{i=1}^r R[\cG] v_i$. There exists an open normal  subgroup $U$ of $\cG$ such that $v_1, \ldots, v_m\in V^U$. Since $U$ is normal we have $gv_i \subset V^U$ for $i=1,\ldots, m$ and $g\in \cG$ hence also $V=V^U$ and therefore $V\in \Mod_{R,f}^{\sm}(\cG)$. The converse implication is obvious.
\end{proof}

It follows that a discrete $R[\cG]$-module $V$ is locally admissible if and only if $V$ is the union of its submodules that are finitely generated as $R$-modules. For an arbitrary discrete $R[\cG]$-module $V$ let $V_{\ladm}$ denote its maximal locally admissible $R[\cG]$-submodule, i.e.\
\begin{equation*}
\label{ladmW}
V_{\ladm} \, =\, \sum_{V'\in \sX} V' \, =\, \bigcup_{V'\in \sX} V'
\end{equation*} 
where $\sX$ be the collection of its $R[\cG]$-submodules $V'\subseteq V$ with $V'\in \Mod_{R,f}^{\sm}(\cG)$. 

\begin{lemma}
\label{lemma:locadm1}
Assume that $\cG$ satisfies condition \ref{assum:normalsubgroup}.
\medskip

\noi (a) The category $\Mod_R^{\ladm}(\cG)$ is a Serre subcategory of $\Mod_R^{\sm}(\cG)$. It is in particular an abelian category. 
\medskip

\noi (b) The category $\Mod_R^{\ladm}(\cG)$ has enough injectives. 
\end{lemma}

\begin{proof} (a) Let $0 \lra V'  \stackrel{\alpha}{\lra}V \stackrel{\beta}{\lra} V'' \lra 0$ be a short exact sequence in $\Mod_R^{\sm}(\cG)$. We have to show that $V\in \Mod_R^{\ladm}(\cG)$ if and only if $V', V''\in \Mod_R^{\ladm}(\cG)$. Firstly, assume that $V$ is locally admissible. Then $V'$ is obviously locally admissible. That $V''$ is locally admissible follows immediately from the fact that quotients of objects in $\Mod_{R,f}^{\sm}(\cG)$ lie again in $\Mod_{R,f}^{\sm}(\cG)$. Conversely, assume that $V', V''\in \Mod_R^{\ladm}(\cG)$. Let $v\in V$ and put $V_0 = R[\cG]v$, $V_0' = \alpha^{-1}(V_0)$ and $V_0''= \beta(V_0)$. Let $U$ be an open normal subgroup of $\cG$ with $v\in V^U$. Then $V_0'$ is a submodule of the $R[\cG/U]$-module $V_0$, hence it is finitely generated as an $R[\cG/U]$-module (since the ring $R[\cG/U]$ is noetherian). Therefore as an $R[\cG]$-module $V_0'$ is finitely generated and it is admissible (as a submodule of $V'$). Also, as submodule of $V''$, the $R[\cG]$-module $V_0''= R[\cG] \beta(v)$ is admissible. Hence the extension $V_0$ of $V_0''$ by $V_0'$ is admissible, too. 

(b) Note that the functor 
\begin{equation}
\label{ladmW2}
\Gamma_{\ladm}: \Mod_R^{\sm}(\cG)\lra \Mod_R^{\ladm}(\cG), \,\, V \mapsto V_{\ladm}
\end{equation} 
is right adjoint to the embedding $\Mod_R^{\ladm}(\cG)\hra \Mod_R^{\sm}(\cG)$ hence it maps injectives objects to injectives objects. Let $V\in \Mod_R^{\ladm}(\cG)$. There exists a monomorphism of $R[\cG]$-modules $\mu: V\hra I$ with $I$ an injective object in $\Mod_R^{\sm}(\cG)$. Then $\mu$ factors in the form $V\hra I_{\ladm}\hra I$ and $I_{\ladm}$ is an injective object in $\Mod_R^{\ladm}(\barT)$.
\end{proof}

\subsection{Locally finite $R[X, X^{-1}]$-modules} 
\label{subsection:lfin}

We now consider the case $\cG =\bZ$, so that the group ring $R[\cG]$ can be identified with the localization $R[X, X^{-1}] = S^{-1}R[X]$ of the polynomial ring $R[X]$ with respect to $S= \{X^n \mid n \ge 1\}$. We will write $\Mod_{R[X,X^{-1}]}^{\lfin}$ instead of $\Mod_R^{\ladm}(\bZ)$ and call an object $M$ of $\Mod_{R[X,X^{-1}]}^{\lfin}$ a {\it locally finite} $R[X, X^{-1}]$-module. Note that an $R[X, X^{-1}]$-module $M$ is locally finite if $R[X, X^{-1}]m$ is finitely generated as an $R$-module for every $m\in M$. For 
$M\in \Mod_{R[X,X^{-1}]}^{\lfin}$ we denote by $M_{\lfin}$ its maximal locally $R[X, X^{-1}]$-finite submodule. 

For $\cG=\bZ$ we will study the left exact functor \eqref{ladmW2}, i.e.\ the functor 
\begin{equation}
\label{lfin}
\Gamma_{\lfin}: \Mod_{R[X,X^{-1}]}\lra \Mod_{R[X,X^{-1}]}^{\lfin}, \quad M \mapsto M_{\lfin}.
\end{equation} 
Let $S_1$ (resp.\ $S_2$) be the multiplicative subset of $R[X,X^{-1}]$ consisting of polynomials $a_n X^n + a_{n-1} X^{n-1} + \ldots + a_0\ne 0$ (resp.\ $a_n X^{-n} + a_{n-1} X^{-n+1} + \ldots + a_0\ne 0$) with leading coefficient $a_n=1$.

\begin{lemma}
\label{lemma:lfinchar}
For $M \in \Mod_{R[X,X^{-1}]}$ we have 
\begin{equation*}
\label{lfin2}
M_{\lfin}\, =\, \ker\left(M\lra S_1^{-1} M \oplus S_2^{-1} M, \, m\mapsto \left(\frac{m}{1}, \frac{m}{1}\right)\right).
\end{equation*} 
\end{lemma}

\begin{proof} For $m\in M$ we have
\begin{eqnarray*}
\label{lfin3}
m\in M_{\lfin} & \Leftrightarrow & R[X, X^{-1}]m = R[X]m + R[X^{-1}]m \in \Mod_{R, f}\\
& \Leftrightarrow & R[X]m ,\, R[X^{-1}] m \in \Mod_{R, f}\nonumber\\
& \Leftrightarrow & \exists \, n\ge 1: \,\, \sum_{i=0}^{n-1} R\, X^i m = R[X] m \quad\mbox{and} \quad \sum_{i=0}^{n-1} R\, X^{-i} m = R[X^{-1}] m\nonumber \\
& \Leftrightarrow & \exists \,s_1\in S_1, s_2\in S_2: \,\, s_1\cdot m \, =\, 0\, =\, s_2 \cdot m.\nonumber
\end{eqnarray*} 
\end{proof}

As a consequence of the above characterization of the functor $\Gamma_{\lfin}$ we get

\begin{prop}
\label{prop:rnlfin}
(a) The embedding $\Mod_{R[X,X^{-1}]}^{\lfin}\hra \Mod_{R[X,X^{-1}]}$ preserves injectives. 
\medskip

\noi (b) If $M\in \Mod_{R[X,X^{-1}]}^{\lfin}$ then $M$ is $\Gamma_{\lfin}$-acyclic, i.e.\ we have $(R^n \Gamma_{\lfin})(M)=0$ for $n\ge 1$.
\end{prop}

We need the following 

\begin{lemma}
\label{lemma:injlocs}
Let $R$ be a noetherian ring, let $S$ be a multiplicative subset of $R$ consisting of non-zero divisors and let $I$ be an injective $R$-module. We have 
\medskip

\noi (a) $J:= \ker\left(I\to S^{-1} I, \, x\mapsto \frac{x}{1}\right)$ is an injective $R$-module. 
\medskip

\noi (b) The sequence of $R$-modules 
\begin{equation}
\label{injlocs2}
\begin{CD}
0@>>> J @> \incl >> I@> x\mapsto \frac{x}{1}>> S^{-1} I@>>> 0
\end{CD}
\end{equation} 
is exact and splits. 
\end{lemma}

\begin{proof} (a) It suffices to see that for any ideal $\fb\subseteq R$ and any homomorphism $\varphi: \fb \to J$ there exists a homomorphism $\psi: R \to J$
extending $\varphi$. Since $\fb$ is finitely generated and $\varphi(\fb)$ lies in $J$ there exists $s\in S$ such that $s \cdot \varphi(\fb)=0$. Let $\fa= R\cdot s$ and put $\Gamma_{\fa}(I)= \{x\in I \mid \exists\, n\ge 1:\, \fa^n x=0\} =\{x\in I\mid \exists\, n\ge 1:\, s^n\cdot x=0\}$. Note that $\varphi(\fb)\subseteq \Gamma_{\fa}(I) \subseteq J$. By (\cite{hartshorne}, Ch.\ III, Lemma 3.2), $\Gamma_{\fa}(I)$ is an injective $R$-module. Hence there exists a homomorphism $\psi: R \to \Gamma_{\fa}(I)$ that extends $\varphi$. Since $\Gamma_{\fa}(I)\subseteq J$ we conclude that $J$ is an injective $R$-module.

(b) Since $s\in S$ is not a zero divisor the map $I\to I, x\mapsto s\cdot x$ is surjective. 
Hence the second map in \eqref{injlocs2} is surjective. That the sequence splits follows from (a). 
\end{proof}

\begin{proof}[Proof of Prop.\ \ref{prop:rnlfin}] (a) Let $I$ be an injective $R[X,X^{-1}]$-module. By applying Lemma \ref{lemma:injlocs} twice we see that $I_{\lfin}$ is an injective $R[X,X^{-1}]$-module as well. Indeed firstly $I_0:=\ker(I \to S_1^{-1}I)$ is an injective $R[X,X^{-1}]$-module, so by Lemma \ref{lemma:lfinchar} the $R[X,X^{-1}]$-module $I_{\lfin}=\ker(I_0 \to S_2^{-1}I_0)$ is injective, too.

Now let $J$ be an injective object in $\Mod_{R[X,X^{-1}]}^{\lfin}$ and let $\mu: J\hra I$ be a monomorphism into an injective $R[X, X^{-1}]$-module $I$. Then $\mu$ factors as $J\hra I_{\lfin} \hra I$. Since $J$ is in $\Mod_{R[X,X^{-1}]}^{\lfin}$ it is a direct summand of $I_{\lfin}$. Since the latter is an injective $R[X,X^{-1}]$-module the same holds for $J$. 

(b) follows immediately from (a) and the fact that $\Gamma_{\lfin}$ is the identical functor when restricted to $\Mod_{R[X,X^{-1}]}^{\lfin}$. 
\end{proof}

Now assume that $R=\cO$ is a complete noetherian local ring with maximal ideal $\fm$ and residue field $k$. Let $M$ be an $\cO[X, X^{-1}]$-module that -- as an $\cO$-module -- has support contained in $\{\fm\}$ (i.e.\ for every $m\in M$, $m\ne 0$ we have $\sqrt{\Ann_{\cO}(m)} =\fm$). In this case the submodule $M_{\lfin}$ admits a simpler description than the one given in Lemma \ref{lemma:lfinchar}. For that consider the multiplicative subset 
\begin{equation*}
\label{lin4}
S\,=\, \{x\in \cO[X,X^{-1}]\mid \bar{x}\ne 0\}
\end{equation*}
of $\cO[X,X^{-1}]$ where $x\mapsto \bar{x}$ denotes the residue map mod $\fm$, i.e.\ the map
\begin{equation*}
\label{modfm}
\cO[X,X^{-1}]\lra k[X,X^{-1}],\,\, x\mapsto \bar{x}:=x\mod \fm.
\end{equation*} 

\begin{lemma}
\label{lemma:lfinchar2}
Let $M$ be an $\cO[X,X^{-1}]$-module whose support -- as an $\cO$-module -- is contained in $\{\fm\}$. We have 
\begin{equation}
\label{lfin5}
M_{\lfin}  \, =\, \ker\left(M\to S^{-1} M, \, m\mapsto \frac{m}{1}\right)
\end{equation} 
\end{lemma}

\begin{proof} Since $S_1, S_2\subseteq S$ we have $M_{\lfin}\subseteq \ker\left(M\to S^{-1} M, \, m\mapsto \frac{m}{1}\right)$ by Lemma \ref{lemma:lfinchar}. Conversely, let $m\in M$ so that there exists $s\in S$ with $s\cdot m=0$. By Lemma \ref{lemma:lfinchar} we have to show that there exists $s_1\in S_1$, $s_2\in S_2$ with $s_1\cdot m=0=s_2\cdot m$. The assumption on $M$ implies that there exists a positive integer $n$ such that $\fm^n \cdot m = 0$. 

After multiplying $s$ with a high enough power of $X$ we may assume $s\in \cO[X]$. We decompose $\bar{s}\in k[X]$ in the form $\bar{s} = g_0 \cdot X^n$ where $n\ge 0$ and $g_0 \in k[X]$ is a polynomial with non-zero constant term. By Hensel's Lemma there exists polynomials $g, h\in \cO[X]$ with $s=g\cdot h$, $\bar{g}=g_0$, $\deg(g) = \deg(g_0)=:d$ and $\bar{h} = X^n$. The element $h$ is -- modulo $\fm^n$ -- a unit (i.e.\ there exists $h'\in \cO[X,X^{-1}]$ such that $h h'-1\in \fm^n \cO[X,X^{-1}]$) since $\bar{h}=X^n\in k[X,X^{-1}]$ is a unit and since the ideal $\fm \cO[X,X^{-1}]/\fm^n\cO[X,X^{-1}]$ of $(\cO/\fm^n)[X,X^{-1}]$ is nilpotent. Hence $g\cdot m=h'\cdot g\cdot h\cdot m= h'\cdot s \cdot m= 0$. We write $g = a_d X^d + \ldots + a_1 X + a_0$. The assumptions $\deg(g)=\deg(g_0)$, $\bar{g} =g_0$ and $g_0(0)\ne 0$ imply that $a_d, a_0$ are units in $\cO$. Hence $s_1:=a_d^{-1} g\in S_1$, $s_2:=a_0^{-1} T^{-d} g\in S_2$ and $s_1\cdot m=0=s_2\cdot m$. This proves \eqref{lfin5}.
\end{proof}

For the right derived functors of \eqref{lfin} we deduce the following 

\begin{lemma}
\label{lemma:lfinchar3}
Let $\cO$ be a complete noetherian local ring with maximal ideal $\fm$ and residue field $k$ and let $M$ be an $\cO[X,X^{-1}]$-module. Assume that $\dim(\cO) \le 1$ and that the support of $M$ -- as an $\cO$-module -- is contained in $\{\fm\}$. Then for $n\ge 0$ we have 
\begin{equation}
\label{lfin6}
R^n\Gamma_{\lfin} M \, =\, \left\{ \begin{array}{cc}\ker\left(M\to S^{-1} M, \, m\mapsto \frac{m}{1}\right) & \mbox{if $n=0$,}\\
\coker\left(M\to S^{-1} M, \, m\mapsto \frac{m}{1}\right) & \mbox{if $n=1$,}\\
0 & \mbox{if $n\ge 2$.}
\end{array}\right.
\end{equation} 
\end{lemma}

\begin{proof} Firstly, assume that $\dim(\cO) =0$, i.e.\ $\cO$ is Artinian. 
Let $0 \to M \to I^{\bu}$ be an injective resolution. By Lemmas \ref{lemma:lfinchar2} and \ref{lemma:injlocs} the sequence of complexes $0\to I^{\bu}_{\lfin}\to I^{\bu}\to S^{-1} I^{\bu}\to 0$
is exact. Passing to the long exact cohomology sequence and yields \eqref{lfin6}.

Now assume that $\dim(\cO) =1$ and let $S_0$ be the multiplicative set of non-zero divisors of $\cO$ so that $S_0^{-1} \cO = \prod_{\fp\in \Spec \cO, \fp\ne \fm} \cO_{\fp}$. Therefore an arbitrary $\cO$-module $N$ has support contained in $\{\fm\}$ if and only if $S_0^{-1} N = 0$. Let $0 \to M \to I^{\bu}$ be an injective resolution and put $J^{\bu} = \ker(I^{\bu}\to S_0^{-1} I^{\bu})$. Since, by assumption, we have $S_0^{-1} M = 0$ , Lemma \ref{lemma:injlocs} implies that $0\to M\to J^{\bu}$ is an injective resolution as well. Moreover for every $n\ge 0$ the support of $J^n$ -- as an $\cO$-module -- is contained in $\{\fm\}$. Again by Lemmas \ref{lemma:lfinchar2} and \ref{lemma:injlocs} the sequence of complexes $0\to J^{\bu}_{\lfin}\to J^{\bu}\to S^{-1} J^{\bu}\to 0$ is exact and the associated long exact cohomology sequence and yields \eqref{lfin6}.
\end{proof}

Let $R$ be again an arbitrary noetherian ring. We will now study certain properties of the right derived functors of the right adjoint of the forgetful functor 
\begin{equation}
\label{lfinxinv}
\Mod_{R[X,X^{-1}]}^{\lfin} \lra \Mod_{R[X]}, \quad N \mapsto N.
\end{equation} 
Our results will be used in section \ref{subsection:derord} and when we the study the derived functors of Emerton's functor $\Ord$ and in section \ref{subsection:principal} where we consider ordinary cohomology.

We need the following simple

\begin{lemma}
\label{lemma:adlfinxinv}
(a) The functor 
\begin{equation}
\label{homxinv}
\Mod_{R[X]}\lra \Mod_{R[X,X^{-1}]}, \quad M \mapsto \Hom_{R[X]}(R[X,X^{-1}], M)
\end{equation} 
is right adjoint to the forgetful functor $\Mod_{R[X,X^{-1}]} \lra \Mod_{R[X]}$.
\medskip

\noi (b) Let $M\in \Mod_{R[X]}$ and $q\ge 0$. We have 
\[
\Ext_{R[X]}^q(R[X,X^{-1}], M) \, =\, \prolimn^{(q)} (\ldots  \lra M \stackrel{X\wcdot}{\lra} M  \lra \ldots \stackrel{X\wcdot}{\lra} M).
\]
In particular we have $\Ext_{R[X]}^q(R[X,X^{-1}], M)=0$ if $q\ge 2$.
\medskip

\noi (c) Let $\ldots \to M_{+1} \to M_n \to \ldots \to M_2 \to M_1$ be an inverse system of $R[X]$-modules such that $X^n M_n =0$ for all $n\ge 1$ and let $M= \prolimn M_n$. Then we have 
\[
\Ext_{R[X]}^q(R[X,X^{-1}], M) \, =\, 0
\]
for every $q\ge 0$.
\end{lemma}

\begin{proof} (a) is obvious. For (b) we remark that 
\[
\Hom_{R[X]}(R[X,X^{-1}], M) \,=\, \prolimn\, (\ldots \stackrel{X\wcdot}{\lra}  M\stackrel{X\wcdot}{\lra} M \stackrel{X\wcdot}{\lra}  \ldots \stackrel{X\wcdot}{\lra} M)
\]
i.e.\ the functor \eqref{homxinv} is the composite of the exact functor 
\begin{equation*}
\label{extrt0rt0inv}
\Mod_{R[X]} \lra (\Mod_{R[X]})^{\bN},\quad M \mapsto  (\ldots \lra M\stackrel{X\wcdot}{\lra} M\lra   \ldots \stackrel{X\wcdot}{\lra} M)
\end{equation*}
with $\prolimn: (\Mod_{R[X]})^{\bN}\to \Mod_{R[X]}$. Here $(\Mod_{R[X]})^{\bN}$ denotes the category of inverse systems $ \ldots \to M_n \to M_{n-1} \to \ldots \to M_0$ of $R[X]$-modules.  Note that an injective $R[X]$-module $I$ is mapped onto a $\prolimn$-acyclic objects in $\Mod_{R[T]}^{\bN}$ since the map $I\to I, m\mapsto X\cdot m$ is surjective (see \cite{weibel}, Cor.\ 3.5.4 and Prop.\ 3.5.7). The assertion follows. 

(c) The group $\Ext_{R[X]}^q(R[X,X^{-1}], M_n)$ vanishes for every $q\ge 0$ and $n\ge 1$ since $X^n$ acts as an isomorphism and annihilates it. It follows 
\[
\Ext_{R[X]}^q(R[X,X^{-1}], M) \,= \, \prolimn \Ext_{R[X]}^q(R[X,X^{-1}], M_n)\, =\, 0
\]
for every $q\ge 0$.  
\end{proof}

Now we discuss the existence and basic properties of the right adjoint of \eqref{lfinxinv}.

\begin{lemma}
\label{lemma:adlfinxinv2}
(a) The functor 
\begin{equation}
\label{adlfintinv3}
\Gamma_R^{X\sord}: \Mod_{R[X]}\lra \Mod_{R[X,X^{-1}]}^{\lfin}, \,\, M \mapsto \Gamma^{X\sord}(M) = 
\Hom_{R[X]}(R[X,X^{-1}], M)_{\lfin}
\end{equation} 
is right adjoint to \eqref{lfinxinv}.
\medskip

\noi (b) The functor \eqref{adlfintinv3} commutes with direct limits. 
\medskip

\noi (c) There exists a spectral sequence
\begin{equation}
\label{adlfintinvss}
E_2^{rs} = R^r \Gamma_{\lfin} \Ext_{R[X]}^s(R[X,X^{-1}], M)\, \Longrightarrow\, E^{r+s} = (R^{r+s} \Gamma_R^{X\sord})(M)
\end{equation}
for every $M\in \Mod_{R[X]}$.
\end{lemma}

\begin{proof} (a) follows immediately from Lemma \ref{lemma:adlfinxinv} (a). (b) is essentially (\cite{emerton2}, Lemma 3.2.2 (2)). The proof given there can be easily adapted to our slightly different framework. For (c) note that \eqref{adlfintinv3} can be factored as 
\begin{equation}
\label{adlfintinv4}
\begin{CD}
\Mod_{R[X]}@>\Hom_{R[X]}(R[X,X^{-1}],\wcdot ) >> \Mod_{R[X,X^{-1}]} @> \Gamma_{\lfin} >>  \Mod_{R[X, X^{-1}]}^{\lfin}.
\end{CD}
\end{equation} 
By Lemma \ref{lemma:adlfinxinv} (a) the first functor preserves injectives. So there exists a Grothendieck spectral sequence \eqref{adlfintinvss} associated to the decomposition \eqref{adlfintinv4}. 
\end{proof}

\begin{remark}
\label{remark:rnXordlim} 
\rm More generally the $q$-th right derived functor of $\Gamma^{X\sord}$ commutes with direct limits for every $q\ge 0$. We only need a special case namely if $(M_n)_{n\ge 1} = (M_1\to M_2 \to \ldots \to M_n \to \ldots)$ is a direct system of $R[X]$-modules then we have 
\begin{equation}
\label{rnXordlim2}
\dlim_n R^q \Gamma^{X\sord} M_n \, \cong \, R^q \Gamma^{X\sord} \dlim_n M_n
\end{equation}
for ever $q\ge 0$. In fact it is easy to construct a sequence of direct systems of $R[X]$-modules
\begin{equation}
\label{liminjres}
0 \to (M_n)_{n\ge 1}\to (I_n^0)_{n\ge 1}\to (I_n^1)_{n\ge 1}\to \ldots 
\end{equation}
such that $0\to M_n \to I_n^0 \to I_n^1\to \ldots $ is an injective resolution of $M_n$ for every $n\ge 1$. 
Since $R[X]$ is noetherian, passing to direct limits in \eqref{liminjres} yields an injective resolution of $\dlim_n M_n$. By applying $\Gamma^{X\sord}$ to it and using then Lemma \ref{lemma:adlfinxinv2} (b) yields \eqref{rnXordlim2}.
\enddemo
\end{remark} 

Let $\varphi: R'\to R$ be a ring homomorphism between noetherian rings. It induces a homomorphism $R'[X] \to R[X]$, $R'[X,X^{-1}]\to R[X,X^{-1}]$ which -- by abuse of notation -- will be denoted by $\varphi$ as well. Note that for $M\in \Mod_{R[X]}$ we have $\Gamma_{R}^{X\sord}(M)_{\varphi}= \Gamma_{R'}^{X\sord}(M_{\varphi})$. 
Since $M\mapsto (R^{\bu} \Gamma_{R}^{X\sord})(M)_{\varphi}$ and $M\mapsto (R^{\bu} \Gamma_{R'}^{X\sord})(M_{\varphi})$ are both $\delta$-functors and the first one is universal there exists canonical homomorphisms 
\begin{equation}
\label{adlfinder1}
(R^n \Gamma_{R}^{X\sord})(M)_{\varphi}\, \lra \, (R^n \Gamma_{R'}^{X\sord})(M_{\varphi})
\end{equation}
for every $M\in \Mod_{R[X]}$ and $n\ge 0$. 

\begin{prop}
\label{prop:artinadlfinder}
Let $\varphi: R'=\cO'\to R=\cO$ be an epimorphism of complete noetherian local rings with maximal ideals $\fm'$ and $\fm$ respectively and let $M\in \Mod_{\cO[X]}$. Assume that $\dim(\cO') \le 1$ and that the support of $M$ -- as an $\cO$-module -- is contained in $\{\fm\}$. Then the homomorphism \eqref{adlfinder1} is an isomorphism for every $n\ge 0$. In particular if $M$ is $\Gamma_{\cO}^{X\sord}$-acyclic then $M_{\varphi}$ is $\Gamma_{\cO'}^{X\sord}$-acyclic.
\end{prop}

\begin{proof} By Remark \ref{remark:rnXordlim}, since $M= \dlim_n M[\fm^n]$ where $M[\fm^n]= \Hom_{\cO}(\cO/\fm^n, M)$, it suffices to prove the assertion for each $\cO[X]$-module $M[\fm^n]$. Thus we may assume $\fm^n M= 0$ for some $n\ge 0$. 

It is easy to see that there exists a morphism of spectral sequences 
\begin{eqnarray*}
&& \left(E_2^{rs}=R^r \Gamma_{\lfin} \Ext_{\cO[X]}^s(\cO[X,X^{-1}], M)_{\varphi} \, \Longrightarrow\, E^{r+s} = (R^{r+s} \Gamma_{\cO}^{X\sord})(M)_{\varphi} \right) \lra \\
&& \hspace{2cm} \left(E_2^{rs}=R^r \Gamma_{\lfin} \Ext_{\cO'[X]}^s(\cO'[X,X^{-1}], M_{\varphi}) \, \Longrightarrow\, E^{r+s} = (R^{r+s} \Gamma_{\cO'}^{X\sord})(M_{\varphi})\right)
\end{eqnarray*}
where the induced maps on the limit terms are the maps \eqref{adlfinder1} (compare with the proof of Prop.\ \ref{prop:hsetsing} in section \ref{subsection:hxmf}). So it suffices to see that the induced maps on the $E_2$-terms are isomorphisms. This follows immediately from Lemmas \ref{lemma:lfinchar3} and Lemma \ref{lemma:adlfinxinv} (b).
\end{proof}

\begin{prop}
\label{prop:finomodacyclic}
Let $R=\cO$ be a complete noetherian local ring and let $M$ be an $\cO[X]$-module that is finitely generated as an $\cO$-module. For $n\ge 0$ we have 
\begin{equation*}
\label{finomod1}
(R^n\Gamma_{\cO}^{X\sord})(M) \, =\, \left\{ \begin{array}{cc}\bigcap_{i\ge 1} X^i M & \mbox{if $n=0$,}\\
0 & \mbox{if $n\ge 1$.}
\end{array}\right.
\end{equation*} 
\end{prop}

We need two Lemmas

\begin{lemma}
\label{lemma:ext0os}
(a) The $\cO[X]$-submodule $M^{X\sord} :=\bigcap_{q\ge 1} X^q M$ of $M$ is a direct summand of $M$, i.e.\ there exists a (unique) $\cO[X]$-submodule $M'$ of $M$ such that  
\begin{equation*}
\label{ord3a}
M\, =\, M^{X\sord} \oplus M'.
\end{equation*}
Moreover the $\cO[X]$-action on $M^{X\sord}$ extends to a $\cO[X,X^{-1}]$-action.
\medskip

\noi (b) The map $\ev:\Hom_{\cO[X]}( \cO[X,X^{-1}], M)\to M^{X\sord}, \Psi\mapsto \Psi(1)$ 
is injective with image $M^{X\sord}$. It induces an isomorphism of $\cO[X, X^{-1}]$-modules
\begin{equation}
\label{eval}
\Hom_{\cO[X]}( \cO[X,X^{-1}], M)\lra M^{X\sord}
\end{equation}
In particular we have $\Hom_{\cO[X]}( \cO[X,X^{-1}], M) \in \Mod_{\cO[X,X^{-1}]}^{\lfin}$.
\end{lemma}

\begin{proof} (a) This is essentially (\cite{emerton2}, Lemma 3.1.5). For completeness we recall argument. Put $h:=X\wcdot : M\to M$ and let $A$ be the image of $\cO[X]$ in $\End_{\cO}(M)$. Then $A$ is a finite commutative $\cO$-algebra so it is a product of finitely many local $\cO$-algebras $A = \prod_{i\in I} A_i$. We decompose $h= (h_i)_{i\in I}\in \prod_{i\in I} A_i$ accordingly. Let $J\subseteq I$ be the subset of indices $i\in I$ such that $h_i$ is a unit in $A_i$. We have a corresponding decomposition of $M$, namely $M= \prod_{i\in I} M_i$ and it is easy to see that $M^{X\sord}$ is the submodule $\prod_{i\in J} M_i$ of $M$. Thus we get $M = M^{X\sord} \oplus M'$ with $M'= \prod_{i\in I\setminus J} M_i$. The proof shows that $M'$ is uniquely determined. Since the restriction of $h|_{M^{X\sord}}: M^{X\sord} \to M^{X\sord}$ is an isomorphism the $\cO[X]$-action extends to an $ \cO[X,X^{-1}]$-action. Thus we have $M^{X\sord}\in \Mod_{\cO[X,X^{-1}]}^{\lfin}$.

(b) Let $\Psi\in \Hom_{\cO[X]}( \cO[X,X^{-1}], M)$. Since $\Psi(1)= X^n \Psi( X^{-n})$ for every $n\ge 1$ we have $\Psi(1) \in M^{X\sord}$ so the image of \eqref{eval} is contained in $M^{X\sord}$. To prove surjectivity of \eqref{eval} let $m\in M^{X\sord}$ and define $\Psi\in  \Hom_{\cO[X]}( \cO[X,X^{-1}], M)$ by $\Psi(x) = x\cdot m$ for $x\in \cO[X,X^{-1}]$ (note that $\Psi$ is well-defined since $M^{X\sord}$ is an $\cO[X,X^{-1}]$-module by (a)). We have $\ev(\Psi) = 1\cdot m = m$.  

If $\Psi\in  \Hom_{\cO[X]}( \cO[X,X^{-1}], M)$ lies in the kernel of \eqref{eval} then we have $\Psi(x) = x\cdot \Psi(1) =0$ for $x\in \cO[X]$. Since $X^n \Psi(X^{-n}) = \Psi(1) =0$, we get $\Psi(X^{-n}) = X^{-n}\Psi(1) =0$ for every $n\ge 1$, whence $\Psi(x)=0$ for all $x\in \cO[X,X^{-1}]$. Thus $\Psi=0$ and \eqref{eval} is bijective. The last assertion follows from the fact that $M^{X\sord}$ is -- as a submodule of $M$ -- finitely generated as an $\cO$-module. 
\end{proof}

\begin{lemma}
\label{lemma:lim1profinite}
Let 
\[
 \ldots \lra M_{n+1} \stackrel{\alpha_{n+1}}{\lra} M_n \stackrel{\alpha_n}{\lra} \ldots \lra M_1\stackrel{\alpha_1}{\lra} M_0
\]
be an inverse system of finitely generated $\cO$-modules. Then we have $\prolimn^{(1)} M_n =0$.
\end{lemma}

\begin{proof} Let $\fm$ be the maximal ideal of $\cO$ and put $M:= \prod_{n\ge 0} M_n$. Recall that $\prolimn^{(1)} M_n$ is the cokernel of the map 
\[
\Delta:  M \lra M,\,\, (m_n)_{n\ge 0} \mapsto (m_n- \alpha_{n+1}(m_{n+1}))_{n\ge 0}.
\]
By (\cite{matsumura}, Thm.\ 8.7) we have $M_n = \prolim_k (M_n \otimes \cO/\fm^k)$ for every $n\ge 0$ hence also $M = \prolim_k (M \otimes \cO/\fm^k)$. Note that the sequence 
\begin{equation}
\label{MLmodfm}
\begin{CD} 
0 @>>> \prolimn\, (M_n \otimes \cO/\fm^k) @>>> M \otimes \cO/\fm^k @> \Delta \otimes_{\cO} \id_{\cO/\fm^k} >> M \otimes \cO/\fm^k @>>> 0
\end{CD}
\end{equation}
is exact for every $k\ge 1$ since $(M_n\otimes_{\cO} \cO/\fm^k, \alpha_n\otimes_{\cO} \id_{\cO/\fm^k})_n$ is an inverse system of $\cO$-modules of finite lenght so it satisfies the Mittag-Leffler condition. Note also that the transition maps in
\[
\ldots \lra \prolimn\, (M_n \otimes \cO/\fm^{k+1})\lra \prolimn\, (M_n \otimes \cO/\fm^k) \lra \ldots \lra \prolimn\, (M_n \otimes \cO/\fm)
\]
are surjective so $\prolim^{(1)}$ of this inverse system vanishes. Therefore by 
passing in \eqref{MLmodfm} to the inverse limit over all $k\ge 1$ we obtain a short exact sequence
\[
0 \lra\prolimk \, \left(\prolimn\, (M_n \otimes \cO/\fm^k) \right) \lra M \stackrel{\Delta}{\lra} M \lra 0,
\]
hence $\prolimn^{(1)} M_n = \coker(\Delta) = 0$. 
\end{proof}

\begin{proof}[Proof of Prop.\ \ref{prop:finomodacyclic}] By Lemmas \ref{lemma:adlfinxinv} (b), \ref{lemma:ext0os} (b) and \ref{lemma:lim1profinite} we have 
\begin{equation*}
\label{extttplus}
\Ext_{\cO[X]}^s(\cO[X,X^{-1}], M) \, =\, \left\{\begin{array}{cc} M^{X\sord} & \mbox{if $s=0$,}\\
0 & \mbox{if $s>0$.}
\end{array}\right.
\end{equation*}
Thus the spectral sequence \eqref{adlfintinvss} degenerates and we have 
\begin{equation*}
\label{extttplus2}
R^n \Gamma_{\cO}^{X\sord}(M)\, =\, R^n \Gamma_{\lfin}(M^{X\sord}) \, =\, \left\{\begin{array}{cc} M^{X\sord} & \mbox{if $n=0$,}\\
0 & \mbox{if $n>0$.}
\end{array}\right.
\end{equation*}
for every $n\ge 0$ by Prop.\ \ref{prop:rnlfin} (b).
\end{proof}

\section{Representations of $\bG_m$ and $\PGL_2$ of a $p$-adic field}
\label{section:borelind}

Throughout this chapter $F$ denotes a $p$-adic field, i.e.\ a finite extension of $\bQ_p$. We let 
$v_F: F\to \bZ\cap\{+\infty\}$ denote the normalized valuation of $F$, $\cO_F$ its valuation ring and $\fp$ its valuation ideal. 
We also put $U^{(n)}=U_F^{(n)} =1 +\fp_F^n$ and $U_F = U_F^{(0)} = \cO_F^*$.

Let $G=\PGL_2(F) = \GL_2(F)/Z$ and let $\pr: \GL_2(F)\to G$ be the projection. Let $B$ be the standard Borel subgroup of $G$ and let $\barB$ be the opposite Borel. Let $B = T N$ be the Levi decomposition of $B$, so $T$ consists of the diagonal matrices (modulo the center $Z$ of $\GL_2(F)$) and $N$ consists of the upper triangular matrices with $1$ as diagonal entries (mod $Z$). In the following we are often going to identify $T$ with the onedimensional split torus $F^*$ via the isomorphism 
\begin{equation}
\label{splittorus}
\delta: F^*\lra T, \,\, x\mapsto \delta(x)=\left(\begin{matrix} x & \\ & 1\end{matrix}\right) \mod Z
\end{equation}
and $N$ with the additive group $F$ via the isomorphism
\begin{equation}
\label{unipotent}
n: F\lra N, \,\, y\mapsto n(y)=\left(\begin{matrix} 1 & y\\  0 & 1\end{matrix}\right) \mod Z.
\end{equation}
Any $b\in B$ can be written uniquely in the form $b=\delta(x)\cdot n(y)$ with $x\in F^*$, $y\in F$. We also remark that the map
\begin{equation}
\label{semidirect}
F^* \ltimes F\lra B, \,\, (x,y)\mapsto \left(\begin{matrix} x & y\\  0 & 1\end{matrix}\right) \mod Z.
\end{equation}
is an isomorphism. Here we let the multiplicative group $F^*$ act on the additive group $F$ by multiplication.

We let $N_0$ be the subgroup of $N$ hat corresponds to the $\cO_F$ under the isomorphism \eqref{unipotent}, i.e.\ $N_0=\{ n(y)\mid y\in \cO_F\}$ and put $T^+= \{t\in T\mid N_0^t= t N_0 t^{-1} \subseteq N_0\}$, $T^-= \{t\in T\mid N_0^t \supseteq N_0\}$ and $T^0= T^+\cap T^-$. Thus $T^+$ (resp.\ $T^0$) corresponds to the monoid $\cO_F-\{0\}$ (resp.\ the group $U_F$) under the isomorphism \eqref{splittorus}. 

For a positive integer $n$ we let $K(\fp^n) = \ker(\GL_2(\cO_F) \to \GL_2(\cO_F/\fp^n))$
and put 
\begin{eqnarray*}
\label{congruence}
&& K_1(\fp^n)\, =\, \{ A\in  \GL_2(\cO_F)\mid A \equiv \begin{pmatrix} 1 & \star \\ 0 & 1 \end{pmatrix} \mod \fp^n \, \} \qquad \mbox{resp.}\\ 
&& K_0(\fp^n)\, =\, \{ A\in  \GL_2(\cO_F)\mid A \equiv \begin{pmatrix} \star & \star \\ 0 & \star \end{pmatrix} \mod \fp^n \, \}.
\end{eqnarray*} 
We define the compact open subgroups $K(n)$ resp.\ $K_1(n)$ resp.\ $K_0(n)$ of $G$ as the image of the group $K(\fp^n)$ resp.\ $K_1(\fp^n)$ resp.\ $K_0(\fp^n)$ under the projection  $\pr:\GL_2(F)\to G$, i.e.\ we have
\begin{equation*}
\label{congruence2}
K(n) \, =\, K(\fp^n)Z(F)/Z(F) \quad \mbox{and} \quad K_{?}(n) \, =\, K_{?}(\fp^n)Z(F)/Z(F)
\end{equation*} 
for $?\in \{0, 1\}$. We also put $K(0) = \PGL_2(\cO_F)$. 

More generally for a closed subgroup $H$ of $T^0$ we put 
\begin{equation*}
\label{congruence3}
K_H(n) \, =\, K(n) H N_0.
\end{equation*}
For example if $H = 1$ (resp.\ $H= \delta(U_F)$) then $K_H(n) = K_1(n)$ (resp.\ $K_H(n) = K_0(n)$). 

\subsection{Finiteness properties for $\Ext$-groups of admissible $R[T]$-modules.}
\label{subsection:ladmext}

Let $R$ be a noetherian ring. In this section we establish some basic properties of the category of locally admissible $R[T]$-modules and then -- under certain assumptions on $R$ -- prove some finiteness properties for the $\Ext$-groups $\Ext_{R, T}^{\bu}(W_1, W_2)$ of two admissible $R[T]$-modules $W_1$ and $W_2$. In fact we work here -- and throughout the paper -- in a slightly more general framework, namely we fix a compact subgroup $H$ of $T^0$ and we put 
$\barT = T/H$, $\barT^0 = T^0/H$, $\barT^+=T^+/H$, $\barT^-=T^-/H$. We denote by $\bar{\delta}$ the composite map 
\begin{equation}
\label{splittorus2}
\bar{\delta}: F^*\stackrel{\eqref{splittorus}}{\lra} T\stackrel{\pr}{\lra} \barT.
\end{equation}
Note that $\barT$ satisfies the condition \ref{assum:normalsubgroup} so that $\Mod_R^{\ladm}(\barT)$ is a Serre subcategory of $\Mod_R^{\sm}(\barT)$ with enough injectives. 

Let $W$ be a $R[\barT]$-module and let $w\in W$. We fix a prime element $\pi\in \cO_F$ and put 
$t_0=\bar{\delta}(\pi)\in \barT^+$. Recall (see \cite{emerton2}, Def.\ 2.3.1) that $w\in W$ is called $t_0^{\pm 1}$-finite if $R[t_0^{\pm 1}]w=R[t_0, t_0^{-1}]w$ is a finitely generated $R$-module. The subset $W_{t_0^{\pm 1}-\fin}$ of all $t_0^{\pm 1}$-finite elements of $W$ is an $R[\barT]$-submodule of $W$. If $W= W_{t_0^{\pm 1}-\fin}$ then $W$ is called locally $t_0^{\pm 1}$-finite. 

\begin{lemma}
\label{lemma:loctfin}
Let $W$ be a discrete $R[\barT]$-module. Then we have 
\begin{equation}
\label{tfin2}
W_{\ladm}\, = \, W_{t_0^{\pm 1}-\fin}.
\end{equation} 
In particular $W$ is locally admissible if and only if it is locally $t_0^{\pm 1}$-finite. 
\end{lemma}

\begin{proof} The inclusion $W_{\ladm}\subseteq W_{t_0^{\pm 1}-\fin}$ is obvious. Conversely if $w\in W_{t_0^{\pm 1}-\fin}$ and if  $U=\Stab_{\barT^0}(w)$ then $R[\barT] w= R[\barT^0/U](R[t_0^{\pm 1}] w)$
is a finitely generated $R$-module since $R[\barT^0/U]$ is a finite $R$-algebra.
\end{proof}

\begin{prop}
\label{prop:rnladm}
(a) The embedding $\Mod_R^{\ladm}(\barT)\hra \Mod_R^{\sm}(\barT)$ maps injectives to injectives. 
\medskip

\noi (b) Let $W$ be a locally admissible $R[\barT]$-module. Then $W$ is $\Gamma_{\ladm}$-acyclic, i.e.\ we have $(R^n \Gamma_{\ladm})(W)\break =0$ for every $n\ge 1$.
\end{prop}

We need the following 

\begin{lemma}
\label{lemma:injcriterion}
Let $I\in \Mod_R^{\sm}(\barT)$ (resp.\ $I\in \Mod_R^{\ladm}(\barT)$). Then $I$ is an injective object of $\Mod_R^{\sm}(\barT)$ (resp.\ of $\Mod_R^{\ladm}(\barT)$) if and only if $I^U$ is an injective object of $\Mod_R^{\sm}(\barT/U) = \Mod_{R[\barT/U]}$ (resp.\ of $\Mod_R^{\ladm}(\barT/U)$) for every open subgroup $U$ of $\barT^0$.
\end{lemma}

\begin{proof} We give the proof only for the category $\Mod_R^{\sm}(\barT)$ (the proof in the other case is similar). For an open subgroup $U$ of $\barT^0$ the functor $\Mod_R^{\sm}(\barT) \to \Mod_{R[\barT/U]}, W\mapsto W^U$ maps injectives to injectives since it has an exact left adjoint, namely the obvious embedding $\Mod_{R[\barT/U]}\hra \Mod_R^{\sm}(\barT)$.
Conversely, assume that $I^U$ is an injective $R[\barT/U]$-module for every open subgroup $U$ of $\barT^0$. Let $W\in \Mod_R^{\sm}(\barT)$, let 
$W'$ be an $R[\barT]$-submodule of $W$ and let $\varphi: W'\to I$ be a homomorphism of $R[\barT]$-modules. We have to show that $\varphi$ extends to homomorphism $\psi: W\to I$. Using a standard argument involving Zorn's Lemma it suffices to consider the case where $W = W' + R[\barT] w_0$ for some $w_0\in W$. Let $U$ be a compact open subgroup of $\barT$ with $w_0\in W^U$. Since $W_0'= W'\cap R[\barT/U]w_0$ is an $R[\barT/U]$-submodule of $R[\barT]w_0 = R[\barT/U]w_0$ we see that the image of $\varphi_0:= \varphi|_{W_0'}: W_0' \to I$ is contained in $I^U$, i.e.\ we can view $\varphi_0$ as a homomorphism of $R[\barT/U]$-modules $W_0'\to I^U$. Because of the assumption $\varphi_0$ extends to a homomorphism $\psi_0: R[\barT/U]w_0 \to I^U$. If we define $\psi: W = W' + R[\barT] w_0\to I$ by $\psi(w) = \varphi(w_1) + \psi_0(w_2)$ if $w= w_1 + w_2$ with $w_1\in W'$ and $w_2\in R[\barT/U]w_0$ then $\psi$ is well-defined and extends $\varphi$. 
\end{proof} 

\begin{proof}[Proof of Prop.\ \ref{prop:rnladm}] The second statement follows immediately from the first. For (a) let $I$ be an injective object in $\Mod_R^{\ladm}(\barT)$ and let $U$ be an open subgroup of $\barT^0$. By Lemma \ref{lemma:injcriterion} we have $I^U\in \Mod_R^{\ladm}(\barT/U)$ and we have to show that $I^U$ is an injective $R[\barT/U]$-module. However this follows from Prop.\ \ref{prop:rnlfin} (a) since the ring $R'= R[\barT^0/U]$ is noetherian and we have $R[\barT/U]= R'[t_0, t_0^{-1}]$.
\end{proof}

\begin{remark}
\label{remark:dliminj} \rm The Lemma implies that a direct limit of injective objects in $\Mod_R^{\sm}(\barT)$ is again injective. This follows from the fact that the ring $R[\barT/U]$ is noetherian if $U\subseteq \barT^0$ is open so any direct limit of injective $R[\barT/U]$-modules is again injective. \enddemo
\end{remark} 

We now turn to the Ext groups. Recall that for 
$W_1, W_2\in \Mod_R^{\sm}(\barT)$ and $n\ge 0$ we have defined $\Ext_{R, \barT}^n(W_1, W_2) : = \Ext_{\Mod_R^{\sm}(\barT)}^n(W_1, W_2)$.
By Prop.\ \ref{prop:rnladm} we have 
\begin{equation*}
\label{extladm2}
\Ext_{R, \barT}^{\bu}(W_1, W_2) \,=\, \Ext_{\Mod_R^{\ladm}(\barT)}^{\bu}(W_1, W_2)
\end{equation*} 
if $W_1, W_2\in  \Mod_R^{\ladm}(\barT)$. Now assume that $W_1$ is finitely generated as an $R$-module and that $W_2$ is admissible. Then there exists an open subgroup $U$ of $\barT^0$ so that $W_1 = W_1^U$, hence 
\[
\Ext_{R, \barT}^0(W_1, W_2)\, = \,\Hom_{R[\barT/U]}(W_1, W_2^U)
\]
is a finitely generated $R$-module. Thus one might expect $\Ext_{R, \barT}^n(W_1, W_2)\in \Mod_{R,f}$ for every $n\ge 0$. We prove this in two cases namely (i) when $\barT$ is discrete (i.e.\ the subgroup $H$ of $T^0$ is open) and (ii) when $R$ is an Artinian local ring with finite residue field of characteristic $p$. 

\begin{lemma}
\label{lemma:extadmfin1}
If  $\,\barT^0$ is discrete then $\Ext_{R, \barT}^n(W_1, W_2)$ is a finitely generated $R$-module for every $n\ge 0$.
\end{lemma}

\begin{proof} We fix $n\ge 0$. The assumption implies that $R[\barT]$ is noetherian, that $W_2$ is also finitely generated as an $R$-module and that $\Ext_{R, \barT}^n(W_1, W_2) = \Ext_{R[\barT]}^n(W_1, W_2)$. Thus if $\fa= \Ann_{R[\barT]}(W_1)$ then $\Ext_{R, \barT}^n(W_1, W_2)$ if a finitely generated $R[\barT]/\fa$-module. Since $R[\barT]/\fa$ is a finite $R$-algebra the assertion follows. 
\end{proof}

In the second case the following holds

\begin{prop}
\label{prop:extadmfin2}
Let $R=A$ be an Artinian local ring with maximal ideal $\fm$ and finite residue field $k$ of characteristic $p$. Let $W_1 \in \Mod_{A, f}^{\sm}(\barT)$ and $W_2 \in \Mod_A^{\adm}(\barT)$. Then 
$\Ext_{A, \barT}^n(W_1, W_2)$ is a finitely generated $A$-module for every $n\ge 0$.
\end{prop}

\begin{proof} The assumptions imply that $W_1$ has finite length. We may thus assume that $W_1$ is a simple object of $\Mod_{A, f}^{\sm}(\barT)$. This implies then that the maximal pro-$p$-subgroup of $\barT^0$ acts trivially on $W_1$ (because $W_1$ is a $p$-group) and also $\fm W_1=0$. 
Thus if $A'= A[\Delta]$ where $\Delta$ is the prime-to-$p$-part of $\barT^0$, we can view $W_1$ as a simple object in the category of $\Mod_{A'[t_0^{\pm}]}$. These can be easily classified: there exists a finite extension $\ell/k$, a character $\chi_0: \Delta\to \ell^*$ and an element $\alpha\in \ell$, $\alpha\ne 0$ so that 
$W_1 = \ell$ with $\Delta$-action given by $\chi_0$ and so that $t_0$ acts by multiplication with $\alpha$. It follows that there exists a resolution in $\Mod_A^{\sm}(\barT)$
\begin{equation*}
\label{resolutioncind}
\begin{CD}
0 @>>> \cInd_{\barT^0}^{\barT} \ell(\chi_0) @> t_0 - a >>  \cInd_{\barT^0}^{\barT} \ell(\chi_0) @>>> W_1 @>>> 0.
\end{CD}
\end{equation*} 
Thus by Lemma \ref{lemma:extadmfin1} it suffices to show that the $A$-module
\[
\Ext_{\Mod_A^{\sm}(\barT)}^n(\cInd_{\barT^0}^{\barT} \ell(\chi_0), W_2) \, =\, \Ext_{\Mod_A^{\sm}(\barT^0)}^n(\ell(\chi_0), W_2) 
\]
is finitely generated for every $n\ge 0$. This follows immediately from (\cite{emerton3}, Prop.\ 2.1.9). 
\end{proof}

We are now going to study the groups $\Ext_{R, \barT}^n(W_1, W_2)$ when both $R[\barT]$-modules $W_1$ and $W_2$ are admissible. For that we recall the notion of an {\it augmented $R[T]$-module} (see \cite{emerton2}, \S 2). 
For an open subgroup $U$ of $\barT$ we consider the $R$-algebra
\begin{equation}
\label{iwasawa}
\Lambda_{R}(U) \, =\, \prolim_{U'} R[U/U']
\end{equation} 
where $U'$ runs over all compact open subgroups of $U$. If $U$ itself is compact then $\La_R(U)= R[\![U]\!]$ is the usual completed $R$-group algebra of $U$. Note that we have 
\begin{equation*}
\label{iwasawa2}
\Lambda_{R}(\barT) \, =\,R[\![\barT^0]\!][t_0^{\pm 1}].
\end{equation*} 
In particular the ring $\Lambda_{R}(\barT)$ is noetherian (more generally $\Lambda_{R}(U)$ is noetherian for any open subgroup $U$ of $\barT$). 

Following (\cite{emerton2}, \S 2) a $\Lambda_{R}(\barT)$-module will be called an {\it augmented $R[\barT]$-module}. 
The category of augmented $R[\barT]$-modules will be denote by $\Mod_R^{\aug}(\barT)$. Let $W$ be a discrete $R[\barT]$-module. The fact that every element $w\in W$ has an open stabilizer in $\barT$ implies that the $R[\barT]$-action on $W$ extends naturally to a $\Lambda_{R}(\barT)$-action. Thus we have a natural fully faithful embedding $\Mod_R^{\sm}(\barT)\hra \Mod_R^{\aug}(\barT)$. An augmented $R[\barT]$-module $L$ is called finitely generated if $L$ is finitely generated as $\Lambda_{R}(U)$-modules for some (equivalently, any) compact open subgroup $U$ of $\barT$. The full subcategory of $\Mod_R^{\aug}(\barT)$ of finitely generated augmented $R[\barT]$-modules will be denoted by $\Mod_R^{\fgaug}(\barT)$. If $R=\cO$ is a complete, noetherian local ring with finite residue field of characteristic $p$ then there exists a natural (profinite) topology on every objects $L$ of $\Mod_{\cO}^{\fgaug}(\barT)$ (see \cite{emerton2}, Prop.\ 2.1.3) such that the action $\Lambda_{\cO}(T) \times L \to L$ is continuous. This topology is called the {\it canonical topology}.

Assume that $R=A$ is an Artinian local ring as in Prop.\ \ref{prop:extadmfin2} and let $W_1, W_2 \in \Mod_A^{\adm}(\barT)$. By (\cite{emerton2}, Lemma 2.2.7), Pontrjagin duality induces an anti-equivalence of categories
\begin{equation*}
\label{pontrjagin1}
\Mod_A^{\adm}(T)\lra \Mod_A^{\fgaug}(T), \,\, W\mapsto D(W)=\Hom_A(W, \bQ_p/\bZ_p)
\end{equation*} 
hence 
\[
\Ext_{A, \barT}^0(W_1, W_2)\, = \,\Hom_{A[\barT]}(W_1, W_2)\, \cong\, \Hom_{\Lambda_{A}(\barT)}(D(W_2), D(W_1))
\]
is a finitely generated augmented $A[\barT]$-module. More generally we have

\begin{prop}
\label{prop:extadmfin3}
Let $R=A$ be an Artinian local ring with finite residue field $k$ of characteristic $p$ and let $W_1, W_2 \in \Mod_A^{\adm}(\barT)$. Then $\Ext_{A, \barT}^n(W_1, W_2)$ is a finitely generated augmented $A[\barT]$-module for every $n\ge 0$. 
\end{prop}

\begin{proof} Let $U$ be an open torsionfree subgroup of $\barT^0$ and put $\La=\La_A(U)$.
Note that $\La$ is a complete noetherian local $A$-algebra with residue field $k$. 
For $m\ge 1$ we put $W_1^{(m)}= W_1^{\bar{\delta}(U_F^{(m)})}$. We have $W_1= \bigcup_{m\ge 1} W_1^{(m)}$ and $W_1^{(m)} \in \Mod_{A, f}(\barT)$ for every $m\ge1$. There exists a short exact sequence of augmented $A[\barT]$-modules
\begin{equation*}
\label{extlim}
0\lra \prolimm^{(1)} \Ext_{A, \barT}^{n-1}(W_1^{(m)}, W_2) \lra \Ext_{A, \barT}^n(W_1, W_2)\lra \prolimm \Ext_{A, \barT}^n(W_1^{(m)}, W_2)\lra 0
\end{equation*}
(this can be deduced from (\cite{weibel}, Thm.\ 3.5.8) in the same way as Application 3.5.10 in loc.\ cit.).
Since, by Prop.\ \ref{prop:extadmfin2} (a), the groups $\Ext_{A, \barT}^n(W_1^{(m)}, W_2)$ are finite for all $m, n\ge 0$ the $\prolimm^{(1)}$-term vanishes so we get 
\begin{equation*}
\label{extlim2}
\Ext_{A, \barT}^n(W_1, W_2)\,\cong\, \prolimm \Ext_{A, \barT}^n(W_1^{(m)}, W_2)
\end{equation*}
for every $n\ge 0$. In particular we see that $\Ext_{A, \barT}^n(W_1, W_2)$ carries a profinite topology thus providing it with the structure of a compact $\La$-module. 

Let $\cO$ denote the image of $\La$ in $\End_{A[\barT]}(W_1)$. Note that $\cO$ is -- as a quotient of $\La$ -- a complete noetherian local $A$-algebra with residue field $k$. To prove that $\Ext_{A, \barT}^n(W_1, W_2)$ is a finitely generated $\La$-module for every $n\ge 0$ we use induction over $\dim(\cO)$. If $\dim(\cO) = 0$ then $\cO$ is an Artinian local $A$-algebra. Hence the $\La$-action on $W_1$ factors over $A[U/U']$ for some open subgroup $U'$ of $U$. Thus $U'$ acts trivially on $W_1$ and we get $W_1\in \Mod_{A, f}(\barT)$ (since $W_1$ is admissible) and $\Ext_{A, \barT}^n(W_1, W_2)\in \Mod_{A, f}$ by Prop.\ \ref{prop:extadmfin2} (a). 

If $\dim(\cO)>0$ we choose $\la$ in the maximal ideal of $\cO$ such that $\dim(\cO/\la \cO) =\dim(\cO) -1$. Put $W_1[\la]=\ker(W_1\stackrel{\la\cdot }{\lra} W_1)$ and $W_1'=\image(W_1\stackrel{\la\cdot }{\lra} W_1)$. Since the images of $\La$ in $\End_{A[\barT]}(W_1[\la])$ and $\End_{A[\barT]}(W_1/W_1')$ are quotients of $\cO/\la \cO$, the assertion holds for $W_1[\la]$ and $W_1/W_1'$ by the induction hypothesis and the fact that the quotient of an admissible $A[\barT]$-module is again admissible (\cite{emerton2}, Prop.\ 2.2.13; here we use the fact that $A$ is Artinian). By considering the long exact $\Ext$-sequences associated to the two short exact sequences 
\[
0\lra W_1[\la] \lra W_1 \lra W_1' \lra 0, \qquad 0\lra W_1' \lra W_1 \lra W_1/W_1' \lra 0
\]
we conclude that kernel and cokernel of the map $\Ext_{A, \barT}^n(W_1, W_2)\to \Ext_{A, \barT}^n(W_1, W_2), x\mapsto \la \cdot x$ are finitely generated $\La$-modules. Since $\Ext_{A, \barT}^n(W_1, W_2)\otimes_{\La} k$ is a quotient of the cokernel we conclude 
\begin{equation*}
\label{extlim3}
\dim_k \Ext_{A, \barT}^n(W_1, W_2)\otimes_{\La} k<\infty
\end{equation*}
for every $n\ge 0$. By the topological Nakayama Lemma (see e.g.\ \cite{ballister-howson}) this implies that $\Ext_{A, \barT}^n(W_1, \break W_2)$ is a finitely generated $\La$-module.
\end{proof}

\subsection{The functor $\Gamma^{\ord}$} 
\label{subsection:gammaord}

In this section we assume that $R=\cO$ is a complete noetherian local ring with maximal ideal $\fm$ and residue field $k$ (the ring $\cO$ has no connection to the valuation ring $\cO_F$ in $F$). We fix again closed subgroup $H$ of $T^0$ and define $\barT^+$ ($\barT^0$ etc.) as before as the  image of the monoid $\cO_F\setminus\{0\}$ (the group $U_F$ etc.) under the epimorphism \eqref{splittorus2}. Note that the monoid $\barT^+$ satisfies condition \ref{assum:normalsubgroup} and that the $\cO$-algebra $\cO[\barT^+]$ is isomorphic to a polynomial ring in one variable over the group ring $\cO[\barT^0]$. An $\cO[\barT^+]$-module $V$ is discrete (resp.\ admissible) if and only if and only if $V$ is discrete (resp.\ admissible) as an $\cO[\barT^0]$-module. 

\begin{df}
\label{df:smallord}
Let $V$ be a discrete $\cO[\barT^+]$-module. We define the $\cO[\barT^+]$-submodule $V^{\ord}$ of $V$ by 
\begin{equation*}
\label{ord1}
V^{\ord} \, =\, \bigcap_{m\ge 0} t_0^m V.
\end{equation*}
\end{df}

We consider the map 
\begin{equation}
\label{ord2}
\ev_1: \Hom_{\cO[\barT^+]}( \cO[\barT], V)\lra V^{\ord}, \,\, \Psi\mapsto \Psi(1).
\end{equation}
The image is contained in $V^{\ord}$ since $\Psi(1) = t_0^m \Psi(t_0^{-m})\in t_0^m V$
for all $m\ge 0$. 

\begin{lemma}
\label{lemma:smallord}
Let $V\in \Mod_{\cO, f}^{\sm}[{\barT}^+]$. 
\medskip 

\noi (a) The $\cO[\barT^+]$-submodule $V^{\ord}$ is a direct summand of $V$, i.e.\ there exists a (unique) $\cO[\barT^+]$-submodule $V'$ of $V$ such that  
\begin{equation*}
\label{ord3b}
V\, =\, V^{\ord} \oplus V'.
\end{equation*}
Moreover the $\barT^+$-action on $V^{\ord}$ extends to a $\barT$-action.
\medskip 

\noi (b) The map \eqref{ord2} is an isomorphism of $\cO[\barT]$-modules. Thus we have 
\[
\Hom_{\cO[\barT^+]}( \cO[\barT], V)\cong V^{\ord} \in \Mod_{\cO, f}^{\sm}(\barT).
\] 
\end{lemma}

\begin{proof} (a) can be proved as (\cite{emerton2}, Lemma 3.1.5) or it can be deduced from Lemma \ref{lemma:ext0os} (a). 

(b) To prove surjectivity of \eqref{ord2} let $v\in V^{\ord}$ and define $\Psi\in \Hom_{\cO[\barT^+]}( \cO[\barT], V)$ by $\Psi(x) = x\cdot v$ for $x\in \cO[\barT]$ (note that $\Psi$ is well-defined because $V^{\ord}$ is a $\cO[\barT]$-module by (a)). We have $\ev_1(\Psi) = 1\cdot v = v$.  

If $\Psi\in \Hom_{\cO[\barT^+]}( \cO[\barT], V)$ lies in the kernel of \eqref{ord2} then we have $\Psi(x) = x\cdot \Psi(1) =0$ for $x\in \cO[\barT^+]$. For $t\in \barT^{-}$ we have $t^{-1} \Psi(t) = \Psi(1) =0$ hence $\Psi(t)=0$ for all $t\in \barT^{-}$ as well. Thus $\Psi=0$ and \eqref{ord2} is  injective. The last assertion follows from the fact that $\Hom_{\cO[\barT^+]}( \cO[\barT], V)\cong V^{\ord}$ is finitely generated as an $\cO$-module. 
\end{proof}

\begin{df}
\label{df:gammaord}
For $V\in \Mod_{\cO}^{\sm}(\barT^+)$ we set
\[
\Gamma_{\cO}^{\ord}(V)\, :=\, \Hom_{\cO[\barT^+]}( \cO[\barT], V)_{t_0^{\pm 1}-\fin}.
\]
\end{df}

When we are dealing with a fixed coefficient ring $\cO$ we will often drop it from the notation, i.e.\ we write $\Gamma^{\ord}(V)$ instead of $\Gamma_{\cO}^{\ord}(V)$. 

\begin{prop}
\label{prop:adTTplus}
(a) For $V\in \Mod_{\cO, f}^{\sm}(\barT^+)$ we have 
\[
\Gamma^{\ord}(V) \, =\, V^{\ord}.
\]
In particular we get $\Gamma^{\ord}(V) \in \Mod_{\cO, f}^{\sm}(\barT)$. 
\medskip

\noi (b) For any $V\in \Mod_{\cO}^{\sm}(\barT^+)$ the $\cO[\barT]$-module $\Gamma^{\ord}(V)$ is locally admissible. 
\medskip

\noi (c) The functor 
\begin{equation}
\label{adTTplus2}
\Gamma^{\ord}: \Mod_{\cO}^{\sm}(\barT^+)\lra \Mod_{\cO}^{\ladm}(\barT), \,\, V \mapsto \Gamma^{\ord}(V)
\end{equation} 
is right adjoint to the forgetful functor $\Mod_{\cO}^{\ladm}(\barT)\lra \Mod_{\cO}^{\sm}(\barT^+)$.
\medskip

\noi (d) The functor \eqref{adTTplus2} commutes with direct limits. 
\end{prop}

\begin{proof} (a) By Lemma \ref{lemma:smallord} we have 
\[
\Gamma^{\ord}(V)\, =\, \Hom_{\cO[\barT^+]}( \cO[\barT], V)_{t_0^{\pm 1}-\fin}\,\cong \,(V^{\ord})_{t_0^{\pm 1}-\fin}\, =\, V^{\ord}.
\]
The last equality follows from the fact that $V^{\ord}$ is -- as direct summand of $V$ -- finitely generated as an $\cO$-module. 

(b) Let $\Psi \in \Gamma^{\ord}(V)$ and let $n\ge 1$ such that $\{t_0^i\cdot \Psi\mid -n\le i\le n\}$ generates $\cO[t_0^{\pm 1}] \Psi$. Let $U$ be an open subgroup of $\barT^0$ such that $\Psi(t_0^i)\in V^U$ for all $i\in \bZ$ with $-n\le i\le n$ and put 
\[
V_0\, : =\,  \sum_{i=-n}^n \cO[\barT^0/U] \Psi(t_0^i) \, =\, \Psi\left(\sum_{i=-n}^n \cO[\barT^0] \, t_0^i\right).
\]
Note that $V_0\in \Mod_{\cO,f}$ and $V_0\subseteq \image(\Psi)$. We claim $\image(\Psi) = V_0$. For that it suffices to see that $\Psi(t)\in V_0$ for every $t\in T$. Fix $t\in T$ and let $m\in \bZ$, $u\in T^0$ with $t= t_0^m u$. There exists $a_i\in \cO$ for $i\in \bZ$, $-n\le i \le n$ such that $t_0^m\cdot \Psi \, =\, \sum_{i=-n}^n a_i (t_0^i\cdot \Psi)$. It follows 
\[
\Psi(t) \, = \, u \Psi(t_0^m) \, =\, u (t_0^m\cdot \Psi)(1)\, =\, \sum_{i=-n}^n a_i u (t_0^i\cdot \Psi)(1)\in V_0.
\]
The equality $\image(\Psi) = V_0$ implies that $V_0$ is a $\cO[\barT^+]$-submodule of $V$ and that 
$\Psi$ lies in the image of the canonical map $V_0^{\ord}=\Gamma^{\ord}(V_0) \hra \Gamma^{\ord}(V)$. 
It follows that $\cO[\barT] \Psi\subseteq \image(\Gamma^{\ord}(V_0) \hra \Gamma^{\ord}(V))$, hence it is a finitely generated as $\cO$-module by (a). Therefore $\cO[\barT] \Psi$ is admissible and $\Gamma^{\ord}(V)$ is locally admissible. 

For (c) let $W\in \Mod_{\cO}^{\ladm}(\barT)$ and let $V\in \Mod_{\cO}^{\sm}(\barT^+)$. Since $W$ is $t_0^{\pm 1}$-finite we have
\begin{eqnarray*}
\Hom_{\cO[\barT^+]}(W, V)& \cong &\Hom_{\cO[\barT]}(W,  \Hom_{\cO[\barT^+]}( \cO[\barT], V)) \\
& \cong & \Hom_{\cO[\barT]}(W,  \Hom_{\cO[\barT^+]}( \cO[\barT], V)_{t_0^{\pm 1}-\fin})
\end{eqnarray*}
so the assertion follows from (a).

For (d) see (\cite{emerton2}, Lemma 3.2.2 (2)). Note that the proof given there remains valid if the coefficient ring $\cO$ is a noetherian local ring.
\end{proof}

We finish this section with the following result which is helpful in studying the derived functors of the functor $\Ord$ in section \ref{section:ordinary}.

\begin{prop}
\label{prop:rgammaord}
(a) The restriction of the functor $\Gamma^{\ord}$ to the subcategory $\Mod_{\cO}^{\ladm}(\barT^+)$, i.e.\ the functor 
\begin{equation*}
\label{ladTTplus2}
\Mod_{\cO}^{\ladm}(\barT^+)\lra \Mod_{\cO}^{\ladm}(\barT), \,\, V \mapsto \Hom_{\cO[\barT^+]}( \cO[\barT], V)_{\ladm}
\end{equation*} 
is exact.
\medskip

\noi (b) Any locally admissible $\cO[\barT^+]$-module is $\Gamma^{\ord}$-acyclic. 
\end{prop}

\begin{proof} (a) We first consider the restriction of $\Gamma^{\ord}$ to the subcategory $\Mod_{\cO, f}(\barT^+/U)$ of $\Mod_{\cO}^{\ladm}(\barT^+)$ where $U$ is an open subgroup of $\barT^0$. 
The $\cO$-algebra $\wcO := \cO[\barT^0/U]$ is finite, so it is a product of a finite number of local $\cO$-algebras $\wcO = \prod_{i\in I} \wcO_i$. We decompose every object $V=\prod_{i\in I} V_i$ in $\Mod_{\cO, f}(\barT^+/U)$ accordingly. We then have 
\begin{equation*}
\label{extttplus2a}
\Gamma^{\ord}(V) \, =\, \prod_{i\in I} \Gamma_{\cO_i}^{t_0\sord}(V_i).
\end{equation*}
Together with Prop.\ \ref{prop:finomodacyclic} it follows that the restriction of $\Gamma^{\ord}$ to $\Mod_{\cO, f}^{\sm}(\barT^+/U)$ is exact. 

Let $0 \lra V'  \stackrel{\alpha}{\lra}V \stackrel{\beta}{\lra} V'' \lra 0$ be a short exact sequence in $\Mod_{\cO}^{\ladm}(\barT^+)$. 
We can write $V$ as a direct limit $V=\dlim_{j\in J} V_j$ where each $V_j$ is a $\cO[\barT^+]$-submodule of $V$ that is finitely generated as an $\cO$-module. If for $j\in J$ we put $V_j' :=\alpha^{-1}(V_j)$ and $V_j'':=\beta(V_j)$ then the sequence $0 \lra V'  \stackrel{\alpha|_{V_j'}}{\lra}V \stackrel{\beta|_{V_j}}{\lra} V'' \lra 0$ is exact and lies in $\Mod_{\cO, f}(\barT^+/U)$ for some open subgroup $U$ of $\barT^0$. It follows (as we have already seen) that 
\[
0 \lra \Gamma^{\ord}(V_j') \lra \Gamma^{\ord}(V_j) \lra \Gamma^{\ord}(V_j'') \lra 0
\]
is exact. Passing to the direct limit over all $j\in J$ using Prop.\ \ref{prop:adTTplus} (d) yields the exactness of $0 \to \Gamma^{\ord}(V') \to \Gamma^{\ord}(V) \to \Gamma^{\ord}(V) \to 0$. 

(b) follows from (a) and Prop.\ \ref{prop:rnladm} (a).
\end{proof}

\subsection{Parabolic induction and the functor $\Ord$}
\label{subsection:parind}

Let $R$ be a noetherian ring. As in the last section $H$ denotes a compact subgroup of $T^0\cong U_F$ and we put $\barT = T/H$, $\barT^+ = T^+/H$. 
We consider the parabolic induction functor 
\begin{equation}
\label{parind1}
\Ind_{B}^{G}: \Mod_R^{\sm}(\barT)\lra \Mod_R^{\sm}(G), \,\, W\mapsto \Ind_{B}^{G} W.
\end{equation} 
It is the composite of the natural embedding $\Mod_R^{\sm}(\barT)\hra \Mod_R^{\sm}(T)$, the inflation functor $\Infl_{T}^{B}: \Mod_R^{\sm}(T)\lra \Mod_R^{\sm}(B)$ -- sending a discrete $R[\barT]$-module $W$ to the $R[B]$-module $\Infl_{T}^{B}(W) = W$ with $B$-action structure induced by the projection $B \to T, b = t\cdot n \mapsto t$ -- and the (smooth) induction functor $\Ind_{B}^{G}: \Mod_R^{\sm}(B)\lra \Mod_R^{\sm}(G)$. By abuse of notation the functor will also be denote by $\Ind_{B}^{G}$. Recall that the elements of $\Ind_{B}^{G} W$ are functions $\Phi: G\to W$ such that there exists $m\ge 0$ so that $\Phi(tngk) = t \Phi(g)$ for $t\in T$, $n\in N$, $g\in G$ and $k\in K(m)$. The $G$-action on $\Ind_{B}^{G} W$ is induced by right multiplication, i.e.\ we have 
\begin{equation*}
\label{parind1a}
(h\cdot \Phi)(g) \, =\, \Phi(gh)
\end{equation*} 
for $g, h\in G$ and $\Phi\in \Ind_{B}^{G} W$. The functor \eqref{parind1} is exact, commutes with inductive limits and with base change of the ring $R$. More precisely we have 

\begin{lemma}
\label{lemma:parindbasechange}
(a) Let $\varphi: R\to R'$ be a ring homomorphism and let $W\in \Mod_R^{\sm}(\barT)$. The canonical map 
\begin{equation*}
\label{parind2}
\left(\Ind_{B}^{G} W\right)_{R'} \, \lra \, \Ind_{B}^{G} W_{R'} 
\end{equation*}
is an isomorphism of discrete $R'[G]$-modules.
\medskip

\noi (b)  Let $\{W_i\}_{i\in I}$ be an inductive system of discrete $R[\barT]$-modules and let $W= \dlim_{i\in I} W_i$. The canonical map of $R[G]$-modules
\begin{equation*}
\label{parind3}
\dlim_{i\in I} \Ind_{B}^{G} W_i \, \cong \, \Ind_{B}^{G} W 
\end{equation*}
is an isomorphism.
\end{lemma}

\begin{proof} (a) If $s: B\backslash G \to G$ is a continuous section of $G\to B\backslash G, g\mapsto Bg$ then the map 
\begin{equation}
\label{parind4}
\Ind_{B}^{G} W\, \lra\, C(G/B, W),\,\, \Phi\mapsto \Phi\circ s
\end{equation}
is an isomorphism of $R$-modules (see \cite{vigneras}, \S 4). Here 
$C(G/B, W)$ is the $R$-module of locally constant maps $G/B\to W$. Thus it suffices to show that the canonical map
\[
C(G/B, W)_{R'} \,\lra\, C(G/B, W_{R'}).
\]
is an isomorphism. For that note that we have 
\begin{equation}
\label{parind5}
C(G/B, W) \, =\, \dlim_{n} C(K(n)\backslash G/B, W)
\end{equation}
and $C(K(n)\backslash G/B, W)_{R'}= C(K(n)\backslash G/B, W_{R'})$ since $K(n)\backslash G/B$ finite for every $n\ge 0$. 

(b) can be proved similarly. Using the isomorphisms \eqref{parind4},  \eqref{parind5} the assertion follows from the fact that the functor $W\mapsto C(K(n)\backslash G/B, W)$ commutes with inductive limits
for every $n\ge 0$.
\end{proof}

\begin{remarks}
\label{remarks:parind}
\rm (a) By (\cite{vigneras}, Lemma 4.7) a discrete $R[\barT]$-module is $W$ is admissible if and only if the $R[G]$-module $\Ind_{B}^{G} W$ is admissible. 
\medskip

\noi (b) Let $W\in \Mod_R^{\sm}(\barT)$. Since $W=\bigcup_{U} W^{U}$ ($U$ ranges over all open subgroups of $\barT^0$) we have 
\begin{equation*}
\label{parind3a}
\Ind_{B}^{G} W  \, = \, \dlim_U \Ind_{B}^{G} W^U
\end{equation*}
by Lemma \ref{lemma:parindbasechange} (b) above.  \enddemo
\end{remarks}

\paragraph{The functor of ordinary parts} We now consider the restriction of the parabolic induction functor \eqref{parind1} to the full subcategory $\Mod_R^{\ladm}(\barT)$ of $\Mod_R^{\sm}(\barT)$, i.e.\ we consider the functor
\begin{equation}
\label{parindladm1}
\Ind_{B}^{G}: \Mod_R^{\ladm}(\barT)\lra \Mod_R^{\sm}(G), \,\, W\mapsto \Ind_{B}^{G} W.
\end{equation} 
By (\cite{vigneras}, Prop.\ 4.3) if $H=1$ then the functor \eqref{parind1} admits a right adjoint $R_B^G:  \Mod_R^{\sm}(G)\to  \Mod_R^{\sm}(T)$. Therefore the functor
\begin{equation}
\label{adparindladm1}
\Mod_R^{\sm}(G) \lra \Mod_R^{\ladm}(T), \,\, V\mapsto R_B^G(V)_{\ladm}
 \end{equation} 
is right adjoint to \eqref{parindladm1}. If $H\ne 1$ then it follows that composing \eqref{adparindladm1} with the functor $W \mapsto W^H$, i.e.\ the functor 
\begin{equation}
\label{adparindladm2}
\Mod_R^{\sm}(G) \lra \Mod_R^{\ladm}(\barT), \,\, V\mapsto \left(R_B^G(V)_{\ladm}\right)^H
 \end{equation} 
is right adjoint to \eqref{parindladm1}.

By (\cite{emerton2}; see also \cite{vigneras}) the functors \eqref{adparindladm1} and \eqref{adparindladm2} admit a rather explicit description which we are going to review. Let $V\in \Mod_R^{\sm}(G)$. Recall that $N$ denotes the unipotent radical of $B$. Given two compact open subgroups $N_1\subseteq N_2$ of $N$ the homomorphism $h_{N_2, N_1}: V^{N_1} \to V^{N_2}$ is defined by $h_{N_2, N_1}(v) = \sum_{n\in N_2/N_1} n v$. If $N_3$ is another compact open subgroup of $N$ then $h_{N_3, N_2}\circ h_{N_2, N_1} = h_{N_3, N_1}$. Recall that $N_0$ denotes the image of $\cO_F$ under the isomorphism \eqref{unipotent}. The submodule $V^{N_0}$ of $N_0$-invariant elements of $V$ carries a canonical $R[T^+]$-module structure given as follows: for $t\in T^+$ define $h_{N_0, t}: V^{N_0} \to V^{N_0}$ by $h_{N_0, t}(v) = h_{N_0, N_0^t}( t v)$ (where $N_0^t = t N_0 t^{-1}$). If $t\in T^0$ then the map $h_{N_0, t}$ is just given by the action of $t$ on $V^{N_0}$. For $t_1, t_2\in T^+$ we have $h_{N_0, t_1} \circ h_{N_0, t_2}= h_{N_0, t_1t_2}$. Thus the maps $h_{N_0, t}$ induce an action of the monoid $T^+$ on $V^{N_0}$ -- i.e.\ an $R[T^+]$-module structure -- that extends the $T^0$-action. 

We now assume that $R=\cO$ is a complete noetherian local ring with maximal ideal $\fm$ and residue field $k$. Emerton's functor of ordinary parts 
\begin{equation*}
\label{Ord}
\Ord=\Ord_{\cO}: \Mod_{\cO}^{\sm}(G)\lra \Mod_{\cO}^{\ladm}(T).
\end{equation*}
is defined as the composite of $\Mod_{\cO}^{\sm}(G)\to \Mod_{\cO}^{\sm}(T^+), V\mapsto V^{N_0}$ with the functor \eqref{adTTplus2}, i.e.\ for $V\in  \Mod_{\cO}^{\sm}(G)$ we have 
\begin{equation*}
\label{Ord2}
\Ord_{\cO}(V) \, =\, \Hom_{\cO[T^+]}( \cO[T], V^{N_0})_{t_0^{\pm 1}-\fin}\, =\, \Gamma_{\cO}^{\ord}(V^{N_0}).
\end{equation*} 
If $H\ne 1$ then we define 
\begin{equation*}
\label{Ord3H}
\Ord_{\cO}^H=\Ord_{\cO}^H: \Mod_{\cO}^{\sm}(G) \lra \Mod_{\cO}^{\ladm}(\barT), \,\, V\mapsto \Ord(V)^H.
\end{equation*} 
It is easy to see (using e.g.\ Lemma \ref{lemma:parindbasechange} and Prop.\ \ref{prop:ordadj} below) that $\Ord_{\cO}^H$ does not depend on the coefficient ring $\cO$ in the following sense: if $\varphi: \cO\to \cO'$ is an epimorphism between
complete noetherian local rings then $\Ord_{\cO}^H(V_{\varphi})= \Ord_{\cO}^H(V_{\varphi})$ for every $V\in \Mod_{\cO'}^{\sm}(G)$. Therefore we will often drop the ring $\cO$ from the notation.

For an $\cO[T]$-module $W$ we let $W^{\io}$ be the $\cO[T]$-module $W^{\io}=W$ but with the new $T$-action given by $t\cdot w = t^{-1} w$. We have 

\begin{prop}
\label{prop:ordadj}
(a) The functor \eqref{adparindladm1} (in the case for $H=1$) is isomorphic to the functor $V\mapsto \Ord(V)^{\io}$, i.e.\ there are canonical isomorphisms 
\begin{equation}
\label{adjunctionord}
\Hom_{\cO[G]}(\Ind_{B}^{G} W, V) \, \cong \, \Hom_{\cO[T]}(W, \Ord(V)^{\io})
\end{equation}
for every $W\in \Mod_{\cO}^{\ladm}(T)$ and $V\in \Mod_{\cO}^{\sm}(G)$.
\medskip

\noi (b) More generally, the functor $\Mod_{\cO}^{\sm}(G) \lra \Mod_{\cO}^{\ladm}(\barT), V\mapsto \Ord^H(V)^{\io}$ is right adjoint to \eqref{parindladm1}. 
\end{prop}

\begin{proof} (a) is proved in (\cite{emerton2}, Thm.\ 4.4.6) in the case where $\cO$ is an Artinian local ring with finite residue field of characteristic $p$ and in (\cite{vigneras}, Cor.\ 7.3 and Remark 7.6) for arbitrary (noetherian) $\cO$. Note that in loc.\ cit.\ the adjunction property \eqref{adjunctionord} is stated for the functors $W\mapsto \Ind_{B}^{G} W$ and $V\mapsto \Ord_{\barB}(V)$ where $\barB\subseteq G$ denotes the opposite Borel subgroup to $B$. However if $w_0\in G$ denotes the element $w_0=\left(\begin{matrix} 0 & 1\\  1 & 0\end{matrix}\right)\!\! \mod \! Z$, then the map $V\to V, v\mapsto w_0\cdot v$ induces an isomorphism $\Ord_{\barB}(V) \cong \Ord_B(V)^{\io}$.

(b) follows immediately from (a).
\end{proof}

\begin{remarks}
\label{remarks:emertonord}
\rm We list a few properties of the functors $\Ord$ and $\Ord^H$. 
\medskip

\noi (a) By Prop.\ \ref{prop:adTTplus} (c) (see also \cite{emerton2}, Prop.\ 3.2.4) $\Ord^H$ commutes with direct limits. 
\medskip

\noi (b) For $V\in  \Mod_{\cO}^{\sm}(G)$ we have (see \cite{vigneras}, Remark 8.1)
\begin{equation*}
\label{Ord3Ha}
\Ord^H(V) \, =\, \Hom_{\cO[\barT^+]}( \cO[\barT], V^{H N_0})_{t_0^{\pm 1}-\fin}.
\end{equation*} 
\noi (c) If $V\in \Mod_{\cO}^{\sm}(G)$ is admissible then $\Ord^H(V)$ is admissible as well by (\cite{emerton2}, Prop.\ 3.3.3) and (\cite{vigneras}, Thm.\ 8.1).
\medskip

\noi (d) Let $V$ be a finitely generated $\cO$-module equipped with the trivial $G$-action. If the residue field $k$ of $\cO$ is of characteristic $p$ then we have $\Ord(V) =0=\Ord^H(V)$. Indeed, in this case $t_0$ acts on $(V^{N_0})=V$ by multiplication with $\Norm(\fp) = [\cO_F:\fp]$ so we have by Prop.\ \ref{prop:adTTplus}
\[
\Ord(V) \, =\, V^{\ord} \, =\, 0.
\]
\noi (e) Assume that $\cO=E$ is a field of characteristic $0$ and let $\mu: T \to \bQ^*\subseteq E^*$ be the modulus quasicharacter of $T$ given by $\mu\left(\left( \begin{smallmatrix} x_1 & \\ & x_2\end{smallmatrix} \right)\!\!\mod Z\right) = \Norm(\fp)^{-v_F(x_1) + v_F(x_2)}$ (where $\Norm(\fp)$ is the order of the residue field of $\cO_F$). It follows immediately from (\cite{casselman}, Cor.\ 4.2.5) that
\begin{equation*}
\label{OrdJacquet}
\Ord_E(V) \, \cong \, J(V)(\mu^{-1})
\end{equation*}
for any admissible $E[G]$-module $V$ where $J(V) = V/\langle n\cdot v- v\mid v\in V, n\in N\}\rangle$ denotes the Jacquet module of $V$. This allows us to describe $\Ord_E(V)$ explicitly if $V$ is an irreducible admissible $E[G]$-module such that $\dim_E V =\infty$.

If $V$ is a principal series, i.e.\ if $V\cong \Ind_B^G E(\chi)$ for a quasicharacter $\chi: T\to E^*$ with $\chi^2\ne 1, \mu^2$ then we have (see e.g.\ \cite{bushhen}, \S 9)
\begin{equation*}
\label{OrdJacquet2}
\Ord_E(\Ind_B^G E(\chi)) \, \cong \,  E(\chi^{-1}) \oplus E(\chi \mu^{-1}).
\end{equation*}
If $V$ is a special series representation, i.e.\ if it is isomorphic to the unique $\infty$-dimensional quotient $\sigma(\chi)$ of $\Ind_B^G E(\chi)$ where $\chi: T\to E^*$ is a quasicharacter with $\chi^2=1$ then
\begin{equation*}
\label{OrdJacquet3}
\Ord_E(\sigma(\chi)) \, \cong \,  E(\chi).
\end{equation*}
Finally, if $V$ is supercuspidal then we have $\Ord_E(V) =0$.  \enddemo
\end{remarks}

\begin{prop}
\label{prop:ordHcomp}
Let $V$ be a discrete $\cO[G]$-module. Assume that $H$ is an open subgroup of $T^0$ and let $n\ge 1$ with $\delta(U_F^{(n)})\subseteq H$. Then we have
\medskip

\noi (a) $V^{K_H(n)}$ is an $\cO[\barT^+]$-submodule of $V^{H N_0}$. 
\medskip

\noi (b) The functor 
\begin{equation*}
\label{ordfinlevel}
\Mod_{\cO}^{\sm}(G) \lra \Mod_{\cO}^{\ladm}(\barT), \,\, V\mapsto \Gamma^{\ord}(V^{K_H(n)})
\end{equation*} 
is isomorphic to $\Ord^H$. 
\end{prop}

\begin{proof} (a) This is statement 2 in the proof of (\cite{vigneras}, Prop.\ 8.2). It follows from the fact that for $t\in T^+$ a system of representatives $\{n_i\}_{i\in I}$ of $N_0/N_0^{t}$ is also a system of representatives of $K_H(n)/K_H(n)\cap K_H(n)^t$ so that the restriction of $t\cdot : V^{HN_0} \to V^{HN_0}$ to $V^{K_H(n)}$ is given by the Hecke operator $[K_H(n)tK_H(n)]: V^{K_H(n)}\to V^{K_H(n)}, v\mapsto \sum_{i\in I} n_i t v$.

(b) By Remark \ref{remarks:emertonord} (b) it suffices to see that the map
\begin{equation}
\label{finlevel}
\Gamma^{\ord}(V^{K_H(n)}) \lra \ \Gamma^{\ord}(V^{H N_0})\, =\, \Ord^H(V)
\end{equation}
induced by the inclusion $V^{K_H(n)}\hra V^{H N_0}$ is an isomorphism. Firstly, note that \eqref{finlevel} is injective since $\Gamma^{\ord}$ is left exact. 

For the surjectivity put $N_1= N_0^{t_0}$ so $N_1$ is the image of the maximal ideal $\fp\subseteq \cO_F$ under the isomorphism \eqref{unipotent}. Note that $K_H(n)\cap K_H(n)^{t_0^{-1}} = K_H(n+1)$ so that the map 
\begin{equation}
\label{heckeup}
h:=(h_{N_0, t_0})|_{V^{K_H(n)}}: V^{K_H(n)}\to V^{K_H(n)}
\end{equation} 
decomposes in the form 
\begin{equation*}
\label{heckeup2}
\begin{CD}
h: V^{K_H(n)} @> \incl >> V^{K_H(n+1)}  @> \tau_{n+1} >> V^{K_H(n)}.
\end{CD}
\end{equation*}
Here $\tau_{n+1}$ is the composite map
\[
\begin{CD}
\tau_{n+1}: V^{K_H(n+1)} @> v\mapsto t_0 \cdot v >> V^{K_H(n)^0} @> v\mapsto h_{N_1, N_0}(v) >> V^{K_H(n)}
\end{CD}
\]
where $K_H(n)^0 = K_H(n)^{t_0}\cap K_H(n)= K(n) H N_1$. Note that
\[
\begin{CD} 
V^{K_H(n+1)} @>\incl >> V^{K_H(n+2)} \\
@VV \tau_{n+1} V @VV \tau_{n+2} V \\
V^{K_H(n)} @>\incl >> V^{K_H(n+1)}
\end{CD}
\]
is a commutative diagram of $\cO[\barT^+]$-modules. By applying the functor $\Hom_{\cO[\barT^+]}( \cO[\barT],\wcdot )$ we obtain a commutative diagram of $\cO[\barT]$-modules
\[
\begin{CD} 
\Hom_{\cO[\barT^+]}( \cO[\barT], V^{K_H(n+1)} )@>\incl_*>> \Hom_{\cO[\barT^+]}( \cO[\barT], V^{K_H(n+2)} ) \\
@VV (\tau_{n+1})_* V @VV (\tau_{n+2})_* V \\
\Hom_{\cO[\barT^+]}( \cO[\barT], V^{K_H(n)} ) @>\incl_* >> \Hom_{\cO[\barT^+]}( \cO[\barT], V^{K_H(n+1)} ).
\end{CD}
\]
Since the action of $t_0$ on $\Hom_{\cO[\barT^+]}( \cO[\barT], V^{K_H(n+1)} )$ is the map $\incl_* \circ (\tau_{n+1})_*= (\tau_{n+2})_* \circ \incl_* = h_*$, it follows that $h_*$ is an isomorphism (because $t_0$ is a unit in $\cO[\barT]$). Hence 
\[
\incl_*: \Hom_{\cO[\barT^+]}( \cO[\barT], V^{K_H(n)} ) \lra \Hom_{\cO[\barT^+]}( \cO[\barT], V^{K_H(n+1)} )
\]
is surjective, whence bijective. Therefore $\incl_*: \Gamma^{\ord}(V^{K_H(n)}) \to \Gamma^{\ord}(V^{K_H(n+1)})$ is an isomorphism as well. Thus the transition maps in the inductive system
\[
\begin{CD} 
\Gamma^{\ord}(V^{K_H(n)})  @>\incl_* >> \Gamma^{\ord}(V^{K_H(n+1)}) @>\incl_* >>\Gamma^{\ord}(V^{K_H(n+2)}) @>>> \ldots 
\end{CD}
\]
are isomorphisms. Since $\bigcap_{m\ge n} K_H(m)= H N_0$, we have $V^{N_0 H} = \bigcup_{m\ge n} V^{K_H(m)}$ and we conclude with Prop.\ \ref{prop:adTTplus} (d)
\[
\Gamma^{\ord}(V^{K_H(n)}) \, \cong \, \dlim_{m\ge n} \Gamma^{\ord}(V^{K_H(m)})\, \cong\, \Gamma^{\ord}(V^{H N_0}) \, =\, \Ord^H(V).
\]
\end{proof}

\begin{remark}
\label{remark:ordHfinlevel}
\rm If $V$ is an admissible $\cO[T]$-module then the fact that \eqref{finlevel} is an isomorphism has already been pointed out in the proof of (\cite{vigneras}, Thm.\ 8.1).  \enddemo
\end{remark}

We now assume that $\cO$ denotes a valuation ring of a $p$-adic field $E$, i.e.\ $\cO$ is a complete discrete valuation ring with finite residue field $k$ of characteristic $p$ and quotient field $E$ of characteristic $0$. The following result will be used in the proof of Theorem \ref{theorem:reciprocity} in section \ref{subsection:shimuracurve2}. 

\begin{prop}
\label{prop:classifyordjacq}
Let $V$ be an irreducible admissible $E[G]$-module. We assume that there exists an open subgroup $H$ of $T^0$, an integer $n\ge 1$ and an $\cO[T^+/H]$-submodule $M$ of $V^{K_H(n)}$ with the following properties:

\noi (i) $\delta(U_F^{(n)})\subseteq H$ and $V^{K_H(n)} \ne 0$;

\noi (ii) $M$ is finitely generated as an $\cO$-module and we have $M^{\ord} \ne 0$. 

\noi Then there exists a quasicharacter $\chi: T\to \cO^*$ with $H\subseteq \ker(\chi)$ such that 
\begin{equation*}
\label{Ordchar0}
V_E \, \cong \, \left\{ \begin{array}{cc} \Ind_B^G E(\chi^{-1}) & \mbox{if $\chi^2\ne 1$,}\\
                                                      \sigma(\chi) & \mbox{if $\chi^2 = 1$.}\end{array}\right.
\end{equation*}
Moreover in this case we have $M^{\ord} = \cO(\chi)$ as an $\cO[T]$-module.
\end{prop}

\begin{proof} It is easy to see that the inclusion $M\hra V^{K_H(n)}$ induces an $T/H$-equivariant monomorphism 
\begin{equation}
\label{Ordbasechange}
M^{\ord} \, \lra \, (V^{K_H(n)})^{\ord} \, = \,\Ord_E^H(V) \, \subseteq \, \Ord_E(V).
\end{equation}
In particular (ii) implies that $\Ord_E^H(V)\ne 0$ hence $V$ is either a principal series or a special representation by Remark \ref{remarks:emertonord} (e). 

Firstly, consider the case when $V$ is a principal series representation. Then there exists a quasicharacter $\chi: T\to E^*$ with $\chi^2 \ne 1, \mu^2$ so that 
$V \cong \Ind_B^G E(\chi^{-1})$, hence $\Ord_E(V) \cong E(\chi) \oplus E((\chi\mu)^{-1})$ by Remark \ref{remarks:emertonord} (e). Since $\Ord_E^H(V)\ne 0 $ we have $E(\chi)^H\ne 0$ or $E((\chi\mu)^{-1})^H \ne 0$ hence $H\subseteq \ker(\chi)$ and  $\Ord_E^H(V)= \Ord_E(V)$. The injectivity of \eqref{Ordbasechange} implies that we can view $M^{\ord}$ as a non-trivial $T$-stable finitely generated $\cO$-submodule of $E(\chi) \oplus E((\chi\mu)^{-1})$. This is impossible except if either $\chi$ or $\chi\mu$ has values in $\cO^*$. Since $\Ind_B^G E(\chi^{-1})\cong \Ind_B^G E(\chi\mu)$ we may assume that $\chi$ has values in $\cO^*$. Hence $(\chi \mu)(T) \not\subseteq \cO^*$ and therefore $M^{\ord} \cap E((\chi\mu)^{-1})=0$. It follows $M^{\ord} \subseteq E(\chi)$, whence $M^{\ord} \cong \cO(\chi)$.

Now assume that $V$ is a special representation, so that $V\cong \sigma(\chi)$ for a quasicharacter $\chi: T\to E^*$ with $\chi^2 = 1$. Note that the values of $\chi$ lie in  $\{\pm 1\} \subseteq \cO^*$ and we have $0\ne M^{\ord} \subseteq \Ord_E^H(V) \subseteq \Ord_E(V) \cong E(\chi)$ by Remark \ref{remarks:emertonord} (e). As above, we get $H\subseteq \ker(\chi)$ and $M^{\ord} \cong \cO(\chi)$. 
\end{proof}

\subsection{The right derived functors of $\Ord$}
\label{subsection:derord}

We again assume that $\cO$ is a complete noetherian local ring with maximal ideal $\fm$. We will denote the $n$-th right derivative of the functor $\Ord_{\cO}$ by $\Ord_{\cO}^n(\wcdot)$, i.e.\ we put
\begin{equation*}
\label{ordderived}
\Ord_{\cO}^n(V)\, =\, (R^n \Ord_{\cO})(V)
\end{equation*}
for $V\in \Mod_{\cO}^{\sm}(G)$. 
More generally for a closed subgroup $H$ of $T^0$ the $n$-th right derivative of the functor $\Ord_{\cO}^H$ will be denoted by $\Ord_{\cO}^{H,n}$. Again we often drop the coefficient ring $\cO$ from the notation.

Since $\Ord^H$ is the composite of the functors $V\mapsto V^{N_0H}$ and $\Gamma^{\ord}$ one may ask whether there exists a corresponding Grothendieck spectral sequence. For the applications in sections \ref{section:ordinary} and \ref{section:cohomqhmsv} it suffices to establish 
part (b) of the following

\begin{prop}
\label{prop:rnordcomp}
Let $V\in \Mod_{\cO}^{\sm}(G)$.
\medskip

\noi (a) We have $\Ord_{\cO}^{H, n}(V) = \dlim_{U \open,\, H\le U \le T^0} \Ord_{\cO}^{U, n}(V)$
for every $n\ge 0$. 
\medskip

\noi (b) Assume that $H$ be an open subgroup of $T^0$ and that $\dim(\cO) \le 1$. Let $n$ be a positive integer with $\delta(U_F^{(n)})\subseteq H$. Then there exists a spectral sequence 
\begin{equation}
\label{Ord5ss}
E_2^{rs} = R^r \Gamma^{\ord}(H^s(K_H(n), V)) \, \Longrightarrow\, \Ord^{H, r+s}(V).
\end{equation}
\end{prop}

For the proof we need some preparation. As before we put $\barT = T/H$, $\barT^0 = T^0/H$, $\barT^+ = T^+/H$ and $t_0=\bar{\delta}(\varpi)\in \barT^+$ where $\varpi\in \cO_F$ is a prime element. Let $W$ be an $\cO[\barT^0]$-module. We can view $W$ as a $\cO[K_0(n)]$-module via the projection $K_0(n)\to K_0(n)/K_H(n)\cong T^0/H=\barT^0$ so we can consider the discrete $\cO[G]$-module 
\begin{equation*}
\label{cindhecke}
\cInd_{K_0(n)}^G W \, =\, \{ \Phi: G\to W\mid \supp(\Phi) \, \mbox{is compact}\, \mbox{and}\,\Phi(kg) = k\Phi(g)\, \forall\, k\in K_0(n), g\in G\}.
\end{equation*}
The $G$-action on $\cInd_{K_0(n)}^G W$ is induced by right multiplication. Moreover $\cInd_{K_0(n)}^G W$ is equipped with an $\cO[\barT^0]$-action -- denote by $\star$ -- induced by the $\cO[\barT^0]$-action on $W$ and with a Hecke action, i.e.\ with an action of $\cO[\barT^+]$ defined as follows: for $t\in \barT^+$ and $\Phi\in \cInd_{K_0(n)}^G W$ the function $t\cdot \Phi$ is given by
\begin{equation*}
\label{cindhecke2}
(t\cdot \Phi)(g)\, =\, \sum_{i\in I}\Phi(t^{-1} n_i^{-1} g) \qquad \forall g\in G
\end{equation*}
where $\{n_i\}_{i\in I}$ is a system of representatives of $N_0/N_0^{t}$. For $t\in \barT^0$ we have $t\cdot \Phi = t\star \Phi$. 

We now consider the case where $\cO=k$ is a field and $W$ is one-dimensional as a $k$-vector space. Thus $\barT^0$ acts on $W$ via a character $\chi: \barT^0\to k^*$. By abuse of notation we denote the character $K_0(n)\stackrel{\pr}{\lra} K_0(n)/K_H(n)\cong T^0/H=\barT^0\to k^*$ also by $\chi$ 
 so we have $W\cong k(\chi)$ as $k[K_0(n)]$-modules. We provide $\cInd_{K_0(n)}^G k(\chi)$ with an $k[X]$-module sturcture by letting $X$ acts via multiplication with $t_0$. 

\begin{lemma}
\label{lemma:noprimehtor}
(a) The $k[X]$-module $\cInd_{K_0(n)}^G k(\chi)/\left(\cInd_{K_0(n)}^G k(\chi)\right)\!\{X\}$ is torsionfree (hence flat). 
\medskip

\noi (b) If $I$ is an injective object of $\Mod_k^{\sm}(G)$ then the $k[X]$-module $\Hom_{k[G]}(\cInd_{K_0(n)}^G k(\chi), I)$ is $\Gamma_k^{X\sord}$-acyclic.
\end{lemma}

Here for a $k[X]$-module $M$, i.e.\ we put $M\{X\} :=\bigcup_{i\ge 1} M[X^i]$ where $M[X^i] := \ker(M\stackrel{X^i\wcdot}{\lra} M)$ for $i\ge 1$.   

\begin{proof} (a) We recall some basic facts about the Bruhat-Tits tree $\cT$ of $G$. Its set of vertices $\cV$ is the set of homothety classes of lattices in $F^2$. For $n\ge 1$ we call a sequence of vertices $c=(v_0, v_1, \ldots, v_n)$ an $n$-path in $\cT$ if $v_i$ is adjacent to $v_{i+1}$ and if and $v_{i+1} \ne v_{i-1}$ for $i=1, \ldots, n-1$. We denote the set of all $n$-paths in $\cT$ by $\oE^{(n)}$. If $n\ge 2$ (resp.\ $n=1$) and if $c=(v_0, v_1, \ldots, v_n), c=(v_0', v_1', \ldots, v_n')\in \oE^{(n)}$ then we say that $c'$ is a successor of $c$ if $v_{i+1} = v_i'$ for $i=0, \ldots, n-1$ (resp.\ if $v_1 = v_0'$ and $v_1'\ne v_0$). Thus if $n=1$ then $\oE^{(1)}=\oE$ is the set of oriented edges of $\cT$. 
 
We use the height function $h: \cV\to \bZ$ introduced in \cite{barthellivne} (see also \cite{spiess1}, \S 3.3). It is defined as follows. For $v\in \cV$ the geodesic ray from $v$ to $\infty$ has a non-empty intersection with the standard apartment $A=\{v^*_m\mid m\in \bZ\}$.\footnote{The vertex $v^*_m$ is the homothety class of the lattice $\cO_F\oplus \fp^m$; see \cite{spiess1}, \S 3.3} If $v^*_m$ is any element of the intersection then we have $h(v) = m- d(v, v^*_m)$ where $d(v, v_m^*)$ is the distance between $v$ and $v_m^*$. For $c=(v_0, v_1, \ldots, v_n)\in \oE^{(n)}$ we define
$h(c) := h(v_n)$. 

For $i\in \{0, \ldots, n-1\}$ we set 
\[
\oE^{(n)}_i\, = \, \{(v_0, v_1, \ldots, v_n)\in \oE^{(n)} \mid h(v_i) > h(v_{i+1}) > \ldots > h(v_n)\}.
\]
and $(\oE^{(n)})^0 : = \, \oE^{(n)}\setminus \oE^{(n)}_{n-1}$. By (\cite{spiess1}, Lemma 3.6) we have
\[
(\oE^{(n)})^0 \, =\, \{(v_0, v_1, \ldots, v_n)\in \oE^{(n)} \mid h(v_0) < h(v_1) < \ldots < h(v_n)\}.
\]
Also if $c\in \oE^{(n)}_i$ and $i\ge 1$ then every successor $c'$ of $c$ lies in $\oE^{(n)}_{i-1}$. If $c\in \oE^{(n)}_0$ then all its successors $c'$ lie in $\oE^{(n)}_0$ and we have $h(c') = h(c) -1$. Finally if $c\in (\oE^{(n)})^0$ then there exists a unique successor -- denoted by $s(c)$ -- that lies again in $(\oE^{(n)})^0$ and we have $h(s(v)) = h(c) +1$. All other successors lie in $\oE^{(n)}_{n-1}$. 

Moreover for $m\in \bZ$ we define an increasing sequence of subsets of $\oE^{(n)}_0$ and $(\oE^{(n)})^0$ by $\oE^{(n)}_{0, m} = \{c\in \oE^{(n)}\mid h(c) \le m\}$ and $(\oE^{(n)})^0_m= \{c\in (\oE^{(n)})^0\mid h(c) \ge - m\}$
respectively. If $c\in \oE^{(n)}_{0, m}$ then its successors lie in $\oE^{(n)}_{0, m-1}$ whereas for $c\in (\oE^{(n)})_m^0$ we have $s(c) \in (\oE^{(n)})_{m-1}^0$.

Put $\cM= \cInd_{K_0(n)}^G k(\chi)$. Note that the support of any function $\Phi: G\to k$ in $\cM$ is the union of finitely many left $K_0(n)$-cosets. Since $G$ acts transitively on $\oE^{(n)}$ and the stabilizer of the $n$-path $c^*=(v_n^*, v_{n-1}^*, \ldots, v_1^*, v_0^*)$ is the group $K_0(n)$, this fact may be rephrased as follows: the support of any $\Phi\in \cM$ is of the form $\nu^{-1}(Y)$ where 
\begin{equation*}
\label{prinbundle}
\nu: G \lra \oE^{(n)}, \,\, g\mapsto g^{-1} c^*
\end{equation*}
and $Y\subseteq \oE^{(n)}$ is finite. A simple computation shows that for $\Phi\in \cM$ with $\supp(\Phi)= \nu^{-1}(c)$ for some $c\in \oE^{(n)}$ we have 
\begin{equation}
\label{prinbundle2}
\supp(X\cdot \Phi)\, =\,  \nu^{-1}(\{c'\in \oE^{(n)}\mid \mbox{$c'$ is successor of $c$}\}).
\end{equation}
We set $\cM^0=\{\Phi\in \cM\mid \supp(\Phi) \subseteq \nu^{-1}((\oE^{(n)})^0)\}$ and 
\[
\cM_i\, =\, \{\Phi\in \cM\mid \supp(\Phi) \subseteq \nu^{-1}(\oE^{(n)}_i)\}
\]
for $i\in \{0, \ldots, n-1\}$. We have $\cM_0\subseteq \cM_1 \subseteq \ldots \subseteq \cM_{n-1}$ and $\cM= \cM_{n-1}\oplus \cM^0$. It follows from \eqref{prinbundle2} that the subspaces $\cM_0, \cM_1, \ldots, \cM_{n-1}$ are 
$k[X]$-submodules of $\cM$ that satisfy
\[
X\cdot \cM_i \subseteq \cM_{i-1}
\]
for every $i\in \{1, \ldots, n-1\}$. The subspace $\cM^0$ is not a $k[X]$-submodule of $\cM$ but by identifying $\cM^0$ with the quotient $\cM/\cM_{n-1}$ we can provide $\cM^0$ with a $k[X]$-module structure as well. 

Furthermore for $m\in \bZ$ we set 
\begin{eqnarray*}
\cM_{0, m} & = & \{\Phi\in \cM\mid \supp(\Phi) \subseteq \nu^{-1}(\oE^{(n)}_{0, m})\} \\
\cM_m^0 & = &\{\Phi\in \cM\mid \supp(\Phi) \subseteq \nu^{-1}((\oE^{(n)})_m^0)\}
\end{eqnarray*}
so that $\{\cM_{0, m}\}_{m\in \bZ}$ and $\{\cM^0_{m}\}_{m\in \bZ}$ are increasing, separating and exhausting sequence of $k[X]$-submodules of $\cM_0$ and $\cM^0$ respectively. By \eqref{prinbundle2} we have
\[
X\cdot \cM_{0, m} \subseteq \cM_{0, m-1} \qquad \mbox{and} \qquad X\cdot \cM_m^0 \subseteq \cM_{m-1}^0
\]
for every $m\in \bZ$.

To show that the $k[X]$-module $\cM/\cM\{X\}$ is torsionfree it suffices to show that if $P\in k[X]$ is a polynomial with constant term $1$ and if $\Phi\in \cM$ with $P\cdot \Phi=0$ then we have $\Phi=0$. Note that $P$ acts as the identity on each of the quotients 
$\cM_{0, m}/\cM_{0, m-1}$, $\cM_m^0/\cM_{m-1}^0$ and $\cM_i/\cM_{i-1}$ since $X$ annihilates them. So if $\overline{\Phi}:=\Phi\! \mod \! \cM_{n-1} \in \cM/\cM_{n-1}\cong \cM^0$ lies in $\cM_m^0$ for some $m\in \bZ$ then it lies in $\cM_{m-1}^0$ hence in $\cM_{m-2}^0$ etc. This shows 
$\overline{\Phi}\in \bigcap_m \cM_m^0=\{0\}$, i.e.\ $\Phi\in \cM_{n-1}$. Similarly we deduce successively $\Phi\in \cM_{n-2}$, $\Phi\in \cM_{n-3}$ etc.\ hence $\Phi\in \cM_0$. Finally if $\Phi\in \cM_{0, m}$ for some $m\in \bZ$ then we get  $\Phi\in \bigcap_m \cM_{0, m}=\{0\}$, so $\Phi=0$.

(b) Since $I$ is injective the sequence of $k[X]$-modules 
\begin{equation}
\label{ordacyclic}
0\lra \Hom_{k[G]}(\cM/\cM\{X\}, I) \lra \Hom_{k[G]}(\cM, I)\lra \Hom_{k[G]}(\cM\{X\}, I)\lra 0
\end{equation}
is exact. The first $k[X]$-module is injective since, by (a), the left adjoint of 
$\Hom_{k[G]}(\cM/\cM\{X\}, \wcdot ): \Mod_k^{\sm}(G)\to \Mod_{k[X]}$, namely the functor 
\[
\Mod_{k[X]}\to \Mod_k^{\sm}(G),\quad M \mapsto \cM/\cM\{X\}\otimes_{k[X]} M,
\] 
is exact. For the $k[X]$-module $\Hom_{k[G]}(\cM\{X\}, I)=\prolimn \Hom_{k[G]}(\cM[X^n], I)$ 
we get 
\[
R^q\Gamma_k^{X\sord}(\Hom_{k[G]}(\cM\{X\}, I)) \,=\, 0
\] 
for every $q\ge 0$ by Lemmas \ref{lemma:adlfinxinv} (c) and \ref{lemma:adlfinxinv2} (b). It follows
\[
R^q\Gamma_k^{X\sord}(\Hom_{k[G]}(\cM, I))\, \cong \, R^q\Gamma_k^{X\sord}(\Hom_{k[G]}(\cM/\cM\{X\}, I))\, =\, 0
\]
for every $q\ge 1$.
\end{proof}

As a consequence of the previous Lemma we obtain

\begin{lemma} 
\label{lemma:injacyclic}
Assume that $\dim(\cO)\le 1$. Let $I$ be an injective object of $\Mod_{\cO}^{\sm}(G)$ and let $W\in \Mod_{\cO[\barT^0], f}$. Then the $\cO[\barT^+]$-module $\Hom_{\cO[G]}(\cInd_{K_0(n)}^G W, I)$ is $\Gamma^{\ord}$-acyclic.
\end{lemma}

\begin{proof} Firstly, some preliminary observations. The finite $\cO$-algebra $\cO[\barT^0]$ is a product of a finite number of complete noetherian local  $\cO$-algebras $\cO[\barT^0]= \prod_{j\in J} \wcO_j$. Every $W\in \Mod_{\cO[\barT^0]}$, $V\in  \Mod_{\cO[\barT^+]}$ as well as the ring $\cO[\barT^+]$ decomposes accordingly
\[
W = \prod_{j\in J} W_j, \qquad V = \prod_{j\in J} V_j, \qquad \cO[\barT^+] = \prod_{j\in J} \wcO_j[t_0].
\]
Note that for $V\in\Mod_{\cO}^{\sm}(\barT^+)=\Mod_{\cO[\barT^+]}$ we have  
\begin{equation}
\label{injacyclic2}
R^q \Gamma_{\cO}^{\ord} V \, =\, \prod_{j\in J} R^q \Gamma_{\wcO_j}^{t_0\sord}(V_j)\qquad \forall \, q\ge 0.
\end{equation}
For $W\in \Mod_{\cO[\barT^0]}$ we have
\[
\Hom_{\cO[G]}(\cInd_{K_0(n)}^G W, I) \cong \prod_{j\in J} \Hom_{\cO[G]}(\cInd_{K_0(n)}^G W_j, I)
 \cong  \prod_{j\in J} \Hom_{\wcO_j[G]}(\cInd_{K_0(n)}^G W_j, I_j)
\]
where $I_j: = \Hom_\cO(\wcO_j, I)$. Note that $I_j$ is again an injective 
object of $\Mod_{\wcO_j}^{\sm}(G)$ since the functor $\Hom_{\cO}(\wcO_j, \wcdot): \Mod_{\cO}^{\sm}(G)\to \Mod_{\wcO_j}^{\sm}(G)$ has an exact left adjoint, namely the forgetful functor $\Mod_{\wcO_j}^{\sm}(G)\to \Mod_{\cO}^{\sm}(G)$. Thus by \eqref{injacyclic2} 
we obtain 
\begin{equation}
\label{injacyclic3}
R^q \Gamma^{\ord}(\Hom_{\cO[G]}(\cInd_{K_0(n)}^G W, I)) \, =\, \prod_{j\in J} R^q \Gamma_{\wcO_i}^{t_0\sord}(\Hom_{\wcO_j[G]}(\cInd_{K_0(n)}^G W_j, I_j))
\end{equation}
for every $q\ge 0$. 

Let $\varphi: \cO \to \cO'$ be an epimorphism of complete noetherian local rings of dimension $\le 1$, let $W\in \Mod_{\cO'[\barT^0]}$ and let $I$ is an injective object of $\Mod_{\cO}^{\sm}(G)$. Put $I': = \Hom_{\cO}(\cO', I)= I[\fa]$ (note that $I'$ is an injective object of  $\Mod_{\cO'}^{\sm}(G)$). We claim that 
\begin{eqnarray}
\label{injacyclic4}
&& \Hom_{\cO'[G]}(\cInd_{K_0(n)}^G W, I')\quad  \mbox{ is $\Gamma_{\cO'}^{\ord}$-acyclic} \\
&&\hspace{2cm} \quad \Longrightarrow \quad 
\Hom_{\cO[G]}(\cInd_{K_0(n)}^G W_{\varphi}, I)\quad \mbox{is $\Gamma_{\cO}^{\ord}$-acyclic}\nonumber
\end{eqnarray}
Indeed, if $\wcO_j': = \wcO_j/\fa \wcO_j$ where $\fa = \ker(\varphi)$ then the decomposition of $\cO'[\barT^0]$ into a product of complete local rings is given by $\cO'[\barT^0]= \prod_{j\in J} \wcO_j'$.
Thus if $W = \prod_{j\in J} W_j$ is the corresponding decomposition of $W$ then we have $(W_{\varphi})_j = (W_j)_{\varphi}$. Since $I'_j = \Hom_{\cO'}(\wcO_j', I') =\Hom_{\cO_j}(\wcO_j', I_j)$ we get
\[ 
\Hom_{\wcO_j[G]}(\cInd_{K_0(n)}^G (W_{\varphi})_j, I_j)\, =\, \Hom_{\wcO_j'[G]}(\cInd_{K_0(n)}^G W_j, I_j')_{\varphi}.
\]
Thus \eqref{injacyclic3} together with Prop.\ \ref{prop:artinadlfinder} implies \eqref{injacyclic4}.

Let $S_0$ denote the set of non-zero divisors of $\cO$. For the assertion it suffices to prove 
\medskip

\begin{claim}
\label{claim:ordacyclic} The $\cO[\barT^+]$-module $\Hom_{\cO[G]}(\cInd_{K_0(n)}^G W, I)$ is in the following cases $\Gamma_{\cO}^{\ord}$-acyclic:

\noi (i) $\dim(\cO) =0$ and $W\in \Mod_{\cO[\barT^0], f}$.

\noi (ii) $\dim(\cO)=1$, $W\in \Mod_{\cO[\barT^0]}$ such that $S_0^{-1} W= 0$ and $W[\fm^m]\in \Mod_{\cO, f}$ for all $m\ge 1$. 

\noi (iii) Let $\dim(\cO)=1$ and $W\in \Mod_{(S_0^{-1}\cO)[\barT^0], f}$.

\noi (iv) Let $\dim(\cO)=1$ and $W\in \Mod_{\cO[\barT^0], f}$.
\end{claim}

In case (i) note that $\dim(\cO) =0$ implies that each of rings $\wcO_j$, $j\in J$ is an Artinian local ring. By \eqref{injacyclic3} it suffices to see that for $j\in J$ and an $\wcO_j$-module $W_j$ of finite length we have 
\begin{equation*}
\label{injacyclic5}
R^q \Gamma_{\wcO_j}^{t_0\sord}(\Hom_{\wcO_j[G]}(\cInd_{K_0(n)}^G W_j, I_j))\, =\, 0 \qquad \forall \, q\ge 1.
\end{equation*}
Since the functor $W_j \mapsto \Hom_{\wcO_j[G]}(\cInd_{K_0(n)}^G W_j, I_j)$ is exact it is enough to prove \eqref{injacyclic4} if $W_j$ has lenght 1. In this case it follows from Lemma \ref{lemma:noprimehtor} (b) and Prop.\ \ref{prop:artinadlfinder}.

For (ii), firstly we remark that the assumption $S_0^{-1} W=0$ implies $W= \dlim_{m\ge 1} W[\fm^m]$. For $m\ge 1$ we have 
\[
\Hom_{\cO[G]}(\cInd_{K_0(n)}^G W[\fm^m], I)\, \cong \, \Hom_{\cO_m[G]}(\cInd_{K_0(n)}^G W[\fm^m], I[\fm^m])
 \]
where $\cO_m = \cO/\fm^m$. Since $I[\fm^m] = \Hom_{\cO}(\cO_m, I)$ is an injective object of $\Mod_{\cO_m}^{\sm}(G)$, case (i) together with \eqref{injacyclic4} implies that 
$\Hom_{\cO_m[G]}(\cInd_{K_0(n)}^G W[\fm^m], I[\fm^m])$ is $\Gamma_{\cO}^{\ord}$-acyclic for every $m\ge 1$. Also since the functor $W'\mapsto \Hom_{\cO[G]}(\cInd_{K_0(n)}^G W', I)$ is exact we see that the transition map
\[
\Hom_{\cO[G]}(\cInd_{K_0(n)}^G W[\fm^{m+1}], I) \to \Hom_{\cO[G]}(\cInd_{K_0(n)}^G W[\fm^m], I)
\]
is surjective for every $m\ge 1$ and that the kernel is isomorphic to 
\[
\Hom_{\cO[G]}(\cInd_{K_0(n)}^G W[\fm^{m+1}]/W[\fm^m], I)\, \cong\, \Hom_{\cO_1[G]}(\cInd_{K_0(n)}^G W[\fm^{m+1}]/W[\fm^m], I[\fm])
\]
hence it is $\Gamma^{\ord}$-acyclic as well. Therefore we can apply Lemma \ref{lemma:derfunclim2} of the appendix 
to deduce that
\[
\Hom_{\cO[G]}(\cInd_{K_0(n)}^G W, I)\, \cong\, \prolimm \Hom_{\cO[G]}(\cInd_{K_0(n)}^G W[\fm^m], I)
\]
is $\Gamma^{\ord}$-acyclic as well.

For the case (iii) note that $(S_0^{-1}\cO)[\barT^0]$ has dimension $0$, so it is a finite product of Artinian local rings $S_0^{-1}\cO[\barT^0] = \prod_{\xi\in \Xi} A_{\xi}$. We denote the residue field of $A_{\xi}$ by $E_{\xi}$. Any $W\in \Mod_{(S_0^{-1}\cO)[\barT^0], f}$ admits a corresponding decomposition 
$W= \prod_{\xi\in \Xi} W_{\xi}$. Again the exactness of the functor $W \mapsto \Hom_{\cO[G]}(\cInd_{K_0(n)}^G W, I)$ implies that it suffices to consider the case when $W$ is a $(S_0^{-1}\cO)[\barT^0]$-module of length 1. In this case there exists $\xi\in \Xi$ such that $W_\xi\cong E_{\xi}$ and $W_{\xi'} =0$ for all $\xi'\ne \xi$. If we put $E:=E_{\xi}$ then we have 
\[
\Hom_{\cO[G]}(\cInd_{K_0(n)}^G W, I)\, =\, \Hom_{E[G]}(\cInd_{K_0(n)}^G E, I_0)
\]
where $I_0:=\Hom_{\cO}(E, I)$. It is an injective object of $\Mod_E^{\sm}(G)$. There exists $j\in J$ such that the canonical map $\cO[\barT^0] = \prod_{j\in J} \wcO_j\to (S_0^{-1}\cO)[\barT^0]\stackrel{\pr}{\lra} A_{\xi}\stackrel{\pr}{\lra} E$ 
factors in the form $\cO[\barT^0]\stackrel{\pr}{\lra} \wcO_j \to E$. Let $\wcO'= \image(\wcO_j\to E)$. Note that $E$ is the quotient field of $\wcO'$. 
By Prop.\ \ref{prop:artinadlfinder} it suffices to show 
\begin{equation}
\label{injacyclic6}
R^q \Gamma_{\wcO'}^{X\sord}(\Hom_{E[G]}(\cInd_{K_0(n)}^G E, I_0))\, =\, 0 \qquad \forall \, q\ge 1
\end{equation}
where $\cM= \Hom_{E[G]}(\cInd_{K_0(n)}^G E, I_0)$. The proof of Lemma \ref{lemma:noprimehtor} (b) shows (see \eqref{ordacyclic}) that there exists an exact sequence of $E[X]$-modules 
\begin{equation}
\label{injacyclic7}
0 \lra \cM \lra \Hom_{E[G]}(\cInd_{K_0(n)}^G E, I_0)\lra \prolimm \cN_m \lra 0
\end{equation}
where $\cM$ is an injective $E[X]$-module and where $\cN_m$ is an $E[X]$-module with $X^m \cN_m$ for every $m\ge 1$. 
Since $E[X]$ is a flat $\wcO'[X]$-algebra, $\cM$ is also injective as an $\wcO'[X]$-module hence $R^q \Gamma_{\wcO'}^{X\sord} \cM=0$ for every $q\ge 0$. Moreover we have $R^q \Gamma_{\wcO'}^{X\sord} \cN=0$ for every $q\ge 0$ by Lemma \ref{lemma:adlfinxinv} (c) and \ref{lemma:adlfinxinv2} (b). Thus by applying the $\delta$-functor $R^{\bu} \Gamma_{\wcO'}$ to \eqref{injacyclic7} yields \eqref{injacyclic6}.

For (iv) let $W\in \Mod_{\cO[\barT^0], f}$ and put $W_1= \ker(W\to S_0^{-1} W)$, $W_2= \image(W\to S_0^{-1} W)$, $W_3
= S_0^{-1} W$ and $W_4= \coker(W\to S_0^{-1} W)$ so that there exists short exact sequences 
\begin{equation}
\label{injacyclic8}
0 \to W_1\to W \to W_2 \to 0, \qquad 0 \to W_2\to W_3 \to W_4 \to 0.
\end{equation}
Clearly, we have $S_0^{-1} W_1 =0= S_0^{-1} W_4$, and $W_3\in \Mod_{(S_0^{-1}\cO)[\barT^0], f}$. 
Also since $W\in \Mod_{\cO, f}$ we have $W_1, W_2\in \Mod_{\cO, f}$ hence $W_1[\fm^m]\in \Mod_{\cO, f}$ and  $W_4[\fm^m]\cong \Ext^1_{\cO}(\cO_m, W_2) \in\Mod_{\cO, f}$ for every $m\ge 1$. By (ii) the claim holds for $W_1$ and $W_4$ and it holds for $W_3$ by (iii). Thus the sequences \eqref{injacyclic8} yield that the claim holds for $W_2$ and $W$ as well. 
\end{proof}

\begin{proof}[Proof of Prop.\ \ref{prop:rnordcomp}] (a) Let $\cU$ denote the set of all open subgroups $U$ of $T^0$ containing $H$. We have
\[
\Ord^H(V) \, = \, \Ord(V)^H \, =\, \dlim_{U\in \cU} \Ord^U(V)
\]
for every $V\in \Mod_{\cO}^{\sm}(G)$. Thus if $0\lra V\lra I^{\bu}$ is an injective resolution then we get
\begin{eqnarray*}
\Ord^{H, n} (V)& = & H^n(\Ord^H(I^{\bu}))\,= \,H^n(\dlim_{U\in \cU} \Ord^U(I^{\bu}))\, =\, \dlim_{U\in \cU} H^n(\Ord^U(I^{\bu}))\\
& = &  \dlim_{U\in \cS} \Ord^{U, n}(V).
\end{eqnarray*}

For (b) let $I$ be an injective object of $\Mod_{\cO}^{\sm}(G)$.
We have 
\[
I^{K_H(n)}\, \cong \, \Hom_{\cO[G]}(\cInd_{K_H(n)}^G \cO, I)\, \cong \, 
\Hom_{\cO[G]}(\cInd_{K_0(n)}^G \Ind_{K_H(n)}^{K_0(n)} \cO, I)
\]
so by Lemma \ref{lemma:injacyclic} (applied to $W= \Ind_{K_H(n)}^{K_0(n)} \cO$) the 
$\cO[\barT^+]$-module $I^{K_H(n)}$ is $\Gamma^{\ord}$-acyclic. Thus by Prop.\ \ref{prop:ordHcomp} (b) there exists a Grothendieck spectral sequence \eqref{Ord5ss} associated to the decomposition
\begin{equation*}
\label{ordfinlevel1a}
\begin{CD}
\Mod_{\cO}^{\sm}(G) @> V\mapsto V^{K_H(n)} >> \Mod_{\cO}(\barT^+)@> \Gamma^{\ord} >> \Mod_{\cO}^{\ladm}(\barT)
\end{CD}
\end{equation*}
of the functor $\Ord^H$. 
\end{proof}

\begin{remark} 
\label{remark:Xtor} 
\rm Let $\cO=k$ be a field of characteristic $p$. The $k[X]$-module (and $G$-representation) $\cInd_{K_0(1)}^G k$ can be identified with $C_c(\oE, k)$, the set of maps $\Phi: \oE\to k$ with finite support. If $e=(v_0, v_1)$ is an oriented edge of $\cT$ then we put $o(e) =v_0$, $t(e) =v_1$ and $\bar{e} = (v_1, v_0)$. For $e\in \oE$ we let $1_e\in C_c(\oE, k)$ be the function $1_e(e') = 1$ if $e'=e$ and $1_e(e')=0$ otherwise. Fix $v\in \cV$ and define $\Phi\in C_c(\oE, k)$ by $\Phi:= \sum_{e\in \oE, t(e) =v} 1_e$. We have 
\[
X\cdot \Phi \, =\, \sum_{e\in \oE, t(e) =v} X\cdot 1_e\, =\, \sum_{e\in \oE, t(e) =v} \sum_{e'\in \oE, o(e') = v, e'\ne \bar{e}} 1_{e'} \, =\, [\cO_F: \fp] \sum_{e\in \oE, t(e) =v} 1_e \, =\, 0.
\]
This example shows that the $X$-primary torsion submodule of $\cInd_{K_0(n)}^G k(\chi)$ considered in Lemma \ref{lemma:noprimehtor} (a) above can be non-trivial. Thus, in general, the functor $\Mod_{\cO}^{\sm}(G) \to  \Mod_{\cO}(\barT^+), V\mapsto V^{K_H(n)}$ does not have an exact left adjoint.
\end{remark}

As a first consequence of the existence of the spectral sequence \ref{Ord5ss} we remark

\begin{corollary}
\label{corollary:derordncoeff}
Let $\varphi:\cO\to \cO'$ be an epimorphism of complete noetherian local rings. If $\dim(\cO)\le 1$ and $\dim(\cO') =0$ then the canonical homomorphism\footnote{The definition of \eqref{derordncoeff1} is similar to the definition of \eqref{adlfinder1}.}
\begin{equation}
\label{derordncoeff1}
\Ord_{\cO}^{H, n}(V)_{\varphi} \lra \Ord_{\cO}^{H, n}(V_{\varphi})
\end{equation}
is an isomorphism for every $V\in \Mod_{\cO'}^{\sm}(G)$ and $n\ge 0$.
\end{corollary}

\begin{proof} By Prop.\ \ref{prop:rnordcomp} (a) it suffices to consider the case where $H$ is an open subgroup of $T^0$. 
Choose $m\ge 1$ with $\delta(U_F^{(m)})\subseteq H$. As in the proof of Lemma \ref{lemma:injacyclic} we write
$\cO[\barT^0]$ (with $\barT^0=T^0/H$) as a product of complete noetherian local rings $\cO[\barT^0]= \prod_{j\in J} \wcO_j$ and decompose $\cO'[\barT^0]= \prod_{j\in J} \wcO_j'$ accordingly. By \eqref{injacyclic2} and Prop.\ \ref{prop:artinadlfinder} we have for $W\in \Mod_{\cO'}(\barT^+)$
\begin{equation}
\label{derordncoeff2}
R^n \Gamma_{\cO}^{\ord} W_{\varphi} \, =\, \prod_{j\in J} R^n \Gamma_{\wcO_j}^{t_0\sord}((W_j)_{\varphi_j})\,\cong\,  \prod_{j\in J} R^n \Gamma_{\wcO_j}^{t_0\sord}((W_j)_{\varphi_j})\, =\, R^n \Gamma_{\cO'}^{\ord} W
\end{equation}
where $\varphi_j: \wcO_j\to \wcO_j'$ denotes the $j$-component of the map $\cO[\barT^0]\to \cO'[\barT^0]$ induced by $\varphi$. 

Since $\Ord_{\cO'}^H(V)_{\varphi} = \Ord_{\cO}^H(V_{\varphi})$ for every $V\in \Mod_{\cO'}^{\sm}(G)$ and since $V\mapsto \Ord_{\cO}^{H, \bu}(V_{\varphi})$ is a $\delta$-functor it suffices to see that for an injective object 
$I\in \Mod_{\cO'}^{\sm}(G)$ we have $\Ord_{\cO}^{H, n}(I_{\varphi})= 0$ for every $n\ge 1$. Consider the spectral sequence \eqref{Ord5ss} for $I_{\varphi}$
\[
E_2^{rs}=R^r \Gamma_{\cO}^{\ord}(H^s(K_H(m), I)_{\varphi}) \, \Longrightarrow\, \Ord^{H, r+s}(I_{\varphi})
\]
Since $V\mapsto H^n(K_H(m), V)$ is the $n$-th right derived functor of $\Mod_{\cO'}^{\sm}(G)\to \Mod_{\cO'}(\barT^+), V\mapsto V^{K_H(m)}$ we have $E_2^{rs}=0$ for $s\ge 0$. Together with \eqref{derordncoeff2} and Lemma \ref{lemma:injacyclic} it follows 
\[
\Ord_{\cO}^{H, n}(I_{\varphi})\, \cong\, R^n \Gamma_{\cO}^{\ord}(I^{K_H(m)}_{\varphi})\, \cong \, R^n \Gamma_{\cO'}^{\ord}(I^{K_H(m)}) \, =\, 0
\]
for every $n\ge 1$. 
\end{proof}

Assume now that $\cO=A$ is an Artinian local ring with finite residue field of characteristic $p$. 
If $H$ is an arbitrary closed subgroup of $T^0$ then the fact that $\Ord^H$ is the composition of the functors $V \mapsto V^{HN_0}$ and $\Gamma^{\ord}$ together with (\cite{emerton3}, Prop.\ 2.1.11) implies that there are canonical homomorphisms 
\begin{equation}
\label{ordderived3}
\Ord^{H, n}(V)\, \lra \, \Gamma^{\ord}(H^n(HN_0, V)) 
\end{equation} 
for every $V\in \Mod_{\cO}^{\sm}(G)$ and $n\ge 0$. 

\begin{corollary}
\label{corollary:adordn}
If $V$ is a locally admissible $A[G]$-module then the homomorphism \eqref{ordderived3} is an isomorphism for every $n\ge 0$. Moreover if $V$ is an admissible $A[G]$-module then $\Ord^{H, n}(V)$ is an admissible $A[\barT]$-module for every $n\ge 0$. 
\end{corollary}

\begin{proof} As a first step we consider the case where $H$ is an open subgroup of $T^0$ and $V\in \Mod_{\cO}^{\ladm}(G)$. Let $m\ge 0$ with $\delta(U^{(m)}) \subseteq H$. By Prop.\ \ref{prop:ordHcomp} (b) there exists a spectral sequence 
\begin{equation}
\label{Ord5ssa}
E_2^{rs} = R^r \Gamma^{\ord}(H^s(K_H(m), V)) \, \Longrightarrow\, E^{r+s}=\Ord^{H, r+s}(V).
\end{equation}
By (\cite{emerton3}, Lemma 3.4.4 and Thm.\ 3.4.7 (a)) the $A[\barT^+]$-modules $H^s(K_H(m), V)$ are admissible (resp.\ locally admissible). By Prop.\ \ref{prop:rgammaord} (b) we have $E_2^{rs}=0$ if $r>0$. Hence the edge map $E^n\to E_2^{0n}$
is an isomorphism for every $n\ge 0$. Moreover if $V$ is admissible then we obtain 
using Prop.\ \ref{prop:adTTplus} (a) that 
\begin{equation*}
\label{Ord5ssedge}
\Ord^{H,n}(V) \, \cong \, H^n(K_H(m), V)^{\ord}
\end{equation*}
is finitely generated as an $\cO$-module, hence admissible as an $\cO[\barT]$-module. 

To simplify the notation we now consider only the other extreme case $H=1$ (the proof for arbitrary $H$ is essentially the same). Consider the diagram 
 \begin{equation*}
\label{Ordlim}
 \begin{CD}
 \dlim_m\Ord^{\delta(U^{(m)}), n}(V) @>\cong >> \dlim_m \Gamma^{\ord}(H^n(K_1(m), V))\\
 @VV\cong V @VV \cong V \\
\Ord^n(V)@>\eqref{ordderived3}>> \Gamma^{\ord}(H^n(N_0, V)).
\end{CD}
\end{equation*}
The upper horizontal map is the limit of edge morphisms of the degenerating spectral sequence \eqref{Ord5ssa} (for $H=\delta(U^{(m)})$) so it is an isomorphism. The left vertical map is an isomorphism by Prop.\ \ref{prop:ordHcomp} (a).
The right vertical map is induced by the restriction maps $H^n(K_1(m), V)\to H^n(N_0, V)$. It is an isomorphism by Prop.\ \ref{prop:adTTplus} (d) since we have $H^n(N_0, V)= \dlim_m H^n(K_1(m), V)$ (see e.g.\ \cite{emerton3}, Lemma 3.4.3). 

It remains to prove that $\Ord^n(V)\in \Mod_A^{\adm}(T)$ if $V\in \Mod_A^{\adm}(G)$. We will prove this by induction on $n$. The assertion holds for $n=0$ by Remark \ref{remarks:emertonord} (c). Assume that $n\ge 1$ and that $\Ord^r(V)$ is admissible for every $r<n$. We fix a compact open subgroup $U$ of $T^0$. There is a Grothendieck spectral sequence
\begin{equation*}
\label{Ord6ss}
E_2^{rs} = H^r(U, \Ord^s(V)) \, \Longrightarrow\, \Ord^{U, r+s}(V).
\end{equation*}
By (\cite{emerton3}, Lemma 3.4.4) we have $E_2^{rs} =H^r(U, \Ord^s(V))\in \Mod_{A, f}$ for every $r\ge 0$ and $s<n$. It follows that $E_m^{rs}\in \Mod_{A,f}$ for every $r\ge 0$, $s<n$ and $m\in \{2, 3, \ldots, \infty\}$ as they are subquotient of $E_2^{rs}$. Since $E^n=\Ord^{U, n}(V)\in \Mod_{A, f}$ as shown above and since $E_{\infty}^{0n} = E_{n+1}^{0n}$ is a quotient of $E^n$ we have $E_{n+1}^{0n}\in \Mod_{A, f}$. The 
exact sequences 
\[
0\lra E_{m+1}^{0n}\lra E_{m}^{0n} \stackrel{d_{m}^{0n}}{\lra} E_{m}^{m(n-m+1)} 
\]
for $m=n, n-1, \ldots, 2$ now allow us to deduce that the $A$-modules $E_n^{0n}, E_{n-1}^{0n}, \ldots, E_2^{0n}$ are finitely generated. Hence $E_2^{0n}= (\Ord^n(V))^U \in \Mod_{A, f}$. Since this holds for every open subgroup $U$ of $T^0$ we conclude that the $A[T]$-module $\Ord^n(V)$ is admissible.
\end{proof}

\begin{remarks}
\label{remarks:derivedord}
\rm (a) The canonical homomorphism \eqref{ordderived3} also exists if $A=E$ is a field of characteristic 0 (instead of (\cite{emerton3}, Prop.\ 2.1.11) one can use here the fact the functor $\cInd_{HN_0}^G: \Mod_E^{\sm}(HN_0)\to \Mod_E^{\sm}(G)$ is exact, so that the restriction $\Mod_E^{\sm}(G)\to \Mod_E^{\sm}(HN_0)$ preserves injectives). Also Cor.\ \ref{corollary:adordn} holds in this case as well and its proof requires very little modification. Indeed, in this case the spectral sequence \eqref{Ord5ssa} degenerates, i.e.\ we have $E_2^{rs}= 0$ if $s\ge 0$. Since the functor $\Mod_E^{\sm}(G)\to \Mod_E^{\sm}(T^+), V\mapsto V^{N_0}$ is exact it thus follows 
\begin{equation*}
\label{ordderivedchar0}
\Ord_E^n(V)\, =\, 0
\end{equation*} 
for every $V\in \Mod_E^{\ladm}(G)$ and $n\ge 0$.
\medskip

\noi (b) In (\cite{emerton3}, Conj.\ 3.7.2) Emerton has given the following conjectural description of the right derived functors of $\Ord$ when restricted to the subcategory $\Mod_{\cO}^{\ladm}(G)$: if $\wOrd$ denotes this restriction then we have
\begin{equation}
\label{wordderived}
R^n \wOrd(V)\, \cong \, \Gamma^{\ord}(H^n(N_0, V)) 
\end{equation}
for every $V\in \Mod_{A}^{\ladm}(G)$ (so together with Cor.\ \ref{corollary:adordn} this implies $R^n \wOrd(V)= R^n \Ord(V)$).

This can be deduced from Prop.\ \ref{prop:rgammaord} (a) as follows.
Firstly, note that for $V\in \Mod_{A}^{\adm}(G)$ the $A[T^+]$-module $V^{N_0}= \bigcup_{m\ge 1} V^{K_1(m)}$ is locally admissible (since $V^{K_1(m)}\in \Mod_{A, f}^{\sm}(T^+)$ for every $m\ge 1$). Since every locally admissible $A[G]$-module is the union of its admissible 
submodules we also have $V^{N_0}\in \Mod_A^{\ladm}(T^+)$ if $V\in\Mod_{A}^{\ladm}(G)$, i.e.\ the functor $\wOrd$ factors as 
\[
\begin{CD}
\Mod_{A}^{\ladm}(G) @> V\mapsto V^{N_0} >> \Mod_A^{\ladm}(T^+) @> \Gamma^{\ord} >> \Mod_A^{\ladm}(T).
\end{CD}
\]
Hence Prop.\ \ref{prop:rgammaord} (a) together with (\cite{emerton3}, Prop.\ 2.1.11) implies \eqref{wordderived}. 
A proof of (\cite{emerton3}, Conj.\ 3.7.2) for $G = \GL_2(F)$ (and implicitly also for $G = \PGL_2(F)$) 
by completely different means has been given in \cite{emerton-paskunas}.
In \cite{spiess3} the derived functors of the functor of ordinary parts will be studied for general $p$-adic reductive groups. 
\enddemo
\end{remarks}

Let $d=[F:\bQ_p]$ and let $V$ be an admissible $A[G]$-module. By Cor.\ \ref{corollary:adordn} and (\cite{emerton3}, 3.6.1, 3.6.2) we have 
\begin{equation*}
\label{ordderived2}
\Ord^n(V)\, = \,\left\{\begin{array}{cc} V_N(\alpha) & \mbox{if $n=d$,}\\
0 &  \mbox{if $n>d$.}
\end{array}\right.
\end{equation*}
Here $V_N$ denote the $N$-coinvariants of $V$ and $\alpha$ the character 
\begin{equation}
\label{grosschar}
\begin{CD}
\alpha: T@>\delta^{-1}>> F^*@>\Norm_{F/\bQ_p}>>\bQ_p^* @>x\mapsto x/p^{v_p(x)}>> \bZ_p^*@>>> A^*
\end{CD}
\end{equation}
(where $v_p:\bQ_p^* \to \bZ\cup\{\infty\}$ is the normalized valuation of $\bQ_p$). In particular we obtain

\begin{corollary}
\label{corollary:ordind}
For $W\in \Mod_A^{\adm}(T)$ we have 
\begin{equation*}
\label{ordderived2b}
\Ord^n(\Ind_B^G(W))\, = \,\left\{\begin{array}{cc} W^{\io} & \mbox{if $n=0$,}\\
W(\alpha) &  \mbox{if $n=d$,}\\
0 & \mbox{if $n>d$.}
\end{array}\right.
\end{equation*}
\end{corollary}

\begin{proof} For $n=0$ this follows from (\cite{emerton2}, Prop.\ 4.3.4). For $n=d$ it suffices to show that the evaluation map
\begin{equation}
\label{res1}
\Ind_{B}^{G} W\,\lra\, W, \,\,\Phi \mapsto \Phi(1).
\end{equation}
induces an isomorphism $(\Ind_{B}^{G} W)_N\cong W$. Note that \eqref{res1} is an epimorphism of $A[B]$-modules. Hence it suffices to see that the $N$-coinvariants of the kernel vanish. 

Let $U$ be this kernel. It consists of those $\Phi\in \Ind_{B}^{G} W$ that have support in the big cell $B \omega N$ (where $\omega = \left(\begin{smallmatrix} 0 & 1\\ 1 & 0 \end{smallmatrix}\right)$). So by mapping $\Phi\in U$ to the function $F\to W, y\mapsto \Phi(\omega n(y))$ the module $U$ can be identified with the $A[N]$-module $C_c(F, W)$ of locally constant maps $f: F\to W$ with compact support. The $N$-action is given by $(n(y)\cdot f)(x) = f(x+y)$ for $x, y\in F$, $f\in C_c(F, W)$. As an $A[N]$-module $C_c(F, W)$ is generated by functions 
$1_{\fp^m} \cdot w$ for $m\in \bZ$, $w\in W$ (where $1_{\fp^m}\in C_c(F, \bZ)$ denotes the characteristic function of $\fp^m$). 
Let $q$ be the order of the residue field of $F$ and choose an integer $N\ge 1$ so that $q^N=0$ in $A$. Since $1_{\fp^m} \cdot w = \sum_y n(y) \cdot 1_{\fp^{m+1}} \cdot w$ (where $y$ ranges over a system of representatives of $\fp^m/\fp^{m+1}$) we see that we have 
\[
1_{\fp^m} \cdot w \, \equiv \, 1_{\fp^{m+1}} \cdot (q w) \,  \equiv \, 1_{\fp^{m+2}} \cdot (q^2 w)\, \equiv \, \ldots \,  \equiv \, 1_{\fp^{m+N}} \cdot (q^N w)\, =\, 0
\]
in $U_N$. Hence $U_N=0$. 
\end{proof}

\subsection{$\varpi$-adically admissible and Banach space representations of $T$}
\label{subsection:admbanach}

In this section $R=\cO$ denotes the valuation ring of a $p$-adic field $E$, i.e.\ $\cO$ is a complete discrete valuation ring with finite residue field $k$ of characteristic $p$ and quotient field $E$ of characteristic $0$. We denote the maximal ideal of $\cO$ by $\fm$ and fix a prime element $\varpi\in \fm$. For $m\ge 1$ we put $\cO_m =\cO/\fm^m$. We denote the normalized absolute value on $E$ by $| \wcdot |: E\to \bR$ (so $|\varpi| = \# k^{-1}$). 

For an $\cO$-module $N$ we denote its torsion submodule by $N_{\tor} = N[\fm^{\infty}]$ and the maximal torsionfree quotient by $N_{\fl}$ so that there exists a short exact sequence 
\begin{equation*}
\label{torfl}
0 \lra N_{\tor} \lra N \lra N_{\fl} \lra 0.
\end{equation*} 
We also set $N_m= N\otimes_{\cO}\cO_m$ and $N[\fm^m] = \Hom_{\cO}(\cO_{m}, N)$. The 
kernel-cokernel exact sequence for the maps $N\stackrel{\varpi}{\lra} N  \stackrel{\varpi^m}{\lra} N$ is
\begin{equation}
\label{kercoker}
0 \lra N[\fm]\lra N[\fm^{m+1}]\lra N[\fm^m]\lra N_1\lra N_{m+1}\lra N_m \lra 0 
\end{equation} 
for $m\ge 1$. 

As before let $H$ be a closed subgroup of $T^0$ and put $\barT= T/H$. We are going to review the notion of $\varpi$-adically continuous and $\varpi$-adically admissible $\cO[\barT]$-modules (see \cite{emerton2}, \S 2.4). Firstly, an $\cO[\barT]$-module $W$ is called $\varpi$-adically continuous if 
(i) $W$ is $\varpi$-adically complete and separated, (ii) $W_{\tor}$ is of bounded exponent, i.e.\ $\varpi^m W=0$ for $m\ge 1$ sufficiently large and (iii) $W_m$ is a discrete $\cO_m[\barT]$-module for every $m\ge 1$.  

A $\varpi$-adically admissible $\cO[\barT]$-module $W$ is a $\varpi$-adically continuous $\cO[\barT]$-module $W$ such that $W_1 \in \Mod_k^{\adm}(\barT)$. The exactness of the sequence \eqref{kercoker} implies then that we have $W_m \in \Mod_{\cO_m}^{\adm}(\barT)$ for every $m\ge 1$. The full subcategory of the category of $\cO[\barT]$-modules consisting of $\varpi$-adically continuous (resp.\ $\varpi$-adically admissible) $\cO[\barT]$-modules will be denoted by 
$\Mod_{\cO}^{\varpi-\cont}(\barT)$ (resp.\ $\Mod_{\cO}^{\varpi-\adm}(\barT)$). The category $\Mod_{\cO}^{\varpi-\adm}(\barT)$ is closed under the formation of kernels, images and cokernels. It is in particular an abelian category. It is also closed under extensions in the category $\Mod_{\cO}^{\varpi-\cont}(\barT)$.

For an $\cO[\barT]$-module $W$ define $\sD(W) : = \Hom_{\cO}(W, \cO)$. By (\cite{emerton2}, Prop.\ 2.4.10) the assignment $W\mapsto \sD(W)$ induces an anti-equivalence of categories
\begin{equation}
\label{pontrajagin3}
\sD: \Mod_{\cO}^{\varpi-\adm}(\barT)_{\fl}\lra \Mod_{\cO}^{\fgaug}(\barT)_{\fl}, \,\, W\mapsto \sD(W).
\end{equation} 
Here for an $\cO$-linear abelian category $\sA$ we denote by $\sA_{\fl}$ the subcategory of torsionfree objects $A\in \sA$ (i.e.\ objects so that $\varpi\wcdot: A\to A$ is a monomorphism). 

\begin{remarks}
\label{remarks:padicadm}
\rm (a) Let $W\in \Mod_{\cO}^{\varpi-\adm}(\barT)$. The fact that $W_m$ is discrete when viewed as an $\cO[\barT]$-module implies that the $\cO[\barT]$-action on $W\cong {\prolim} \,_m\, W_m$ extends naturally to a $\La_{\cO}(\barT)$-action, i.e.\ any $\varpi$-adically admissible $\cO[\barT]$-module $W$ can be viewed as an augmented $\cO[\barT]$-module. In particular the category $\Mod_{\cO}^{\varpi-\adm}(\barT)$ as well as the functor \eqref{pontrajagin3} are $\La_{\cO}(\barT)$-linear.

\noi (b) Let $W\in \Mod_{\cO}^{\varpi-\adm}(\barT)$ and assume that $W$ is torsionfree (i.e.\ $W_{\tor} =0$). Using the exact sequence \eqref{kercoker} one can easily see that the $\cO_m$-module $W_m$ is free for every $m\ge 1$.

\noi (c) Simple examples of $\varpi$-adically admissible $\cO[T]$-modules are those associated to continuous characters. For that let $\cA$ be an $\cO$-algebra that is $\varpi$-adically complete and separated and let $\chi: F^*\to \cA^*$ be a continuous character (i.e.\ $\chi$ is continuous with respect to the $\varpi$-adic topology on $\cA$). We 
attach to $\chi$ the following $\cO[T]$-module $\cA(\chi)$. For that we identify $T$ with $F^*$ via the isomorphism \eqref{splittorus}, so that we can view $\chi$ as a character of $T$. For $t\in T$ and $a\in \cA$ define
\begin{equation}
\label{chitwist}
t\cdot a : = \chi(t) a.
\end{equation}
Then $\cA(\chi)$ is the $\cO$-module $\cA$ equipped with the $T$-action given by \eqref{chitwist}.
If $A$ is a finite $\cO$-algebra then we have $\cA(\chi)\in \Mod_{\cO}^{\varpi-\adm}(T)$. Moreover if $\chi$ is a quasicharacter (i.e.\ if $\ker(\chi)$ is open in $F^*$) then $\cA(\chi)\in \Mod_{\cO}^{\adm}(T)$.
If $\chi=1: F^*\to \cA^*$ is the trivial character then we will also write $\cA(0)$ rather than $\cA(\chi)$. 

\noi (d) Let $\psi: F^*\to \cO$ (i.e.\ we have $\psi(xy) = \psi(x) + \psi(y)$ for all $x,y\in \cO$) and let
$\wcO = \cO[\vep] = \cO[X]/(X^2)$, $\vep := X + (X^2)$ be the $\cO$-algebra of dual numbers. The character 
\begin{equation*}
\label{extstein1}
\Theta_{\psi}: F^*\lra \wcO^*,\,\, \Theta_{\psi}(x)=  1+ \psi(x)\ep.
\end{equation*}
is obviously continuous so as a special case of the above construction we obtain the $\varpi$-adically admissible $\cO[T]$-module $\wcO(\Theta_{\psi})$. It is an admissible $\cO[T]$-module if and only if $\ker(\psi)$ is open in $F^*$ (note that this holds if and only if $\psi$ is some multiple $c\cdot v_F$, $c\in \cO$ of the normalized valuation $v_F$ of $F$).  \enddemo
\end{remarks}

\paragraph{Banach space representations} Recall that an $E$-Banach space representation of $\barT$ is an $E$-Banach space $V=(V, \| \cdot \|)$ together with a continuous $E$-linear action $\barT\times V \to V, (t,v)\mapsto t\cdot v$. An $E$-Banach space representation $V$ of $\barT$ is called admissible if there exists an open and bounded $\cO[\barT]$-submodule $W\subseteq V$ such that the $U$-invariant $(W/V)^U$ of the quotient $V/W$ are an $\cO$-module of cofinite type for every compact open subgroup $U$ of $\barT$ (i.e.\ the Pontrjagin dual $D((W/V)^U)=\Hom((W/V)^U, E/\cO)$ is a finitely generated $\cO$-module). \footnote{Note that this condition implies that we can choose the norm $\|\wcdot \|$ on $V$ so that $(V, \| \cdot \|)$ is a unitary Banach space representation of $\barT$, i.e.\ we have $\|t\cdot v\|=\|v\|$ for every $t\in \barT$ and $v\in V$.}

The category of admissible $E$-Banach space representations of $\barT$ will be denoted by $\Ban_E^{\adm}(\barT)$. It is an abelian category. This can be easily deduced from the duality theorem (\cite{schickhof}, \cite{schneider-teitelbaum}, Thm.\ 3.5) or from the fact that $\Ban_E^{\adm}(\barT)$ is equivalent to the localized category $\Mod_{\cO}^{\varpi-\adm}(\barT)_E$ (see Lemma \ref{lemma:locbanach} below). 

We first review the duality theorem. By slight abuse of notation, a $\Lambda_{\cO}(\barT)_E$-module $M$ will be called an augmented $E[\barT]$-module. Again, $\Mod_E^{\aug}(\barT)$ denotes the category of augmented $E[\barT]$-modules. An augmented $E[\barT]$-module $M$ is called finitely generated if there exists a $\La(\barT)$-submodule $L$ with $L\in \Mod_{\cO}^{\fgaug}(\barT)$ and $L_E = M$. We equip $M$ with the topology induced by the canonical topology on $L$, i.e.\ $M$ is a topological vector space and the inclusion is $L\hra M$ is open and continuous (hence $M$ is locally compact). 
As before $\Mod_E^{\fgaug}(\barT)$ denotes the full category of $\Mod_E^{\aug}(\barT)$ of finitely generated augmented $E[\barT]$-modules. 

There is a canonical contravariant functor
\begin{equation}
\label{pontrajagin4}
\cD: \Ban_E^{\adm}(\barT)\lra \Mod_E^{\fgaug}(\barT), \,\, V\mapsto \cD(V)= V'.
\end{equation} 
Here for an $E$-Banach space $V=(V, \| \wcdot\|)$, $\cD(V) = V'$ denotes the dual Banach space equipped with the weak topology. 

\begin{lemma}
\label{lemma:locdual}
The functor \eqref{pontrajagin4} is an anti-equivalence of $E$-linear categories. 
Its quasi-inverse is given by
\begin{equation}
\label{pontrajagin5}
\Mod_E^{\fgaug}(\barT) \lra  \Ban_E^{\adm}(\barT), \, \,\, M\mapsto \Hom_{E, \cont}(M, E).
\end{equation} 
\end{lemma}

\begin{proof} This follows immediately from (\cite{schneider-teitelbaum}, Thm.\ 3.5).
\end{proof}

\begin{remark}
\label{remark:locdualpply}
\rm Let $V_1$, $V_2\in \Ban_E^{\adm}(\barT)$ and assume that $\dim_E(V_1) <\infty$. The Lemma implies that the $E$-vector space $\Hom_{E[\barT]}(V_1, V_2)$ is finite-dimensional as well. In fact any $\barT$-equivariant homomorphism $V_1\to V_2$ is automatically continuous so we have 
\begin{equation}
\label{pontrajagin5a}
\Hom_{E[\barT]}(V_1, V_2)\, \cong \, \Hom_{\Mod_E^{\fgaug}(\barT)}(\cD(V_2), \cD(V_1)).
\end{equation}
If $U$ is a compact open subgroup of $\barT$ then the $\Lambda_{\cO}(U)_E$-module \eqref{pontrajagin5a} is finitely generated, hence it has finite length (since $\cD(V_1)$ has finite length because of $\dim_E(\cD(V_1))= \dim_E(V_1) <\infty$) and is therefore finite-dimensional as an $E$-vector space. \enddemo
\end{remark}

Recall that for an $\cO$-linear additive category $\cA$ the $E$-linear additive category $\cA_E$ has the same objects as $\cA$ whereas the morphisms are given by 
$\Hom_{\cA_E}(A, B) = \Hom_{\cA}(A, B)\otimes_{\cO} E$. If $\cA$ is abelian then $\cA_E$ is abelian as well (see Lemma \ref{lemma:locab2} of the appendix). 

Note that $L_{\tor} = L[\fm^{\infty}]$ is of bounded exponent for any object $L\in \Mod_{\cO}^{\fgaug}(\barT)$. 
Indeed, since $\La_{\cO}(\barT)$ is noetherian and $L$ is finitely generated as $\La_{\cO}(\barT)$-module the sequence of submodules $L[\fm]\subseteq L[\fm^2]\subseteq \ldots \subseteq L[\fm^m]\subseteq \ldots$ becomes stationary. Therefore the induced functor 
\begin{equation*}
\label{loclambda}
(\Mod_{\cO}^{\fgaug}(\barT))_E\lra \Mod_E^{\fgaug}(\barT), \,\, L\mapsto L_E
\end{equation*} 
is an equivalence of categories. Similarly, the categories $\Mod_{\cO}^{\varpi-\adm}(\barT)_E$ and $\Ban_E^{\adm}(\barT)$
are equivalent.

\begin{lemma}
\label{lemma:locbanach}
For $W\in \Mod_{\cO}^{\varpi-\adm}(\barT)$ let $\|\wcdot\|$ be the norm on $V=W_E$ so that  
$\,\image(W\to V, w\mapsto w\otimes 1)$ is the unit ball $\{v\in V\mid \|v\|\le 1\}$ in $V$. Then $(V, \|\wcdot \|)$ is an admissible $E$-Banach space representation of $\barT$. The induced functor 
\begin{equation}
\label{locbanach2}\Mod_{\cO}^{\varpi-\adm}(\barT)_E\lra \Ban_E^{\adm}(\barT)
\end{equation} 
is an equivalence of categories. Moreover for $W\in \Mod_{\cO}^{\varpi-\adm}(\barT)$ and $V= W_E$ we have 
\begin{equation*}
\label{pontrajagin6}
\sD(W)\otimes_{\cO} E\,\cong\, \cD(V)
\end{equation*} 
\end{lemma}

\begin{proof} Since $W_E = (W_{\fl})_E$ we may assume that $W_{\tor} =0$, i.e.\ $W\subseteq V$
(note that $W_{\fl}\in \Mod_{\cO}^{\varpi-\adm}(\barT)$ by \cite{emerton2}, Cor.\ 2.4.13). Let $U$ be a compact open subgroup of $\barT$. The fact that $W_1 \in \Mod_k^{\adm}(\barT)$ implies that $\dim_k D((V/W)^U)\otimes_{\cO} k = \dim_k \Hom_k(W_1^U, k)< \infty$, hence 
$D((V/W)^U)\in \Mod_{\cO,f}$ by the topological Nakayama Lemma. This proves that $V$ is admissible.

To show that the induced functor \eqref{pontrajagin5} is an equivalence we have to see that every $V\in \Ban_E^{\adm}(\barT)$ lies in its essential image. For that we choose a norm $\|\wcdot\|$ so that $\|t\cdot v\|=\|v\|$ for every $t\in \barT$ and $v\in V$. Let $W$ be the unit ball in $V$. Then $(W/V)^U$ is of cofinite type for every compact open subgroup $U$ of $\barT$ (as shown in \cite{schneider-teitelbaum} prior to Thm.\ 3.5). Passing in the sequence $0\to W_1 \to V/W\stackrel{\varpi\cdot }{\lra} V/W$ to $U$-invariants and Pontrajagin duals implies that $\dim_k D(W_1^U) = \dim_k D((V/W)^U)\otimes_{\cO} k< \infty$ for every compact open subgroup $U$ of $\barT$. Hence $W_1\in \Mod_k^{\varpi-\adm}(\barT)$ and $W\in \Mod_{\cO}^{\varpi-\adm}(\barT)$. The last assertion is obvious.
\end{proof}
 
\begin{remarks}
\label{remarks:lamod}
\rm (a) The fact that the unit ball $W$ of an admissible $E$-Banach space representation $V$ of $\barT$ is $\varpi$-adically admissible implies that it carries a canonical $\La_{\cO}(\barT)$-action (see Remark \ref{remarks:padicadm} above). Hence the $E[\barT]$-action on $V$ extends naturally to a $\La_{\cO}(\barT)_E$-action and the functor \eqref{pontrajagin4} is $\La_{\cO}(\barT)_E$-linear.

\noi (b) Let $A$ be a finite-dimensional commutative $E$-algebra equipped with the canonical topology (induced by any choice of a norm $\| \wcdot \|: A\to \bR_{\ge 0}$). Let $\chi: F^*\to A^*$ be a continuous character. As in Remark \ref{remarks:padicadm} (c) we define the $E$-Banach space representation $A(\chi)$ of $T$ by $A(\chi) = A$ with $T$-action given by \eqref{chitwist}. We have $A(\chi)\in \Ban_E^{\adm}(\barT)$ if and only if the image of $\chi$ is bounded. Indeed if the latter holds then there exists an $\cO$-subalgebra $\cA$ of $A$ that is a lattice in $A$ so that $\chi(F^*)\subseteq \cA^*$. In this case $A(\chi)$ is the image of $\cA(\chi)$ under the functor \eqref{locbanach2}. 

\noi (c) Let $\psi: F^*\to E$ be a continuous character and let
$\wE = E[\vep]$ be the $E$-algebra of dual numbers. Since the image of $\psi$ is bounded in $E$ the image of the character
\begin{equation*}
\label{extstein1a}
\Theta_{\psi}: F^*\lra \wE^*,\,\, \Theta_{\psi}(x)=  1+ \psi(x)\vep.
\end{equation*}
is bounded in $\wE$. Therefore similar to Remark \ref{remarks:padicadm} (d) we obtain an admissible $E$-Banach space representation $\wE(\Theta_{\psi})$ of $T$ (even if $\psi(F^*)$ is not contained in $\cO$).  \enddemo
\end{remarks}

\section{Cohomology of $\sS_{G, \sK}$-spaces and schemes}
\label{section:sspaces}

\subsection{Notation and preliminary remarks}
\label{subsection:notpreparation}

Let $R$ be a ring, let $G$ be a locally profinite group and let $H$ be a closed subgroup of $G$. The  smooth induction functor
\begin{equation*}
\label{coind}
\Ind_H^G : \Mod_R^{\sm}(H) \lra \Mod_R^{\sm}(G)
\end{equation*}
is the right adjoint of the forgetful functor $\Res_H^G: \Mod_R^{\sm}(G) \to \Mod_R^{\sm}(H)$. Recall that $\Ind_H^G M$ for $M\in \Mod_R^{\sm}(H)$ consists of maps $\Phi: G\to M$ satisfying $\Phi(hg) = h\Phi(g)$ for all $h\in H$, $g\in G$ and such that there exists an open subgroup $K$ of $G$ with $\Phi(gk) = \Phi(g)$ for all $g\in G$, $k\in K$. The $G$-action is induced by right multiplication. Also if $H=K$ is an open subgroup of $G$ then $\Res_K^G$ has a left adjoint
\begin{equation*}
\label{indcoind2}
\cInd_K^G: \Mod_R^{\sm}(K) \lra \Mod_R^{\sm}(G)
\end{equation*}
For $M\in \Mod_R^{\sm}(K)$, $\cInd_K^G M$ is the $R[G]$-submodule of $\Ind_K^G N$ consisting of maps $\Phi\in \Ind_K^G N$ that have compact support (modulo $H$). 

Let $\cC$ be a site, i.e.\ $\cC$ is a category equipped with a (Grothendieck) topology (see \cite{artin}, Def.\ 1.2). We denote by $\PSh(\cC)$ the category of presheaves (of sets) on $\cC$ and by $\Sh(\cC)$ the category of sheaves on $\cC$. For a ring $R$ we let $\PSh(\cC, R)$ resp.\ $\Sh(\cC, R)$ denote the category of presheaves (resp.\ sheaves) of $R$-modules on $\cC$. By $\PSh(\cC, R)\to \Sh(\cC, R), \sF\mapsto \sF^{\sharp}$ we denote the sheafification functor. An object of $\Sh(\cC, R)$ will be called an $R$-sheaf for short. Note that if $\cO$ denotes the sheafification of the constant presheaf of rings $U\in \cC\mapsto R$ then the category $\Sh(\cC, R)$ can be identified with the category of sheaves of $\cO$-modules on $\cC$ (see \cite{stacksproject}, Tag 03CY). 

Recall that an ($R$-)topos is the category of $R$-sheaves on a site. By a morphism of ($R$-)topoi 
\[
(f^*, f_*): \Sh(\cC, R)\lra \Sh(\cD, R)
\]
(where $\cC$ and $\cD$ are sites) we mean a pair of adjoint $R$-linear functors $f^*:\Sh(\cD, R)\lra \Sh(\cC, R)$ and $f_*:\Sh(\cC, R)\lra \Sh(\cD, R)$. Note that we do not require $f^*$ to be exact. 
For an object $X\in \cC$ and $\cF\in \Sh(\cC, R)$ we denote by $H^{\bu}(X, \cF)$ the cohomology groups of $\cF$ over $X$. Also if $\cC/X$ denotes the localization of $\cC$ at $X$ (see \cite{stacksproject}, Tag 00Y0) and if $\cF\in \Sh(\cC/X, R)$ then we define $H^{\bu}(X, \cF)$ as the cohomology of $\cF$ over the object $X\stackrel{1_X}{\lra} X$ of $\cC/X$. If $\cF= \sF|_X$ is the restriction of an $R$-sheaf
on $\sF$ on $\cC$ then we have $H^{\bu}(X, \cF)= H^{\bu}(X, \sF)$ by (\cite{stacksproject}, Tag 03F3).

Let $u: \cC \to \cD$ be a continuous functor giving rise to a morphism of topoi $(u^*, u_*):= (u_s, u^s): \Sh(\cD, R) \to \Sh(\cC, R)$ (i.e.\ we assume that $u_s$ is exact). Let $X\in \cC$ and let $Y= u(X)$. The functor $u':\cC/X \to \cD/Y, (\varphi: U\to X)\mapsto (u(\varphi): u(U) \to Y)$ induces a morphism of topoi as well (see \cite{stacksproject}, Tag 03CF). By abuse of notation it will be denoted by $(u^*, u_*): \Sh(\cD/Y, R) \to \Sh(\cC/X, R)$ too. For $\cF\in  \Sh(\cC/X, R)$ there exists a canonical homomorphism 
\begin{equation*}
\label{basechangecohomcl}
H^n(X, \cF) \lra H^n(Y, u^*(\cF))
\end{equation*}
for every $n\ge 0$ defined as the the composite of the morphism $H^n(X, \cF)\to H^n(X, u_*u^*(\cF))$ (induced by the unit of adjunction) with edge morphism $H^n(X, u_*u^*(\cF))\to H^n(Y, u^*(\cF))$ of the Leray spectral sequence.

Let $f: X\to Y$ be a morphism in a site $\cC$. Then the site $\cC/X$ can be identified with the localization of the site $\cC/Y$ at $(X\stackrel{f}{\lra} Y)\in \cC/Y$ (\cite{stacksproject}, Tag 04BB). So there exists a morphism of topoi 
\begin{equation}
\label{topos1}
(f^*, f_*): \Sh(\cC/X, R) \to \Sh(\cC/Y, R).
\end{equation}
Here $f^* = \widetilde{f}^s$ where $\widetilde{f}$ denotes the functor $\cC/X\to \cC/Y, (U\to X) \mapsto (U\to X\stackrel{f}{\lra} Y)$. By (\cite{stacksproject}, Tag 03DI) the functor $f^*$ is exact as it has a left adjoint as well. 
If $\cC$ has fibre products then the base change functor with respect to $f:X\to Y$
\begin{equation*}
\label{basechange}
v: \cC/Y\lra \cC/X, \quad (V\to Y) \mapsto (V\times_Y X\to X).
\end{equation*}
is continuous and right adjoint to $\widetilde{f}$. By (\cite{stacksproject}, Tag 00XY) this implies $f_* = v^s$, i.e.\ for $\cF\in \Sh(\cC/X, R)$ we have 
\begin{equation}
\label{basechange2}
f_*(\cF)(V\to Y) \, =\, \cF(V\times_Y X\to X)
\end{equation}
for all $(V\to Y)\in \cC/Y$. Moreover we remark that given a commutative diagram in $\cC$ 
\begin{equation*}
\label{lss2}
\begin{CD} 
X_1 @> f_1 >> Y_1\\
@VV g V@VV h V\\
X_2 @> f_2 >> Y_2
\end{CD}
\end{equation*}
and $\cF\in \Sh(X_2, R)$ there exists a canonical morphism of Leray spectral sequences (see \cite{stacksproject}, Tag 0E46)
\begin{eqnarray*}
\label{lss3}
&& \left( \, E_2^{rs} = H^r(Y_2, R^s (f_2)_* \cF) \, \Longrightarrow\, H^{r+s}(X_2, \cF) \, \right) \, \lra \, \\
&& \hspace{3cm} \left(\, E_2^{rs} = H^r(Y_1, R^s (f_1)_* g^* \cF) \, \Longrightarrow\, H^{r+s}(X_1, g^*\cF)\, \right).
\nonumber
\end{eqnarray*}

\subsection{Sites of discrete $G$-sets}
\label{subsection:gsets}

Let $G$ be a locally profinite group and let $R$ be a ring. We denote by $\sS_G$ the category of discrete left $G$-sets and $G$-equivariant maps. Recall that a left $G$-set $S$ is discrete if the stabilizer $\Stab_G(s)$ of any $s\in S$ is open in $G$. Note that the category $\sS_G$ admits fibre products. We equip $\sS_G$ with the canonical topology (\cite{artin}, 1.1.2), i.e.\ we consider $\sS_G$ as a site. Recall (\cite{artin}, 1.1.5) that a family of morphisms $\left\{ S_i \stackrel{\rho_i}{\lra} S\right\}$ in $\sS_G$ is a covering if $\bigcup_i \rho_i(S_i) = S$. A set $S\in \sS_G$ will be called {\it connected} if $G$ acts transitively on $S$, i.e.\ if $S\cong G/K$ for some open subgroup $K$ of $G$. Note that if $K_1$, $K_2$ are open subgroups of $G$ then every $G$-equivariant map $\rho: G/K_1\to G/K_2$ is necessarily of the form $\rho(hK_1) = hgK_2$ for some element $g\in G$ with $K_1\subseteq K_2^g$.

The category of sheaves of $R$-modules on $\sS_G$ is equivalent to the category of discrete $R[G]$-modules and we will sometimes identify the two. Indeed, the functors 
\begin{eqnarray*}
\label{sheafgmod1}
&& \Mod_R^{\sm}(G) \lra \Sh(\sS_G, R),\quad M\mapsto \Maps_G(\wcdot, M)\\
&& \Sh(\sS_G, R)\lra \Mod_R^{\sm}(G),\quad\cM \mapsto \dlim_{U\le G,\, U\open} \cM(G/U)
\label{sheafgmod2}
\end{eqnarray*}
are mutually quasi-inverse equivalences of categories. 

Let $\varphi: G \to G'$ be a continuous homomorphism of locally profinite groups. Then any discrete left $G'$-set $S$ becomes a left $G$-set $S_{\varphi}$ via $\varphi$, so we obtain a continuous functor $\wvphi: \sS_{G'}\to \sS_{G}, S\mapsto S_{\varphi}$ between sites. It induces a morphism of $R$-topoi
\begin{equation*}
\label{grouphomtopos}
(\varphi^*, \varphi_*): \Sh(\sS_{G}, R) \lra \Sh(\sS_{G'}, R)
\end{equation*}
where $\varphi^*=\wvphi_s$ and $\varphi_*=\wvphi^s$.

Under the identifications $\Mod_R^{\sm}(G)=\Sh(\sS_G, R)$ and $\Mod_R^{\sm}(G')=\Sh(\sS_{G'}, R)$ the functors $\varphi_*$ and $\varphi^*$ are given as follows. For $M\in \Mod_R^{\sm}(G)$ we have 
\begin{equation*}
\label{sheafgmodfunc2}
\varphi_*(M) = \left\{\Phi\in \Maps(G', M)| \exists U'\le G' \open \forall  g\in G, g'\in G', u'\in U' : \Phi(\varphi(g)g' u') =g \Phi(g')\,\right\}
\end{equation*}
with $G'$-action induced by right multiplication. For $M'\in \Mod_R^{\sm}(G')$ the pull-back $\varphi^*(M')$ is the $G$-module $M'_{\varphi}$. 

If $H$ is a closed subgroup of $G$ and $\varphi=\iota: H\hra G$ is the inclusion then $\iota_*= \Ind_H^G: \Mod_R^{\sm}(H)\to \Mod_R^{\sm}(G)$ is the (smooth) induction functor. If $H$ is a closed normal subgroup of $G$ and $\varphi=\pi: G\to G/H$ is the projection then $\pi^*$ is the inflation functor $\Infl_{G/H}^G: \Mod_R^{\sm}(G/H)\to \Mod_R^{\sm}(G)$ and $\pi_*$ is given by $\pi_*(M) =M^H$ for $M\in \Mod_R^{\sm}(G)$.

Now we assume that $H=K$ is an open subgroup of $G$ and $\varphi=\iota: H\hra G$ is again the inclusion. In this case the functor $\wio: \sS_{G}\to \sS_{K}, S\mapsto S_{\iota}$ has a left adjoint defined as follows. Given $S'\in \sS_K$ we define a right $K$-action on $G\times S'$ by $k\cdot (g, s') = (gk, k^{-1} s')$ and a left $G$-action by $h \cdot (g, s') = (hg, s')$. The latter induces a discrete left $G$-action on $G\times_K S' : =  K\backslash (G\times S')$. If for example $S' = K/K'$ for some open subgroup $K'$ of $K$ then $G\times_K S'\cong G/K'$. 
Thus we obtain a functor 
\begin{equation}
\label{fibergk}
u: \sS_K \lra \sS_{G}, \qquad S' \mapsto G\times_K S'.
\end{equation}
that is easily seen to be cocontinuous. By (\cite{stacksproject}, Tag 00XY) we get $\io^* =\wio_s= u^s$. 

Now we consider certain full subcategories of $\sS_G$ associated to a cofinal system of open subgroups of $G$. More precisely we consider a subset $\sK$ of the set of open subgroups of the locally profinite group $G$ satisfying the following

\begin{assum}
\label{assum:kgood}
\rm (i) $\sK$ is cofinal in the set of all open subgroups of $G$, i.e.\ for every open subgroup $U$ of $G$ there exists a $K\in \sK$ with $K\subseteq U$.

\noi (ii) $\sK$ is closed under conjugation, i.e.\ for every $K\in \sK$, $g\in G$ we have $K^g\in \sK$.

\noi (iii) For $K\in \sK$ and an open subgroup $K'$ of $K$ we have $K'\in \sK$.
\end{assum}

We denote by $\sS_{G, \sK}$ the full subcategory of $S\in \sS_G$ such that $\Stab_G(s)\in \sK$. Condition (ii) guarantees that for $K\in \sK$ the left $G$-set $G/K$ lies in $\sS_{G, \sK}$. Condition (iii) implies that if $S'\to S$ is a morphism in $\sS_G$ then $S\in \sS_{G, \sK}$ implies $S'\in \sS_{G, \sK}$ as well. It follows that $\sS_{G, \sK}$ has fibre products and equalizers. 

We equip $\sS_{G, \sK}$ with the canonical topology. Condition (i) implies that any $S\in \sS_G$ admits a covering by objects in $\sS_{G, \sK}$. Therefore the inclusion $u: \sS_{G, \sK} \hra \sS_G$ induces mutually quasi-inverse equivalences of categories $u^s: \Sh(\sS_{G, \sK}, R)\to \Sh(\sS_G, R)$ and $_s u: \Sh(\sS_{G, \sK}, R)\to \Sh(\sS_G, R)$. Recall (\cite{stacksproject}, Tag 00XH) that $(_s u)(\cF)$ is given by
\begin{equation*}
\label{extendsheaf}
(_s u)(\cF)(S) \, =\, \prolim_{(S', \rho)} \cF(S')
\end{equation*}
where the limit is taken over the category of pairs $(S', \rho)$ where $S'\in \sS_{G, \sK}$ and $\rho:S'\to S$ is a morphism in $\sS_G$.

Hence we can identify $\Sh(\sS_{G, \sK}, R)$ also with the category of discrete $R[G]$-modules via the equivalences of categories
\begin{eqnarray}
\label{sheafgmod3}
&& \Mod_R^{\sm}(G) \lra \Sh(\sS_{G, \sK}, R),\quad M \mapsto\Maps_G(\wcdot, M)\nonumber\\
&& \Sh(\sS_{G, \sK}, R)\lra \Mod_R^{\sm}(G),\quad \cM \mapsto \dlim_{K\in \sK^{\opp}} \cM(G/K).
\label{sheafgmod4}
\end{eqnarray}

For $S\in \sS_{G, \sK}$ consider the localization morphism 
\begin{equation*}
\label{sitegsetloc}
j_S: \sS_{G, \sK}/S \lra \sS_{G, \sK}, \quad (\rho: S'\to S) \mapsto S'.
\end{equation*}
It induces a morphism of topoi
\begin{equation*}
\label{sheafgsetloc}
(j_S^*, (j_S)_*): \Sh(\sS_{G, \sK}/S, R)\lra \Sh(\sS_{G, \sK}, R) .
\end{equation*}
The functor $j_S^*:= j_S^s$ is exact. Indeed, by (\cite{stacksproject}, Tag 03DI) it has also a left adjoint $(j_S)_!$. 

Now assume that $S=G/K$ where $K\in \sK$. In this case the localization $\sS_{G, \sK}/S$ can be identified with the site $\sS_K$ of discrete left $K$-sets. Indeed, the functors
\begin{eqnarray}
\label{sitegsetlock1}
&& \sS_{G, \sK}/S \lra \sS_K, \quad (\rho: S'\to G/K) \mapsto \rho^{-1}(1\cdot K)\nonumber\\
&& \sS_K \lra \sS_{G, \sK}/S, \quad S' \mapsto (G\times_K S' \to S)
\label{sitegsetlock2}
\end{eqnarray}
are mutually quasi-inverse equivalences of sites. Here for $S'\in \sS_K$ the morphism 
$G\times_K S' \to S$ is induced by the unique map $S'\to \pt=K/K$. It follows that 
the category $\Sh(\sS/S, R)$ can be identified with the category $\Sh(\sS_K, R)$ (and hence 
with $\Mod_R^{\sm}(K)$.

\begin{lemma}
\label{lemma:reslocal} 
Let $K\in \sK$, put $S=G/K$ and let $\io: K\to G$ be the inclusion. Under the identifications $\Sh(\sS_K, R)=\Sh(\sS/S, R)$ we have an equality of morphisms of topoi
\begin{equation*}
\label{reslocal1}
(j_S^*, (j_S)_*) \, =\, (\io^*, \io_*): \Sh(\sS/S, R)=\Sh(\sS_K, R)\lra \Sh(\sS, R).
\end{equation*}
\end{lemma}

\begin{proof} Since the composition of \eqref{sitegsetlock2} with $j_S$ is the functor \eqref{fibergk} we get $j_S^*= j_S^s = u^s = \io^*$. The equality $\io_* = (j_S)_*$ follows from the adjointness property. 
\end{proof}

\begin{remarks}
\label{remarks:locresexact}
\rm (a) Let $K\in \sK$ and put $S= G/K$. By identifying the category $\Sh(\sS_{G, \sK}/S, R)$ with $\Mod_R^{\sm}(K)$ and $\Sh(\sS_{G, \sK}, R)$ with $\Mod_R^{\sm}(G)$ we see that the functors $(j_S)_!$, ${j_S}^*$ and $(j_S)_*$ and correspond to compact induction, the restriction and the (smooth) induction functor  
\begin{eqnarray*}
\label{sheafgsetloca}
&& \cInd_K^G: \Mod_R^{\sm}(K) \lra\Mod_R^{\sm}(G), \\
&& \Res_K^G: \Mod_R^{\sm}(G) \lra \Mod_R^{\sm}(K), \\
&& \Ind_K^G: \Mod_R^{\sm}(K) \lra \Mod_R^{\sm}(G)
\end{eqnarray*}
respectively. 

\noi (b) Let $\varphi: G \to G'$ be a continuous homomorphism of locally profinite groups and let
$\sK'$ be a subset of the set of open subgroups of $G'$ that satisfies the hypotheses \ref{assum:kgood}. The fact that the topoi $\Sh(\sS_{G, \sK}, R)$ and $\Sh(\sS_{G', \sK'}, R)$ can be identified with $\Sh(\sS_G, R)$ and $\Sh(\sS_{G'}, R)$ respectively, implies that 
$\varphi$ induces a morphism of topoi
\begin{equation}
\label{grouphomtopos2}
(\varphi^*, \varphi_*): \Sh(\sS_{G, \sK}, R) \lra \Sh(\sS_{G', \sK'}, R).
\end{equation}
In general the morphism \eqref{grouphomtopos2} is not induced by a continuous functor $\sS_{G', \sK'}\to \sS_{G, \sK}$. The functor $\varphi_*$ has the following concrete description. For $S'\in \Sh(\sS_{G', \sK'}, R)$ let $_{S'} \cJ$ denote the category of pairs $(S, \xi)$ where $S\in \sS_{G, \sK}$ and $\xi: S\to S'_{\varphi}$ is a $G$-equivariant map. A morphism between two objects $(S_1, \xi_1)$, $(S_2, \xi_2)$ in $_{S'} \cJ$ is a $G$-equivariant map $\rho: S_1\to S_2$ such that $\xi_2\circ \rho = \xi_1$. For $\cF\in \Sh(\sS_{G, \sK}, R)$ we define a functor 
\begin{equation*}
\label{grouphomtopos3}
_{S'} \cF: {_{S'}\cJ}^{\opp} \lra \Mod_R, \qquad (S, \xi)\mapsto \cF(S).
\end{equation*}
We then have 
\begin{equation*}
\label{grouphomtopos3a}
(\varphi_*)(\cF)(S') \, =\, \prolim_{_{S'} \cJ^{\opp}}  {_{S'} \cF}.
\end{equation*} \enddemo
\end{remarks}

\subsection{The topos associated to an $\sS$-object in a site}
\label{subsection:ssschemes}

Let $G$ be a locally profinite group. As in the last section we fix a subset $\sK$ of the set of open subgroups of $G$ satisfying the conditions (i)--(iii) of assumption \ref{assum:kgood} and we denote by $\sS=\sS_{G,\sK}$ the associated site of discrete left $G$-sets whose stabilizers contained in $\sK$. Let $\cC$ be a site. Associated to an object $U$ of $\cC$ are the two functors (co- and contravariant respectively) $h^U= \Hom_{\cC}(U, \wcdot): \cC\to \sets$ and $h_U= \Hom_{\cC}(\wcdot, U): \cC\to \sets$. We make the following 

\begin{assum}
\label{assum:sitegood}
$\cC$ has fibre products, equalizers and coproducts. 
\end{assum}

The first two conditions imply that arbitrary connected finite limits exists in $\cC$ (\cite{stacksproject}, Tag 04AT). The last condition implies in particular that $\cC$ has an initial object. 

\begin{df}
\label{df:Sspace}
(a) A continuous functor $X: \sS\to \cC$ that commutes with coproducts will be called an $\sS$-object in $\cC$.\footnote{For 
the definition of a continuous functor between sites see (\cite{stacksproject}, Tag 00WV).}
\medskip

\noi (b) An $\sS$-object $X$ in $\cC$ is called (left) exact if it commutes with equalizers. 
\end{df}

Let $X$ be an $\sS$-object in $\cC$. For $S\in \sS$ we will write $X_S$ rather than $X(S)$. Also for $K\in \sK$ we write $X_K$ instead of $X_{G/K}$. If $\rho:S_1\to S_2$ is a morphism in $\sS$ then the induced morphism will be denoted simply by $\rho: X_{S_1}\to X_{S_2}$. Since every morphism in $\sS$ is part of a covering the functor $X$ commutes with fibre products. So by (\cite{stacksproject}, Tag 04AT), $X$ is exact if and only if $X$ commutes with connected finite limits. If $G\in \sK$ then $\sS=\sS_{G, \sK} = \sS_G$ has a final object $\pt=G/G$. In this case every $\sS$-object in $\cC$ is exact.

\begin{remark}
\label{remark:Sspacealt}
\rm Let $\sS^c$ denote the full subcategory $\sS^c$ of $\sS$ with set of objects $\{G/K\mid K\in \sK\}$, let
$X$ be an $\sS$-object in $\cC$ and let 
\begin{equation}
\label{ssclassical}
\sS^c\to \cC, \quad G/K\mapsto X_K , \quad (\rho: G/K\to G/K)\mapsto (\rho: X_{K}\to X_{L})
\end{equation}
be the restriction of $X$ to $\sS^c$. The fact that $X$ commutes with coproducts implies that the latter completely determines $X$. The functor \eqref{ssclassical} has the following properties
\medskip

\noi (i) $\rho: X_K\to X_L$ is a covering for every morphism $\rho: G/K\to G/L$ in $\sS^c$.
\medskip

\noi (ii) Let $G/K_1\to G/L$ and $G/K_2\to G/L$ be morphism in $\sS^c$ and let
$G/K_1\times_{G/L} G/K_2= \coprod_{i\in I} G/M_i$ be their fibre product in $\sS$. Then the induced diagram 
\[
\begin{CD}
\coprod_{i\in I} X_{M_i} @>>> X_{K_1}\\
@VVV @VVV \\
X_{K_2} @>>> X_L
\end{CD}
\]
is cartesian.
\medskip

If $X$ is exact then we also have 
\medskip

\noi (iii) For two different morphisms $\rho_1, \rho_2: G/K\to G/L$ in $\sS^c$ the equalizer of 
$\rho_1, \rho_2: X_K\to X_L$ is the initial object in $\cC$. 
\medskip

Conversely, a functor \eqref{ssclassical} with the properties (i) and (ii) extends to an $\sS$-object $X$ in $\cC$ (i.e.\ the restriction of $X$ to $\sS^c$ is \eqref{ssclassical}). Moreover if additionally (iii) holds then $X$ is exact.\enddemo
\end{remark}

Let $X$ be an $\sS$-object in $\cC$. For $U\in \cC$ we consider the colimit of the functor $h^U\circ X:\sS\to \sets$, i.e.\ the set 
\begin{equation*}
\label{sspace2}
\colim h^U\circ X= \colim_{S\in \sS} \Hom_{\cC}(U, X_S).
\end{equation*}
A morphism $\phi: U_1 \to U_2$ in $\cC$ induces a morphism of functors $h^{U_2}\to h^{U_1}$ hence a map 
\begin{equation*}
\label{sspace2a}
\phi^*: \colim h^{U_2}\circ X \to \colim h^{U_1}\circ X.
\end{equation*}

\begin{df}
\label{df:CuebercolimX} 
Let $X$ be an $\sS$-object in $\cC$. 
\medskip

\noi (a) We introduce the site $\cC/\colim X$. 
\medskip

\noi (i) Objects of $\cC/\colim X$ are pairs $(U, \alpha)$ where $U\in \cC$ and $\alpha\in \colim h^U\circ X$. 
\medskip

\noi (ii) A morphisms in $\cC/\colim X$
\[
\phi: (U_1, \alpha_1) \lra (U_2, \alpha_2)
\] 
is a morphism $\phi: U_1 \to U_2$ in $\cC$ such that $\phi^*(\alpha_2) = \alpha_1$. 
\medskip

\noi (iii) A family of morphisms $\left\{ \phi_i: (U_i, \alpha_i)\to (U, \alpha)\right\}$ in $\cC/\colim X$ is a covering if $\left\{ \phi_i: U_i\to U\right\}$ is a covering in $\cC$. 
\medskip

\noi (b) The category $\Sh(X, R)$ of $R$-sheaves $\sF$ on $X$ is defined as the category of $R$-sheaves of the site $\cC/\colim X$. 
\end{df}

\begin{remarks}
\label{remarks:colimexplizit} 
\rm (a) Note that the category $\cC/\colim X$ also satisfies assumption \ref{assum:sitegood}. 
\medskip

\noi (b) For every object $(U, \alpha)$ of $\cC/\colim X$ the obvious functor between localizations 
\begin{equation}
\label{sspace3}
\left(\cC/\colim X\right)\!/(U, \alpha)\lra \cC/U, \,\,\, \left(\phi: (U', \alpha')\to (U, \alpha)\right) \mapsto (\phi: U'\to U)
\end{equation}
is an isomorphism of sites. 
\medskip

\noi (c) If $G\in \sK$ then we can identify $\colim h^U\circ X$ with $\Hom_{\cC}(U, X_G)$, the site $\cC/\colim X$ with the localization $\cC/X_G$ and the category $\Sh(X, R)$ with $\Sh(\cC/X_G, R)$. 
\medskip

\noi (d) Let $U\in \cC$ and let $\alpha\in \colim h^U\circ X$. Let $\cI = \cI_{(U,\alpha)}$ be the category of pairs $(S, \psi)$ where $S\in \sS$ and $\psi: U \to X_S$ is a morphism representing $\alpha$, i.e.\ the image of $\psi$ under the canonical map $\Hom(U, X_S) \to \colim_{S'\in \sS} \Hom_{\cC}(U, X_{S'})$ is $\alpha$. A morphism $\tau: (S_1, \psi_1)\to (S_2, \psi_2)$ in $\cI$ is a morphism $\tau: S_1\to S_2$ in $\sS$ such that 
\begin{equation*}
\label{colimexplizit1}
{\xymatrix@-0.5pc{& \ar[dl]_{\psi_1} U \ar[dr]^{\psi_2}\\
X_{S_1'}\ar[rr]^{\tau} && X_{S_2'}
}}
\end{equation*}
commutes. Objects of $\cI_{(U,\alpha)}$ are called representatives of $(U,\alpha)$.

Note that the category $\cI_{(U,\alpha)}$ has fibre products and is connected, i.e.\ any two objects $(S, \psi)$, $(S', \psi')$ can be linked 
by a chain of morphisms
\[
(S, \psi) =(S_0, \psi_0)\stackrel{\tau_1}{\lra}(S_1, \psi_1) \stackrel{\tau_2}{\lla}(S_2, \psi_2) \stackrel{\tau_3}{\lra} \ldots \stackrel{\tau_{n-1}}{\lla}(S_n, \psi_n)=(S', \psi').
\]
By replacing every pair of morphisms $(S_{i-1}, \psi_{i-1}) \stackrel{\tau_{i-1}}{\lra}(S_i, \psi_i) \stackrel{\tau_i}{\lla}(S_{i+1}, \psi_{i+1})$ in the chain successively by its fibre product we see that every two objects $(S, \psi)$, $(S', \psi')$ can be linked by a diagram of the form 
\begin{equation*}
\label{colimexplizit2}
(S, \psi)\stackrel{\tau}{\lla}(S'', \psi'') \stackrel{\tau'}{\lra}(S', \psi'). 
\end{equation*} 
Moreover, if $X$ is exact then the category $\cI_{(U,\alpha)}^{\opp}$ is filtered. For that we are left to show that given two morphism $\tau_1, \tau_2: (S, \psi)\to (S', \psi')$ in $\cI_{(U,\alpha)}$ there exists a morphism $\xi: (S'', \psi'') \to (S, \psi)$ with $\tau_1\circ \xi = \tau_2\circ \xi$. We define $\xi: S''\to S$ to be the equalizer of $\tau_1, \tau_2: S\to S'$. The fact that $\tau_1\circ \psi = \psi'=\tau_2\circ \psi$ and that $X$ commute with equalizers implies that there exists a morphism $\psi'': U\to X_{S''}$ with $\xi \circ \psi''=\psi$.
\enddemo
\end{remarks}
\medskip

Let $X$ be an $\sS$-object in $\cC$, let $\cD$ be another site that satisfies assumption \ref{assum:sitegood} and let $u: \cC \to \cD$ be a continuous functor commuting with fibre products, equalizers and coproducts. Then $u(X):= u\circ X$ is an $\sS$-object in $\cD$. 
For $U\in \cC$ the collection of maps $\Hom_{\cC}(U, X_S) \to \Hom_{\cD}(u(U), u(X)_S), \psi\mapsto u(\psi)$ for $S\in \sS$ defines a morphism of functors $h^U \circ X \to h^{u(U)} \circ v(X)$ hence a map
\begin{equation*}
\label{changesite}
\colim h^U\circ X\lra \colim h^{u(U)}\circ u(X), \,\,\, \alpha\mapsto u(\alpha)
\end{equation*}
Thus $u$ induces a continuous functor between sites 
\begin{equation}
\label{changesite2}
\wu: \cC/\colim X\lra \cD/\colim u(X),\quad (U, \alpha) \mapsto (u(U), u(\alpha))
\end{equation}
that commutes with fibre products, equalizers and coproducts. It induces a morphism of topoi 
\begin{equation}
\label{changesite3}
(u^*, u_*):=(\wu_s, \wu^s): \Sh(u(X), R) \lra \Sh(X, R).
\end{equation}

\begin{prop}
\label{prop:pullbackexact}
The functor $u^*:=\wu_s: \Sh(X, R) \to \Sh(u(X), R)$ is exact. 
\end{prop}

\begin{proof} For an object $(V, \beta)\in \cD/\colim u(X)$ let $\cI=\cI_{(V, \beta)}^u$ be the category of triples $(U, \alpha, \zeta)$ where $(U, \alpha)\in \cC/\colim X$ and $\zeta: V \to u(U)$ is a morphism in $\cD$ such that $\zeta^*(u(\alpha)) = \beta$. A morphism $\eta: (U_1, \alpha_1, \zeta_1)\to (U_2, \alpha_2, \zeta_2)$ in $\cI_{(V, \beta)}^u$ is a morphism $\eta: (U_1, \alpha_1)\to (U_2, \alpha_2)$ in $\cC/\colim X$ such that $u(\eta)\circ \zeta_1=\zeta_2$. By (\cite{stacksproject}, Tag 00X5) it suffices to show that the category $\cI^{\opp}$ is filtered. 

Since $\wu$ commutes with fibre products and equalizers, the category $\cI$ has the following properties (see \cite{stacksproject}, Tag 00X4)
\medskip

\noi (i) For every pair of morphism $\eta_1: j_1 \to i$, $\eta_2: j_2\to i$ in $\cI$ there exists an object $k\in \cI$ and morphisms $\xi_1: k\to  j_1$, $\xi_2: k\to j_2$ in $\cI$ such that $\eta_1\circ \xi_1=\eta_2\circ \xi_2$.
\medskip

\noi (ii) For every $i, j \in \cI$ and morphisms $\eta_1: j\to i$, $\eta_2: j \to i$ in $\cI$ there exists $k\in \cI$ and a morphism $\xi: k\to j$ such that $\eta_1 \circ \xi= \eta_2 \circ \xi: k \to i$. 
\medskip

It remains to show 
\medskip

\noi (iii) For every $i_1, i_2\in \cI$ there exists $i\in \cI$ and morphisms $i\to i_1$, $i\to i_2$ in $\cI$.
\medskip

For that let $i_1=(U_1, \alpha_1, \zeta_1)$, $i_2=(U_2, \alpha_2, \zeta_2)$ and let $(S_1, \varphi_1)$ and $(S_2, \varphi_2)$ be representative of $(U_1, \alpha_1)$ and $(U_2, \alpha_2)$ respectively. Then $(S_1, u(\varphi_1)\circ \zeta_1)$, $(S_2, u(\varphi_2)\circ \zeta_2)$ are representatives of $(V, \beta)$. By Remark \ref{remarks:colimexplizit} (d) there exists a representative $(S, \varphi)$ of 
$(V, \beta)$ and morphisms $\tau_1: (S, \varphi)\to (S_1, u(\varphi_1)\circ \zeta_1)$ and $\tau_2: (S, \varphi)\to (S_2, u(\varphi_2)\circ \zeta_2)$.
Note that $j:=(X_S, \alpha_S, \varphi)$, $j_1:=(X_{S_1}, \alpha_{S_1}, u(\varphi_1)\circ \zeta_1)$, $j_2:=(X_{S_2}, \alpha_{S_2}, u(\varphi_2)\circ \zeta_2)$, are objects of $\cI$. So we have a diagram in $\cI$
\begin{equation*}
\label{resjshrieka}
{\xymatrix@-0.1pc{i_1 \ar[dr]^{\varphi_1} && \ar[dl]_{\tau_1} j \ar[dr]^{\tau_2} && \ar[dl]_{\varphi_2} i_2\\
& j_1 && j_2 
}}
\end{equation*}
Thus applying property (i) three times yields (iii). 
\end{proof}

The continuous functor $X:\sS\to \cC$ can be ``lifted'' to a functor $\sS\to \cC/\colim X$. For that let $S\in \sS$ and let $\alpha_S\in \colim h^{X_S}\circ X$ be the image of $1_{X_S}\in \Hom_{\cC}(X_S, X_S)$ under the canonical map 
\begin{equation*}
\label{sspace4}
\Hom_{\cC}(X_S, X_S)\lra  \colim_{S'\in \sS} \Hom_{\cC}(X_S, X_{S'})=\colim h^{X_S}\circ X.
\end{equation*}
The assignment 
\begin{equation}
\label{sspace5}
\wX: \sS\lra \cC/\colim X, \quad S\mapsto (X_S, \alpha_S).
\end{equation}
is a continuous functor between sites, lifting the functor $X: \sS\to \cC$ (i.e.\ the composition of $\wsX$ with the forgetful functor $\cC/\colim X\to \cC, (U, \alpha)\mapsto U$ is the functor $X$). If $G\in \sK$ then $\wX$ can be identified with the functor 
\begin{equation}
\label{sspace5b}
\sS\lra \cC/X_G, \quad S\mapsto (X_S\to X_G)
\end{equation}
(the morphism $X_S\to X_G$ is the image under $X$ of the unique map $S\to \pt$). 

The functor \eqref{sspace5} induces a morphism of topoi 
\begin{equation}
\label{sspace5a}
(\wX_s, \wX^s):\Sh(X, R) \lra \Sh(\sS, R).
\end{equation}
As in the previous section we can identify the categories $\Sh(\sS, R)$ and $\Mod_R^{\sm}(G)$. Under this identification we denote the first component of \eqref{sspace5b} by 
\begin{equation*}
\label{gmodsheafx}
\Mod_R^{\sm}(G) \lra \Sh(X, R), \,\, M\mapsto M_{X, G}
\end{equation*}
If $M\in \Mod_R^{\sm}(G)$ is equipped with the trivial $G$-action then we will write $M_X$ instead of $M_{X, G}$. 
The second component of  \eqref{sspace5b} will be denoted by and the second by
\begin{equation}
\label{sheafxgmod}
H^0(X_{\infty}, \wcdot): \Sh(X, R) \lra \Mod_R^{\sm}(G), \,\, \sF\mapsto H^0(X_{\infty}, \sF).
\end{equation}
So there exists a natural isomorphisms
\begin{equation*}
\label{sheafxgmod2}
\Hom_{R[G]}(M, H^0(X_{\infty}, \sF))\, \cong\, \Hom_{\Sh(X, R)}(M_{X,G}, \sF)
\end{equation*}
for every $\sF\in \Sh(X, R)$ and every $M\in \Mod_R^{\sm}(G)$. 

\begin{remark}
\label{remark:sheavestogmod}
\rm For $\sF\in \Sh(X, R)$ we have
\begin{equation*}
\label{sheafxgmod3}
H^0(X_{\infty}, \sF)\, = \, \dlim_{K\in \sK^{\opp}} \sF_K(X_K).
\end{equation*}
This follows immediately from the definition of $\wX^s$ and of the isomorphism \eqref{sheafgmod4}. Moreover for $K\in \sK$ the $K$-invariants of $H^0(X_{\infty}, \sF)$ can be identified with $\sF_K(X_K)$. \enddemo
\end{remark}

We now assume that $X$ is exact. Then \eqref{sspace5} can be viewed as a special case of the functor \eqref{changesite2}. Namely for $u:=X: \sS\to \cC$ and $X=\id_{\sS}$ the identical $\sS$-object we obtain a continuous functor
\begin{equation*}
\label{sspace5c}
\wu: \sS/\colim \id_{\sS} \lra \cC/\colim X.
\end{equation*}
Since the obvious functor $\sS\to \sS/\colim \id_{\sS}$ is an equivalence of sites the functor $\wX$ can be identified with $\wu$. So by Prop.\ \ref{prop:pullbackexact} we obtain

\begin{corollary}
\label{corollary:gmodsheafxexact}
Let $X$ be an exact $\sS$-object in $\cC$.Then the functor $\Mod_R^{\sm}(G) \to \Sh(X, R), \,M\mapsto M_{X, G}$ is exact. 
\end{corollary}

By \eqref{sspace3}, for $S\in \sS$ the localization of $\cC/\colim X$ at $(X_S, \alpha_S)$ can be identified with the localization of $\cC$ at $X_S$.  Therefore the (continuous and cocontinuous) functor 
\begin{equation*}
\label{sspace6}
j_{X_S}: \cC/X_S \lra \cC/\colim X, \quad (\psi: U\to X_S) \mapsto (U, \psi^*(\alpha_S))
\end{equation*}
induces a morphism of topoi 
\begin{equation*}
\label{sspace7}
(j_{X_S}^*, (j_{X_S})_*): \Sh(\cC/X_S, R)\lra \Sh(X, R).
\end{equation*}
The functor $j_{X_S}^*$ will be denote by
\begin{equation*}
\label{sspace7a}
\res_S: \Sh(X, R)\lra \Sh(\cC/X_S, R), \quad \sF\mapsto \sF_S, \quad (u:\sF \to \sG)\mapsto (u_S: \sF_S\to \sG_S).
\end{equation*}
Also for $K\in \sK$ we put $\sF_K:=\sF_{G/K}\in \Sh(\cC/X_K, R)$. For an $R$-sheaf $\sF$ on $X$ the sheaves $\sF_S\in \Sh(\cC/X_S, R)$, $S\in \sS$ will be called the layers of $\sF$. 

The right adjoint $(j_{X_S})_*$ of $\res_S$ is called direct image functor. By (\cite{stacksproject}, Tags 03DI and 03DJ) $\res_S$ also has an exact left adjoint 
\begin{equation*}
\label{sspace8}
(j_{X_S})_! : = (j_{X_S})_s: \Sh(\cC/X_S, R)\lra \Sh(X, R). 
\end{equation*}
This implies part (a) of the following Lemma.

\begin{lemma}
\label{lemma:resexactfaith} 
Let  $X$ be an $\sS$-object in $\cC$ and let $S\in \sS$ be non-empty. 
\medskip

\noi (a) The functor $\res_S$ is exact and preserves injectives. 
\medskip

\noi (b) The functor $\res_S$ is faithful.
\end{lemma}

\begin{proof}[Proof of (b)] It suffices to show for each $\sF\in \Sh(X, R)$
\begin{equation}
\label{resfaithful}
\sF_S \,= \,0 \quad \Longrightarrow \quad \sF\,=\, 0.
\end{equation}
Since every object of $\cC/\colim X$ lies in the image the functor $j_{X_{S'}}$ for some $S'\in\sS$ it is enough to see that $\sF_{S}=0$ implies $\sF_{S'}= 0$ for every $S'\in \sS$. We first consider the case when there exists a morphism $\rho: S'\to S$ in $\sS$. Let $\widetilde{\rho}: \cC/X_{S'}\to \cC/X_S$ denote the continuous functor $(U\to X_{S'})\mapsto (U\to X_{S'} \stackrel{\rho}{\lra} X_S)$. Since the diagram 
\begin{equation*}
\label{resjshriekb}
{\xymatrix@-0.5pc{\cC/X_{S'}\ar[dd]^{\widetilde{\rho}} \ar[rrd]^{j_{X_{S'}}}\\
&& \cC/\colim X\\
\cC/X_S \ar[rru]^{j_{X_S}}
}}
\end{equation*}
commutes we have $\res_{S'} = \rho^*\circ \res_S$ hence $\sF_{S'} = \rho^*(\sF_S) =0$. 

Now assume that $S'$ is an arbitrary discrete left $G$-set in $\sS$. It is easy to see that there exists a diagram in $\sS$ 
\begin{equation*}
\label{scovering}
{\xymatrix@+0.1pc{ & \ar[dl] S'' \ar[rd]^{\rho} & \\
S && S'
}}
\end{equation*}
such that $\rho$ is a covering (i.e.\ a surjective map). We have already seen that $\sF_{S''} =0$ and we want to deduce $\sF_{S'}=0$. Let $U\to X_{S'}$ be an object of $\cC/X_{S'}$. We have  
\begin{eqnarray*}
&& \sF_{S'}(X_{S''}\times_{X_S'} U \to X_{S'}) \, =\, \sF(j_{X_{S'}}(X_{S''}\times_{X_S'} U\to X_{S'}))\\ 
&& =\, \sF(j_{X_{S''}}(X_{S''}\times_{X_S'} U\to X_{S''}))\, = \, \sF_{S''}(X_{S''}\times_{X_S'} U \to X_{S''})\, =\, 0
\end{eqnarray*}
hence also $\sF_{S'}(U\to X_{S'})=0$ since $X_{S''}\times_{X_S'} U \to U$ is (as a base change of $\rho: X_{S''} \to X_{S'}$) a covering 
in $\cC$. 
\end{proof}

A consequence of the previous Lemma is that exactness of a sequence in $\Sh(X, R)$ can be checked  
 at one layer. 

\begin{lemma}
\label{lemma:exactlayer} 
Let  $X$ be an $\sS$-object in $\cC$ and let $S\in \sS$ be non-empty. 
\medskip

\noi (a) Let
\begin{equation}
\label{exactlayers}
\sF'\stackrel{u}{\lra} \sF \stackrel{v}{\lra}\sF'
\end{equation}
be a sequence of morphism in $\Sh(X, R)$. If the sequence $\sF_S'\stackrel{u_S}{\lra} \sF_S \stackrel{v_S}{\lra}\sF_S''$ is exact then \eqref{exactlayers} is exact as well. 
\medskip

\noi (b) Let $\cA$ be an $R$-linear abelian category and let $F: \cA\to \Sh(X, R)$ be an $R$-linear functor. If the functor $\res_S\circ F: \cA\to \Sh(\cC/X_S, R)$ is exact (resp.\ $=0$) then $F$ is exact (resp.\ $=0$) as well. 
\end{lemma}

\begin{proof} (b) follows immediately from (a) (resp.\ Lemma \ref{lemma:resexactfaith} (b)). For (a) note that $(v\circ u)_S = v_S\circ u_S=0$, hence part (b) of Lemma \ref{lemma:resexactfaith} implies that $v\circ u=0$. Let $\sH$ denote the homology of \eqref{exactlayers}. We have $\sH_S = \ker(u_S)/\image(v_S)=0$ by Lemma \ref{lemma:resexactfaith} (a) hence $\sH =0$ by \eqref{resfaithful}.
\end{proof} 

\begin{remarks}
\label{remarks:Xloccomp}
\rm (a) Let $X$ be an $\sS$-object in $\cC$ and let $S\in \sS$ be non-empty. The continuous functor $\wX/S: \sS/S\to \cC/X_S$ induces a morphism of topoi
\begin{equation*}
\label{Xloccomp1}
((\wX/S)_s, (\wX/S)^s): \Sh(\cC/X_S, R)\lra \Sh(\sS/S, R).
\end{equation*}
By (\cite{stacksproject}, Tag 03CF) the diagram of topoi 
\begin{equation*}
\label{Xloccomp2}
\begin{CD} \Sh(\cC/X_S, R) @> (j_{X_S}^*, (j_{X_S})_*) >> \Sh(X, R)\\
@VV ((\wX/S)_s, (\wX/S)^s) V @VV (\wX_s, \wX^s) V \\
\Sh(\sS/S, R) @> (j_S^*, (j_S)_*) >> \Sh(\sS, R) 
\end{CD}
\end{equation*}
commutes. 
\medskip

\noi (b) Assume that $X$ is exact. An $R$-sheaf on $X$ can be defined in more concrete terms (i.e.\ without introducing the site $\cC/\colim X$). Namely, we could define an $R$-sheaf $\sF'$ on $X$ as 
\medskip

\noi $\bu$ A collection of sheaves $\sF'_S\in \Sh(\cC/X_S, R)$ for each $S\in \sS$ and
\medskip

\noi $\bu$ A collection of isomorphisms 
\begin{equation}
\label{sheaf1}
\varphi_{\rho}: \rho^*\sF'_{S_2}\lra \sF'_{S_1} 
\end{equation}
for every morphism $\rho: S_1\to S_2$ in $\sS$ that satisfy the cocycle condition 
\begin{equation*}
\label{sheaf2}
\varphi_{\rho_1}\circ \rho_1^*(\varphi_{\rho_2})\, =\, \varphi_{\rho_2\circ \rho_1}: (\rho_2\circ \rho_1)^*\sF'_{S_3}\lra \sF'_{S_1} 
\end{equation*}
for all morphisms $\rho_1: S_1\to S_2$, $\rho_2: S_2\to S_3$ in $\sS$.
\medskip

Given two $R$-sheaves $\sF'$, $\sG'$ on $X$ as defined above a morphism $\alpha:\sF' \to \sG'$ consists of collection of morphisms $\alpha_S: \sF'_S \to \sG'_S$ for every $S\in \sS$ that are compatible with the isomorphisms \eqref{sheaf1}, i.e.\ the diagram 
\begin{equation*}
\label{sheaf3}
\begin{CD}
\rho^*\sF'_{S_2}@> \varphi_{\rho} >> \sF'_{S_1}\\
@VV \rho^* \alpha_{S_2} V @VV  \alpha_{S_1} V \\
\rho^*\sG'_{S_2}@> \varphi_{\rho} >> \sG'_{S_1}
\end{CD}
\end{equation*}
commutes for all morphisms $\rho: S_1\to S_2$ in $\sS$. 

Let $\Sh'(X, R)$ denote the category of $R$-sheaves $\sF'$ in this sense.
To explain why this notion is equivalent to the one given in Def.\ \ref{df:CuebercolimX} (b) 
we define mutually quasi-inverse equivalences of categories 
\begin{equation*}
\label{sheaf4}
\Sh(X, R) \lra \Sh'(X, R),\quad \Sh'(X, R)\lra \Sh(X, R)
\end{equation*}
Given $\sF\in \Sh(X, R)$ we define $\sF'\in \Sh'(X, R)$ by $\sF_S' :=\res_S \sF$ for every $S\in \sS$. For a morphism $\rho: S_1\to S_2$ there exists also canonical isomorphism $\rho^*(\res_{S_2}(\sF))\to \res_{S_1}(\sF)$. Indeed, if
$\widetilde{\rho}: \cC/X_{S_1}\to \cC/X_{S_2}$ denotes the continuous functor $(U\to X_{S_1})\mapsto (U\to X_{S_1} \stackrel{\rho}{\lra} X_{S_2})$ then we have $\rho^*\circ \res_{S_2} \cong \res_{S_1}$ (see the proof of Lemma \ref{lemma:resexactfaith}). 

Conversely, given an $R$-sheaf $\sF'$ in the sense above we define $\sF\in \Sh(X, R)$ as follows. For an object $(U, \alpha)$ of $ \cC/\colim X$ define a functor 
\begin{equation}
\label{sheaf6}
\cI_{(U, \alpha)}^{\opp}\lra \Mod_R
\end{equation}
by mapping $(S, \psi)\in \cI_{(U, \alpha)}$ to $\sF'_{S}(U\stackrel{\psi}{\lra} X_{S})$ and a morphism 
$\tau: (S_1, \psi_1)\to (S_2, \psi_2)$ to the isomorphism
\begin{equation}
\label{sheaf7}
\sF'_{S_2}(U\stackrel{\psi_2}{\lra} X_{S_2}) \, =\, \sF'_{S_2}(U\stackrel{\tau\circ \psi_1}{\lra} X_{S_2})
\, =\, \tau^*(\sF'_{S_2})(U\stackrel{\psi_1}{\lra} X_{S_1}) \stackrel{\varphi_{\tau}}{\lra}\sF'_{S_1}(U\stackrel{\psi_1}{\lra} X_{S_1})
\end{equation}
(the category $\cI_{(U, \alpha)}$ has been defined in Remark \ref{remarks:colimexplizit} (d); it is filtered since $X$ is exact; note that in order for \eqref{sheaf7} to define a functor we need the cocycle condition \eqref{sheaf2}).
We define $\sF((U, \alpha))$ as the colimit of the \eqref{sheaf6}. It is then easy to see that this construction defines a sheaf on the site $\cC/\colim X$. \enddemo
\end{remarks}

To conclude this section we show that the topos $\Sh(X, R)$ does not depend on the choice of $\sK$. More precisely we show that shrinking $\sK$ does not effect $\Sh(X, R)$. So let $\sK'$ be a subset of $\sK$ that satisfies assumption \ref{assum:kgood} and let $X'$ denote the restriction of $X$ to the subcategory $\sS'= \sS_{G, \sK'}$ of $\sS_{G, \sK}$.
Consider the continuous and cocontinuous functor 
\begin{equation}
\label{shrinkk1}
u: \cC/\colim X'\lra \cC/\colim X, \quad (U, \alpha') \mapsto (U, \alpha)
\end{equation}
where $\alpha$ denotes the image of $\alpha'$ under the canonical map 
\begin{equation}
\label{shrinkk2}
\colim h^U\circ X' \lra \colim h^U\circ X. 
\end{equation}
By (\cite{stacksproject}, Tags 09W7, 03DI and 03DJ) the functor $u$ yields a morphism of topoi 
\begin{equation*}
\label{shrinkk3}
(g^*, g_*): \Sh(X', R)\lra \Sh(X, R) 
\end{equation*}
with $g^* = u^s$. The functor $g^*$ is exact since it has a left adjoint $u_s$.

\begin{prop}
\label{prop:toposcan}
The functors $g^*$ and $g_*$ are mutually quasi-inverse equivalences of categories.
\end{prop}

For the proof we need 

\begin{lemma}
\label{lemma:ufullyfaith} 
The functor \eqref{shrinkk1} is fully faithful.
\end{lemma}

\begin{proof} Firstly, we show that the map \eqref{shrinkk2} is injective. For that let $\alpha', \alpha''\in \colim h^U\circ X'$ with the same image $\alpha$ under \eqref{shrinkk2}. Let $(S', \psi')$ and $(S'', \psi'')$ be representatives of the objects $(U, \alpha')$ and $(U, \alpha'')$ of $\cC/\colim X'$ respectively. Then $(S', \psi')$ and  $(S'', \psi'')$ can be viewed as representatives of the object $(U, \alpha)$ of $\cC/\colim X$. By Remark \ref{remarks:colimexplizit} (d) there exists a diagram of the form
\begin{equation*}
\label{colimexplizit2a}
(S', \psi')\stackrel{\tau}{\lla}(S, \psi) \stackrel{\tau'}{\lra}(S'', \psi'')
\end{equation*}
in $\cI_{(U, \alpha)}$. The fact that $\tau$ is a morphism with target in $\sS'$ implies that its source lies in $\sS'$ as well, i.e.\ $S\in \sS'$. It follows $\alpha'= (\psi')^*(\alpha_{S'})= (\psi)^*(\alpha_{S})= (\psi'')^*(\alpha_{S''})=\alpha''$.  

Let $(U_1,\alpha_1')$, $(U_2,\alpha_2')$ be objects of $\cC/\colim X'$ and let $(U_1,\alpha_1)$, $(U_2,\alpha_2)$ be their images under $u$. We have to see that the natural map
\begin{equation*}
\label{shrinkk4}
\Hom_{\cC/\colim X'}((U_1,\alpha_1'), (U_2,\alpha_2'))\lra \Hom_{\cC/\colim X}((U_1,\alpha_1), (U_2,\alpha_2))
\end{equation*}
is bijective, i.e.\ for a morphism $\phi :U_1 \to U_2$ in $\cC$ we have to show that $\phi^*(\alpha_2) = \alpha_1$ implies
$\phi^*(\alpha'_2) = \alpha'_1$. Let $(S', \varphi_2')$ be a representative of $(U_2,\alpha_2')$. If $\phi^*(\alpha_2) = \alpha_1$ then $(U_1, \varphi_2' \circ \phi)$ is a representatives $(U_1,\alpha_1)$ hence also of $(U_1,\alpha_1')$ because of the injectivity of \eqref{shrinkk2}. It follows $\phi^*(\alpha'_2) = \alpha'_1$.
\end{proof} 

\begin{proof}[Proof of Prop.\ \ref{prop:toposcan}] By (\cite{stacksproject}, Tag 00XT) and Lemma \ref{lemma:ufullyfaith} the counit of adjunction $\varep:g^*\circ g_*\to 1_{\Sh(X, R)}$ is an isomorphism. We have to show that the unit $\eta: 1_{\Sh(X, R)}\to g_*\circ g^*$ is an isomorphism as well. Since $g^*(\eta) = g_*(\varep)^{-1}$ is an isomorphism, it suffices to see that $g^*$ is faithful. 
For that choose a non-empty $S'\in \sS$ and consider the diagram 
\begin{equation*}
\label{resjshriekc}
{\xymatrix@-0.5pc{\Sh(X, R)\ar[dd]^{g^*}\ar[drr]^{\res_{S'}}\\
&& \Sh(\cC/X_{S'}, R)\\
\Sh(X', R)\ar[urr]^{\res_{S'}}
}}
\end{equation*}
By Lemma \ref{lemma:resexactfaith} (b) the vertical functors are faithful. Hence $g^*$ is faithful as well.
\end{proof}

\subsection{Morphisms between $\sS$-objects}
\label{subsection:ssmorphism}

We keep the notation and assumptions of the last section. 

\begin{df}
\label{df:Smorph}
Let $X$, $Y$ be $\sS$-objects in $\cC$. A morphism $f: X\to Y$ of $\sS$-objects in $\cC$ is a morphism of functors such that the diagram
\begin{equation}
\label{sspacemor1}
\begin{CD} X_{S_1} @> f_{S_1} >> Y_{S_1}\\
@VV \rho V @VV \rho V \\
X_{S_2} @> f_{S_2} >> Y_{S_2}
 \end{CD}
\end{equation}
is cartesian for every morphism $\rho: S_1\to S_2$ in $\sS$.
\end{df}

Let $f: X\to Y$ be a morphism of $\cC$-objects in $\cC'$. It induces two functors 
\begin{eqnarray}
&& u: \cC/\colim X\lra \cC/\colim Y,\quad (U, \alpha) \mapsto (U, f\circ \alpha),\nonumber \\
&& v: \cC/\colim Y\lra \cC/\colim X,\quad (V, \beta) \mapsto (V, \beta)\times_Y X
\label{sspacemor3}
\end{eqnarray}
For the first note that $f$ induces an obvious morphism of functors $h^U\circ X\to h^U\circ Y$ hence a map between colimits $\colim h^U\circ X\to \colim h^U\circ Y, \alpha\mapsto f\circ \alpha$. 
The functor \eqref{sspacemor3} is defined a follows. Let $(V, \beta)\in  \cC/\colim Y$ and let $(S, \psi: V\to Y_S)$ be a representative of $(V, \beta)$. Then we defined $(V, \beta)\times_Y X): = (X_S\times_{Y_S} V, \pr_{X_S}^*(\alpha_S))$, i.e.\ the pair $(S, \pr: X_S\times_{Y_S} V\to X_S)$ is a representative of $(V, \beta)\times_Y X$. That this definition is independent of the choice of the representative $(S, \psi)$ follows easily from Remark \ref{remarks:colimexplizit} (d)
and the fact that the diagrams \eqref{sspacemor1} are cartesian. One easily verifies that $u$ is continuous and cocontinuous and that $v$ is continuous and right adjoint to $u$. By (\cite{stacksproject}, Tags 00XO, 00XR and 00XY) the functors $u$ and $v$ induces the same morphism of topoi 
\begin{equation*}
\label{sspacemor4}
(f^*, f_*):= (u^s, \, _s u) = (v_s, v^s) : \Sh(X, R) \lra \Sh(Y, R).
\end{equation*}

\begin{prop}
\label{prop:sspacemortop}
Let $f:X\to Y$ be a morphism of $\sS$-Objects in $\cC$. 
\medskip

\noi (a) The functor $f^*: \Sh(Y, R)\lra \Sh(X, R)$ is exact. 
\medskip

\noi (b) For $S\in \sS$ the diagram of topoi 
\begin{equation*}
\label{sspacemor5}
\begin{CD} 
\Sh(\cC/X_S, R) @> (j_{X_S}^*, (j_{X_S})_*) >> \Sh(X, R)\\
@VV (f_S^*, (f_S)_*) V @VV (f^*, f_*) V \\
\Sh(\cC/Y_S, R) @> (j_{Y_S}^*, (j_{Y_S})_*) >> \Sh(Y, R) 
\end{CD}
\end{equation*}
commutes. 
\medskip

\noi (c) The diagram of topoi
\begin{equation*}
\label{sspacemor6}
{\xymatrix@+0.1pc{\Sh(X, R) \ar[rr]^{(f^*, f_*)}\ar[rdd]_{(\wX_s, \wX^s)} && \Sh(Y, R)\ar[ldd]^{(\wY_s, \wY^s)}\\\\
& \Sh(\sS, R)
}}
\end{equation*}
commutes, i.e.\ we have 
\begin{equation*}
\label{sspacemor6a}
H^0(Y_{\infty}, f_*(\sF))\, = \, H^0(X_{\infty}, \sF) \qquad\mbox{and}\qquad f^*(M_{Y, G}) \, = \, M_{X, G}
\end{equation*}
for every $M\in \Mod_R^{\sm}(G)$ and $\sF\in \Sh(X, R)$.
\end{prop}

\begin{proof} (a) is a consequence of the fact that $f^*$ has a left adjoint (namely $u_s$) and
(b) follows from the commutativity of the diagram 
\begin{equation*}
\label{sspacemor5a}
\begin{CD} \cC/X_S @> j_{X_S} >> \cC/\colim X\\
@VVV @VV u V \\
\cC/Y_S@> j_{Y_S} >> \cC/\colim Y 
\end{CD}
\end{equation*}
together with (\cite{stacksproject}, Tag 03L5) where the left vertical functors is given by $(U\to X_S) \mapsto (U\to X_S\stackrel{f_S}{\lra} Y_S)$. Finally (c) is a consequence of the commutativity of the diagram
\begin{equation*}
\label{sspacemor7a}
{\xymatrix@-0.5pc{ & \sS\ar[ldd]_{\wY}\ar[rdd]^{\wX}\\\\
\cC/\colim Y  \ar[rr]^v && \cC/\colim X.
}}
\end{equation*}
\end{proof}

\subsection{Change of the group $G$}
\label{subsection:functinG}

As in the previous section $\cC$ denotes a site that has fibre products. In this section we consider several functorial properties of $\sS$-objects in $\cC$ and their associated topoi with respect to changing the group $G$. To begin with let $G'$ be a second locally profinite groups and let $\varphi: G\to G'$ be a continuous homomorphism. Let $\sK$ (resp.\ $\sK'$) be a subset of the set of open subgroups of $G$ (resp.\ $G'$) that satisfies assumption \ref{assum:kgood} and put $\sS=\sS_{G, \sK}$ and $\sS'=\sS_{G', \sK'}$.

\begin{df}
\label{df:Smorphphi}
Let $X$ be an $\sS$-objects and let $Y$ be an $\sS'$-object in $\cC$. A morphism $f: X\to Y$ compatible with $\varphi$ consists of a collection of morphism in $\cC$
\begin{equation*}
\label{sspacemorgen1}
f_{\xi}: X_S\lra Y_{S'}
\end{equation*}
for every $G$-equivariant morphism $\xi: S\to S'_{\varphi}$ with $S\in \sS$ and $S'\in \sS'$. We require that the following two conditions hold
\medskip

\noi (i) If 
\begin{equation}
\label{sspacemorgen2}
\begin{CD} 
S_1 @>  \rho  >> S_2\\
@VV \xi_1 V @VV \xi_2 V \\
S_1' @>  \rho' >> S_2'
 \end{CD}
\end{equation}
is commutative where $\rho$ and $\rho'$ are morphisms in $\sS$ and $\sS'$ respectively and where $\xi_1$ and $\xi_2$ are $G$-equivariant maps then the diagram 
\begin{equation}
\label{sspacemorgen3}
\begin{CD} 
X_{S_1} @> \rho >> X_{S_2}\\
@VV  f_{\xi_1}V @VV f_{\xi_2} V \\
Y_{S_1'} @>  \rho' >> Y_{S_2'}
 \end{CD}
\end{equation}
commutes as well.
\medskip

\noi (ii) If the diagram \eqref{sspacemorgen2} is cartesian (as a diagram in the category $\sS_G$) then \eqref{sspacemorgen3} is cartesian as well.
\end{df}

\begin{remark}
\label{remark:ssmoroldnew} 
\rm If $G=G'$, $\sK=\sK'$ and if $\varphi=\id_G:G \to G$ then it is easy to see that the notion of a morphism compatible 
$\id_G$ is equivalent to the notion of a morphism between $\sS$-objects given before. Indeed if $f':X\to Y$ is a morphism in the sense of Def.\ \ref{df:Smorph} then we obtain a morphism $f: X\to Y$as above by defining $f_{\xi}: X_S\to Y_{S'}$ as the composite $X_S\stackrel{\xi}{\lra} X_{S'} \stackrel{f_{S'}}{\lra}Y_{S'}$ where $\xi:S\to S'$ is a morphism in $\sS$. Conversely, if $f: X\to Y$ is a morphism in the sense of Def.\ \ref{df:Smorphphi}
then the collection of morphism $f'_S := f_{\id_S}: X_S\to Y_S$, $S\in \sS$ defines a morphism between $\sS$-objects. \enddemo
\end{remark}

We are going to show that a morphism $f: X\to Y$ compatible with $\varphi$ induces a morphism of topoi 
\begin{equation}
\label{smorgen1}
(f^*, f_*): \Sh(X, R) \lra \Sh(Y, R).
\end{equation}
For that we may assume that every $K\in \sK$ is contained in $\varphi^{-1}(K')$ for some $K'\in \sK$. Indeed, by Prop.\ \ref{prop:toposcan} replacing $\sK$ with the subset 
\begin{equation*}
\label{shrinksk}
\sK_0=\{K\in \sK\mid \exists K'\in \sK'\, \mbox{such that } K\subseteq \varphi^{-1}(K')\}.
\end{equation*}
does not change the topos $\Sh(X, R)$.

The above hypothesis on $\sK$ implies that for every object $S\in \sS$ there exists an object $S'\in \sS'$ and a $G$-equivariant map $\xi: S\to S'_{\varphi}$. Hence for $U\in \cC$ the morphism $f_{\xi}: X_S\to Y_{S'}$ induces a map
\begin{equation}
\label{smorgen2}
\begin{CD}
h^U(X_S)@> \psi\mapsto f_{\xi}\circ \psi >> h^U(Y_{S'}) @>\can >> \colim h^U \circ Y
\end{CD}
\end{equation}
It is easy to see that \eqref{smorgen2} is independent of the choice of $\xi$. Thus the collection of maps \eqref{smorgen2} for varying $S\in \sS$ induce a map
\begin{equation}
\label{smorgen3}
\colim h^U\circ X \lra \colim h^U \circ Y.
\end{equation}
Hence we obtain a functor 
\begin{equation*}
\label{smorgen4}
u: \cC/\colim X\lra \cC/\colim Y, \qquad (U, \alpha) \mapsto (U, f\circ \alpha)
\end{equation*}
Here $f\circ \alpha$ denotes the image of $\alpha$ under \eqref{smorgen3}. Clearly, $u$ is continuous and cocontinuous so it induces the morphism of topoi \eqref{smorgen1} with $(f^*, f_*)= (u^s, \, _s u)$. 

\begin{prop}
\label{prop:smorcomp}
Let $X$ be an $\sS_{G, \sK}$-object, let $Y$ be an $\sS_{G', \sK'}$-object in $\cC$ and let 
$f:X\to Y$ be a morphism compatible with $\varphi:G \to G'$. 
\medskip

\noi (a) The functor $f^*: \Sh(Y, R)\lra \Sh(X, R)$ is exact. 
\medskip

\noi (b) For $S\in \sS$, $S'\in \sS'$ and a $G$-equivariant map $\xi:S\to S'_{\varphi}$ the diagram of topoi 
\begin{equation*}
\label{smorcomp2}
\begin{CD} 
\Sh(\cC/X_S, R) @> (j_{X_S}^*, (j_{X_S})_*) >> \Sh(X, R)\\
@VV (f_{\xi}^*, (f_{\xi})_*) V @VV (f^*, f_*) V \\
\Sh(\cC/Y_{S'}, R) @> (j_{Y_{S'}}^*, (j_{Y_{S'}})_*) >> \Sh(Y, R) 
\end{CD}
\end{equation*}
commutes. 
\medskip

\noi (c) The diagram of morphism of topoi
\begin{equation*}
\label{smorcomp3}
\begin{CD} 
\Sh(X, R) @> (f^*, f_*) >> \Sh(Y, R)\\
@VV (\wX_s, \wX^s) V @VV (\wY_s, \wY^s) V \\
\Sh(\sS, R) @> (\varphi^*, \varphi_*) >> \Sh(\sS', R) 
\end{CD}
\end{equation*}
commutes.
\end{prop}

\begin{proof} The proofs of (a) and (b) are similar to the proofs of the corresponding statement in Prop.\ \ref{prop:sspacemortop}. 

For (c) note that it suffices to prove one of the equalities $\wY^s \circ f_* \, =\, \wY^s$, $f^* \circ \wY_s= \wX_s$. 
Firstly, we consider the case $G\in \sK$, $G'\in \sS'$ so that $\sS=\sS_G$, $\sS'=\sS_{G'}$ and so that the morphism $(f^*, f_*):\Sh(X, R)\to \Sh(Y, R)$ can be identified with $(f_{\xi}^*, (f_{\xi})_*): \Sh(\cC/X_G, R)\to \Sh(\cC/Y_{G'}, R)$ where $\xi$ is the unique map $G/G=\pt\to G'/G'=\pt$. Consider the following diagram of continuous functors 
\begin{equation*}
\label{smorcomp3a}
\begin{CD} 
\sS_{G'} @> \wvphi >> \sS_G\\
@VV S'\mapsto (Y_{S'}\to Y_{G'}) V @VV S\mapsto (Y_S\to X_G) V \\
\cC/Y_{G'} @> v >> \cC/X_G\\
\end{CD}
\end{equation*}
where $\wvphi$ is given by $S'\mapsto S'_{\varphi}$ and $v$ is the base change functor with respect to $f_{\xi}: X_G\to Y_{G'}$. It commutes since condition (ii) of Def.\ \ref{df:Smorphphi} above implies that 
\[
v(Y_{S'}\to Y_{G'}) = (Y_{S'}\times_{Y_{G'}} X_G \to X_G) = (X_{S'_{\varphi}\times_{\pt} \pt} \to X_{\pt}) =(X_{S_{\varphi}} \to X_G).
\]
By \eqref{basechange2} we have $v^s=(f_{\xi})_*$ so we obtain
\[
\wY^s \circ f_* = \wY^s\circ (f_{\xi})_*= \wY^s\circ v^s= \wvphi^s \circ \wX^s = \varphi_*\circ \wX^s.
\]
In the general case we choose $K\in \sK$, $K'\in \sK'$ with $\varphi(K) \subseteq K'$. Put $S=G/K$, $S'=G'/K'$ and let $\xi: S\to S', gK \mapsto \varphi(g) K'$. Let $\io: K\to G$ and $\io': K'\to G'$ denote the inclusions and let $\varphi_0: K\to K'$ be the restriction of $\varphi$ to $K$. We have already seen that 
\[
(f_{\xi})^* \circ (\widetilde{Y|_{K'}})_s\,= \,(\widetilde{X|_K})_s \circ (\varphi_0)^*.
\]
Together with (b) and Prop.\ \ref{prop:restiota} we get
\begin{eqnarray*}
&& \res_S\circ f^* \circ \wY_s \, = \, j_{X_K}^* \circ f^* \circ \wY_s\, =\, (f_{\xi})^* \circ j_{Y_{K'}}^* \circ \wY_s\, =\, (f_{\xi})^* \circ (\widetilde{Y|_{K'}})_s\circ (\io')^*\\
 && \, =\,(\widetilde{X|_K})_s \circ (\varphi_0)^* \circ (\io')^* 
  \, =\,(\widetilde{X|_K})_s \circ \io^* \circ \varphi^*   \, =\,  j_{X_K}^* \circ \wX_s \circ \varphi^*   
 \,=\,  \res_S \circ \wX_s \circ \varphi^*
   \end{eqnarray*}
hence $f^* \circ \wY_s=\wX_s\circ \varphi^*$ by Lemma \ref{lemma:exactlayer} (b).
\end{proof}

\paragraph{Push-out of an $\sS$-object} We keep the assumption of the previous paragraph and assume additionally that 
\begin{equation*}
\label{kkphicomp}
\varphi^{-1}(K') \in \sK 
\end{equation*}
holds for every $K'\in \sK'$. This implies that 
\begin{equation*}
\label{kkphicomp2}
\wvphi: \sS'\to \sS, \qquad S'\to S'_{\varphi}
\end{equation*} 
is a continuous functor between sites. Let $X: \sS\to \cC$ be an $\sS$-object in $\cC$. Then $X_{\varphi}: =X\circ \wvphi : \sS'\to \cC$ is an $\sS'$-object in $\cC$. It will be called the {\it push-out of $X$ with respect to $\varphi$}. We define a (continuous and cocontinuous) functor 
\begin{equation*}
\label{phipushout1}
\wPhi: \cC/\colim X_{\varphi}\lra \cC/\colim X,\quad  (U, \alpha)\mapsto (U, \beta)
\end{equation*}
where $\beta$ is the image of $\alpha$ under the canonical map (see (\cite{stacksproject}, Tag 002K))
\begin{equation*}
\label{phipushou2}
\colim h^U\circ X_{\varphi} \lra \colim h^U\circ X.
\end{equation*}
We denote the image of an $R$-sheaf $\sF$ on $X$ under $\wPhi^s$ by $\sF_{X_{\varphi}}$ so that $\wPhi^s$ can be written as 
\begin{equation*}
\label{phipushou3}
\Sh(X, R) \lra \Sh(X_{\varphi}, R), \quad \sF\mapsto \sF_{X_{\varphi}}.
\end{equation*}

\begin{lemma}
\label{lemma:pushoutopos} 
The functor \eqref{phipushou3} is an equivalence of categories.
\end{lemma}

\begin{proof} Assume first that for every $K\in \sK$ there exists $K'\in \sK'$ such that $K\subseteq \varphi^{-1}(K')$. Then the category $\sS'$ is cofinal in $\sS$, so the map \eqref{phipushou2} is bijective (see \cite{stacksproject}, Tag 04E7) and $\wPhi$ is an equivalence of sites. In the general case we may replace $\sK$ by the subset 
\begin{equation*}
\label{shrinksk2}
\sK_0=\{K\in \sK\mid \exists K'\in \sK'\, \mbox{such that } K\subseteq \varphi^{-1}(K')\}.
\end{equation*}
By Prop.\ \ref{prop:toposcan} this does not effect the category $\Sh(X, R)$. 
\end{proof}

Note that the diagram of continuous functors
\begin{equation*}
\label{phipushout4}
\begin{CD}
\cC/\colim X_{\varphi}@> \wPhi>> \cC/\colim X\\
@V \widetilde{X_{\varphi}} VV @V \wX VV\\
\sS' @> \wphi >> \sS
\end{CD}
\end{equation*}
commutes. This implies that the diagram of topoi 
\begin{equation}
\label{phipushout5}
{\xymatrix@-0.5pc{ & \ar@/_/[ddl]_{(\wX_s, \wX^s)}\Sh(X, R) \cong \Sh(X_{\varphi}, R)\ar@/^/[ddr]^{((\widetilde{X_{\varphi}})_s, \widetilde{X_{\varphi}}^s)}\\\\
\Sh(\sS, R) \ar[rr]^{(\varphi^*, \varphi_*)}&& \Sh(\sS', R)}}
\end{equation}

\paragraph{Induced $\sS$-objects} Now we consider a special instance of the previous construction namely the case where $\varphi=\io: H\to G$ is the inclusion of a closed subgroup $H$ of $G$. We let $\sK^H$ be the collection of open subgroups of $H$ of the form $K\cap H$ for $K\in \sK$. Let $Y$ be an $\sS_{H, \sK^H}$-object in $\cC$. We define the induced $\sS$-object $\caInd_H^G \,Y$ as the push-out of $Y$ with respect to $\io$, so it is the functor 
\begin{equation}
\label{indssob}
\caInd_H^G\, Y: \sS\lra \cC, \quad S\mapsto Y_S
\end{equation}
where we view a left $G$-set $S$ as an $H$-set via $\io$. The commutativity of \eqref{phipushout5} yields

\begin{lemma}
\label{lemma:indssho} 
Let $Y$ be an $\sS_{H, \sK^H}$-object in $\cC$ and put $X=\caInd_H^G \,Y$. We have
\begin{equation*}
\label{indcompres}
H^0(X_{\infty}, \sF_X)\, \cong \, \Ind_H^G H^0(Y_{\infty}, \sF) \quad \mbox{and}\quad M_{X,G}\, \cong \, ((\Res_H^G M)_{Y, H})_X
\end{equation*}
for every $\sF\in \Sh(Y, R)$ and $M\in \Mod_R^{\sm}(G)$.
\end{lemma}

\paragraph{Restriction of $\sS$-objects} Let $G$ be a locally profinite group, let 
$\sK$ be a subset of the set of open subgroups of $G$ that satisfies hypothesis \ref{assum:kgood} and put $\sS=\sS_{G, \sK}$. 

\begin{df}
\label{df:Sspaceres}
Let $X$ be an $\sS$-object in $\cC$ and let $K\in \sK$. The restriction $X|_K$ of $X$ to $\sS_K$ is defined as the composite of the continuous functor \eqref{fibergk} with $X$, i.e.\ $X|_K$ is defined as 
\begin{equation}
\label{sspaceres1}
X|_K: \sS_K \lra \cC, \quad S'\mapsto X_{G\times_K S'}.
\end{equation}
\end{df}

Note that we have $G\times_K S'\in \sS$ for every $S'\in \sS_K$ so that \eqref{sspaceres1} is well-defined. Note also that we have $(X|_K)_{K'} = X_{K'}$ for every open subgroup $K'$ of $K$. Moreover since the category $\sS_K$ has a final object we can identify the topos $\Sh(X|_K, R)$ with $\Sh(\cC/X_K, R)$. 

\begin{lemma}
\label{lemma:reslocal2} 
Under the identifications $\Sh(\sS_K, R)=\Sh(\sS/S, R)$ and $\Sh(X|_K, R) = \Sh(\cC/X_K, R)$ we have an equality of morphisms of functors
\begin{equation*}
\label{reslocal2}
((\widetilde{X|_K})_s, (\widetilde{X|_K})^s) \, =\, ((\wX/S)_s, (\wX/S)^s): \Sh(X|_K, R)\lra \Sh(\sS_K, R).
\end{equation*}
\end{lemma}

\begin{proof} By the definition of $X|_K$ the diagram of continuous functors between sites
\begin{equation*}
\label{reslocal2a}
{\xymatrix@-0.5pc{\sS_K \ar[rr]^{\eqref{sitegsetlock2}}\ar[rdd]_{\widetilde{X|_K}} && \ar[ldd]^{\wX/S}\sS/S\\\\
& \cC/X_K
}}
\end{equation*}
commutes. The assertion follows. 
\end{proof}

As a consequence of Lemmas \ref{lemma:reslocal}, \ref{lemma:reslocal2} and Remark \ref{remarks:Xloccomp} (a) we obtain

\begin{prop}
\label{prop:restiota} 
The following diagram of topoi commutes 
\begin{equation*}
\label{resiota1}
\begin{CD} 
\Sh(\cC/X_K , R) @> (j_{X_K}^*, (j_{X_K})_*) >> \Sh(X, R)\\
@VV ((\widetilde{X|_K})_s, (\widetilde{X|_K})^s) V @VV (\wX_s, \wX^s) V \\
\Sh(\sS_K, R) @> (\io^*, \io_*) >> \Sh(\sS, R) 
\end{CD}
\end{equation*}
commutes. Here $\io: K \hra G$ denotes the inclusion. 
\end{prop}

\begin{remark}
\label{remark:resxinft} 
\rm Since $j_{X_K}\circ \widetilde{X|_K} = \wX \circ u: \sS_K \to \cC/\colim X$ where $u$ is the functor \eqref{fibergk} we have $\widetilde{X|_K}^s \circ \res_S = u^s\circ \wX^s = \io^*\circ \wX^s$. Therefore we get 
\begin{equation*}
\label{resiota2}
\Res_K^G  H^0(X_{\infty}, \sF)\, = \, H^0((X|_K)_{\infty}, \sF_K).
\end{equation*}
for every $\sF\in \Sh(X, R)$.  \enddemo
\end{remark}

\subsection{Cohomology of $\sS$-objects}
\label{subsection:cohomology}

As in the previous section $G$ denotes a locally profinite group, $\sK$ a subset of the set of open subgroups of $G$ satisfying hypothesis \ref{assum:kgood} and $\sS=\sS_{G,\sK}$ the associated site. Let $\cC$ be a site that satisfies hypothesis \ref{assum:sitegood}.
 
\begin{df}
\label{df:etcohom}
Let $X$ be an $\sS$-object in $\cC$ and let $n\ge 0$. 
\medskip

\noi (a) We denote the $n$-th right derived functor of  \eqref{sheafxgmod} by 
\begin{equation*}
\label{sheafxgmodder}
H^n(X_{\infty}, \wcdot): \Sh(X, R) \lra \Mod_R^{\sm}(G).
\end{equation*}
\medskip

\noi (b) The $n$-th right derived functor of  
\begin{equation}
\label{Xsection}
H^0(X, \wcdot): \Sh(X, R) \lra \Mod_R, \,\, \sF\mapsto H^0(X, \sF):=H^0(X_{\infty}, \sF)^G
\end{equation}
will be denoted by 
\begin{equation*}
\label{Xcohom}
H^n(X, \wcdot): \Sh(X, R) \lra \Mod_R.
\end{equation*}
\end{df}

\begin{remark}
\label{remark:gcompactcoh}
\rm If $G\in \sK$ then $\sS=\sS_G$ and the topos $\Sh(X, R)$ can be identified with $\Sh(\cC/X_G, R)$ and the functor \eqref{Xsection} with $H^0(X_G, \wcdot)$. So in this case we have
$H^n(X, \sF) = H^n(X_G, \sF_G)$
for every $n\ge 0$. \enddemo
\end{remark}

For $n\ge 0$ the $n$-th right derived functors of the functor $\Mod_R^{\sm}(G) \to \Mod_R, M \mapsto M^G$ will be denoted by $H^n(G, \wcdot)$. 

\begin{prop}
\label{prop:hs}
For $\sF\in \Sh(X, R)$ we have 
\medskip

\noi (a) $H^n(X_{\infty}, \sF) \, =\, \, \dlim_{K\in \sK^{\opp}} H^n(X_K, \sF_K)$.
\medskip

\noi (b) For $K\in \sK$ we have $\Res_K^G H^n(X_{\infty}, \sF)= H^n((X|_K)_{\infty}, \sF_K)$.
\medskip

\noi (c) (Covering spectral sequence) Assume that $X$ is exact. Then there exists a spectral sequence
\begin{equation}
\label{hss}
E_2^{rs} = H^r(G, H^s(X_{\infty}, \sF)) \, \Longrightarrow\, H^{r+s}(X, \sF).
\end{equation}
\end{prop}

\begin{proof} (a) Let $0\to \sF \to \sI^{\bu}$ be an injective resolution. By Lemma \ref{lemma:resexactfaith} (a) the sequence $0\to \sF_K \to \sI_K^{\bu}$ is an injective resolution of $\sF_K$ for every $K\in \sS$. Together with Remark \ref{remark:sheavestogmod} (a) it follows
\[
H^n(X_{\infty}, \sF) \, =\, H^n( H^0(X_{\infty}, \sI^{\bu})) \, =\,  \dlim_K H^n( H^0(X_K, \sI_K^{\bu})) \, =\,  \dlim_K H^n(X_K, \sF_K).
\]
(b) follows from Remark \ref{remark:resxinft}.

\noi (c) The functor $H^0(X, \wcdot)$ is the composite of two left exact functors 
\[
H^0(X_{\infty}, \wcdot) = H^0(G, \wcdot) \circ H^0(X_{\infty}, \wcdot): \Sh(X, R) \lra \Mod_R^{\sm}(G) \lra \Mod_R.
\]
By Cor.\ \ref{corollary:gmodsheafxexact} the functor $H^0(X_{\infty}, \wcdot)$ has an exact left adjoint. Hence the corresponding Grothendieck spectral sequence \eqref{hss} exists.
\end{proof}

\begin{remarks}
\label{remarks:gcompactcoh2}
\rm (a) By Prop.\ \ref{prop:hs} (a) and (\cite{stacksproject}, Tag 03FD) the cohomology groups with coefficient in an $R$-sheaf $\sF$ introduced in Def.\ \ref{df:etcohom} (a) and (b) depend only on the underlying sheaf of abelian groups on $X$. 
\medskip

\noi (b) If $G$ is compact (hence profinite) and if $G\in \sK$ then by Remark \ref{remark:gcompactcoh} the limit terms of \eqref{hss} can be identified with the cohomology groups $H^\bu(X_G, \cF)$ where $\cF=\sF_{G}$. In this case the spectral sequence is of the form 
\begin{equation*}
\label{hsscomp}
E_2^{rs} = H^r(G, H^s(X_{\infty}, \sF)) \, \Longrightarrow\, H^{r+s}(X_G, \sF_G).
\end{equation*}
so we recover the Hochschild-Serre spectral sequence of (\cite{artin}, 3.4.7). \enddemo
\end{remarks}

By using the spectral sequence \eqref{hss} we get the following criterion for the cohomology groups $H^{\bu}(X_{\infty}, \sF)$ to be admissible $R[G]$-modules.

\begin{corollary}
\label{corollary:gadmissible}
Let $X$ be an exact $\sS$-object in $\cC$ and let $\sF\in \Sh(X, R)$. Assume that \\
(i) $R=A$ is an Artinian local ring with finite residue field $k$ of characteristic $p$.\\
(ii) $G$ is a $p$-adic (i.e.\ the $\bQ_p$-valued points of a) reductive group. \\
(iii) $H^n(X_K, \sF_K)\in \Mod_{A, f}$ for every $K\in \sK$ and $n\ge 0$. \\
Then the discrete $A[G]$-module $H^n(X_{\infty}, \sF)$ is admissible for every $n\ge 0$.
\end{corollary}

\begin{proof} We first remark by applying Prop.\ \ref{prop:hs} (b), (c) to the restriction $X|_K: \sS_K \to \cC$ we obtain a covering spectral sequence 
\begin{equation}
\label{hsscomp2}
E_2^{rs} = H^r(K, H^s(X_{\infty}, \sF)) \, \Longrightarrow\, H^{r+s}(X_K, \sF_K)
\end{equation}
for every $K\in \sK$. We prove the assertion by induction on $n$. Assume that $n\ge 0$ and that $H^r(X_{\infty}, \sF)$ is admissible for every $r<n$. For $K\in \sK$ we have 
$H^r(K, H^s(X_{\infty}, \sF))\in \Mod_{A,f}$ by (\cite{emerton3}, Prop.\ 2.1.9. and Lemma 3.4.4) for every $r\ge 0$ and $s<n$. It follows that the terms $E_m^{rs}$ of the spectral sequence \eqref{hsscomp2} are finitely generated $A$-modules for every $r\ge 0$, $s<n$ and $m\in \{2, 3, \ldots, \infty\}$ as they are subquotient of $E_2^{rs}$. Since $E^n\in \Mod_{A,f}$ by assumption and since $E_{\infty}^{0n} = E_{n+1}^{0n}$ is a quotient of $E^n$ we have $E_{n+1}^{0n}\in \Mod_{A,f}$. The 
exact sequences 
\[
0\lra E_{m+1}^{0n}\lra E_{m}^{0n} \stackrel{d_{m}^{0n}}{\lra} E_{m}^{m(n-m+1)} 
\]
for $m=n, n-1, \ldots, 2$ imply that the $A$-modules $E_n^{0n}, E_{n-1}^{0n}, \ldots, E_2^{0n}$ are finitely generated.
Hence $E_2^{0n}= H^n(X_{\infty}, \sF)^K\in \Mod_{A,f}$. Since this holds for every $K\in \sK$ we conclude $H^n(X_{\infty}, \sF)\in \Mod_A^{\adm}(G)$.
\end{proof}

Now assume that $H$ is a closed subgroup of $G$ and put $\sK^H = \{K\cap H\mid K\in \sK\}$. Regarding the cohomology of induced $\sS$-objects we have

\begin{lemma}
\label{lemma:indsshn} 
Let $Y$ be an exact $\sS_{H, \sK^H}$-object in $\cC$ and let $X=\caInd_H^G \,Y$. For $\sF\in \Sh(Y, R)$ there exists a spectral sequence 
\begin{equation}
\label{coindcohom}
E_2^{rs} = R^s \Ind_H^G H^s(Y_{\infty}, \sF) \, \Longrightarrow\, H^{r+s}(X_{\infty}, \sF_X).
\end{equation}
In particular if the functor $\Ind_H^G: \Mod_R^{\sm}(H)\to  \Mod_R^{\sm}(G)$ is exact then we have 
\begin{equation*}
\label{coindcohom2}
H^n(X_{\infty}, \sF_X)\,=\, \Ind_H^G H^n(Y_{\infty}, \sF)
\end{equation*}
for every $n\ge 0$.
\end{lemma}

\begin{proof} By Lemma \eqref{lemma:indssho} the functor $\Sh(Y, R) \to \Mod_R^{\sm}(G), \,\, \sF\mapsto H^0(X_{\infty}, \sF_X)$ can be factored as 
\[
\begin{CD}
\Sh(Y, R)@> H^0(Y_{\infty}, \wcdot) >> \Mod_R^{\sm}(H) @> \Ind_H^G >> \Mod_R^{\sm}(G).
\end{CD}
\]
By Cor.\ \ref{corollary:gmodsheafxexact} the functor $H^0(Y_{\infty}, \wcdot)$ has an exact left adjoint so the associated Grothendieck spectral 
sequence \eqref{coindcohom} exists.
\end{proof} 

\begin{remark}
\label{remark:coindexact}
\rm In general the functor $\Ind_H^G: \Mod_R^{\sm}(H)\to  \Mod_R^{\sm}(G)$ is not exact. It is exact however if the quotient space $G/H$ is compact or if $R$ is a field of characteristic $0$. \enddemo
\end{remark}

\paragraph{Double edge homomorphism} Let $G'$ be a second locally profinite group, let $\varphi: G\to G'$ be a continuous homomorphism and let $f: X\to Y$ be a $\varphi$-compatible homomorphism in the sense of Def.\ \ref{df:Smorphphi} with $X$ exact. By Prop.\ \ref{prop:smorcomp} (a) there exists a canonical homomorphism 
\begin{equation}
\label{twoedge}
H^n(Y_{\infty}, f_* \sF) \lra \varphi_* H^n(X_{\infty}, \sF)
\end{equation}
for every $\sF\in \Sh(X, R)$ and $n\ge 0$. It is defined as the composite of the two edge morphisms
\[
H^n(Y_{\infty}, f_* \sF) \lra R^n(\wY^s\circ f_*)( \sF)\, =\, R^n(\varphi_* \circ \wX^s)( \sF)\lra \varphi_* H^n(X_{\infty}, \sF).
\]
We are interested in the case when $G'$ is a quotient of $G$. 

So let us assume that $H$ is a closed normal subgroup of $G$, put $\barG= G/H$ and let $\pi: G\to \barG$ denote the projection. We set $\barsK= \{\pi(K) \mid K\in \sK\}$ and $\barsS =\sS_{\barG, \barsK}$. Note that the set $\barsK$ of open subgroups of $\barG$ satisfies hypothesis \ref{assum:kgood}. Also note that for $S\in \sS$ the set of $H$-orbits $\barS : = H\backslash S$ becomes a discrete left $\barG$-set and in fact $\barS$ lies in $\barsS$. Thus we obtain a functor 
\begin{equation*}
\label{sSmodH} 
\sS\to \barsS, \qquad S\mapsto \barS=H\backslash S, \qquad (\rho:S_1\to S_2)\mapsto (\barrho: \barS_1\to \barS_2).
\end{equation*}

Let $f: X\to Y$ be a $\pi$-compatible morphism between an exact $\sS$-object $X$ and a $\barsS$-object $Y$ in $\cC$. For $\sF\in \Sh(X, R)$ the target of \eqref{twoedge} are the $H$-invariant elements of $H^n(X_{\infty}, \sF)$. Thus \eqref{twoedge} can be viewed as a $G$-equivariant map
\begin{equation}
\label{twoedge2}
H^n(Y_{\infty}, f_* \sF) \lra H^n(X_{\infty}, \sF).
\end{equation}
For $S\in \sS$ let $\xi_S: S\to \barS$ denote the map $s\mapsto H\cdot s$ and define 
\begin{equation*}
\label{twoedge3}
f_S:= f_{\xi_S}: X_S\lra Y_{\barS}.
\end{equation*}
Also for $K\in \sK$ we put $\barK=\pi(K)$, $\xi_K= \xi_{G/K}: G/K\to \barG/\barK$ and $f_K:=f_{G/K}: X_K \to Y_{\barK}$. 

By Prop.\ \ref{prop:smorcomp} (b) associated to $\xi_K$ there is a base change homomorphism
\begin{equation*}
\label{basechange1a}
\BC_{\sF, K}= \BC_{\sF, \xi_K}: (f_*(\sF))_{\barK} \lra (f_K)_*(\sF_K). 
\end{equation*}
for every $\sF\in \Sh(X, R)$. Also for $K_1, K_2\in \sK$ with $K_2\subseteq K_1$ associated to the diagram 
\begin{equation*}
\label{sspacemork}
\begin{CD} 
X_{K_2} @> \rho >> X_{K_1}\\
@VV  f_{K_2}V @VV f_{K_1} V \\
Y_{\barK_2} @>\barrho >> Y_{\barK_1}
 \end{CD}
\end{equation*}
(with $\rho: G/K_2\to G/K_1$ and $\barrho:\barG/\barK_1\to \barG/\barK_2$ being the obvious maps) there exists a 
base change morphism (see \cite{stacksproject}, Tag 0736)
\begin{equation}
\label{basechangecohom}
\BC_{\sF, K_1, K_2}^n: \barrho^* R^n (f_{K_1})_*(\sF_{K_1})\lra R^n (f_{K_2})_*(\sF_{K_2}).
\end{equation}
The following criterion for \eqref{twoedge2} to be an isomorphism will be applied in section \ref{subsection:borelserre} to compute the cohomology of the boundary of Hilbert modular $\sS$-schemes. 

\begin{prop}
\label{prop:basechangeiso}
Let $\sF\in \Sh(X, R)$ and assume that for every $K\in \sK$ we have

\noi (i) $\BC_{\sF, K}$ is an isomorphism.

\noi (ii) There exists $K'\in \sK$, $K'\subseteq K$ such that $\BC_{\sF, K, K'}^n = 0$ for all $n\ge 1$.

\noi Then the homomorphism \eqref{twoedge2} is an isomorphism for every $n\ge 0$ (in particular $H$ acts trivially on $H^n(X_{\infty}, \sF)$).
\end{prop}

\begin{proof} More generally than \eqref{basechangecohom} given $K_1, K_2\in \sK$ with $K_2\subseteq K_1$ there exists a base change morphism in the derived category $D^+(Y_{\barK_1}, R):=D^+(\Sh(\cC/Y_{\barK_1}, R))$
\begin{equation}
\label{basechangecohom2}
\BC_{K_1, K_2}: \barrho^* R (f_{K_1})_*(\sF_{K_1})\lra R (f_{K_2})_*(\sF_{K_2}).
\end{equation}
hence a morphisms 
\begin{equation*}
\label{basechangecohom3}
\BC_{K_1, K_2, [m,n]}: \barrho^* t_{[m,n]} R (f_{K_1})_*(\sF_{K_1})\lra  t_{[m,n]} R (f_{K_2})_*(\sF_{K_2}).
\end{equation*}
for every pair of integers $(m,n)$ with $0\le m\le n$ (where $t_{[m,n]} = t_{\ge m}\circ t_{\le n}$). So for $m=n$ we have in particular $\BC_{K_1, K_2, [n,n]}[n]=\BC_{K_1, K_2}^n$. We claim that for every $K\in\sK$ and $n\ge 1$ the following holds
\begin{itemize}
\item[(iii)] There exists $K'\in \sK$ with  $K'\subseteq K$ such that $\BC_{K, K', [1,n]} =0$ and $\BC_{K, K'}^m = 0$ for all $m\ge n+1$.
\end{itemize}
We prove this by induction. For $n=1$ this is assumption (ii). 
Assume that $n\ge 2$ and that the assertion holds for $n-1$. 
We fix $K\in \sK$. By the induction hypothesis there exists $K_1\subseteq K$ that satisfies (iii) for $K$ and $n-1$ and also $K_2\subseteq K_1$ that satisfies (iii) for $K_1$ and $n-1$. 
We claim that (iii) holds for $K_2\subseteq K$ and $n$. Consider the diagram of morphisms between triangles (where $\rho: G/K_2\to G/K$ and $\rho_2: G/K_1\to G/K_2$ are the obvious maps)
\[
\footnotesize{
\begin{CD}
\rho^*t_{[1,n-1]} R (f_K)_*  \sF_K @>>> \rho^*t_{[1,n]} R (f_K)_*  \sF_K@>>>  \rho^*R^n (f_K)_*  \sF_K[-n]@>>> \rho^*t_{[1, n-1]} R (f_K)_* \sF_K[1] \\
@VV\rho_2^*(\BC_{K, K_1, [1,n-1]}) V @VV\rho_2^*(\BC_{K, K_1, [1,n]}) V @VV \rho_2^*(\BC_{K, K_1}^n)[-n] V @VVV\\
\rho_2^*t_{[1,n-1]} R (f_{K_1})_*  \sF_{K_1} @>\alpha >> \rho_2^*t_{[1,n]} R (f_{K_1})_*  \sF_{K_1}@>\beta >>  \rho_2^*R^n (f_{K_1})_*  \sF_{K_1}[-n]@>>> \rho_2^*t_{[1, n-1]} R (f_{K_1})_* \sF_{K_1}[1] \\
@VV\BC_{K_1, K_2, [1,n-1]} V @VV\BC_{K_1, K_2, [1,n]} V @VV \BC_{K_1, K_2}^n[-n] V @VVV\\
t_{[1, n-1]} R (f_{K_2})_* \sF_{K_2} @>\gamma >> t_{[1, n]} R (f_{K_2})_* \sF_{K_2}@>>> R^n (f_{K_2})_* \sF_{K_2}@>>> t_{[1, n-1]} R (f_{K_2})_* \sF_{K_2}[1] \\
\end{CD}
}
\]
The vertical morphisms in the first, third and fourth column vanish. The vanishing of $\BC_{K, K_1}^n$ implies $\beta \circ \rho_2^*(\BC_{K, K_1, [1,n]})= 0$ hence there exists $\mu:  \rho^*t_{[1,n]} R (f_K)_*  \sF_K\to \rho_2^*t_{[1,n-1]} R (f_{K_1})_*  \sF_{K_1} $ such that 
$\rho_2^*(\BC_{\rho_1, [1,n]})=\alpha\circ \mu$. It follows 
\[
\BC_{K, K_2, [1,n]}= \BC_{K_1, K_2, [1,n]} \circ \rho_2^*(\BC_{K, K_1, [1,n]})= \BC_{K_1, K_2, [1,n]} \circ \alpha \circ \mu = \gamma \circ \BC_{K_1, K_2, [1,n-1]} \circ \mu =0.
 \]
This proves (iii).

Now for $K\in \sK$ consider the following distinguished triangle in $D^+(Y_{\barS}, R)$
\begin{equation}
\label{basechangecohom5}
(f_K)_*\sF_K \lra R (f_K)_*\sF_K \lra t_{\ge 1}R (f_K)_*\sF_K\lra \sF_K[1].
\end{equation} 
Because of assumption (i) the first term can be identified with $(f_* \sF)_{\barK}$. By applying the cohomological functor $H^0(Y_{\barK}, \wcdot): D^+(Y_{\barK}, R)\to \Mod_R$ to the triangle \eqref{basechangecohom5}  we obtain the long exact sequence 
\begin{equation*}
\label{basechangecohom6}
\ldots \to H^n(Y_{\barK}, (f_* \sF)_{\barK}) \to H^n(X_K, \sF_K) \to H^n(Y_{\barK}, t_{\ge 1}R (f_K)_* \sF_K) \to H^{n+1}(Y_{\barK}, (f_*\sF)_{\barK})\to \ldots
\end{equation*}
By passing to the direct limit over all $K$ yields 
\begin{equation*}
\label{basechangecohom11}
\begin{CD}
\ldots H^n(Y_{\infty}, f_* \sF)\stackrel{\eqref{basechangecohom2}}{\lra} H^n(X_{\infty}, \sF)
\to \dlim_K H^n(Y_{\barK}, t_{\ge 1}R (f_K)_* \sF_K) \to H^{n+1}(Y_{\infty}, f_* \sF)\ldots
\end{CD}
\end{equation*}
Thus it suffices to show 
\begin{equation}
\label{basechangecohom12}
\dlim_K H^n(Y_{\barK}, t_{\ge 1}R (f_K)_* \sF_K) \, =\, 0
\end{equation}
for every $n\ge 0$. Note that 
\[
H^n(Y_{\barK}, t_{\ge 1}R (f_K)_* \sF_K)\, =\, H^n(Y_{\barK}, t_{[1, n]}R (f_K)_* \sF_K)
\]
for $n\ge 0$ and $K\in \sK$. Therefore \eqref{basechangecohom12} follows from (iii).
\end{proof}

\subsection{The cohomology groups $H_R^{\bu}(X; M, \sF)$}
\label{subsection:hxmf}

As in the last section $G$ denotes a locally profinite group, $\sK$ a subset of the set of open subgroups of $G$ satisfying hypothesis \ref{assum:kgood}, $\cC$ a site that satisfies \ref{assum:sitegood} and $R$ an arbitrary ring. 
Let $\sS=\sS_{G,\sK}$ be the associated site of discrete left $G$-sets and let $X$ be an $\sS$-object in $\cC$. 
For $M\in \Mod_R^{\sm}(G)$ and $\sF\in \Sh(X, R)$ we put 
\[
H_R^0(X; M, \sF) : = \, \Hom_{R[G]}(M, H^0(X_{\infty}, \sF))\, \cong\, \Hom_{\Sh(X, R)}(M_{X,G}, \sF).
\]
Thus we obtain a bifunctor 
\begin{equation*}
\label{ext}
 H_R^0(X; \wcdot, \wcdot): \Mod_R^{\sm}(G)^{\opp} \times \Sh(X, R) \lra R.
\end{equation*}
We denote the $n$-th right derived functor of $H_R^0(X; M, \wcdot)$ (for fixed $M\in \Mod_R^{\sm}(G)$) by
\begin{equation*}
\label{Xcohomm}
H^n(X; M, \wcdot): \Sh(X, R) \lra \Mod_R.
\end{equation*}
Thus we have 
\begin{equation}
\label{ext2}
H_R^n(X; M, \sF) \, = \, \Ext_{\Sh(X, R)}^n(M_{X,G}, \sF) \qquad n\ge 0.
\end{equation}
If $X$ then by Cor.\ \ref{corollary:gmodsheafxexact} the bifunctor 
\begin{equation*}
\label{ext3}
H_R^{\bu}(X; \wcdot, \wcdot): \Mod_R^{\sm}(G)^{\opp} \times \Sh(X, R) \lra R, \quad (M, \sF)\mapsto H_R^{\bu}(X; M, \sF)
\end{equation*}
is a $\delta$-functor in both arguments $M$ and $\sF$.

\begin{remark}
\label{remark:extcohom}
\rm If $M=R$ is equipped with the trivial $G$-action then we have 
\[
H_R^0(X; R, \sF)\, =\, H^0(X_{\infty}, \sF)^G \, =\, H^0(X, \sF)
\]
hence $H_R^n(X; R, \sF) =H^n(X, \sF)$ for every $n\ge 0$. \enddemo
\end{remark}

The Ext groups for the category $\Mod_R^{\sm}(G)$ will be denoted by $\Ext_{R, G}^{\bu}(\wcdot , \wcdot)$.

\begin{prop}
\label{prop:extcohom}
(a) Let $K\in \sK$, $M\in  \Mod_R^{\sm}(K)$ and $\sF\in \Sh(X,R)$. We have 
\[
H_R^{\bu}(X; \cInd_K^G M, \sF) \, = \, H_R^{\bu}(X|_K; M, \sF_K).
\]
In particular if $M=R$ is equipped with the trivial $K$-action then we get
\[
H_R^{\bu}(X; \cInd_K^G R, \sF) \, \cong \, H^{\bu}(X_K, \sF_K).
\]
(b) Assume that $X$ is exact. For $M\in \Mod_R^{\sm}(G)$ and $\sF\in \Sh(X, R)$ there exists a spectral sequence
\begin{equation}
\label{exthss}
E_2^{rs} = \Ext_{R, G}^r(M, H^s(X_{\infty}, \sF)) \, \Longrightarrow\, H_R^{r+s}(X; M, \sF).
\end{equation}
(c) Let $f: X\to Y$ be a morphism of $\sS$-objects in $\cC$, let $M\in \Mod_R^{\sm}(G)$ and let $\sF\in \Sh(X, R)$. 
There exists a spectral sequence
\begin{equation}
\label{extlerayss}
E_2^{rs} = H_R^r(Y; M, R^s f_* \sF)\, \Longrightarrow\, H_R^{r+s}(X; M, \sF).
\end{equation}
\end{prop}

\begin{proof} (a) We have 
\begin{eqnarray*}
&& \Hom_{\Sh(X, R)}((\cInd_K^G M)_{X,G}, \sF) \, \cong \, \Hom_{R[G]}(\cInd_K^G M, H^0(X_{\infty}, \sF))\\
&& \, \cong \, \Hom_{R[K]}(M, H^0((X|_K)_{\infty}, \sF_K)) 
 \,=\, H_R^0(X|_K; M, \sF_K)
 \end{eqnarray*}
for all $\sF\in \Sh(X, R)$ and therefore $(\cInd_K^G M)_{X,G}\cong \cInd_K^G M_{X_K, K}$. The assertion follows now from Lemma \ref{lemma:resexactfaith} (a).

(b) The functor $H_R^0(X; M, \wcdot): \Sh(X,R)\to \Mod_R$ can be factored as 
\[
\begin{CD}
\Sh(X,R)@> H^0(X_{\infty}, \wcdot) >> \Mod_R^{\sm}(G) @> \Hom_{R[G]}(M, \wcdot)>> \Mod_R
\end{CD}
\]
and \eqref{exthss} is the associated Grothendieck spectral sequence.

(c) The functor $H_R^0(X; M, \wcdot)$ can be factored as 
\[
\begin{CD}
\Sh(X, R)@> f_* >> \Sh(Y, R) @> H_R^0(Y; M, \wcdot) >> \Mod_R
\end{CD}
\]
by Prop.\ \ref{prop:sspacemortop} (c) so we obtain a Grothendieck spectral sequence \eqref{extlerayss}.
\end{proof}

\begin{remark} 
\label{remark:extpullback}
\rm Let $f: X\to Y$ be a morphism of $\sS$-objects in $\cC$, let $M\in \Mod_R^{\sm}(G)$ and let $\sF\in \Sh(Y, R)$. 
There exists a canonical homomorphism 
\begin{equation}
\label{extpullback1}
H_R^n(Y; M, \sF)\lra H_R^n(X; M, f^*\sF)
\end{equation}
for $n\ge 0$ defined as the composite 
\[
H_R^n(Y; M, \sF) \lra H_R^n(Y; M, f_* f^* \sF)\lra H_R^n(X; M, f^*\sF).
\]
The first map is induced by the unit of adjunction $\sF\to f_* f^* \sF$ and the second is an edge morphism 
of the spectral sequence \eqref{extlerayss}. \enddemo
\end{remark}

\paragraph{Change of the coefficient ring} Let $\varphi:R\to R'$ be a ringhomomorphism. A discrete $R'[G]$-module can be viewed as a discrete $R[G]$-module via the homomorphism $\varphi:R\to R'$, i.e.\ $\varphi$ induces a homomorphisms of group rings $R[G]\to R'[G]$, hence an exact functor 
\begin{equation*}
\label{homchangering}
\Mod_{R'}^{\sm}(G) \lra \Mod_R^{\sm}(G),\,\, M\mapsto M_{\varphi}.
\end{equation*}
Also given an $R'$-sheaf $\sF$ on $X$ we may view it as an $R$-sheaf via $\varphi$. So we get an exact functor 
\begin{equation}
\label{homchangering2}
\Sh(X,R') \lra \Sh(X,R),\,\, \sF\mapsto \sF_{\varphi}.
\end{equation}
Note that $(M_{\varphi})_{X,G}= (M_{X,G})_{\varphi}$ and $H^0(X_{\infty}, \sF_{\varphi}) =H^0(X_{\infty}, \sF)_{\varphi}$ for every $M\in \Mod_{R'}^{\sm}(G)$ and $\sF\in \Sh(X, R')$. The functor \eqref{homchangering2} has the right adjoint 
\begin{equation}
\label{homchangering3}
\Sh(X,R) \lra \Sh(X,R'),\,\, \sG\mapsto \uHom_R(R'_X, \sG)
\end{equation}
and the left adjoint  
\begin{equation}
\label{homchangering4}
\Sh(X,R) \lra \Sh(X,R'),\,\, \sG\mapsto \sG_{R'} :=\sG\otimes_R R'_X.
\end{equation}

\begin{lemma}
\label{lemma:basechangefree}
Let $\varphi:R\to R'$ be a ringhomomorphism, let $M$ be a discrete $R'[G]$-module and let $\sF\in \Sh(X, R)$. Assume that $R'$ is free and of finite rank as an $R$-module and that $\Hom_R(R', R)$ is free of rank 1 as an $R'$-module. Then the choice of an $R'$-module isomorphism $\psi: R'\cong \Hom_R(R', R)$ induces isomorphisms
\begin{equation*}
\label{drmqch1}
H_R^n(X; M_{\varphi}, \sF) \, \cong \, H_{R'}^n(X; M, \sF_{R'})
\end{equation*}
for every $n\ge 0$ and $\sF\in \Sh(X,R)$. 
\end{lemma}

\begin{proof} The isomorphism $\psi$ induces an isomorphism between the functors \eqref{homchangering3} and \eqref{homchangering4}. Hence we have 
\[
\Hom_{\Sh(X, R)}((M_{\varphi})_{X, G}, \wcdot ) \, \cong \, \Hom_{\Sh(X, R)}((M_{X, G})_{\varphi}, \wcdot ) \cong \, \Hom_{\Sh(X, R)}(M_{X, G},  R'\otimes_R \wcdot ).
\]
so the right derived functors are isomorphic as well.
\end{proof}

Now assume that $\varphi: R\to R'$ is an arbitrary ringhomomorphism. Let $M$ be a discrete $R[G]$-module and put $M_{R'}\colon = M\otimes_{R'} R$. Because of $(M_{R'})_{X, G}\cong (M_{X,G})_{R'}$ the functor $\sF\mapsto \Hom_{\Sh(X, R)}(M_{X, G}, \sF_{\varphi})$ is isomorphic to $\sF\mapsto \Hom_{\Sh(X, R)}((M_{R'})_{X, G}, \sF)$. Note that $\sF \mapsto H_R^{\bu}(X; \break M, \sF_{\varphi})$ is a $\delta$-functor whereas $\sF \mapsto H_{R'}^{\bu}(X; M_{R'}, \sF)$ is a universal a $\delta$-functor. Hence there are canonical homomorphisms 
\begin{equation}
\label{drmqch2}
H_{R'}^n(X; M_{R'}, \sF) \, \lra \, H_R^n(X; M, \sF_{\varphi})
\end{equation}
for every $n\ge 0$. 

\begin{prop}
\label{prop:drmqch}
Assume that $X$ is exact and that $M$ is projective as an $R$-module. Then the map \eqref{drmqch2} is an isomorphism for every $n\ge 0$ and $\sF\in \Sh(X, R')$. 
\end{prop}

We need some preparation. Since $\Hom_{R[G]}(M, N_{\varphi})\,\cong \,\Hom_{R'[G]}(M_{R'}, N)$ for $M\in \Mod_R^{\sm}(G)$ and $N\in \Mod_{R'}^{\sm}(G)$ there exists 
canonical homomorphisms 
\begin{equation}
\label{drmqch3}
\Ext^n_{R', G}(M_{R'}, N) \, \lra \, \Ext_{R, G}^n(M, N_{\varphi})
\end{equation}
for all $n\ge 0$. 

\begin{lemma}
\label{lemma:extind}
Assume that $M\in \Mod_R^{\sm}(G)$ is projective as an $R$-module. 
Then the homomorphisms \eqref{drmqch3} are isomorphisms.
\end{lemma}

\begin{proof} Since $\Ind_1^G: \Mod_R \to \Mod_R^{\sm}(G)$ is exact and preserves injectives we have  
\begin{equation*}
\label{drmqch3a}
\Ext_{R, G}^{\bu}(M, \Ind_1^G Q) \, \cong \, \Ext_{R}^{\bu}(M, Q) 
\end{equation*}
for all $Q\in \Mod_R$. Because $M$ is projective as an $R$-module this implies $\Ext_{R, G}^n(M, \Ind_1^G Q)=0$ for all $n\ge 1$. In particular we obtain 
\begin{equation}
\label{drmqch3b}
\Ext_{R, G}^n(M, (\Ind_1^G Q')_{\varphi}) \, = \, \Ext_{R, G}^n(M, \Ind_1^G (Q'_{\varphi}))\, =\, 0
\end{equation}
for every $Q'\in \Mod_{R'}$ and $n\ge 1$. 

The functor $\Mod_{R'}^{\sm}(G)\to \Mod_R, N \mapsto \Hom_{R[G]}(M, N_{\varphi})$ factors in the form
\begin{equation}
\label{homchangering6}
\begin{CD}
\Mod_{R'}^{\sm}(G) @ > N\mapsto N_{\varphi} >> \Mod_R^{\sm}(G)@>  \Hom_{R[G]}(M, \wcdot) >> \Mod_R
\end{CD}
\end{equation}
The first functor is exact and maps injective objects to $\Hom_{R[G]}(M, \wcdot)$-acyclic objects. 
Indeed, if $I\in \Mod_{R'}^{\sm}(G)$ is injective then it is a direct summand of $\Ind_1^G I$ hence by \eqref{drmqch3b} we get $\Ext_{R, G}^n(M, I_{\varphi}) =  0$
for $n\ge 1$. Therefore there exists a degenerating Grothendieck spectral sequence corresponding to the composite of functors \eqref{homchangering6}. The assertion follows.
\end{proof}

\begin{proof}[Proof of Proposition~\ref{prop:drmqch}] It suffices to show that for an injective object $\sI\in \Sh(X, R')$, the object $\sI_{\varphi}\in \Sh(X,R)$ is $\Hom_{\Sh(X,R)}(M_{X,G}, \wcdot)$-acyclic. The functor $\Sh(X, R') \to \Mod_R, \sF\mapsto \Hom_{\Sh(X,R)}(M_{X,G}, \sF_{\varphi})$ factors in the form 
\begin{equation}
\label{homchangering7}
\begin{CD}
\Sh(X, R') @ > H^0(X_{\infty}, \wcdot) >> \Mod_{R'}^{\sm}(G)@>  N\mapsto N_{\varphi} >> \Mod_R^{\sm}(G)@>  \Hom_{R[G]}(M, \wcdot) >> \Mod_R.
\end{CD}
\end{equation}
The first has an exact left adjoint so it preserves injectives. As shown in the proof of Lemma \ref{lemma:extind} the second functor maps injective objects to $\Hom_{R[G]}(M, \wcdot)$-acyclic objects. If we denote the composition of the first two functors in \eqref{homchangering7} by $T:\Sh(X, R')\to \Mod_R^{\sm}(G)$
then there exists a Grothendieck spectral sequence 
\begin{equation*}
\label{exthss1a}
E_2^{rs} = \Ext_{R, G}^r(M, R^q T(\sF)) \, \Longrightarrow\, H_R^{r+s}(X; M, \sF_{\varphi}).
\end{equation*}
For the injective object $\sI\in \Sh(X, R')$ we have $R^s T(\sI)=H^s(X_{\infty}, \sI)_{\varphi} = 0$ for $s\ge 1$ hence 
\[
H_R^n(X; M, \sI_{\varphi})\, \cong \,  \Ext_{R, G}^n(M, \sI(X_{\infty})_{\varphi}) \, =\, 0
\]
for $n\ge 1$. Hence $\sF \mapsto H_R^{\bu}(X; M, \sF_{\varphi})$ is also a universal $\delta$-functor.\end{proof}

\paragraph{Change of the underlying site} Let $R$ be a ring, let $\cC$ and $\cD$ be two sites that satisfy hypothesis \ref{assum:sitegood} and let $u: \cC \to \cD$ be a continuous functor commuting with fibre products, equalizers and coproducts. Let $X$ be an $\sS$-object in $\cC$, let $u(X):= u\circ X$ be the induced $\sS$-object in $\cD$ and let $\wu: \cC/\colim X \to \cD/\colim u(X)$ be the functor \eqref{changesite2} induced by $u$. Note that the diagram of topoi
\begin{equation*}
\label{resjshriekd}
{\xymatrix@-0.5pc{\Sh(u(X), R)\ar[rr]^{(u^*, u_*)}\ar[ddr]^{\eqref{sspace5a}} && \ar[ddl]_{\eqref{sspace5a}}\Sh(X, R) \\
\\
& \Mod_R^{\sm}(G)
}}
\end{equation*}
commutes. Since $u^*$ is exact by Prop.\ \ref{prop:pullbackexact}, there exists a Grothendieck spectral sequence 
\begin{equation}
\label{leraychangeofsite}
E_2^{rs} = H^s(X_{\infty}, R^r u_* \sG) \, \Longrightarrow\, H^{r+s}(u(X)_{\infty}, \sG)
\end{equation}
for every $\sG\in \Sh(u(X), R)$. Also we get $u^*(M_{X, G}) = M_{u(X), G}$ for $M\in \Mod_R^{\sm}(G)$. Thus the functor $u^*$ induces canonical homomorphisms 
\begin{equation}
\label{pullbackcohomss}
H_R^{\bu}(X; M, \sF) \lra H_R^{\bu}(u(X); M, u^*(\sF))
\end{equation}
for every $\sF\in \Sh(X, R)$. Furthermore there exists canonical homomorphisms 
\begin{equation}
\label{pullbackcohomss2}
H^{\bu}(X_{\infty}, \sF) \lra H^{\bu}(u(X)_{\infty}, u^*(\sF))
\end{equation}
defined as the composite of the map $H^{\bu}(X_{\infty}, \sF)\to H^{\bu}(X_{\infty}, u_*u^*(\sF))$ (induced by the unit of adjuction) with an edge morphism of the spectral sequence \eqref{leraychangeofsite}. We have 

\begin{prop}
\label{prop:hsetsing}
Assume that $X$ is exact. Then there is a morphism of spectral sequences 
\begin{eqnarray}
\label{comexthsspb}
&& \left( \,\Ext_{R, G}^r(M, H^s(X_{\infty}, \sF)) \, \Longrightarrow\, H_R^{r+s}(X; M, \sF) \, \right) \, \lra \, \\
&& \hspace{3cm} \left(\, \Ext_{R, G}^r(M, H^s(u(X)_{\infty}, u^*(\sF))) \, \Longrightarrow\, H_R^{r+s}(u(X); M, u^*(\sF))\, \right).
\nonumber
\end{eqnarray}
\end{prop}

{\em Proof.} Let $0\to \sF \to \sI^{\bu}$ be an injective resolution of $\sF$ in $\Sh(X, R)$. By Prop.\ \ref{prop:pullbackexact} applying $u^*$ yields a resolution $0\to u^*(\sF) \to u^*(\sI^{\bu})$ of $u^*(\sF)$. Let $0\to u^*(\sF) \to \sJ^{\bu}$ be an injective resolution of $u^*(\sF)$ in $\Sh(X, R)$ and let 
$\alpha:  u^*(\sI^{\bu}) \to \sJ^{\bu}$ be a morphism such that 
\[
\begin{CD}
0 @>>> u^*(\sF) @>>> u^*(\sI^{\bu})\\
@.@VV \id V @VV \alpha V\\
0 @>>> u^*(\sF) @>>> \sJ^{\bu}
\end{CD}
\]
commutes. Consider the homomorphism of complexes of discrete $G$-modules 
\begin{equation*}
\label{analytic6f}
\begin{CD} 
\sI^{\bu}(X_{\infty}) @> \eqref{pullbackcohomss2} >> u^*(\sI^{\bu})(u(X)_{\infty}) @> \alpha >> \sJ^{\bu}(u(X)_{\infty}).
\end{CD} 
\end{equation*}
It induces a morphism between $G$-hypercohomology spectral sequences
\begin{eqnarray*}
&&  \left(\, \Ext_{R, G}^r(M, H^s(\sI^{\bu}(X_{\infty})))\, \Longrightarrow\, \Ext_{R, G}^{r+s}(M, \sI^{\bu}(X_{\infty})) \,\right) \, \lra \, \\
&& \hspace{3cm} \left( \, \Ext_{R, G}^r(M, H^s(\sJ^{\bu}(u(X)_{\infty}))) \, \Longrightarrow\, \Ext_{R, G}^{r+s}(M, \sJ^{\bu}(u(X)_{\infty}))\, \right).
\end{eqnarray*}
The first spectral sequence can be identified with the source of \eqref{comexthsspb}
 and the second with the target. 
\enddemo.

\subsection{$\sS$-spaces and $\sS$-schemes}
\label{subsection:ssspacescheme}

As in the last section $G$ denotes a locally profinite group, $\sK$ a subset of the set of open subgroups of $G$ satisfying hypothesis \ref{assum:kgood} and $\sS$ the site $\sS=\sS_{G,\sK}$. In this section we consider $\sS$-objects in certain sites of topological spaces and schemes. We first consider schemes. 

Let $B$ be a reduced noetherian base scheme. By $(\Sm/B)_{\et}$ we denote the following site. The objects of $(\Sm/B)_{\et}$ are smooth $B$-schemes that are locally of finite type. The morphisms are {\'e}tale morphisms of $B$-schemes and the coverings are {\'e}tale coverings. The empty scheme is the initial object of $(\Sm/B)_{\et}$. For $X\in (\Sm/B)_{\et}$ the localization $(\Sm/B)_{\et}/X$ is the small {\'e}tale site of $X$. Clearly, the category $(\Sm/B)_{\et}$ has fibre products and coproducts. We also have 

\begin{lemma}
\label{lemma:etaleequal}
The category $(\Sm/B)_{\et}$ has equalizers.
\end{lemma}

\begin{proof} Let $f,g: X_1 \to X_2$ be two morphism in $(\Sm/B)_{\et}$. It is sufficient to consider the case where $X_1$ is connected. Recall that the equalizer $\iota': Z'\to X_1$ of the pair $(f,g)$ exists in the category of $B$-schemes. Namely, it is the base change of diagonal embedding $\Delta_{X_2/B}: X_2\to X_2\times_B X_2$ with respect to the morphism $X_1\to X_2\times_B X_2$ with components $f$ and $g$. Since $\Delta_{X_2/B}$ is an immersion the same holds for $\iota'$. Therefore $\io'$ can be factored in the form $Z'\stackrel{j_1}{\lra} U \stackrel{j_2}{\lra} X_1$ with $j_1$ a closed and $j_2$ an open immersion. Now if there exists a morphism $h:Y \to X_1$ in $(\Sm/B)_{\et}$ with $Y$ non-empty, then $\io'$ is an open immersion. Indeed, since the image of $h$ is open and non-empty it is dense in $X_1$, hence $Z'$ is dense in $X_1$, whence in $U$, so $j_1(Z') = U$. It follows that  $j_1: Z'\to U$ is an isomorphism because $U$ is reduced. 
Thus the equalizer of $(f, g)$ in $(\Sm/B)_{\et}$ is $\iota': Z'\to X_1$ if it is an open immersion and $\emptyset \hra X_1$ otherwise. 
\end{proof}

\begin{df}
\label{df:Gscheme}
Let $B$ be a reduced noetherian scheme. We call an $\sS$-object $X$ in $(\Sm/B)_{\et}$ an $\sS$-scheme (over $B$). For a ring $R$ the category of $R$-sheaves on $X$ will be denoted by $\Sh(X_{\et}, R)$. 
\end{df}

We call an $\sS$-scheme $X$ over $B$ to be of finite type, if $X_K$ is a $B$-scheme of finite type for every $K\in \sK$. 
An $R$-sheaf $\sF\in \Sh(X_{\et}, R)$ will be called constructible if $\sF_K$ is constructible for every $K\in \sK$. 

\begin{remarks}
\label{remarks:trace} \rm (a) Assume that $f: X\to Y$ is a finite morphism of $\sS$-schemes over a base $B$, i.e.\ we assume that $f_S: X_S \to Y_S$ is finite (and {\'e}tale) for every $S\in \sS$. Then $f_*: \Sh(X_{\et}, R) \to \Sh(Y_{\et}, R)$ is exact and also left adjoint to $f^*$. By Prop.\  \ref{prop:extcohom} (c) the canonical map 
\begin{equation}
\label{trace1}
H_R^n(Y; M, f_*\sG)\lra H_R^n(X; M, \sG)
\end{equation}
is an isomorphism for every $n \ge 0$, $M\in \Mod_R^{\sm}(G)$ and $\sG\in \Sh(X_{\et}, R)$. 

For $\sF\in \Sh(Y_{\et}, R)$ and $n\ge 0$ we define a trace map
\begin{equation}
\label{trace2}
\begin{CD}
H_R^n(X; M, f^*\sF) @> \cong >> H_R^n(Y; M, f_* f^* \sF)@> \tr >> H_R^n(Y; M, \sF)
\end{CD}
\end{equation}
where the first map is the inverse of \eqref{trace1} and the second is induced by the counit of adjunction. 
\medskip 

\noi (b) Let $u: B'\to B$ be a morphism between reduced noetherian schemes. The functor 
\begin{equation}
\label{basechangescheme}
(\Sm/B)_{\et}\to (\Sm/B')_{\et}, \quad X\to X_{B'} = X\times_B B'
\end{equation} 
is continuous and commutes with fibre products and coproducts. 

Now assume that $u$ is a faithfully flat and quasi-compact morphism. Then \eqref{basechangescheme} commutes also with equalizers. This follows from the proof of Lemma \ref{lemma:etaleequal} together with (\cite{ega4}, Prop.\ 2.7.1). 
Let $X$ is an $\sS$-scheme over $B$. We define the base change $X_{B'}$ of $X$ with respect to $u:B'\to B$ as the composite of $X:\sS\to (\Sm/B)_{\et}$ with \eqref{basechangescheme}. 
We obtain an induced morphism of topoi (see \eqref{changesite3})
\begin{equation*}
\label{basechangescheme2}
(u^*, u_*): \Sh((X_{B'})_{\et}, R)\lra \Sh(X_{\et}, R).
\end{equation*} 
We denote the image of $\sF\in \Sh(X, R)$ under $u^*$ also by $\sF\mapsto \sF_{B'}$. 
\medskip 

\noi (c) Let $k'/k$ be an extension of separably closed fields and let $A$ is an Artinian local ring with finite residue field of characteristic $\ne \charac \bar{k}$. Let $X$ be an exact $\sS$-scheme over $\Spec k$, let $\sF\in \Sh(X_{\et}, A)$ and let 
$M\in \Mod_A^{\sm}(G)$. We put $X_{k'}= X_{\Spec k'}$ and $\sF_{k'}= \sF_{\Spec k'}$. Since for any $K\in \sK$ the canonical homomorphisms $H^{\bu}(X_K, \sF_K)\to H^{\bu}((X_K)_{k'}, (\sF_K)_{k'})$ are isomorphisms by the smooth base change theorem, Prop.\ \ref{prop:hs} (a) and Prop.\ \ref{prop:hsetsing} imply that the homomorphism (see \eqref{pullbackcohomss})
\begin{equation*}
\label{pullbackfields}
H_A^{\bu}(X; M, \sF) \lra H_A^{\bu}(X_{k'}; M, \sF_{k'})
\end{equation*}
is an isomorphism as well. 
\enddemo
\end{remarks}

Now we are going to define the notion of an $\sS$-space. We denote by $\spaces_{\cl}$ the following site. Objects of $\spaces_{\cl}$ are arbitrary topological spaces and the morphism in $\spaces_{\cl}$ are local homeomorphisms, i.e.\ continuous maps $f:X\to Y$ such that every $x\in X$ has an open neighbourhood which is mapped homeomorphically 
under $f$ onto a neighbourhood of $f(x)$. Note that the morphism in $\spaces_{\cl}$ are open maps. Coverings in $\spaces_{\cl}$ are families of morphisms $\{f_i: X_i \to X\}$ such that $\bigcup_i f_i(X_i)= X$. The category $\spaces_{\cl}$ has fibre products and coproducts. Also equalizers exists in $\spaces_{\cl}$, Indeed, if $f, g: X\to Y$ are two morphisms in $\spaces_{\cl}$ then the interior (i.e.\ the maximal open subset) of $\{x\in X\mid f(x) = g(x)\}$ is the equalizer of $f$ and $g$. 

\begin{df}
\label{df:Gspace}
An $\sS$-space is an $\sS$-object in $\spaces_{\cl}$. 
\end{df}

We now state and prove a comparison theorem between the cohomology of an $\sS$-scheme over $\bC$ and the associated $\sS$-space. Recall that for a locally algebraic $\bC$-scheme $X$, the set of $\bC$-points $X(\bC)$ carries a canonical topology. We denote $X(\bC)$ with this toplogy by $X^{\an}$. We obtain in particular a functor 
\begin{equation}
\label{analytic}
\Sm/\bC \lra \Man_{\bC},\,\, X\mapsto X^{\an}, \quad f\mapsto f^{\an}
\end{equation}
from the category of smooth, locally algebraic $\bC$-schemes to the category of complex manifolds with holomorphic maps as morphisms. If $f:X\to Y$ is an {\'e}tale morphism between smooth schemes over $\bC$ then $f^{\an}$ is a local homeomorphism. So \eqref{analytic} induces a continuous functor between sites
\begin{equation}
\label{analytic2}
\ep: (\Sm/\bC)_{\et} \lra \spaces_{\cl},\,\, X\mapsto X^{\an}
\end{equation}
It commutes with fibre products, coproducts and equalizers. 

Let $X$ be an $\sS$-scheme over $\bC$ and let $X^{\an}:= \ep(X)$ be the associated $\sS$-space. Then \eqref{analytic2} induces a morphism of topoi $(\ep^*, \ep_*): \Sh(X^{\an}, R) \lra \Sh(X_{\et}, R)$. We denote the image of 
$\sF\in \Sh(X_{\et}, R)$ under $\ep^*$ by $\sF^{\an}$. By \eqref{pullbackcohomss} there exists canonical comparison homomorphisms 
\begin{equation}
\label{analytic3}
H_R^{\bu}(X; M, \sF) \lra H_R^{\bu}(X^{\an}; M, \sF^{\an})
\end{equation}
for every $M\in \Mod_R^{\sm}(G)$ and $\sF\in \Sh(X^{\an}, R)$. Note that $X$ is exact if and only if $X^{\an}$ is exact. 

\begin{prop}
\label{prop:etancomp}
Assume that $X$ is exact and that $X_K$ is of finite type for every $K\in \sK$. Let $A$ be an Artinian local ring with finite residue field, let $M\in \Mod_A^{\sm}(G)$ and let $\sF\in \Sh(X_{\et}, A)$ be constructible. Then the comparison homomorphisms \eqref{analytic3} are isomorphisms.
\end{prop} 

\begin{proof} By Prop.\ \ref{prop:hsetsing} there exists a morphism of spectral sequences 
\begin{eqnarray}
\label{analytichss}
&& \left( \,\Ext_{A, G}^r(M, H^s(X_{\infty}, \sF)) \, \Longrightarrow\, H_A^{r+s}(X; M, \sF) \, \right) \, \lra \, \\
&& \hspace{3cm} \left(\, \Ext_{A, G}^r(M, H^s(X^{\an}_{\infty}, \sF)) \, \Longrightarrow\, H_A^{r+s}(X^{\an}; M, u^*(\sF))\, \right).
\nonumber
\end{eqnarray}
By (\cite{sga4.5}, Arcata, (3.5.1)) the assumptions imply that the canonical maps $H^{\bu}(X_K, \sF_K) \to H^{\bu}(X^{\an}_K, \sF^{\an}_K)$ 
are isomorphisms for every $K\in \sK$. Hence by Prop.\ \ref{prop:hs} (a) the maps between $E_2$-terms of \eqref{analytichss} are isomorphisms and \eqref{analytichss} is an isomorphism of spectral sequences. \end{proof}

\paragraph{The site $\spaces_{\bg}$} In order to compute the cohomology of the boundary of the Borel-Serre compactification of Hilbert modular $\sS$-varieties in section \ref{subsection:borelserre} we need to consider a second site of topological spaces which will be denoted by $\spaces_{\bg}$. It is defined as follows. Objects of $\spaces_{\bg}$ are arbitrary topological spaces and morphisms are open and continuous maps. Coverings in $\spaces_{\bg}$ are again families of local homeomorphisms $\{f_i: X_i \to X\}$ such that $\bigcup_i f_i(X_i)= X$. The category $\spaces_{\bg}$ has fibre products. This follows from the fact that the base change of an open and continuous map is again open and continuous, i.e.\ if $f:X\to Z$, $g:Y\to Z$ are continuous maps between topological spaces such that $g$ is open then the projection $\pr_X: X\times_Z Y=\{(x,y)\in X\times Y\mid f(x)=g(y)\} \to X$ is open as well. Indeed, if $U\subseteq X$ and $V\subseteq Y$ are open then $\pr(U\times V\cap X\times_Z Y) = U\cap f^{-1}(g(V))$ is open in $X$. Also the category $\spaces_{\bg}$ has coproducts and equalizers (the equalizer of two morphisms in $f, g: X\to Y$ in $\spaces_{\bg}$ is again the interior of $\{x\in X\mid f(x) = g(x)\}$). 

The inclusion functor $\io: \spaces_{\cl}\to \spaces_{\bg}$ is continuous and commutes with fibre products, equalizers and coproducts. If $X$ is a topological space then $\io$ induces a 
continuous functor between localizations. It will be denoted by $\io: \spaces_{\cl}/X\to \spaces_{\bg}/X$ as well (note that it is continuous, cocontinuous, fully faithful and commutes with fibre products). We obtain a morphisms of topoi 
\begin{equation*} 
\label{bigsmallsheaf}
(\ext, \res): = (\io_s, \io^s): \Sh(X/\spaces_{\bg}, R)\lra \Sh(X/\spaces_{\cl}, R). 
\end{equation*} 
Recall that $\Sh(X/\spaces_{\cl}, R)$ can be identified with the category of usual sheaves of $R$-modules on $X$. By (see \cite{stacksproject}, Tag 00XU) the functor $\res$ is exact, $\ext$ is fully faithful and we have $\res\circ \ext \cong \id$. 

Let $f: X\to Y$ be an arbitrary continuous map between topological spaces. If $U\to Y$ is an object in $Y/\spaces_{\cl}$ (resp.\ $Y/\spaces_{\bg}$) then $X\times_Y U\to X$ lies in $X/\spaces_{\cl}$ (resp.\ $X/\spaces_{\bg}$) hence base change with respect to $f$ yields continuous functors 
$v: Y/\spaces_{\cl}\to X/\spaces_{\cl}$ and $v: Y/\spaces_{\bg}\to X/\spaces_{\bg}$. For $?\in \{\cl, \bg\}$ the associated morphism of topoi  will be denoted by
\begin{equation*} 
\label{bigsmallsheaf2}
(f^*, f_*)= (v_s, v^s): \Sh(X/\spaces_{?}, R) \lra \Sh(Y/\spaces_{?}, R).
\end{equation*}
By (\cite{stacksproject}, Tag 03CY) the derived functors of $f_*$ commute with $\res$, i.e.\ we have 
\begin{equation*} 
\label{bigsmallsheaf3}
(R^n f_*)(\res(\cF))\, \cong \, \res(R^n f_*(\cF))
\end{equation*}
for every $\cF\in \Sh(X/\spaces_{\bg}, R)$ and $n\ge 0$. If $f:X\to Y$ is open then $f$ induces also a continuous and cocontinuous functor $\wf: \spaces_{\bg}/X\to \spaces_{\bg}/Y, (\phi: U\to X)\mapsto (f\circ \phi: U\to Y)$ that is left adjoint to the base change functor. It follows (see \cite{stacksproject}, Tag 00XY) 
\begin{equation*} 
\label{bigsmallsheaf4}
(f^*, f_*) = (\wf^s, {_s\wf}):  \Sh(X/\spaces_{\bg}, R)\lra \Sh(Y/\spaces_{\bg}, R).
\end{equation*}
Now assume that $X$ is an $\sS$-space. 
Composing $X$ with $\io$ yields an $\sS$-object $\io(X) =\io\circ X$ in $\spaces_{\bg}$ and we define the category of {\it big $R$-sheaves} on $X$ by
\begin{equation*} 
\label{bigsmallsheaf5}
\Sh_{\bg}(X, R) :=  \Sh(\spaces_{\bg}/\colim \io(X), R).
\end{equation*}
Again the functor $\wio: X/\spaces_{\cl}\to X/\spaces_{\bg}$ induced by $\io$ (see \eqref{changesite2}) is continuous, cocontinuous, fully faithful and commutes with fibre products. Thus for the associated morphism of topoi 
\begin{equation*} 
\label{bigsmallsheaf6}
(\ext, \res) := : (\wio_s, \wio^s): \Sh_{\bg}(X, R)\lra \Sh(X, R) 
\end{equation*} 
we have: $\res\circ \ext\cong \id$, $\res$ is exact and $\ext$ is fully faithful. 
The proof of the following comparison result is straight forward.

\begin{prop}
\label{prop:bigsmallcohom}
Let $X$ be an $\sS$-space, let $M\in \Mod_R^{\sm}(G)$ and let $\sF\in \Sh_{\bg}(X, R)$. 
\medskip

\noi (a) We have $H^n(X_{\infty}, \res(\sF)) \cong H^n(\io(X)_{\infty}, \sF)$ for every $n\ge 0$. 
\medskip

\noi (b) If $X$ is exact then $H_R^n(\io(X); M, \sF)\cong H_R^n(X; M, \res(\sF))$ for every $n\ge 0$. 
\end{prop}

\subsection{Galois operation on $H_R^{\bu}(X_{\barF}; M, \sF_{\barF})$}
\label{subsection:galcohomology}

As in the last section $G$ denotes a locally profinite group, $\sK$ a subset of the set of open subgroups of $G$ satisfying hypothesis \ref{assum:kgood} and $\sS$ the site $\sS_{G,\sK}$. In this section we consider an $\sS$-scheme $X$ over a field $F$ (i.e.\ an $\sS$-scheme over $\Spec F$). For simplicity we assume that $F$ is perfect. We fix an algebraic closure $\barF$ of $F$ and let $\fG = \fG_F= \Gal(\barF/F)$ denote the absolute Galois group of $F$. Similarly to usual {\'e}tale cohomology our aim is to show that for $M\in \Mod_R^{\sm}(G)$ and $\sF\in \Sh(X_{\et}, R)$ the cohomology groups $H_R^{\bu}(X_{\barF}; M, \sF_{\barF})$ carry a canonical $\fG_F$-action of $\fG_F$. 

Recall that the small {\'e}tale site $(\Spec F)_{\et}$ can be identified with the site $\sS_{\fG}$ of discrete left $\fG$-sets. We let $\wsK$ be the set of open subgroups of $G\times \fG$ of the form $K\times \fU$ where $K\in \sK$ and 
$\fU$ is an open subgroup of $\fG$. We put $\wsS = \sS_{G\times \fG, \wsK}$. By Remark \ref{remark:Sspacealt} there exists a unique $\wsS$-scheme $X^{\gal}$ characterized by $(X^{\gal})_{S\times S'}= X_S\otimes_F S'$ for every $S\in \sS$ and $S'\in \sS_{\fG}=(\Spec F)_{\et}$. Then $X$ is the push-out of $X^{\gal}$ with respect to the projection $\pi: G\times \fG\to G$, so by Lemma \ref{lemma:pushoutopos} the functor 
\begin{equation*}
\label{galsheaf}
\Sh(X^{\gal}_{\et}, R)\lra \Sh(X_{\et}, R),\quad \sF \mapsto \sF_{\pi}
\end{equation*}
(see \eqref{phipushou3}) is an equivalence of categories. The quasi-inverse will be denoted by 
\begin{equation}
\label{galsheaf2}
\Sh(X_{\et}, R)\lra \Sh(X^{\gal}_{\et}, R), \quad \sF \mapsto \sF^{\gal}.
\end{equation}
If $\sF\in \Sh(X_{\et}, R)$ then the layers of $\sF^{\gal}$ are given by
\[
\sF_{K, F'} := \,\sF^{\gal}_{K\times \fG_{F'}} \, =\,(\sF_K)_{F'}\, =\, (\sF_{F'})_K
\]
where $K\in \sK$ and where $F'/F$ is a finite subextension of $\barF/F$ (here for an {\'e}tale $R$-sheaf $\cF$ on a scheme $Y\in \Sm/F$ we denote by $\cF_{F'}$ the pull-back of $\cF$ with respect $Y_{F'}\to Y$). 

We consider the composite of the functor \eqref{galsheaf2} and the functor $H^0((X^{\gal})_{\infty}, \wcdot)$, i.e.\ we consider the functor 
\begin{equation}
\label{galsheaf2a}
\Sh(X_{\et}, R) \lra \Mod_R^{\sm}(G\times \fG), \,\, \sF \mapsto H^0((X^{\gal})_{\infty}, \sF^{\gal}).
\end{equation}
Its $n$-th right derived functor is therefore
\begin{equation*}
\label{galsheaf2b}
\Sh(X_{\et}, R) \lra \Mod_R^{\sm}(G\times \fG), \,\, \sF \mapsto H^n((X^{\gal})_{\infty}, \sF^{\gal}).
\end{equation*}

\begin{lemma}
\label{lemma:galXinfty}
Let $\sF\in \Sh(X_{\et}, R)$. As $R[G]$-modules we have 
\[
H^n((X^{\gal})_{\infty}, \sF^{\gal})\, \cong \, H^n((X_{\barF})_{\infty}, \sF_{\barF})
\]
for every $n\ge 0$. In particular the discrete $R[G]$-module $H^n((X_{\barF})_{\infty}, \sF_{\barF})$ carries additionally a discrete $\fG$-action (commuting with the $G$-action).
\end{lemma} 

\begin{proof} Using Prop.\ \ref{prop:hs} (a) we obtain
\begin{eqnarray*}
H^n((X^{\gal})_{\infty}, \sF^{\gal}) & = & \dlim_{K, F'} H^n(X_{K, F'}, \sF_{K, F'})\, =\, \dlim_{K} \dlim_{F'} H^n((X_K)_{F'}, (\sF_K)_{F'})\\
 &= &\dlim_{K} H^n((X_K)_{\barF}, (\sF_K)_{\barF}) \, =\, H^n((X_{\barF})_{\infty}, \sF_{\barF}).
\end{eqnarray*}
where $F'$ ranges over all finite extensions of $F$ contained in $\barF$ and $K$ ranges over all subgroups in $\sK$.
 \end{proof}

Now assume that $X$ is exact. Note that $X^{\gal}$ and $X_{\barF}$ are exact as well. We will show that the spectral sequence in Prop.\ \ref{prop:extcohom} (b) for $X_{\barF}$ and an {\'e}tale $R$-sheaf that is the pull-back of an $R$-sheaf on $X$, is Galois-equivariant. To begin with, for $M\in \Mod_R^{\sm}(G)$ we consider the left exact functor
\begin{equation}
\label{galext}
\Sh(X_{\et}, R) \lra \Mod_{R[\fG]}, \quad \sF\mapsto \Hom_{R[G]}(M, H^0((X^{\gal})_{\infty}, \sF^{\gal})).
\end{equation}
Here the $\fG$-action on the image of $\sF\in\Sh(X_{\et}, R)$ under \eqref{galext} is induced by the $\fG$-action on $H^0((X^{\gal})_{\infty}, \sF^{\gal})$.

\begin{prop}
\label{prop:galactionhs}
Assume that $X$ is exact and let $\sF\in \Sh(X_{\et}, R)$ and $M\in \Mod_R^{\sm}(G)$. 
\medskip

\noi (a) The image of $\sF$ under the $n$-th right derived functor of \eqref{galext} is isomorphic to the $R$-module $H^n((X_{\barF})_{\infty}; M, \sF_{\barF})$. In particular the latter is equipped with a canonical $\fG$-action.\medskip

\noi (b) There exists a spectral sequence of $R[\fG]$-modules
\begin{equation*}
\label{galexthss}
E_2^{rs} = \Ext_{R, G}^r(M, H^s((X_{\barF})_{\infty}, \sF_{\barF})) \, \Longrightarrow\, H_R^{r+s}(X_{\barF}; M, \sF_{\barF}).
\end{equation*}
Here the $\fG$-action on the $E^2$-terms is induced by the $\fG$-action on $H^{\bu}((X_{\barF})_{\infty}, \sF_{\barF})$ provided for by Lemma \ref{lemma:galXinfty}.
\end{prop}

We need the following 

\begin{lemma}
\label{lemma:extgfg}
The following diagram commutes 
\[
\begin{CD}
D^+(\Mod_R^{\sm}(G\times \fG)) @> R^+\Hom_{R[G]}(M, \wcdot ) >> D^+(\Mod_{R[\fG]})\\
@VVV @VVV\\
D^+(\Mod_R^{\sm}(G)) @> R^+\Hom_{R[G]}(M, \wcdot )>> D^+(\Mod_R).
\end{CD}
\]
Here the vertical maps are induced by the forgetful functors.
\end{lemma}

Here as usual for a left exact functor $\Phi:\sA\to \sB$ between abelian category with enough injectives $R^+\Phi: D^+(\sA) \to D^+(\sB)$ denotes the (total) right derived functor of $\Phi$. 

\begin{proof} This follows immediately from the fact that the two forgetful functors $\Mod_R^{\sm}(G\times \fG)\to \Mod_R^{\sm}(G)$ and $\Mod_R(\fG)\to \Mod_R$ are exact and admit an exact left adjoint. 
\end{proof}

\begin{proof}[Proof of Prop.\ \ref{prop:galactionhs}] Let $\Phi$ denote the functor \eqref{galext}. We factor it as 
\begin{equation*}
\label{galext2}
\begin{CD}
\Sh(X_{\et}, R)@> \eqref{galsheaf2a} >> \Mod_R^{\sm}(G\times \fG) @> \Hom_{R[G]}(M, \wcdot)>> \Mod_{R[\fG]}.
\end{CD}
\end{equation*}
So there exists a Grothendieck spectral sequence 
(of $R[\fG]$-modules)
\begin{equation}
\label{galext3}
E_2^{rs} = \Ext_{R, G}^r(M, H^s((X^{\gal})_{\infty}, \sF^{\gal})) \, \Longrightarrow\, E^{r+s}=(R^{r+s}\Phi)(\sF).
\end{equation}
By Lemma \ref{lemma:galXinfty} we have 
\begin{eqnarray*}
\Hom_{R[G]}(M, H^0((X^{\gal})_{\infty}, \sF^{\gal})) & \cong & \Hom_{R[G]}(M, \sF_{\barF}((X_{\barF})_{\infty}) \, \cong \, \Hom_{\Sh(X_{\barF}, R)}(M_{X_{\barF}, G}, \sF_{\barF})\\
& = & H^0(X_{\barF}; M, \sF_{\barF})
\end{eqnarray*}
for every $\sF\in \Sh(X_{\et}, R)$. Let $0\to \sF \to \sI^{\bu}$ be an injective resolution of $\sF$ in $\Sh(X_{\et}, R)$ and let $0\to \sF_{\barF} \to \sJ^{\bu}$ be an injective resolution of $\sF_{\barF}\in \Sh((X_{\barF})_{\et}, R)$. 
Since $0\to \sF_{\barF} \to \sI_{\barF}^{\bu}$ is a resolution there exists a morphism 
$\alpha: \sI_{\barF}^{\bu} \to \sJ^{\bu}$ such that 
\[
\begin{CD}
0 @>>> \sF_{\barF}  @>>> \sI_{\barF}^{\bu}\\
@.@VV \id V @VV \alpha V\\
0 @>>> \sF_{\barF}  @>>> \sJ^{\bu}
\end{CD}
\]
commutes. We obtain homomorphisms of complexes of $R[G]$-modules 
\begin{equation*}
\label{galext3a}
\begin{CD} 
(\sI^{\gal})^{\bu})((X^{\gal})_{\infty})\, \cong \, \sI_{\barF}^{\bu}((X_{\barF})_{\infty})
@> \alpha_* >> (\sJ^{\bu})((X_{\barF})_{\infty}).
\end{CD} 
\end{equation*}
It induces a morphism between hypercohomology spectral sequences
\begin{eqnarray}
\label{galext4}
&&  \left(\, \Ext_{G, R}^r(M, H^s((\sI^{\gal})^{\bu}((X^{\gal})_{\infty})))\, \Longrightarrow\, \Ext_{R, G}^{r+s}(M,  (\sI^{\gal})^{\bu}((X^{\gal})_{\infty})) \,\right) \, \lra \, \\
&& \hspace{3cm} \left( \, \Ext_{R, G}^r(M, H^s(\sJ^{\bu}((X_{\barF})_{\infty})))\, \Longrightarrow\, \Ext_{R, G}^{r+s}(M, \sJ^{\bu}((X_{\barF})_{\infty}))\, \right).\nonumber
\end{eqnarray}
By Lemma \ref{lemma:extgfg} 
the first spectral sequence can be identified with \eqref{galext3} (so it is in particular $\fG$-equivariant). The second spectral sequence is the spectral sequence \eqref{exthss} (for 
$X_{\barF}$ and $\sF_{\barF}$). By Lemma \ref{lemma:galXinfty} the morphism of spectral sequences \eqref{galext4} is an isomorphism. This proves (a) and (b).
\end{proof}

\subsection{$\Gamma$-equivariant $\sS$-spaces}
\label{subsection:uniformizable}

As in the previous section $\sK$ denotes a subset of the set of open subgroups of a locally profinite group $G$ satisfying hypothesis \ref{assum:kgood} and $\sS=\sS_{G,\sK}$ the associated site of discrete left $G$-sets. We fix a group $\Gamma$ and a ring $R$. A (left) action of $\Gamma$ on an $\sS$-space $Y$ consists of an action $\Gamma\times Y_S\to Y_S, (\gamma, y)\mapsto \gamma \cdot y$ on $Y_S$ for every $S\in \sS$, such that $\rho: Y_{S_1} \to Y_{S_2}$ is $\Gamma$-equivariant for every morphism $\rho: S_1 \to S_2$ in $\sS$. Equivalently, a $\Gamma$-action on $Y$ is a group homomorphism $\Gamma \to \Aut(Y)$ where $\Aut(Y)$ denotes the group of automorphisms of $Y$. For an $\sS$-space $Y$ equipped with a $\Gamma$-action we consider the following

\begin{assum}
\label{assum:goodact}
(i) The space $Y_S$ is a locally compact Hausdorff space for every $S\in \sS$ with $S\ne \emptyset$.\\
(ii) The group $\Gamma$ acts properly discontinuously on $Y_S$ for every $S\in \sS$ with $S\ne \emptyset$.\footnote{Recall that a group $\Gamma$ acts properly discontinuously on a locally compact Hausdorff space $Y$ if every point $y\in Y$ has a neighbourhood $U$ such that $gU\cap U = \emptyset$ for every $g\in G$, $g\ne 1$.}
\\
(iii) Let $\rho_1, \rho_2: S\to S'\in \sS$ be two morphisms in $\sS$. Assume that $S$ is connected and that $\rho_1\ne \rho_2$. Then the interior of $\{y\in Y_{S}\mid \exists \gamma\in \Gamma : \, \rho_2(y) = \gamma \cdot \rho_1(y)\}$ is empty.
\end{assum}

We equip $\Gamma$ with the discrete topology. We put $\wsS= \sS_{G\times \Gamma, \wsK}$ where $\wsK$ is the of set of open subgroups of $G\times \Gamma$ of the form $K\times H$ with $K\in \sK$ and $H$ is any subgroup of $\Gamma$.  

\begin{lemma}
\label{label:equivariantsspace}
(a) Let $Y$ be an $\sS$-space equipped with a $\Gamma$-action so that \ref{assum:goodact} (i), (ii) holds. 
Then there exists a unique $\wsS$-space $\wY$ such that 
\begin{equation}
\label{equivisspace}
\wY_{S\times S'} \, =\, \Gamma\backslash (Y_S\times S')
\end{equation}
for every $S\in \sS$ and every left $\Gamma$-set $S'$. 
\medskip

\noi (b) If moreover $Y$ is exact and \ref{assum:goodact} (iii) holds then $\wY$ is exact. 
\end{lemma}

\begin{proof} (a) We denote by $\wsS^c$ the full subcategory of $\wsS$ whose objects are $G/K\times \Gamma/H$ with $K\times H\in \wsK$. Consider the functor 
\begin{equation}
\label{equivisspace2}
\wsS^c\to \cC, \quad (G/K\times \Gamma/H) \mapsto \Gamma\backslash (Y_K\times \Gamma/H) \, =\, H\backslash Y_K.
\end{equation}
It is easy to see that \ref{assum:goodact} (i), (ii) implies that \eqref{equivisspace2} satisfies the properties (i) and (ii) of Remark \ref{remark:Sspacealt}. Therefore a unique $\wsS$-space $\wY$ exists that satisfies \eqref{equivisspace} if $S$ and $S'$ are connected. Since $Y$ commutes with coproducts \eqref{equivisspace}, holds for arbitrary $S$ and $S'$.

(b) To show property \ref{remark:Sspacealt} (iii) let $(\rho_1, \tau_1), (\rho_2, \tau_2): G/K\times \Gamma/H\to 
G/K'\times \Gamma/H'$ be two different $G\times \Gamma$-equivariant maps with $K, K'\in \sK$. 
There exists $\gamma_1, \gamma_2\in \Gamma$ with $H^{\gamma_i^{-1}}\subseteq H'$ and $\tau_i(\delta H) = \delta \gamma_i H'$ for $i=1,2$ and every $\delta\in \Gamma$. For $i=1,2$ the pair $(\rho_i, \tau_i)$ induces the map 
\begin{equation}
\label{equivisspace3}
H\backslash Y_K\lra H'\backslash Y_{K'}, \quad H\cdot y\mapsto H'\cdot(\gamma_i ^{-1}\rho_i(y)).
\end{equation}
We denote by $Z$ the subset of points of $H\backslash Y_K$ where the two maps (i.e.\ for $i=1, 2$) coincide. 
If $\rho_1 = \rho_2$ then $\tau_1\ne \tau_2$ hence $\gamma_1 H' \ne \gamma_2 H'$ and \ref{assum:goodact} (ii) yields $Z=\emptyset$. If on the other hand $\rho_1 \ne \rho_2$ then \ref{assum:goodact} (iii) implies that the interior of $Z$ is empty so the equalizer of the two maps \eqref{equivisspace3} in the category $\spaces_{\cl}$ is the empty set. 
\end{proof}

In the following $Y$ denotes an $\sS$-space equipped with a $\Gamma$-action satisfying hypothesis \ref{assum:goodact} (i), (ii) and $\wY$ denotes the associate $\wsS$-space satisfying \eqref{equivisspace}. We denote by $X=:\Gamma \backslash Y$ the push-out of $\wY$ with respect to the projection $q: G\times \Gamma\to G$, so $X$ is the $\sS$-space given by 
\begin{equation*}
\label{equivisspace2a}
X_S \, =\, \wY_{S\times \pt} \, =\, \Gamma\backslash Y_S.
\end{equation*}
for every $S\in \sS$. Condition \ref{assum:goodact} (ii) implies that the canonical projection $\pi_S: Y_S\to X_S, y\mapsto \Gamma\cdot y$ for every $S\in \sS$ defines a morphism of $\sS$-spaces $\pi: Y\to X$. 
If additionally $Y$ is exact and \ref{assum:goodact} (iii) holds then $X$ is exact as well.

By Lemma \ref{lemma:pushoutopos} the functor $\Sh(\wY, R)\to \Sh(X, R), \, \sF \mapsto \sF_{q}$
is an equivalence of categories. We denote its quasi-inverse by 
\begin{equation*}
\label{equisheaf2}
\Sh(X, R) \lra \Sh(\wY, R), \quad\sF \mapsto \wsF. 
\end{equation*}
Consider its composite with $H^0(\wY_{\infty}, \wcdot)$
\begin{equation}
\label{equisheaf3}
\Sh(X, R) \lra \Mod_R^{\sm}(G\times \Gamma), \quad \sF \mapsto H^0(\wY_{\infty}, \wsF).
\end{equation}
Note that we have 
\begin{equation}
\label{equisheaf3a}
H^0(\wY_{\infty}, \wsF)^{\Gamma}\, =\, H^0(X_{\infty}, \sF)
\end{equation}
for every $\sF\in \Sh(X, R)$. We let 
\begin{equation*}
\label{equisheaf4}
\Sh(X, R) \lra \Mod_R^{\sm}(G\times \Gamma), \quad \sF \mapsto H^n(\wY_{\infty}, \wsF).
\end{equation*}
denote the $n$-th right derived functor of \eqref{equisheaf3}. By a simple adaptation of the proof of Lemma \ref{lemma:galXinfty} we obtain

\begin{lemma}
\label{lemma:equiYinfty}
Let $\sF\in \Sh(X, R)$. As $R[G]$-modules we have
\[
H^n(\wY_{\infty}, \wsF)\, \cong \, H^n(Y_{\infty},\pi^*(\sF))
\]
for every $n\ge 0$. In particular $H^n(Y_{\infty},\pi^*(\sF))$ carries additionally a $\Gamma$-action (commuting with the $G$-action).
\end{lemma} 

We now assume that additionally that $Y$ is exact and that \ref{assum:goodact} (iii) holds so that $\wY$ and $X$ are exact as well. For $M\in \Mod_R^{\sm}(G)$ we consider the left exact functor
\begin{equation}
\label{equiext}
\Sh(X, R) \lra \Mod_{R[\Gamma]}, \quad \sF\mapsto \Hom_{R[G]}(M, H^0(\wY_{\infty}, \wsF)).
\end{equation}
Similar to Prop.\ \ref{prop:galactionhs} one shows

\begin{prop}
\label{prop:gammaactionhs}
(a) The $n$-th right derived functor of \eqref{equiext} is isomorphic to the functor $\sF\mapsto H^n(Y_{\infty}; M, \pi^*(\sF))$. In particular for $\sF\in \Sh(X, R)$ the cohomology group $H^n(Y_{\infty}; M, \pi^*(\wcdot))$ is equipped with a canonical $R[\Gamma]$-action.
\medskip

\noi (b) For $\sF\in \Sh(X, R)$ there exists a spectral sequence of $R[\Gamma]$-modules
\begin{equation*}
\label{gammaequivhss}
E_2^{rs} = \Ext_{R, G}^r(M, H^n(Y_{\infty},\pi^*(\sF))) \, \Longrightarrow\, H_R^{r+s}(Y; M, \pi^*(\sF)).
\end{equation*}
Here the $\Gamma$-action on the $E^2$-terms is induced by the $\Gamma$-action on $H^{\bu}(Y_{\infty},\pi^*(\sF))$.
\end{prop}

By \eqref{equisheaf3a} we have for $M\in \Mod_R^{\sm}(G)$ and $\sF\in \Sh(X, R)$ we have 
\begin{eqnarray}
\label{equisheaffix}
H_R^0(Y; M, \pi^*(\sF))^{\Gamma} & = & \Hom_{R[G]}(M, H^0(\wY_{\infty}, \wsF)^{\Gamma})\, =\, \Hom_{R[G]}(M, H^0(X_{\infty}, \sF))\\
& = & H_R^0(X; M, \sF).\nonumber
\end{eqnarray} 

\begin{prop}
\label{prop:gammahsss2}
Assume that $Y$ is exact and that \ref{assum:goodact} (i)--(iii) holds. 
Let $M\in \Mod_R^{\sm}(G)$, let $\sF\in \Sh(X, R)$ and assume that $M$ is flat as an $R$-module. Then there exists a spectral sequence 
\begin{equation}
\label{gammaequivss}
E_2^{rs} = H^r(\Gamma, H_R^s(Y; M, \pi^*(\sF))) \, \Longrightarrow\, H_R^{r+s}(X; M, \sF).
\end{equation}
\end{prop}

\begin{proof} By \eqref{equisheaffix} the functor $\sF\mapsto H_R^0(X; M, \sF)$ factors in the form 
\begin{equation}
\label{gammaequivss2}
\begin{CD}
\Sh(X, R) @>  \sF\mapsto H_R^0(Y; M, \pi^*(\sF)) >> \Mod_R(\Gamma) @> N\mapsto N^{\Gamma} >> \Mod_R.
\end{CD}
\end{equation}
The first functor can be identified with $\sF\mapsto \Hom_{R[G]}(M, H^0(\wY_{\infty}, \wsF))$,
so for $N\in \Mod_R(\Gamma)$ and $\sF\in \Sh(X, R)$ we get
\begin{eqnarray*}
\Hom_{R[\Gamma]}(N, H_R^0(Y; M, \pi^*(\sF))) & = & \Hom_{R[G\times \Gamma]}(M\otimes_R N, H^0(\wY_{\infty}, \wsF))\\
& = &  \Hom_{\Sh(\wY, R)}((M\otimes_R N)_{\wY, G\times \Gamma}, \wsF).
\end{eqnarray*} 
Here the $G\times \Gamma$-action on $M\otimes_R N$ is given by $(g, \gamma)\cdot m\otimes n=(gm) \otimes (\gamma n)$. So the assumption on $Y$ and $M$ implies that the first functor in \eqref{gammaequivss2} has an exact left adjoint and \eqref{gammaequivss} is the Grothendieck spectral sequence associated to the decomposition \eqref{gammaequivss2}.
\end{proof}

\paragraph{Cohomology of uniformizable $\sS$-spaces}
We now consider a particular example of an $\sS$-space $Y$ equipped with an action of group. 
Let $\Gamma$ be a subgroup of $G$ and let $\sX$ be a locally compact Hausdorff space with contractible connected components, equipped with a continuous action of $\Gamma$ (i.e.\ $\gamma\cdot \sX\to \sX$ is continuous for every $\gamma\in \Gamma$). We define an $\sS$-space $Y$ with a $\Gamma$-action as follows. For $S\in \sS$ we set
\begin{equation*}
\label{uni}
Y_S\, =\, S\times \sX.
\end{equation*}
For a morphism $\rho: S_1\to S_2$ in $\sS$ we define 
\begin{equation*}
\label{uni2}
\rho: Y_{S_1}\lra Y_{S_2}, \,\, (s, x) \mapsto (\rho(s), x).
\end{equation*}
The left $\Gamma$-action on $Y_S$, $S\in \sS$ is given by $\gamma \cdot (s, x) := (\gamma \cdot s, \gamma \cdot x)$. We consider the following

\begin{assum}
\label{assum:goodact2}
(i) $\Gamma\cap K$ acts properly discontinuously on $\sX$ for every $K\in \sS$. \\
(ii) The interior of $\{x\in \sX\mid \Stab_{\Gamma}(x)\ne 1\}$ is empty. 
\end{assum}

\begin{lemma}
\label{lemma:uniformsspace}
(a) If \ref{assum:goodact2} (i) holds then 
conditions (i) and (ii) of \ref{assum:goodact} hold for $Y$. 
\medskip

\noi (b) If additionally \ref{assum:goodact2} (ii) holds then condition \ref{assum:goodact} (iii) holds for $Y$ as well.
\end{lemma}

\begin{proof} (a) is clear. For (b) let $\rho_1, \rho_2: S\to S'$, $\rho_1\ne \rho_2$ be morphisms in $\sS$ and assume that $S$ is connected. Then we have $\rho_1(s) \ne \rho_2(s)$ for every $s\in S$. 
Thus if for $x\in \sX$ there exists $\gamma\in \Gamma$ with $\gamma (\rho_1(s), x)=(\rho_2(s), x)$ then $\gamma\in \Stab_{\Gamma}(x)$ and $\gamma\ne 1$. It follows 
\begin{equation*}
\{y\in Y_{S}\mid \exists \gamma\in \Gamma : \, \rho_2(y) = \gamma \cdot \rho_1(y)\}\,\subseteq \, S\times \{x\in \sX\mid \Stab_{\Gamma}(x)\ne 1\}
\end{equation*}
so (ii) implies that the last condition of \ref{assum:goodact} holds as well.
\end{proof}

For the rest of this section we assume that the conditions \ref{assum:goodact2} (i) and (ii) hold. We let $X=:\Gamma \backslash Y$ be the $\sS$-space associated to $Y$ and let $\pi: Y\to X$ be the projection. We will describe the cohomology groups $H_R^n(X; M, \sF)$ in terms of the group cohomology of $\Gamma$ in the case where $\sF = R_X$ is the constant sheaf and $M\in \Mod_R^{\sm}(G)$ is projective as an $R$-module. We write $H_R^n(X; M, R)$ rather than $H_R^n(X; M, R_X)$. Note that $\pi^*(R_X) = R_Y$. 

\begin{lemma}
\label{lemma:wxinf}
There exists a canonical isomorphism of discrete $\Gamma\times G$-modules
\begin{equation*}
\label{uni3}
H^n(Y_{\infty}, R) \, \cong \, \left\{\begin{array}{ll} \Ind_1^G\Maps(\pi_0(\sX), R) &  \mbox{if $n=0$,}\\
0 &  \mbox{otherwise.}
\end{array}\right.
\end{equation*}
\end{lemma}

Recall that $\Ind_1^G = \Ind_1^G: \Mod_R\to \Mod_R^{\sm}(G)$ denotes the right adjoint of the forgetful functor $\Mod_R^{\sm}(G)\to \Mod_R$. The $\Gamma$-action on $\Ind_1^G\Maps(\pi_0(\sX), R)$ is induced by the $\Gamma$-action on $\pi_0(\sX)$.

\begin{proof} For $K\in \sK$ note that that connected component of $Y_K$ are contractible and that 
$\pi_0(Y_K) = G/K\times \pi_0(\sX)$. Hence we have 
\begin{equation*}
\label{uni4}
H^n(Y_K, R) \,=\,  \left\{\begin{array}{cc} \Maps(G/K\times \pi_0(\sX), R) & \mbox{if $n=0$,}\\
0 & \mbox{otherwise}
\end{array}\right.
\end{equation*}
for every $K\in \sK$. It follows 
\begin{eqnarray*}
\label{uni5}
H^n(Y_{\infty}, R) & \cong & \left\{\begin{array}{cc}  \dlim_{K\in \sK^{\opp}} \Maps(G/K\times \pi_0(\wX), R) & \mbox{if $n=0$,}\\
0 & \mbox{otherwise,}
\end{array}\right.\nonumber\\
& = & \left\{\begin{array}{ll} \Ind_1^G\Maps(\pi_0(\sX), R) &  \mbox{if $n=0$,}\\
0 & \mbox{otherwise.}
\end{array}\right.\nonumber
\end{eqnarray*}
\end{proof}

For $M\in \Mod_R^{\sm}(G)$ we define a left $\Gamma$-module $\cC_R(M, R)$ as follows. Elements of $\cC_R(M, R)$ are maps $\psi: M\times \pi_0(\sX) \to R$ that are homomorphisms of $R$-modules in the first component (i.e.\ $\psi(\wcdot , x):M \to R$ is an $R$-linear map for every $x\in \pi_0(\sX)$). The $\Gamma$-action on $\cC_R(M, R)$ is given by $(\gamma \psi)(m, x)= \psi(\gamma^{-1}m, \gamma^{-1} x)$
for $\gamma\in \Gamma$, $\psi\in \cC_R(M, R)$, $m\in M$ and $x\in \pi_0(\sX)$. 

\begin{prop}
\label{prop:uni4}
Let $M\in \Mod_R^{\sm}(G)$ and assume that $M$ is projective as an $R$-module. Then we have: 
\medskip

\noi (a) There exists a canonical isomorphism of discrete $\Gamma\times G$-modules
\begin{equation*}
\label{uni6}
H_R^n(Y; M, R) \, \cong \, \left\{\begin{array}{ll} \cC_R(M, R) &  \mbox{if $n=0$,}\\
0 &  \mbox{otherwise.}
\end{array}\right.
\end{equation*}
\medskip

\noi (b) For $n\ge 0$ we have 
\begin{equation*}
\label{uni7}
H_R^n(X; M, R) \, \cong \, H^n(\Gamma, \cC_R(M, R)).
\end{equation*}
\end{prop}

\begin{proof} (a) Consider the spectral sequence (see Prop.\ \ref{prop:gammaactionhs} (b))
\begin{equation*}
\label{exthssa}
E_2^{rs} = \Ext_{R, G}^r(M, H^s(Y_{\infty}, R)) \, \Longrightarrow\, H_R^{r+s}(Y; M, R).
\end{equation*}
By Lemma \ref{lemma:wxinf} it degenerates and we have
\begin{equation*}
\label{pullbacksd}
\Ext_{R, G}^{\bu}(M, \Ind_1^G\Maps(\pi_0(\sX), R))\, =\, H_R^{\bu}(Y; M, R).
\end{equation*}
Since functor $\Ind_1^G: \Mod_R\to \Mod_R^{\sm}(G)$ is exact and has an exact left adjoint and since the composition 
\[
\begin{CD} 
\Mod_R@> \Ind_1^G >> \Mod_R^{\sm}(G)@> \Hom_{R[G]}(M, \wcdot ) >> \Mod_R
\end{CD}
\]
can be identified with the exact functor $\Hom_R(M, \wcdot): \Mod_R \to \Mod_R$ we obtain
\begin{eqnarray}
\label{uni6e}
H_R^n(Y; M, R) & = &  \Ext_R^n(M,  \Maps(\pi_0(\sX), R)) \\
& = &  \left\{\begin{array}{cc} \Hom_R(M,  \Maps(\pi_0(\sX), R))  & \mbox{if $n=0$,}\\
0 & \mbox{otherwise}
\end{array}\right.\nonumber\\
& = & \left\{\begin{array}{ll} \cC_R(M, R) & \mbox{if $n=0$,}\\
0 & \mbox{otherwise.}
\end{array}\right.\nonumber
\end{eqnarray}
By tracing through the definitions one verifies that \eqref{uni6e} is an isomorphism of $\Gamma$-modules.

To prove (b) we use the spectral sequence \eqref{gammaequivss}. By (a) it degenerates so we get $H_R^n(X; M, R) = H^n(\Gamma, \cC_R(M, R))$.
\end{proof}

\section{Ordinary and $\varpi$-adic cohomology}
\label{section:ordinary}

Throughout this chapter we let $F$ denote a $p$-adic field, i.e.\ $F/\bQ_p$ is finite extension for some fixed prime $p$. We use the notation as in section \ref{section:borelind} so $\cO_F$ denotes the valuation ring and $\fp$ the valuation ideal of $F$. Also recall that for $n\ge 1$ we have $U^{(n)}=U_F^{(n)} =1 +\fp_F^n$ and $U_F = U_F^{(0)} = \cO_F^*$. The group $G$ denotes the group $\PGL_2(F)$, $B$ the standard Borel subgroup of $G$ and $T$ the maximal torus in $B$. If $H$ is a closed subgroup of $T^0$ then we will write $\barT = T/H$ etc.\

We fix a subset $\sK$ of the set of open subgroups of $G$ satisfying the conditions (i)--(iii) of assumption \ref{assum:kgood} and we denote by $\sS=\sS_{G,\sK}$ the associated site of discrete left $G$-sets whose stabilizers contained in $\sK$.

\subsection{Ordinary cohomology}
\label{subsection:principal}

\paragraph{Ordinary cohomology of $\sS$-spaces} Let $X$ be an $\sS$-space and let $\cO$ be a complete noetherian local ring with residue field $k$. For a locally admissible $\cO[\barT]$-module $W$ and $\sF\in \Sh(X, \cO)$ we consider the cohomology groups 
\begin{equation}
\label{parindmod}
\bH_{\cO}^n(X; W, \sF)\,: =\, H_{\cO}^n(X, \Ind_B^{G} W, \sF)
\end{equation}
for $n\ge 0$. Since the functor $\Ind_B^{G}: \Mod_{\cO}(\barT)\to \Mod_{\cO}(G)$ is exact the cohomology groups \eqref{parindmod} define for a fixed $\sF\in \Sh(X, \cO)$ an $\cO$-linear $\delta$-functor
\begin{equation*}
\label{parindmod2}
\Mod_{\cO}(\barT)^{\opp} \lra \Mod_{\cO},\,\,\, W \mapsto \bH_{\cO}^{\bu}(X; W, \sF)\qquad n\ge 0
\end{equation*} 

As explained in section \ref{subsection:admbanach} a locally admissible $\cO[\barT]$-module $W$ can be viewed as an augmented $\cO[\barT]$-module (i.e.\ it carries a natural $\La_{\cO}(\barT)$-module structure; see \eqref{iwasawa} for the definition of $\La_{\cO}(\barT)$). It induces also on $\Hom_{\Sh(X, \cO)}((\Ind_B^{G} W)_{X}, \sF)$ a structure of an augmented $\cO[\barT]$-module for every $\sF\in \Sh(X, \cO)$. Hence we may view the assignment $\sF\mapsto \Hom_{\Sh(X, \cO)}((\Ind_B^{G} W)_{X}, \sF)$ as a functor 
\begin{equation}
\label{parindmod3}
\Hom_{\Sh(X ,\cO)}((\Ind_B^{G} W)_{X}, \wcdot): \Sh(X, \cO)\lra \Mod_{\cO}^{\aug}(\barT). 
\end{equation}
Therefore for fixed $W\in \Mod_{\cO}^{\ladm}(\barT)$ the cohomology groups \eqref{parindmod} can be viewed as the $n$-th right derived functor of \eqref{parindmod3} evaluated at $\sF$ so they are equipped with the structure of an augmented $\cO[\barT]$-module as well.

By Prop.\ \ref{prop:ordadj} (a) we have 
\begin{eqnarray*}
\label{adjoint}
\Hom_{\Sh(X, \cO)}((\Ind_B^{G} W)_{X}, \sF) & \cong & \Hom_{\cO[G]}(\Ind_B^{G} W, \sF(X_{\infty}))\\
& \cong & \Hom_{\cO[\barT]}(W^{\io}, \Ord^H(\sF(X_{\infty})))\nonumber
\end{eqnarray*}
so the functor \eqref{parindmod3} (for fixed $W$) can be decomposed as
\begin{equation*}
\label{parindmod4}
\begin{CD}
\Sh(X, \cO)@> H^0(X_{\infty}, \wcdot ) >> \Mod_{\cO}^{\sm}(G) @> \Ord^H >> \Mod_{\cO}^{\ladm}(\barT) @> \Hom_{\cO[\barT]}(W^{\io}, \wcdot) >> \Mod_{\cO}^{\aug}(\barT).
\end{CD}
\end{equation*}
We denote the composite of the first two functors by 
\begin{equation}
\label{ordsheaf}
\fOrd_{\cO}^H(X, \wcdot): \Sh(X, \cO) \lra \Mod_{\cO}^{\ladm}(\barT).
\end{equation}
It is left exact since it has left adjoint. We define {\it the $n$-th ordinary cohomology} of $X$ -- denoted by $\fOrd_{\cO}^{H, n}(X, \wcdot)$ -- as the $n$-th right derived functor by $\fOrd_{\cO}^H(X, \wcdot)$, i.e.\  for $\sF\in \Sh(X,\cO)$ we define
\begin{equation}
\label{ordsheaf2}
\fOrd_{\cO}^{H, n}(X, \sF)\, =\, (R^n \fOrd_{\cO}^H(X, \wcdot)(\sF)\in \Mod_{\cO}^{\ladm}(\barT).
\end{equation}
If $H=1$ then $H$ will be drop from the notation, i.e.\ we will write $\fOrd_{\cO}^n(X, \sF)$ for the cohomology group \eqref{ordsheaf2}. 

\begin{prop}
\label{prop:ford}
Let $\sF\in \Sh(X, \cO)$ and $W\in \Mod_{\cO}^{\ladm}(\barT)$.
\medskip

\noi (a) We have
\begin{equation*}
\label{fordlim}
\fOrd_{\cO}^{H, n}(X, \sF) \, =\,  \dlim_{U \open,\, H\le U \le T^0} \fOrd_{\cO}^{U, n}(X, \sF)
\end{equation*}
(i.e.\ $U$ ranges over all open subgroups of $T^0$ containing $H$).
\medskip

\noi (b) If $X$ is exact then there exists two spectral sequences 
\begin{eqnarray}
\label{ordinftyss}
&& E_2^{rs} = \Ord_{\cO}^{H,r}(H^s(X_{\infty}, \sF)) \, \Longrightarrow\, \fOrd_{\cO}^{H, r+s}(X, \sF),\\
\label{fordextss}
&& E_2^{rs} = \Ext_{\cO, \barT}^r(W^{\io}, \fOrd_{\cO}^{H, s}(X, \sF)) \, \Longrightarrow\, \bH_{\cO}^{r+s}(X; W, \sF).
\end{eqnarray}
\end{prop}

\begin{proof} (a) Let $\cU$ denote the set of all open subgroups $U$ of $T^0$ containing $H$. We have
\begin{eqnarray*}
\fOrd_{\cO}^H(X, \sG) =   \Ord^H( \sG(X_{\infty})) = \dlim_{U\in \cU} \Ord^U( \sG(X_{\infty}))
  =  \dlim_{U\in \cU} \fOrd_{\cO}^U(X, \sG)
\end{eqnarray*}
for every $\sG\in \Sh(X, \cO)$. Thus if $0\lra \sF\lra \sI^{\bu}$ is an injective resolution then we have
\begin{eqnarray*}
\fOrd_{\cO}^{H, n}(X, \sF) & = & H^n(\fOrd_{\cO}^H(X, \sI^{\bu}))\,= \,H^n(\dlim_{U\in \cU} \fOrd_{\cO}^U(X, \sI^{\bu}))\\
& = & \dlim_{U\in \cU} H^n(\fOrd_{\cO}^U(X, \sI^{\bu}))\,=\,  \dlim_{U\in \cU} \fOrd_{\cO}^{U, n}(X, \sF).
\end{eqnarray*}
(b) As \eqref{ordsheaf} is the composite of two left exact functors the first one preserving injectives by Cor.\ \ref{corollary:gmodsheafxexact}, the corresponding Grothendieck spectral sequence exists. Also we can decompose \eqref{parindmod3} as the composite of \eqref{ordsheaf} and  
$\Hom_{\cO[\barT]}(W^{\io}, \wcdot): \Mod_{\cO}(\barT)^{\ladm}\to \Mod_{\cO}^{\aug}(\barT)$. The exactness of parabolic induction functor implies that $\fOrd_{\cO}^H(X, \wcdot)$ has an exact left adjoint as well. Hence the corresponding Grothendieck spectral sequence \eqref{fordextss} exists, too.
\end{proof}

\begin{corollary}
\label{corollary:ordnxcoeff}
Let $\varphi:\cO\to \cO'$ be an epimorphism of complete noetherian local rings and assume that $\dim(\cO)\le 1$ and $\dim(\cO')= 0$. Then the canonical homomorphism
\begin{equation*}
\label{ordnxcoeff1}
\fOrd_{\cO}^{H, n}(X, \sF)_{\varphi} \lra \fOrd_{\cO}^{H, n}(X, \sF_{\varphi}) 
\end{equation*}
is an isomorphism for every $\sF\in \Sh(X, \cO')$ and $n\ge 0$.
\end{corollary}

\begin{proof} As in the proof of Prop.\ \ref{prop:drmqch} it suffices to see that for an injective $\cO'$-sheaf $\sI$ on $X$ we have $\fOrd_{\cO}^{H, n}(X, \sI_{\varphi})=0$. For that consider the spectral sequence \eqref{ordinftyss} for $\sF=\sI_{\varphi}$. 
By Remark \ref{remarks:gcompactcoh2} (a) we have $E_2^{rs}=0$ for $s>0$. Hence together with Cor.\ \ref{corollary:derordncoeff}
we get 
\[
\fOrd_{\cO}^{H, n}(X, \sI_{\varphi}) =\Ord_{\cO}^{H, n}(\sI(X_{\infty})_{\varphi})=\Ord_{\cO'}^{H, n}(\sI(X_{\infty}))_{\varphi}= \fOrd_{\cO'}^{H, n}(X, \sI)_{\varphi}=0
\]
for $n\ge 1$.
\end{proof}

We will study finiteness properties of the groups \eqref{parindmod}.
For that consider the following conditions for a sheaf $\sF\in \Sh(X, \cO)$.

\begin{assum}
\label{assum:sheaffin}
(i) We have $H^n(X_K, \sF_K)\in \Mod_{\cO, f}$ for every $K\in \sK$ and $n\ge 0$. \\
(ii) There exists $d\in \bN$ such that $H^n(X_K, \sF_K)=0$ for every $K\in \sK$ and $n\ge d+1$.
\end{assum}

\begin{remark}
\label{remark:constrfin}
\rm Condition \ref{assum:sheaffin} (i) and (ii) holds e.g.\ in the case where $X_K$ is a $d$-dimensional compact topological manifold with boundary for each $K\in \sK$ and where $\sF$ is a constant sheaf 
$\sF=N_{X}$ associated to a finitely generated $\cO$-module $N$ (see \cite{iversen}, Ch.\ III, 9.11 and Thm.\ 10.1).
\enddemo
\end{remark}

Condition \ref{assum:sheaffin} (i) yields good finiteness properties for the cohomology groups \eqref{parindmod} in the following cases: (i) $W$ is finitely generated as an $\cO$-module and (ii) $W$ is admissible and $\cO=A$ is Artinian. 

\begin{prop}
\label{prop:finiteness2}
Assume that $X$ is exact and that $\dim(\cO) \le 1$.
\medskip

\noi (a) Assume is an open subgroup of $T^0$ and let $m\ge 1$ such that $\delta(U^{(m)})\subseteq H$. 
Then there exists a spectral sequence
\begin{equation}
\label{fordssfinlev}
E_2^{rs} = R^r\Gamma^{\ord} H^s(X_{K_H(m)}, \sF_{K_H(n)}) \, \Longrightarrow\, \fOrd_{\cO}^{H, r+s}(X, \sF)
\end{equation}
for every $\sF\in \Sh(X, \cO)$. Moreover if $\sF\in \Sh(X, \cO)$ satisfies condition \ref{assum:sheaffin} (i) then 
\begin{equation}
\label{ordsheaf3}
\fOrd_{\cO}^{H, n}(X, \sF)\, =\, H^n(X_{K_H(m)}, \sF_{K_H(m)})^{\ord}.
\end{equation}
In particular in this case $\fOrd_{\cO}^{H, n}(X, \sF)$ is finitely generated as an $\cO$-module for every $n\ge 0$. 
\medskip

\noi (b) Let $W\in \Mod_{\cO,f}(\barT)$ and let $\sF$ be an $\cO$-sheaf satisfying condition \ref{assum:sheaffin} (i). Then the $\cO$-module $\bH_{\cO}^m(X; W, \sF)$ is finitely generated for every $m\ge 0$. 
\end{prop}

Recall that $K_H(m)$ denotes the compact open subgroup $K(\fp^m) H N_0$ of $G$ for any closed subgroup $H$ of $T^0$ and any $m\ge 1$. 

\begin{proof} (a) By Prop.\ \ref{prop:ordHcomp} (b) the functor \eqref{ordsheaf} can be decomposed as
\begin{equation*}
\label{parindmod4a}
\begin{CD}
\Sh(X, \cO)@> H^0(X_{\infty}, \wcdot ) >> \Mod_{\cO}^{\sm}(G) @> H^0(K_H(m), \wcdot) >> 
\Mod_{\cO}(\barT^+) @>  \Gamma^{\ord} >> \Mod_{\cO}^{\ladm}(\barT).
\end{CD}
\end{equation*}
The composite of the first two functors is the functor
\begin{equation}
\label{parindmod5}
\Sh(X, \cO) \lra \Mod_{\cO}(\barT^+),\,\, \sF \mapsto H^0(X_{K_H(m)}, \sF_{K_H(m)}).
\end{equation}
By Cor.\ \ref{corollary:gmodsheafxexact} and Lemma \ref{lemma:injacyclic} it maps injective objects of $\Sh(X, \cO)$ onto $\Gamma^{\ord}$-acyclic objects of $\Mod_{\cO}(\barT^+)$. By Lemma \ref{lemma:resexactfaith} (a) its $n$-th right derived functor is 
\begin{equation*}
\label{parindmod5d}
\Sh(X, \cO) \lra \Mod_{\cO}(\barT^+),\,\, \sF \mapsto H^n(X_{K_H(m)}, \sF_{K_H(m)}).
\end{equation*}
Also by Cor.\ \ref{corollary:gmodsheafxexact} and Lemma \ref{lemma:injacyclic} the functor \eqref{parindmod5} maps injective objects of $\Sh(X, \cO)$ onto $\Gamma^{\ord}$-acyclic objects of $\Mod_{\cO}(\barT^+)$. Hence there exists a Grothendieck spectral sequence \eqref{fordssfinlev}. Now if $\sF$ is an $\cO$-sheaf on $X$ such that $H^n(X_K, \sF_K)\in \Mod_{\cO,f}$ for every $K\in \sK$ and $n\ge 0$ then we have $E_2^{rs}= 0$ if $r\ge 1$ by Prop.\ \ref{prop:rgammaord} (a). Hence the spectral sequence degenerates and \eqref{ordsheaf3} follows from Prop.\ \ref{prop:adTTplus} (a).

(b) If $W\in \Mod_{\cO,f}(\barT)$ then there exists an open subgroup $U$ of $T^0$ with $H\subseteq U$ and $W^U= W$, so that we can view $W$ as a $\cO[T/U]$-module. Lemma \ref{lemma:extadmfin1} together with (a) implies that the $E_2$-terms of the spectral sequence \eqref{fordextss} (for $H=U$) are finitely generated $\cO$-modules. So the limit terms are finitely generated $\cO$-modules as well. 
\end{proof}

\begin{remarks}
\label{remarks:fordhida}
\rm (a) Let $X$ and $\cO$ be as in the Prop.\ \ref{prop:finiteness2}. Assume that $\sF\in \Sh(X, \cO)$ satisfies condition \ref{assum:sheaffin} (i) and that $H$ is an arbitrary closed subgroup of $T^0$. Then $\delta(U^{(1)})\cdot H\supseteq \delta(U^{(2)})\cdot H \supseteq \delta(U^{(3)})\cdot H\supseteq \dots$ is a decreasing sequence of open subgroups of $T^0$ with $\bigcap_{m\ge 1} \delta(U^{(m)})\cdot H = H$. Thus by combining Prop.\ \ref{prop:ford} (a) and Prop.\ \ref{prop:finiteness2} (a) we obtain
\begin{equation}
\label{fordlim2}
\fOrd_{\cO}^{H, n}(X, \sF) \, =\,  \dlim_m H^n(X_{K_H(m)}, \sF_{K_H(m)})^{\ord}.
\end{equation}
This implies in particular for $H=1$ 
\begin{equation}
\label{fordlim3}
\fOrd_{\cO}^n(X, \sF) \, =\,  \dlim_m H^n(X_{K_1(m)}, \sF_{K_1(m)})^{\ord}.
\end{equation}
where $K_1(m)= K_1(\fp^m)$ denotes the congruence subgroup \eqref{congruence} of $G$.
We note that for $n=0$ \eqref{fordlim2} and \eqref{fordlim3} hold even if $X$ is not exact.
\medskip

\noi (b) Note that if additionally condition \ref{assum:sheaffin} (ii) holds then \eqref{fordlim2} implies
\begin{equation*}
\label{fordvanishing}
\fOrd_{\cO}^{H, n}(X, \sF) \, =\,  0
\end{equation*}
for every $n\ge d+1$.  \enddemo
\end{remarks}

We now turn to the case when $\cO=A$ is Artinian and $W$ is admissible. 

\begin{prop}
\label{prop:finiteness1}
Let $X$ be exact and let $A$ be an Artinian local ring with finite residue field $k$ of characteristic $p$. 
Let $W\in \Mod_A^{\adm}(\barT)$, let $\sF\in \Sh(X, A)$ and assume that \ref{assum:sheaffin} (i), (ii) hold. 
\medskip 

\noi (a) The $A[\barT]$-module $\fOrd_A^{H,n}(X, \sF)$ is admissible for every $n\in \bN$. It vanishes if $n\ge d +1$. 
\medskip

\noi (b) The augmented $A[\barT]$-module $\bH_A^n(X; W, \sF)$ is finitely generated for every $n\ge 0$. 
\end{prop}

\begin{proof} (a) The first assertion follows from Cor.\ \ref{corollary:gadmissible} and Cor.\ \ref{corollary:adordn} using the spectral sequence \eqref{ordinftyss} (for $\cO=A$) and the second from Remark \ref{remarks:fordhida} (b) above. 

(b) follow from (a), Prop.\ \ref{prop:extadmfin2} and Prop.\ \ref{prop:extadmfin3} using the spectral sequence \eqref{fordextss}. 
\end{proof}

\paragraph{Ordinary cohomology with compact support} Next we introduce {\it ordinary cohomology with compact support} $\fOrd_{\cO, c}^{\bu}(X, \wcdot)$ using a rather ad hoc definition (we only consider the case $H=1$ so we drop $H$ from the notation). In the following we assume that $\cO$ is a complete noetherian local ring of dimension $\le 1$ with finite residue field of characteristic $p$ and that $X_S$ is a locally compact Hausdorff space for every $S\in \sS$. It is then easy to see that for $\sF\in \Sh(X, \cO)$ and $K, K'\in \sK$ with $K'\subseteq K$ the induced monomorphism $H^0(X_K, \sF_K) \hra H^0(X_{K'}, \sF_{K'})$ maps the submodule $H^0_c(X_K, \sF_K)\subseteq H^0(X_K, \sF_K)$ into $H^0_c(X_{K'}, \sF_{K'})\subseteq H^0(X_{K'}, \sF_{K'})$
so that we can define an $\cO[G]$-submodule 
\[
H^0_c(X_{\infty}, \sF) :=\, \dlim_{K\in \sK}
\]
of $H^0(X_{\infty}, \sF)$. Furthermore it is easy to see that we have $H^0_c(X_{\infty}, \sF)^K=H^0_c(X_K, \sF_K)$. By Prop.\ \ref{prop:ordHcomp} (a) we obtain a left exact functor 
\begin{equation*}
\label{sectcsk1m1}
\Sh(X, \cO) \, \lra \, \Mod_{\cO}^{\sm}(T^+), \,\, \sF\mapsto  H_c^0(X_{K_1(m)}, \sF_{K_1(m)}) \, =\, H_c^0(X_{\infty}, \sF)^{K_1(m)}.
\end{equation*}
for every $m\ge 1$. By Lemma \ref{lemma:resexactfaith} its $n$-th right derived functor is $\sF\mapsto H_c^n(X_{K_1(m)}, \sF_{K_1(m)})$. For $\sF\in \Sh(X,\cO)$ the canonical homomorphism of $\cO[T^+]$-modules
\begin{equation*}
\label{sectcsk1m2}
H_c^0(X_{K_1(m)}, \sF_{K_1(m)}) = H_c^0(X_{\infty}, \sF)^{K_1(m)} \hra H_c^0(X_{\infty}, \sF)^{K_1(m+1)}  = H_c^0(X_{K_1(m+1)}, \sF_{K_1(m+1)})
\end{equation*}
induces a morphism of $\delta$-functors 
\begin{equation}
\label{sectcsk1m3}
 H_c^n(X_{K_1(m)}, \sF_{K_1(m)})\lra H_c^n(X_{K_1(m+1)}, \sF_{K_1(m+1)}), \qquad n\ge 0.
\end{equation}
Let $S_{\fcoh}(X, \cO)$ denote the full subcategory of $\Sh(X, \cO)$ consisting of $\cO$-sheaves $\sF$ such that $H_c^n(X_K,\break \sF_K)\in \Mod_{\cO, f}$ for every $K\in \sK$ and $n\ge 0$. For $\sF\in S_{\fcoh}(X, \cO)$ and $m\ge 0$ we define the $n$-th {\it ordinary cohomology of $X$ with compact support} by 
\begin{equation*}
\label{fordcs1}
\fOrd_{\cO, c}^n(X, \sF) : = \, \dlim_m H_c^n(X_{K_1(m))}, \sF_{K_1(m)})^{\ord}
\end{equation*}
(the transitions maps of the direct system are induced by \eqref{sectcsk1m3}). If the layers $X_K$ of $X$ are compact then $\fOrd_{\cO, c}^n(X, \sF)$ can be identified with $\fOrd_{\cO}^n(X, \sF)$ by Remark \ref{remarks:fordhida} (a). Also the fact that the functor $\Mod_{\cO, f}(T^+)\to \Mod_{\cO, f}, V\mapsto V^{\ord}$ is exact implies that if $0\to \sF_1\to \sF_2\to \sF_3\to 0$ is 
a short exact sequence of $\cO$-sheaves on $X$ all lying in $S_{\fcoh}(X, \cO)$ then there exists an associated long exact $\fOrd_{\cO}^{\bu}(X, \wcdot)$-sequence. 

Now assume that $U$ is an open $\sS$-subspace of $X$, i.e.\ $U_K$ is an open subset of $X_K$ for every $K\in \sK$ and we have $\rho^{-1}(U_{K_2}) = U_{K_1}$ for every morphism $\rho: K_1\to K_2$ in $\sS$. Then $K\mapsto Y_K : = X_K\setminus U_K$ defines a (closed) $\sS$-subspace $Y$ of $X$. For an $\cO$-sheaf $\sF$ on $X$ we put $\sF|_U = j^*(\sF)$ and $\sF|_Y = i^*(\sF)$ where $j:U\hra X$, $i:Y\hra X$ denote the inclusions. 

\begin{lemma}
\label{lemma:csres}
Let $\sF\in S_{\fcoh}(X, \cO)$ and assume that $\sF|_U\in S_{\fcoh}(U, \cO)$ and $\sF|_Y\in S_{\fcoh}(Y, \cO)$. There exists a natural long exact sequence
\begin{equation*}
\label{fordcohomverscs}
\ldots \lra \fOrd_{\cO, c}^n(U, \sF|_U) \lra  \fOrd_{\cO, c}^n(X, \sF) \lra \fOrd_{\cO, c}^n(Y, \sF|_Y) \lra \fOrd_{\cO, c}^{n+1}(U, \sF|_U)\lra\ldots
\end{equation*}
\end{lemma}

\begin{proof} This follows by passing in the long exact sequence
\begin{eqnarray*}
\label{fordcohomverscs2}
&& \ldots \lra H_c^n(U_{K_1(m)}, \sF_{K_1(m)}) \lra H_c^n(X_{K_1(m)}, \sF_{K_1(m)}) \\
&& \hspace{3cm} \lra H_c^n(Y_{K_1(m)}, (\sF_{K_1(m)})|_Y) \lra H_c^{n+1}(U_{K_1(m)}, \sF_{K_1(m)})\lra\ldots
\end{eqnarray*}
to ordinary parts and then passing to the direct limit over all $m\ge 1$. 
\end{proof}

\begin{remark}
\label{remark:fpordcsuniv}
\rm A more natural definition for ordinary cohomology of $X$ with compact support would be to consider the right derived functors of 
\begin{equation*}
\label{ordcohomocsupp2}
\fOrd_{\cO, c}: \Sh(X, \cO) \lra \Mod_{\cO}^{\ladm}(T), \,\,\, \sF\mapsto \Ord(H_c^0(X_{\infty}, \sF)).
\end{equation*}
For $\sF\in S_{\fcoh}(X, \cO)$ it is easy to see that there exists a canonical homomorphism 
\begin{equation}
\label{fordlimcs}
(R^n \fOrd_{\cO, c})(\sF) \lra  \fOrd_{\cO, c}^n(X, \sF) 
\end{equation}
that is an isomorphism if $n=0$. In general it seems however that \eqref{fordlimcs} is not an isomorphism.  \enddemo
\end{remark} 

\paragraph{Ordinary cohomology of $\sS$-schemes} In the following $L$ denotes a subfield of $\bC$, $\barL\subseteq \bC$ the algebraic closure of $L$ and $\fG=\Gal(\overline{L}/L)$ the absolute Galois group of $L$. We consider an exact $\sS$-scheme $X$ of finite type over $L$ and denote by $X^{\an} = (X_{\bC})^{\an}$ the associated $\sS$-space. Let $A$ be an Artinian local ring with finite residue field of characteristic $p$. Otherwise we keep the notation and assumptions from the beginning of this section. As explained in section \ref{subsection:galcohomology} the {\'e}tale cohomology groups 
\begin{equation}
\label{parindmodgal}
\bH_A^n(X_{\barL}; W, \sF)\,: =\, H_A^n(X_{\barL}; (\Ind_{B}^{G} W)_{X_{\barL}, G}, \sF_{\barL})
\end{equation}
for $W\in \Mod_A^{\ladm}(\barT)$, $\sF\in \Sh(X_{\et}, A)$ and $n\ge 0$ carry a canonical $\fG=\fG_L$-action. We will now show that there exists an {\'e}tale variant of $\fp$-ordinary cohomology and a Galois equivariant version of the spectral sequence \eqref{fordextss}. 

Note that the locally profinite group $\barT\times \fG$ satisfies condition \ref{assum:normalsubgroup}. Therefore the category $\Mod_A^{\ladm}(\barT\times \fG)$ is a Serre subcategory of $\Mod_A^{\sm}(\barT\times\fG)$ with enough injectives. Moreover using Lemma \ref{lemma:tadm} one shows that a discrete $A[\barT\times\fG]$-module $W$ is locally admissible if and only if $W$ is locally admissible when viewed solely as an $A[\barT]$-module. 

Let $W\in \Mod_A^{\sm}(\barT\times\fG)$ and let $V\in \Mod_A^{\sm}(G\times\fG)$. It is easy to see that the $A[\barT\times \fG]$-module $\Ord^H(V)$ is locally admissible and (using Lemma \ref{lemma:parindbasechange} (b)) that the $A[G\times \fG]$-module $\Ind_{B}^{G} W$ is discrete. Therefore by Prop.\ \ref{prop:ordadj}
\begin{eqnarray*}
\label{parindgal1}
&& \Mod_A^{\sm}(G\times \fG)\lra \Mod_A^{\ladm}(\barT\times \fG), \,\, V \mapsto \Ord^H(V),\\
&&  \Mod_A^{\ladm}(\barT\times \fG)\lra \Mod_A^{\sm}(G\times \fG),\,\, W\mapsto \Ind_{B}^{G} W^{\iota}
\end{eqnarray*}
is a pair of adjoint functors. 

Consider the composite of functors
\begin{equation}
\label{galparindmod4}
\begin{CD}
\Sh(X_{\et}, A) @> \sF \mapsto H^0((X^{\gal})_{\infty}, \sF^{\gal}) >> \Mod_A^{\sm}(G\times \fG) @> \Ord^H >> \Mod_A^{\ladm}(\barT\times \fG)
\end{CD}
\end{equation}
where the first functor \eqref{galsheaf2} has been considered in section \eqref{subsection:galcohomology}.

\begin{df}
\label{df:pordet}
We denote the $n$-th right derived functor of  \eqref{galparindmod4} by 
\begin{equation*}
\label{pordet2}
\fOrd_A^{H, n}(X_{\barL}, \wcdot): \Sh(X_{\et}, A)\lra \Mod_A^{\ladm}(\barT\times \fG),\,\, \sF \mapsto \fOrd_A^{H, n}(X_{\barL}, \sF).
\end{equation*}
The $A[\barT\times \fG]$-module $\fOrd_A^{H, n}(X_{\barL}, \sF)$ will be called {\it {\'e}tale ordinary cohomology} of $X_{\barL}$ with coefficients in $\sF\in \Sh(X_{\et}, A)$. If $H=1$ then, as before, $H$ will be dropped from the notation.
\end{df}

Similar to Prop.\ \ref{prop:ford} and Prop.\ \ref{prop:finiteness2} we have 

\begin{prop}
\label{prop:fordet}
Let $\sF\in \Sh(X_{\et}, A)$ and let $W\in \Mod_A^{\ladm}(\barT)$.
\medskip

\noi (a) Assume that $\sF$ is constructible and that $H$ is an open subgroup of $T^0$. Let $m\ge 1$ so that $\delta(U^{(m)})\subseteq H$. Then there exists a natural isomorphism of $A[\barT\times \fG]$-modules 
\begin{equation*}
\label{pordet3}
\fOrd_A^{H, n}(X_{\barL}, \sF)\, =\, H^n((X_{K_H(m)})_{\barL}, (\sF_{K_H(m)})_{\barL})^{\ord}.
\end{equation*}
for every $n\ge 0$. In particular $\fOrd_A^{H, n}(X_{\barL}, \sF)$ is finitely generated as an $A$-module. 
\medskip
 
\noi (b) We have
\begin{equation*}
\label{pordetlim}
\fOrd_A^{H, n}(X_{\barL}, \sF) \, =\,  \dlim_{U \open,\, H\le U \le T^0} \fOrd_A^{U, n}(X_{\barL}, \sF)
\end{equation*}

\noi (c) There exists a spectral sequence of $A[\fG]$-modules
\begin{equation*}
\label{pordetextss}
E_2^{rs} = \Ext_{A, \barT}^r(W^{\io}, \fOrd_A^{H, s}(X_{\barL}, \sF)) \, \Longrightarrow\, \bH_A^{r+s}(X_{\barL}; W, \sF).
\end{equation*}
Here the $\fG$-action on the $E_2$-terms is induced by the $\fG$-action on $\fOrd_A^{H, \bu}(X_{\barL}, \sF)$. 
\medskip

\noi (d) If $\sF$ is constructible and if $W\in \Mod_{A, f}^{\sm}(\barT)$ (resp.\ $W\in \Mod_A^{\adm}(T)$) then we have 
\[
\bH_A^n(X_{\barL}; W, \sF_{\barL})\in \Mod_{A, f}(\barT)\qquad  \mbox{(resp.} \quad \bH_A^n(X_{\barL}; W, \sF_{\barL})\in \Mod_{A}^{\fgaug}(\barT))
\]
and the $\fG$-action on $\bH_A^n(X_{\barL}; W, \sF_{\barL})$ is discrete (resp.\ continuous with respect to the canonical topology) for every $n\ge 0$.
\end{prop}

\begin{proof} The proofs of (a) and (b) are simple adaptations of the proofs of \ref{prop:finiteness2} (a) and \ref{prop:ford} (a).
For (c) we use the fact that the functor 
\begin{equation*}
\label{parindet1}
\Sh(X_{\et}, A) \lra \Mod_{A[\fG]}, \,\, \sF\mapsto \Hom_{A[G]}(\Ind_{B}^{G} W, H^0((X^{\gal})_{\infty}, \sF^{\gal}))
\end{equation*}
is isomorphic to the composite 
\begin{equation*}
\label{parindet2}
\begin{CD}
\Sh(X_{\et}, A) @> \eqref{galparindmod4}>>  \Mod_R^{\ladm}(\barT\times \fG) @> \Hom_{A[\barT]}(W^{\io}, \wcdot) >> \Mod_{A[\fG]}
\end{CD}
\end{equation*}
so there exists an associated Grothendieck spectral sequence. By Prop.\ \ref{prop:galactionhs} (a) the limit terms can be identifies with the $R[\fG]$-modules $\bH_A^{r+s}(X_{\barL}; W, \sF_{\barL})$. Using Lemma \ref{lemma:extgfg} (for $G = \barT$) the $E_2^{rs}$-term can be identified with the $A[\fG]$-module $\Ext_{A[\barT]}^r(W^{\io}, \fOrd_A^{H, s}(X_{\barL}, \sF_{\barL}))$. 

For (d) assume first that $W\in \Mod_{A, f}^{\sm}(\barT)$. Then there exists an open subgroup $U$ of $T^0$ with $W^U=W$, i.e.\ we may assume that $H$ is open and that $\barT^0$ is discrete. So the assertion follows from 
(a) and (c). Now assume that $W\in \Mod_A^{\adm}(T)$. By Remark \ref{remarks:trace} (c), Prop.\ \ref{prop:etancomp} and Prop.\ \ref{prop:finiteness1} the cohomology group $\bH_A^n(X_{\barL}; W, \sF_{\barL})\cong \bH_A(X^{\an}, W, \sF^{\an})$ is finitely generated as an augmented $A[\barT]$-module, so it is equipped with a canonical (profinite) topology. To see that the $\fG$-action on $\bH_A^n(X_{\barL}; W, \sF_{\barL})$ is continuous we put 
$W^{(m)}= W^{\bar{\delta}(U_F^{(m)})}$ so that $W= \bigcup_{m\ge 1} W^{(m)}$ and $W^{(m)} \in \Mod_{A, f}(\barT)$ for every $m\ge1$. As in the proof of Prop.\ \ref{prop:extadmfin3} (a) one shows that there exists an exact sequence 
\[
0\lra \prolimm^{(1)} \bH_A^{n-1}(X_{\barL}; W^{(m)}, \sF_{\barL})\lra \bH_A^n(X_{\barL}; W, \sF_{\barL}) \lra \prolimm \bH_A^{n-1}(X_{\barL}; W^{(m)}, \sF_{\barL})\lra 0
\]
The assertion now follows from the fact that $\bH_A^n(X_{\barL}; W^{(m)}, \sF_{\barL})$ is finitely generated as an $A$-module for all $m, n\ge 0$, so the first term vanishes. 
\end{proof}

Finally, we have the following comparison result between {\'e}tale ordinary cohomology of $X_{\barL}$ and ordinary cohomology of $X^{\an}$. 

\begin{prop}
\label{prop:fordetsing}
(a) There exists a natural homomorphism 
\begin{equation}
\label{compordetan1}
 \fOrd_A^{H, n}(X_{\barL}, \sF)\lra  \fOrd_A^{H, n}(X^{\an}, \sF^{\an})
\end{equation}
for every $\sF\in \Sh(X_{\et}, A)$ and $n\ge 0$ (where $\sF^{\an}: =(\sF_{\bC})^{\an}$). If $\sF$ is constructible then \eqref{compordetan1} is an isomorphism. 
\medskip

\noi (b) Let $\sF\in \Sh(X, A)$ and let $W\in \Mod_A^{\ladm}(\barT)$. There exists a natural morphism of spectral sequences 
\begin{eqnarray}
\label{compordetan2}
&& \left( \, \Ext_{A, \barT}^r(W^{\io}, \fOrd_A^{H, s}(X_{\barL}, \sF)) \, \Longrightarrow\, \bH_A^{r+s}(X_{\barL}; W, \sF)\, \right) \lra \, \\
&& \hspace{3cm} \left(\,\Ext_{A, \barT}^r(W^{\io}, \fOrd_A^{H, s}(X^{\an}, \sF^{\an})) \, \Longrightarrow\, \bH_A^{r+s}(X^{\an}; W, \sF^{\an})\, \right).
\nonumber
\end{eqnarray}
It is an isomorphism if $\sF$ is constructible. 
\end{prop}

\begin{proof} (a) For $n=0$ the canonical homomorphism $H^0((X_{\bC})_{\infty}, \sF_{\bC}) \to H^0(X^{\an}_{\infty}, \sF^{\an})$ induces a homomorphism
\begin{eqnarray}
\label{compordetan1deg0}
\fOrd_A^{H, 0}(X_{\barL}, \sF) & = & \Ord^H( H^0((X_{\barL})_{\infty}, \sF_{\barL})) \,\,= \,\, \Ord^H( H^0((X_{\bC})_{\infty}, \sF_{\bC}))\\
& \lra & \Ord^H(H^0(X^{\an}_{\infty}, \sF^{\an})) = \fOrd_A^{H, 0}(X^{\an}, \sF^{\an}).\nonumber
\end{eqnarray}
Since $\sF\mapsto \fOrd_A^{H, n}(X^{\an}, \sF^{\an})$, $n\ge 0$ is a $\delta$-functor the homomorphisms \eqref{compordetan1deg0} extend uniquely to the morphism \eqref{compordetan1} of $\delta$-functors. 
That \eqref{compordetan1} is an isomorphism when $\sF$ is constructible follows immediately from Prop.\ \ref{prop:ford} (a), Prop.\ \ref{prop:finiteness2} (a) and Prop.\ \ref{prop:fordet}.

(b) The existence of the morphism \eqref{compordetan2} can be proved in the same way as in Prop.\ \ref{prop:hsetsing} (a). That it is an isomorphism if $\sF$ is constructible follows then from  (a).
\end{proof}

\begin{remark}
\label{remark:fpordcsuppet}
\rm In section \ref{section:modcurves} we also need {\'e}tale ordinary cohomology with compact support (only in the case $H=1$). Again we will give an ad hoc definition: for a constructible sheaf $\sF\in \Sh(X_{\et}, A)$ and $n\ge 0$ we define the $A[T\times \fG]$-module
\begin{equation*}
\label{ordcohomocsupp4}
\fOrd_{A, c}^n(X_{\barL}, \sF) \, =\,  \dlim_m H_c^n((X_{K_1(m)})_{\barL}, (\sF_{K_1(m)})_{\barL})^{\ord}.
\end{equation*}
That the group $H_c^n((X_{K_1(m)})_{\barL}, (\sF_{K_1(m)})_{\barL})$ lies in $\Mod_{A, f}^{\sm}(T^+)$ (so that it makes sense
to consider its ordinary part) can be seen e.g.\ by using the comparison isomorphism 
\[
H_c^n((X_{K_1(m)})_{\barL}, (\sF_{K_1(m)})_{\barL})\,\cong\, H_c^n(X_{K_1(m)}^{\an}, \sF_{K_1(m)}^{\an})
\]
and by using the fact that $\sF^{\an}$ lies in the category $S_{\fcoh}(X^{\an}, A)$ introduced earlier in this section. The assignment $\sF\mapsto \fOrd_{A, c}^n(X_{\barL}, \sF)$, $n\ge 0$ is a $\delta$-functor from the category of constructible {\'e}tale $A$-sheaves on $X$ to the category $\Mod_A^{\ladm}(T\times\fG)$ and there exists natural comparison isomorphisms of locally admissible $A[T]$-modules 
\begin{equation*}
\label{ordcohomocsupp5}
\fOrd_{A, c}^n(X_{\barL}, \sF) \, =\, \fOrd_{A, c}^n(X^{\an}, \sF^{\an}).
\end{equation*}
\end{remark} 

\subsection{$\varpi$-adic cohomology}
\label{subsection:varpiadic}

As at the end of the previous section $X$ denotes an exact $\sS$-scheme of finite type defined over a subfield $L$ of $\bC$ and $\fG= \Gal(\overline{L}/L)$ the absolute Galois group of $L$. Let $d$ be the dimension of $X$, i.e.\ the dimension of $X_K$ for some (hence all) $K\in \sK$.
Let $\cO$ be a complete discrete valuation ring with maximal ideal $\fm= (\varpi)$, finite residue field $k$ of characteristic $p$ and quotient field $E$ of characteristic $0$. As before for $m\ge 1$ we put $\cO_m = \cO/\fm^m$. We aim to extend the definition of the cohomology groups \eqref{parindmodgal} to the case where $W$ is a $\varpi$-adically admissible $\cO[\barT]$-module $W$ with $W_{\tor} =0$ (compare section \ref{subsection:admbanach}). 
and where the coefficients $\sF$ are {\'e}tale constructible $\varpi$-adic sheaf on $X$. 

As before we put $\cO_m = \cO/\fm^m$ for $m\ge 1$. Also for any $\cO$-module $N$ we put $N_m= N\otimes_{\cO} \cO_m$ for $m\ge 1$. We denote by $\Sh_{\varpi}(X_{\et}, \cO)$ the full subcategory of $\Sh(X_{\et}, \cO)$ consisting of $\cO$-sheaves $\sF$ such that $\varpi^m \cdot \sF=0$ for some $m \ge 0$. By $\Sh_{\varpi, c}(X_{\et}, \cO)$ we denote the full subcategory of sheaves $\sF\in \Sh_{\varpi}(X_{\et}, \cO)$ that are constructible. It follows easily from Lemma \ref{lemma:resexactfaith} that $\Sh_{\varpi, c}(X_{\et}, \cO)$ is an $\cO$-linear abelian subcategory of $\Sh(X_{\et}, \cO)$ that is closed under extensions. Moreover any object of $\Sh_{\varpi, c}(X_{\et}, \cO)$ is noetherian. 

Before we extend the definition of the cohomology groups \eqref{parindmodgal} we introduce some notation and review basic notions and facts about Artin-Rees categories (see \cite{sga5}, Exp.\ V). Let $\cA$ be an $\cO$-linear abelian category. We assume that $\cA$ is noetherian, i.e.\ every object is noetherian. The objects of the Artin-Rees category $\AR \cA$ of $\cA$ are projective systems $A_{\bu} = (A_m)_{m\in \bZ}$ in $\cA$ with $A_m=0$ for $m\ll 0$. Morphisms in $\AR \cA$ are 
\[
\Hom_{\AR \cA}(A_{\bu}, B_{\bu}) \, =\, \dlim_d \Hom(A_{\bu}[d], B_{\bu}) 
\]
where $A_{\bu}[d] := (A_{m+d})$ for $d\in \bZ$ and $\Hom(A_{\bu}[d], B_{\bu})$ denotes the abelian group of morphisms of projective systems. The transition maps of the direct system $(\Hom(A_{\bu}[d], B_{\bu}))_{d\ge 0}$ are induced by the canonical morphism $\io: A_{\bu}[d+1]\to A_{\bu}[d]$. The Artin-Rees category of $\cA$ is again abelian. A projective system $A_{\bu} = (A_m)_{m\in \bZ}$ in $\cA$ is called a strictly $\varpi$-adic system if $A_m=0$ for $m\le 0$, $\varpi^m A_m =0$ for $m\ge 0$ and the natural map $A_{m+1}/\varpi^m A_{m+1} \to A_m$ is an isomorphism for all $m\ge 0$. A projective system $B_{\bu}$ is called $\AR \varpi$-adic if it is $\AR$isomorphic (i.e.\ isomorphic in $\AR \cA$) to a strictly $\varpi$-adic system. The full subcategory of $\AR \varpi$-adic objects of 
$\AR \cA$ will be denoted by $\ARp \cA$. It is a noetherian abelian subcategory of $\AR \cA$ (i.e.\ the inclusion $\ARp \cA\hra \AR \cA$ is exact). If inverse limits exists in $\cA$ then we can consider the functor $A_{\bu} = (A_m)_{m\in \bZ}\mapsto \prolimm A_m$. It induces a functor 
\begin{equation*}
\label{prolimar}
\prolim: \AR \cA \lra \cA, \qquad  A_{\bu}\mapsto \prolimm A_{\bu}.
\end{equation*}
Let $H = (H^n): \cA\to \cB$ be an $\cO$-linear $\delta$-functor between $\cO$-linear noetherian abelian categories. Then $H$ extends to an $\cO$-linear $\delta$-functor $\AR \cA\to \AR \cB$ which -- by abuse of notation -- will be denoted by $H = (H^n)$ as well. For $A_{\bu}\in \AR \cA$ and $n\ge 0$ we have $H^n(A_{\bu}) = (H^n(A_m))_{m\in \bZ}$. 


Returning to our task of extending the definition of the cohomology groups \eqref{parindmodgal}, 
let $W$ be a $\varpi$-adically admissible $\cO[\barT]$-module with $W_{\tor}=0$ and let
$\sF\in \Sh_{\varpi}(X_{\et}, \cO)$. We choose $m\ge 1$ with $\varpi^m \sF=0$, so that $\sF\in \Sh_{\varpi}(X_{\et}, \cO_m)$ and we define
\begin{equation}
\label{extco2}
\bH_{\{\fm\}}^{\bu}(X_{\barL}; W, \sF) \, :=\, \bH_{\cO_m}^{\bu}(X_{\barL}; W_m, \sF).
\end{equation}
By Remark \ref{remarks:padicadm} (a) and Lemma \ref{lemma:parindbasechange} the $\cO_m$-module $\Ind_{B_{\fp}}^{G_{\fp}} W_m$ is projective for every $m\ge 0$. Hence, by Prop.\ \ref{prop:drmqch}, this definition is independent of the choice of $m$. By Prop.\ \ref{prop:finiteness2} and Prop.\ \ref{prop:fordet} (d) the cohomology groups \eqref{extco2} are finitely generated augmented $\cO[\barT]$-modules equipped with a continuous $\fG$-action. Thus we obtain an $\cO$-linear $\delta$-functor
\begin{equation}
\label{laextgal}
\Sh_{\varpi, c}(X_{\et}, \cO) \lra \Mod_{\cO}^{\fgaug}(\barT)(\fG),\quad \sF \mapsto \bH_{\{\fm\}}^n(X_{\barL}; W, \sF) \qquad n\ge 0.
\end{equation} 
Here $\Mod_{\cO}^{\fgaug}(\barT)(\fG)$ denotes the category of pairs $(M, \rho)$ where $M\in \Mod_{\cO}^{\fgaug}(\barT)$ and $\rho: \fG\to \Aut(M)$ is a homomorphism so that the action $\fG\times M \to M, (\sigma, m)\mapsto \rho(\sigma)(m)$ is continuous. We need the following 

\begin{lemma}
\label{lemma:araugfglim}
(a) The functor 
\begin{equation*}
\label{prolimar2}
\prolim: \AR \Mod_{\cO}^{\fgaug}(\barT) \lra \Mod_{\cO}^{\aug}(\barT), \qquad  M_{\bu}\mapsto \prolimm M_m.
\end{equation*}
is exact. 
\medskip

\noi (b) Let $M_{\bu}\in \AR \Mod_{\cO}^{\fgaug}(\barT)$ and assume that $(\prolim M_{\bu})\otimes_{\cO} k\in  \Mod_k^{\fgaug}(\barT)$. Then $\prolim M_{\bu}\in \Mod_{\cO}^{\fgaug}(\barT)$.
\end{lemma}

\begin{proof} (a) follows from Lemma \ref{lemma:lim1profinite} applied to projective systems of finitely generated modules of the complete noetherian local ring $\Lambda_{\cO}(U)$ (where $U$ is any open subgroup of $\barT^0$). For (b) note that $\prolim M_{\bu}$ is a compact $\Lambda_{\cO}(U)$-module so the assertion follows from the topological Nakayama Lemma. 
\end{proof}

\begin{df}
\label{df:constructible}
We define the category $\Sh(X_{\et}, \pad)$ of {\'e}tale constructible $\varpi$-adic sheaves on $X$ by
\[
 \Sh(X_{\et}, \pad):=\, \ARp \Sh_{\varpi, c}(X_{\et}, \cO).
\] 
\end{df}

As explained above \eqref{laextgal} extends to a $\delta$-functor $\AR \Sh_{\varpi, c}(X_{\et}, \cO)  \to \AR \Mod_{\cO}^{\fgaug}(\barT)$. We compose it with the functor \eqref{prolimar2} and restrict it to the subcategory $\Sh(X_{\et}, \pad)$. So for an {\'e}tale constructible $\varpi$-adic sheaf $\sF=\sF_{\bu}$ on $X$ and $n\ge 0$ we consider the group
\begin{equation}
\label{laextgal2}
\bH_{\varpi-\ad}^n(X_{\barL}; W, \sF) : =\, \prolimm \bH_{\{\fm\}}^n(X_{\barL}; W, \sF_m).
\end{equation} 

\begin{prop}
\label{prop:varpiadcohom}
The assignment $\sF \mapsto \bH_{\varpi-\ad}^n(X_{\barL}; W, \sF)$, $n\ge 0$ defines a $\delta$-functor 
\begin{equation}
\label{laextgal3}
\Sh(X_{\et}, \pad) \lra \Mod_{\cO}^{\fgaug}(\barT)(\fG), \quad \left(\sF\mapsto \bH_{\varpi-\ad}^n(X_{\barL}; W, \sF)\right)_{n\ge 0}. 
\end{equation} 
\end{prop}

\begin{proof} By Lemma \ref{lemma:araugfglim} (a) we only have to show $\bH_{\varpi-\ad}^n(X_{\barL}; W, \sF)\in \Mod_{\cO}^{\fgaug}(\barT)$ for $\sF=\sF_{\bu}\in \Sh(X_{\et}, \pad)$ and $n\ge 0$. Firstly, we consider the case $\varpi \sF = 0$, i.e.\ we assume that $\varpi\sF_m =0$ holds for every $m\in \bZ$. There exists a strictly $\varpi$-adic system $\sG = \sG_{\bu}$ in $\Sh_{\varpi, c}(X_{\et}, \cO)$ that is $\AR$isomorphic to $\sF$, i.e.\ there are morphisms of projective systems $\phi: \sF_{\bu}[d_1] \to \sG_{\bu}$, $\psi: \sG_{\bu}[d_2] \to \sF_{\bu}$ for some $d_1, d_2\ge 0$ such that $\psi\circ \phi$ (resp.\ $\phi\circ \psi$) is equal to $\id_{\sF}$ (resp.\ $\id_{\sG}$) in $\AR\Sh_{\varpi, c}(X_{\et}, \cO)$. Thus there exists $d\ge d_1 + d_2$ so that the composite 
\[
\sG_{\bu}[d] \stackrel{\tau}{\lra} \sG_{\bu}[d_1+d_2] \stackrel{\psi}{\lra} \sF_{\bu}[d_1] \stackrel{\phi}{\lra} \sG_{\bu}
\]
is equal to the canonical map $\tau:\sG_{\bu}[d]\to \sG_{\bu}$ induced by the transition maps of the projective system $\sG_{\bu}$. The condition $\varpi \sF_m =0$ for all $m \in \bZ$ thus implies that the image of the transition map $t_{m+d, m}: \sG_{d+m} \to \sG_m$ is annihilated by $\varpi$. On the other hand $t_{m+d, m}$ is surjective since $\sG = \sG_{\bu}$ is a strictly $\varpi$-adic. It follows that $\varpi \sG_m =0$ and $t_{m+1, m}: \sG_{m+1} \to \sG_m$ is an isomorphisms for every $m\ge 1$. Hence we have 
\begin{equation}
\label{laextgal4}
\bH_{\varpi-\ad}^n(X_{\barL}; W, \sF)\, \cong \, \bH_{\varpi-\ad}^n(X_{\barL}; W, \sG)\, \cong \, \bH_{\cO_1}^n(X_{\barL}; W_1, \sG_1)\in \Mod_k^{\fgaug}(\barT).
\end{equation}
Now let $\sF_{\bu}$ be an arbitrary strictly $\varpi$-adic system in $\Sh_{\varpi, c}(X_{\et}, \cO)$. Put 
\[
\sF^{(1)}:= \ker(\sF\stackrel{\varpi\wcdot}{\lra}\sF), \quad \sF^{(2)}:= \image(\sF\stackrel{\varpi\wcdot}{\lra}\sF), \quad \sF^{(3)}:= 
\coker(\sF\stackrel{\varpi\wcdot}{\lra}\sF).
\]
Note that $\sF^{(1)}, \sF^{(2)}, \sF^{(3)}\in \Sh_{\varpi, c}(X_{\et}, \pad)$ and that we have $\varpi \sF^{(1)}=0= \varpi \sF^{(3)}$. 
By considering the long exact sequences for the $\delta$-functor \eqref{laextgal3} associated to the two short exact sequences 
$0\to \sF^{(1)}\to \sF\to \sF^{(2)}\to 0$ and $0\to \sF^{(2)}\to \sF\to \sF^{(3)}\to 0$ and using \eqref{laextgal4} (for the coefficients
$\sF^{(1)}$ and $\sF^{(3)}$) we obtain
\begin{equation*}
\label{laextgal5a}
\coker\left(\bH_{\varpi-\ad}^n(X_{\barL}; W, \sF)\stackrel{\varpi\wcdot}{\lra}\bH_{\varpi-\ad}^n(X_{\barL}; W, \sF)\right)\, \in \Mod_k^{\fgaug}(\barT).
\end{equation*}
for all $n\ge 0$. Together with Lemma \ref{lemma:araugfglim} (b) we conclude $\bH_{\varpi-\ad}^n(X_{\barL}; W, \sF)\in \Mod_{\cO}^{\fgaug}(\barT)$.
\end{proof}

We will now discuss the functorial properties of the cohomology groups \eqref{laextgal2} with respect to $W$. Let $\varphi: W'\to W$ be a homomorphism 
 of $\varpi$-adically admissible $\cO[\barT]$-module with $W_{\tor}=0 = W'_{\tor}$ (i.e.\ we have $W, W'\in \Mod_{\cO}^{\varpi-\adm}(\barT)_{\fl}$).
Then $\varphi$ induces a morphism of the associated $\delta$-functors \eqref{laextgal3}, so for every $\sF\in \Sh(X_{\et}, \pad)$, $n\ge 0$ there exists a canonical homomorphism 
\begin{equation*}
\label{laextgal5b}
\varphi^*:  \bH_{\varpi-\ad}^n(X_{\barL}; W, \sF)\lra \bH_{\varpi-\ad}^n(X_{\barL}; W', \sF)
\end{equation*}
 We have

\begin{prop}
\label{prop:mexact}
Let 
\begin{equation*}
\label{mseq}
\begin{CD}
0@>>> W'@> \varphi >> W @> \psi >> W'' \lra 0
\end{CD}
\end{equation*}
be a short exact sequence in $\Mod_{\cO}^{\varpi-\adm}(\barT)_{\fl}$. Then there exists a $\fG$-equivariant long exact sequence of augmented $\cO[\barT]$-modules
\begin{eqnarray}
\label{laextgal5}
&& \ldots \lra \bH_{\varpi-\ad}^n(X_{\barL}; W'', \sF) \stackrel{\psi^*}{\lra} \bH_{\varpi-\ad}^n(X_{\barL}; W, \sF) \stackrel{\varphi^*}{\lra} \bH_{\varpi-\ad}^n(X_{\barL}; W', \sF)\\
&& \hspace{2cm} \lra \bH_{\varpi-\ad}^{n+1}(X_{\barL}; W'', \sF)\lra  \ldots\nonumber 
\end{eqnarray}
that is functorial in $\sF\in \Sh(X_{\et}, \pad)$.
\end{prop}

\begin{proof} This follows from the exactness of the sequence 
\begin{equation*}
\label{mseq2}
\begin{CD}
0@>>> W'_m@>>> W_m @>>> W''_m \lra 0
\end{CD}
\end{equation*}
for every $m\ge 1$, the fact that parabolic induction is an exact functor and from Lemma \ref{lemma:araugfglim} (a).
\end{proof}

As before let $X^{\an} = (X_{\bC})^{\an}$ be the $\sS$-space associated to $X$. We let $\Sh_{\varpi, c}(X^{\an}, \cO)$ be the category of $\cO$-sheaves of $X^{\an}$ of the form $\sF^{\an}$ with $\sF\in \Sh_{\varpi, c}((X_{\bC})_{\et}, \cO)$ and we put $\Sh(X^{\an}, \pad): = \ARp \Sh_{\varpi, c}(X^{\an}, \cO)$. Similarly to \eqref{extco2} and \eqref{laextgal2}
we define the groups $\bH_{\{\fm\}}^n(X^{\an}; W, \sF)$ and $\bH_{\varpi-\ad}^n(X^{\an}; W, \sF)$ for 
$\sF\in \Sh_{\varpi, c}(X^{\an}, \cO)$ and $\sF\in \Sh(X^{\an}, \pad)$ respectively. Again the latter defines a $\delta$-functor 
\begin{equation*}
\label{laextan3} 
\Sh(X^{\an}, \pad) \lra \Mod_{\cO}^{\fgaug}(\barT),\quad \sF \mapsto \bH_{\varpi-\ad}^n(X^{\an}; W, \sF) \quad n\ge 0
\end{equation*} 
We note that the analog of Prop.\ \ref{prop:mexact} holds for this $\delta$-functor as well. Of course if $\sF=\sF_{\bu}\in \Sh(X_{\et}, \pad)$ and if $\sF^{\an} = \sF_{\bu}^{\an}\in \Sh(X^{\an}, \pad)$ with $\sF_m^{\an} :=((\sF_m)_{\bC})^{\an}$ for $m\in \bZ$ then by Remark \ref{remarks:trace} (c) and Prop.\ \ref{prop:etancomp} we have 
\begin{equation*}
\label{laextan4} 
\bH_{\varpi-\ad}^{\bu}(X_{\barL}; W, \sF)\,\cong \, \bH_{\varpi-\ad}^{\bu}(X^{\an}; W, \sF^{\an}).
\end{equation*} 
Here the groups on the left are equipped with a continuous $\fG$-action. To justify introducing  
the cohomology groups on the right separately we like to point out that they are sometimes equipped with an additional structure defined purely analytically. 

The projective system of constant sheaves $\cO_{\bu}=((\cO_m)_{X})_{m}$ defines an {\'e}tale constructible $\varpi$-adic sheaf on $X$. We put 
\begin{equation}
\label{madiccohomet}
\bH_{\varpi-\ad}^n(X_{\barL}; W, \cO) \,: =\, \bH_{\varpi-\ad}^n(X_{\barL}; W, \cO_{\bu}) \, =\,  \prolimm \bH_{\varpi-\ad}^n(X_{\barL}; W, \cO_m).
\end{equation}
and 
\begin{equation*}
\label{madiccohom}
\bH_{\varpi-\ad}^n(X^{\an}; W, \cO) \,: =\, \bH_{\varpi-\ad}^n(X^{\an}; W, \cO_{\bu}) \, =\,  \prolimm \bH_{\varpi-\ad}^n(X^{\an}; W, \cO_m).
\end{equation*}
Now assume that $W$ is a discrete $\cO[\barT]$ that is free and of finite rank as an $\cO$-module. Then $W$ is also a $\varpi$-adically admissible $\cO[\barT]$-module so we can compare the cohomology groups $\bH_{\cO}^n(X^{\an}; W, \cO)$ and $\bH_{\varpi-\ad}^n(X^{\an}; W, \cO)$. 

\begin{prop}
\label{prop:compwwhat}
Assume that $W\in \Mod_{\cO, f}^{\sm}(\barT)$ with $W_{\tor} =0$. Then the canonical homomorphism 
\begin{equation}
\label{quasichar}
\bH_{\cO}^n(X^{\an}; W, \cO)\, \lra \, \bH_{\varpi-\ad}^n(X^{\an}; W, \cO)
\end{equation}
is an isomorphism of finitely generated augmented $\cO[\barT]$-modules for every $n\ge 0$.
\end{prop}

\begin{proof} The map \eqref{quasichar} is defined as follows: for $m\ge 1$ the canonical projection $\cO\to \cO_m$ induces a homomorphism 
\begin{equation}
\label{quasichar2}
\bH_{\cO}^n(X^{\an}; W, \cO) \to \bH_{\cO}^n(X^{\an}; W, \cO_m) \cong  \bH_{\cO_m}^n(X^{\an}; W_m, \cO_m) = \bH_{\{\fm\}}^n(X^{\an}; W, \cO_m)
\end{equation}
where the second isomorphism follows from Prop.\ \ref{prop:drmqch}. Passing in \eqref{quasichar2} to the projective limit over all $m\ge 1$ yields the map \eqref{quasichar}. For $n\ge 0$ consider the short exact sequence 
\begin{equation}
\label{quasichar3}
0\lra \bH^n(W)\otimes_{\cO}\cO_m \lra \bH_{\cO}^n(X^{\an}; W, \cO_m) \lra \bH^{n+1}(W)[\varpi^m]\lra 0
\end{equation}
induced by $0\to \cO\stackrel{\varpi^m\cdot}{\lra} \cO \to \cO_m \to 0$
where we have abbreviated $\bH^{\bu}(W) = \bH_{\cO}^{\bu}(X^{\an}; W, \cO)$. 
By Prop.\ \ref{prop:finiteness2} (b) we have $\bH^n(W)\in \Mod_{\cO, f}$ for every $n\ge 0$. 
Hence $\bH^n(W)$ is $\varpi$-adically complete with finite torsion, so we get
\[
\bH^n(W)\,\cong \,\prolimm \bH^n(W)\otimes_{\cO}\cO_m \qquad \mbox{and} \qquad \prolimm \bH^n(W)[\varpi^m] \, =\, 0
\]
for every $n\ge 0$. The assertion follows by passing in \eqref{quasichar3} to projective limits.
\end{proof}

\paragraph{The cohomology groups $\bH_{\varpi-\ad}^n(X_{\barL}; V, E)$.} Finally, we show that inverting $\varpi$ allows us to replace $W\in \Mod_{\cO}^{\varpi-\adm}(\barT)_{\fl}$ in \eqref{madiccohomet} with an admissible $E$-Banach space representation of $\barT$. 
Namely, given $V\in \Ban_E^{\adm}(\barT)$ we can choose $W\in \Mod_{\cO}^{\varpi-\adm}(T)$ that is mapped to $V$ under the functor \eqref{locbanach2} and we define 
\begin{equation*}
\label{madiccohomet2}
\bH_{\varpi-\ad}^n(X_{\barL}; V, E) \,: =\, \bH_{\varpi-\ad}^n(X_{\barL}; W_{\fl}, \cO)_E 
\end{equation*}
for $n\ge 0$. Recall that the $\cO$-torsion subgroup $W_{\tor}$ of $W$ is of bounded exponent so the projection $\pr: W\to W_{\fl} = W/W_{\tor}$ becomes an isomorphism in the category $(\Mod_{\cO}^{\varpi-\adm}(\barT))_E$.

\begin{prop}
\label{prop:banachdelta}
The assignment $V \mapsto \bH_{\varpi-\ad}^n(X_{\barL}; V, E)$, $n\ge 0$ defines a contravariant $\delta$-functor 
\begin{equation*}
\label{laextgal7}
\Ban_E^{\adm}(\barT) \lra \Mod_{E}^{\fgaug}(\barT)(\fG),\quad \left(V \mapsto \bH_{\varpi-\ad}^n(X_{\barL}; V, E)\right)_{n\ge 0}
\end{equation*}
\end{prop}

\begin{proof} For $n\ge 0$ the assignment $W\mapsto \bH_{\varpi-\ad}^n(X_{\barL}; W_{\fl}, \cO)_E$ induces contravariant functors  
\begin{equation}
\label{madiccohomet3}
\Mod_{\cO}^{\varpi-\adm}(\barT) \lra \Mod_E^{\fgaug}(\barT)(\fG),\quad W \mapsto \bH_{\varpi-\ad}^n(X_{\barL}; W_{\fl}, \cO)_E\quad n\ge 0.
\end{equation}
We want to show that they give rise to a $\delta$-functor. Let $0\to W_1 \stackrel{\alpha}{\lra}  W_2 \stackrel{\beta}{\lra} W_3 \to 0$ be a short exact sequence in $\Mod_{\cO}^{\varpi-\adm}(\barT)$ and consider the diagram 
\begin{equation}
\label{madiccohomet4}
\begin{CD} 
0 @>>> W_1@>\alpha>>  W_2 @>\beta >> W_3@>>> 0\\
@. @VVV @V\pr VV @VVV \\
0 @>>> W_1' @>> \incl >  (W_2)_{\fl} @>\beta_{\fl} >> (W_3)_{\fl} @>>> 0.
\end{CD} 
\end{equation}
where $W_1':= \ker(\beta_{\fl})$. The first vertical arrow induces a map $\gamma: (W_1)_{\fl} \to W_1'$ whose kernel and cokernel are $\cO$-torsion modules of bounded exponent so the induced maps 
\[
\gamma^*: \bH_{\varpi-\ad}^{\bu}(X_{\barL}; W_1', \cO)_E \lra \bH_{\varpi-\ad}^{\bu}(X_{\barL}; (W_1)_{\fl}, \cO)_E
\]
are isomorphisms. By Prop.\ \ref{prop:mexact} there exists a long exact sequence \eqref{laextgal5} associated to the lower row of \eqref{madiccohomet4} (and $\sF=\cO_{\bu}$). Let  
\[
\delta' : \bH_{\varpi-\ad}^n(X_{\barL}; W_1', \cO)\lra \bH_{\varpi-\ad}^{n+1}(X_{\barL}; (W_3)_{\fl}, \cO)
\] 
be the connecting homomorphism. Then the sequence 
\begin{eqnarray*}
\label{laextgal8}
&& \ldots \lra \bH_{\varpi-\ad}^n(X_{\barL}; (W_3)_{\fl}, \cO)_E \to \bH_{\varpi-\ad}^n(X_{\barL}; (W_2)_{\fl}, \cO)_E \lra \bH_{\varpi-\ad}^n(X_{\barL}; (W_1)_{\fl}, \cO)_E\\
&& \hspace{2cm} \stackrel{\delta}{\lra} \bH_{\varpi-\ad}^{n+1}(X_{\barL}; (W_3)_{\fl}, \cO)\lra  \ldots\nonumber 
\end{eqnarray*}
is exact where $\delta:= \gamma^*\circ (\delta')_E \circ (\gamma^*)^{-1}$. Hence \eqref{madiccohomet3} defines a $\delta$-functor. Now the assertion follows from 
Prop.\ \ref{prop:locdelta2} and Lemma \ref{lemma:locbanach}.
\end{proof}

Similarly, the functors $W\mapsto \bH_{\varpi-\ad}^n(X^{\an}; W_{\fl}, \cO)$ give rise to a $\delta$-functor  
\begin{equation*}
\label{laextan7}
\Ban_E^{\adm}(T) \lra  \Mod_{E}^{\fgaug}(\barT),\quad V \mapsto \bH_{\varpi-\ad}^n(X^{\an}; V, E) \quad n\ge 0.
\end{equation*}

\section{Quaternionic Hilbert modular $\sS$-varieties} 
\label{section:cohomqhmsv}

\subsection{Definitions and basic properties}
\label{subsection:qhmsv}

Let $\barQ$ be the algebraic closure of $\bQ$ in $\bC$ and let $F\subset \barQ$ be a totally real number field. The ring of integers of $F$ will be denoted by $\cO_F$. Let $\cS$ be a finite set of nonarchimedean places of $F$. We denote by $\bA$ (resp.\ $\bA^{\cS}$ resp.\ $\bA_f$ resp.\ $\bA_f^{\cS}$) the ring of adeles of $F$ (resp.\ ring of prime-to-$\cS$ resp.\ finite resp.\ finite prime-to-$\cS$ adeles). We let $\cS_{\infty}$ denote the set of archimedean places of $F$ and put $F_{\infty} = F\otimes_{\bQ} \bR$. For a nonarchimedean place $\fq$ of $F$ we let $\cO_{\fq}$ denote the valuation ring in $F_{\fq}$ and we put $U_{\fq}^{(0)} = U_{\fp} = \cO_{\fq}^*$ and $U_{\fq}^{(n)} =1 +\fq^n\cO_{\fq}$ and for $n\ge 1$. 

Let $D$ be a quaternion algebra over $F$, let $\wG=D^*$ (viewed as an algebraic group over $F$), let $Z\cong \bG_m$ be the center of $\wG$ and put $G = \wG/Z$ (thus $G=\PGL_2/F$ if $D=M_2(F)$ is the algebra of $2\times 2$ matrices). We put $G_{\cS} = G(F_{\cS})$ where $F_{\cS}=\prod_{v\in S} F_v$. Let $\Ram_D$ be the set of (archimedian or nonarchimedian) places of $F$ that are ramified in $D$. We assume that $S$ and $\Ram_D$ are disjoint so that $G_{\cS} \cong \PGL_2(F_{\cS})$. We denote by $\Sigma$ the set of archimedean places of $F$ that split $D$ (i.e.\ $\Sigma = \cS_{\infty}\setminus \Ram_D$) and we put $d=\# \Sigma$. If $d\ge 1$, i.e.\ if $D$ is not totally ramified then we choose an ordering $\Sigma = \{\sigma_1, \ldots, \sigma_d\}$ of the places in $\Sigma$. 

We identify $\cS_{\infty}$ with the set of embeddings $\Hom(F, \bR)$ and assume that $\sigma_1$ is the inclusion (if $d\ge 1$). We let $\fd$ be the ideal of $\cO_F$ that is the product of the primes which are ramified in $D$, so that $\Ram_D = (\cS_{\infty}\setminus \Sigma) \cup \{\fq\mid \fd\}$. 

For $v\in S_{\infty}$ we denote by $G_{v, +}$ the connected component of $1$ in $G_v$. Thus $G_{v, +}=\wG(F_v)_+/Z(F_v)$ where $\wG(F_v)_+$ is the subgroup of elements $g\in \wG(F_v)$ with $\Nrd(g) >0$. We choose a maximal compact subgroup $K_v$ of $G_v$. If $v\in \Sigma$ we choose an isomorphism $G_v\cong \PGL_2(\bR)$ that maps $K_v$ to $\Orth(2)Z(F_v)/Z(F_v)\cong \Orth(2)/\{\pm 1\}$, so that $G_v/K_{v, +}$ can be identified with $\bH^{\pm}$, the union of the upper and lower complex half plane. If $v\in S_{\infty}\cap \Ram_D$ we have $G_v= G_{v, +} = K_v$. Put $K_{\infty, +} = \prod_{v\mid\infty} K_{v, +}$ and $G_{\infty, +}= \prod_{v\mid\infty} G_{v, +}$ so that $G_{\infty, +}/K_{\infty, +}\cong \bH^d$.

We fix a compact open subgroup $K_f^{\cS}$ of $G(\bA_f^{\cS})$ (called a prime-to-$S$ level) and put
\begin{equation*}
\label{quotients}
\sX^{\cS} \, =\, G(\bA^{\cS})/\left(K^{\cS}_f\times K_{\infty, +}\right).
\end{equation*}
For a compact open subgroup $K$ of $G_{\cS}$ and $g\in G(\bA_f^{\cS})$ we put
\begin{equation*}
\label{congrkg}
\Gamma_{K, g} \, =\, G(F) \cap \left( K\times  g K_f^{\cS} g^{-1}\right).
\end{equation*}
The group $\Gamma_{K, g}$ contains only finitely many torsion elements. If it is torsionfree then it acts properly discontinuously on the symmetric space $G_{\infty}/K_{\infty, +}$.

\begin{df}
\label{df:goodsk}
By $\sK$ we denote the set of compact open subgroups $K$ of $G_{\cS}$ such that $\Gamma_{K, g}$ is torsionfree for every $g\in G(\bA_f^{\cS})$. 
\end{df}

Clearly, $\sK$ is closed under conjugation and if $K\in \sK$ and $K'\subseteq K$ is an open subgroup then $K'$ is also contained in $\sK$. Thus $\sK$ satisfies assumption \ref{assum:kgood} because of the following

\begin{prop}
\label{prop:borel}
The set $\sK$ is a cofinal subset of the set of all compact open subgroups of $G_{\cS}$. 
\end{prop}

Note that the reduced norm $\Nrd: \wG\to (\bG_m)_F$ induces homomorphisms 
\begin{equation}
\label{detcenter}
\barnu: G(A) \lra A^*/(A^*)^2 
\end{equation}
for every $F$-algebra $A$. We need the following

\begin{lemma}
\label{lemma:conncomp}
For every compact open subgroup $K$ of $G_{\cS}$ the homomorphism \eqref{detcenter} induces a bijection
\begin{equation*}
\label{conncomp1}
G(F)\backslash G(\bA_f)/(K\times K_f^{\cS})\, \cong\, \coker(F^*\lra  \bA_f^*/(\bA_f^*)^2/\barnu(K\times K_f^{\cS})).
\end{equation*}
In particular there exists a finite number of elements $g_1, \ldots, g_h$ of $G(\bA_f^{\cS})$ such that 
\begin{equation*}
\label{stabdoublecoset}
G(F)\backslash G(\bA_f)/(K\times K_f^{\cS})\, =\, \bigcup_{i=1}^r G(F) (1, g_i) (K\times K_f^{\cS}).
\end{equation*}
holds for every $K$.
\end{lemma}

\begin{proof} The first assertion can be easily deduced from (\cite{deligne}, Th{\'e}or{\`e}me 2.4). The second follows from the fact that the group $\coker(F^*\lra  \left(\bA_f^*/(\bA_f^*)^2 \right)/\barnu(1\times K_f^{\cS}))$ is finite. 
\end{proof} 

\begin{proof}[Proof of Prop.\ \ref{prop:borel}] By Lemma \ref{lemma:conncomp} a compact open subgroup $K$ of $G_{\cS}$ lies in $\sK$ if and only if $\Gamma_{K, g_i}$ is torsionfree for $i=1, \ldots, h$. Since there exists only finitely many torsion elements in each $\Gamma_{K, g_i}$ for a given $K$ we can shrink it so that all groups $\Gamma_{K, g_i}$ become torsionfree.
\end{proof} 

As in section \ref{subsection:gsets} let $\sS=\sS_{G_{\cS}, \sK}$ be the site whose objects are discrete left $G_{\cS}$-sets $S$ such that 
the stabilizer $\Stab_{G_{\cS}}(s)$ of any $s\in S$ lies in $\sK$ and whose morphisms are $G_{\cS}$-equivariant maps. For $S\in \sS$ we define 
\begin{equation}
\label{hbcov}
\wX_S \, =\, S \times \sX^{\cS} \, =\,  S \times G(\bA_f^{\cS})/K^{\cS}_f\times G_{\infty}/K_{\infty, +}. 
\end{equation}
Also for a morphism $\rho\in \Hom_{\sS}(S_1, S_2)$ we define
\begin{equation}
\label{hbcov2}
\rho: \wX_{S_1}^{\an} \lra \wX_{S_2}^{\an}, \,\, (s, x) \mapsto (\rho(s), x).
\end{equation}
The collection of $d$-dimensional complex manifolds \eqref{hbcov} and maps \eqref{hbcov2} defines an $\sS$-space denoted by $\wX$. The embeddings $G(F)\hra G_{\cS}$ and $G(F) \hra G(\bA^{\cS})$ yields a left $G(F)$-action on $\wX$ given by
\begin{equation*}
\label{hbcov3}
\gamma \cdot (s, x) \, =\, (\gamma \cdot s, \gamma \cdot x)
\end{equation*}
for $S\in \sS$, $s\in S$, $x\in \sX^{\cS}$ and $\gamma\in G(F)$. 

\begin{prop}
\label{prop:hbvar}
(a) Passing in \eqref{hbcov} to $G(F)$-orbits yields an $\sS$-space 
\begin{equation*}
\label{hbvar}
X(K_f^{\cS})^{\an} \, =\, G(F)\backslash \wX.
\end{equation*}
(b) If $d\ge 1$ then $X(K_f^{\cS})^{\an}$ is exact.
\end{prop}

\begin{proof} By Lemmas \ref{label:equivariantsspace} and \ref{lemma:uniformsspace} it suffices to verify that the conditions \ref{assum:goodact2} (i) and (if $d\ge 1$) (ii) hold for the $G(F)$-action on $\sX^{\cS}$. The first follows from the fact that $\Gamma_{K, g}$ acts properly discontinuously on $G_{\infty}/K_{\infty, +}$ for every $K\in \sK$ and $g\in G(\bA_f^{\cS})$. For the second note that if $d\ge 1$ then the interior of the set $\{x_{\infty}\in G_{\infty}/K_{\infty, +} \mid \Stab_{G(F)}(x)\ne 1\}$ is empty since $G(F)$ is countable and since every non-trivial element of $G(F)$ has at most two fix points acting on $G_{\infty}/K_{\infty, +}$. Since $\{x\in \sX^{\cS}\mid \Stab_{G(F)}(x)\ne 1\}$ is a subset $G(\bA_f^{\cS})/K^{\cS}_f\times \{x_{\infty}\in G_{\infty}/K_{\infty, +} \mid \Stab_{G(F)}(x)\ne 1\}$ its interior is empty as well. 
\end{proof} 

Unless stated otherwise we assume in the following $d\ge 1$ (the case $d=0$ will be considered again only in section \ref{subsection:shimuracurve}). We will often drop the level $K_f^{\cS}$ from the notation and just write $X^{\an}$ for $X(K_f^{\cS})^{\an}$. Note that $X_S^{\an}$ is a $d$-dimensional complex manifold for every non-empty $S\in \sS$ and that the map $\rho: X_{S_1}^{\an} \to X_{S_2}^{\an}$ induced by a morphism $\rho: S_1\to S_2$ in $\sS$ is holomorphic. As before for $K\in \sK$ we write $X_K^{\an}$ for $X_{G^{\cS}/K}^{\an}=G(F)\backslash \wX_K$. There is a natural action of the group 
\begin{equation*}
\label{pi0ginfty}
\Delta := \, G(F)/G(F)_+ \, \cong \,  G(F_{\infty})/G(F_{\infty})_+ \, \cong \, K_{\infty}/K_{\infty, +}\,
\cong\, \{\pm1\}^d
\end{equation*}
(where $G(F)_+=G(F)\cap G_{\infty, +}$) on $X^{\an}$. It is induced by right multiplication of $K_{\infty}$ on $G_{\infty}/K_{\infty, +}\cong (\bH^{\pm})^d$. As an immediated consequence of Lemma \ref{lemma:conncomp} we obtain the following description of the set of connected components $\pi_0(X_K^{\an})$ of $X_K^{\an}$. 

\begin{lemma}
\label{lemma:conncomp2}
The map \eqref{detcenter} induces a bijection
\begin{equation*}
\label{conncomp2}
\pi_0(X_K^{\an})\, \cong\,  \left(\bA_f^*/(\bA_f^*)^2 \right)/\barnu(K\times K_f^{\cS})
\end{equation*}
for every $K\in \sK$. In particular the set $\pi_0(X_K^{\an})$ stabilizes for small $K$, i.e.\ there exists $K_0\in \sK$ such that the covering $X_K^{\an}\to X_{K_0}^{\an}$ induces a bijection $\pi_0(X_K^{\an})\cong \pi_0(X_{K_0}^{\an})$ for every open subgroup $K$ of $K_0$. 
\end{lemma}

For every $S\in \sS$ there exists a smooth scheme $X_S$ over $\Spec \barQ$ such that $X_S^{\an}$ is the complex manifold associated to $X_S\otimes_{\barQ} \bC$. For $K\in \sK$ the scheme $X_K$ is a 
Shimura variety, the quaternionic Hilbert modular variety of level $K\times K_f^{\cS}$ associated to the reductive group $\Res_{F/\bQ} D^*/F^*$. If $S\in \sS$ has only finitely many $G_{\cS}$-orbits then $X_S$ has a finite number of connected components and it is quasi-projective. The scheme $X_S$ is defined (i.e.\ it admits a canonical model) over the reflex field $L$, i.e.\ the fixed field of the group $\{\tau \in \Gal(\bC/\bQ)\mid \tau \Sigma =\Sigma\}$. It is a subfield of $F^{\gal}$, the normal closure of $F/\bQ$ in $\barQ$. If $\Sigma=\cS_{\infty}$, i.e.\ if $D$ is totally indefinite then we have $L= \bQ$. If $d=1$ (the case of Shimura curves) then $L= F$.

Also for every morphism $\rho: S_1\to S_2$ in $\sS$ the morphism $\rho: X_{S_1}^{\an} \to X_{S_2}^{\an}$
is induced from an {\'e}tale morphism $\rho: X_{S_1} \to X_{S_2}$ of $L$-schemes. If $\rho: S_1\to S_2$ is a surjective then $\rho: X_{S_1} \to X_{S_2}$ is an {\'e}tale covering. Hence $S\mapsto X_S$ is an $\sS$-scheme $X=X(K_f^{\cS})$ defined over $L$. For $K\in \sK$ the variety $X_K$ is projective if and only if $D$ is a skew field, i.e.\ if $\Ram_D\ne \emptyset$. 

Let $R$ be a ring and let $M$ be an $R[G_{\cS}]$-module. We let $\cA_R(M, K_f^{\cS}; R)$ denote the $R[G(F)]$-module of maps 
\[
\psi: M\times G(\bA_f^{\cS})/K_f^{\cS} = M\times \pi_0(\sX^{\cS}) \lra R
\]
that are homomorphism of $R$-modules in the first component.  The $G(F)$-action is given by 
\[
(\gamma \psi)(m, gK_f^{\cS}) \, =\, \psi(\gamma^{-1} m, \gamma^{-1}gK_f^{\cS})
\]
for $\psi\in \cA_R(M, K_f^{\cS}; R)$, $\gamma\in G(F)$, $m\in M$ and $g\in G(\bA_f^{\cS})$. 

The cohomology groups
\begin{equation}
\label{derham}
H^{\bu}(G(F)_+, \cA_R(M, K_f^{\cS}; R)) 
\end{equation}
where considered in \cite{spiess1}. Since $\cA_R(M, K_f^{\cS}; R)$ is an $R[G(F)]$-module the groups \eqref{derham} are naturally $R[\Delta]$-modules. There is also a canonical action of the Hecke algebra $R[G(\bA_f^{\cS})//K_f^{\cS}]$ on \eqref{derham} commuting with the $\Delta$-action. 

\begin{prop}
\label{prop:drmqch3}
Assume that $d\ge 1$. Let $R$ be a noetherian ring and let $M$ be a discrete $R[G_{\cS}]$-module.
\medskip

\noi (a) Suppose that $R=A$ is an Artinian local ring with finite residue field. Then 
there exists a canonical isomorphism 
\begin{equation}
\label{comphbvb}
H_A^n(X_{\barQ}; M, \sF)\, \cong\, H_A^n(X^{\an}; M, \sF^{\an})
\end{equation}
for every constructible sheaf $\sF\in \Sh(X_{\et}, A)$ and every $n\ge 0$. 
\medskip

\noi (b) Assume that $M$ is projective as $R$-module. Then 
there exists a canonical isomorphism of $R[\Delta]$-modules
\begin{equation}
\label{comphbva}
H_R^n(X^{\an}; M, R)\, \cong\, H^n(G(F)_+, \cA_R(M, K_f^{\cS}; R))
\end{equation}
for every $n\ge 0$. 
\end{prop}

\begin{proof} (a) follows from Remark \ref{remarks:trace} (c) and Prop.\ \ref{prop:etancomp}.

For (b) put $\sX_+^{\cS} \, =\, G(\bA_f^{\cS})/K^{\cS}_f\times G_{\infty, +}/K_{\infty, +}$. Note that $\pi_0(\sX_+^{\cS}) = G(\bA_f^{\cS})/K^{\cS}_f$ and that the space $X_K^{\an}$ for $K\in \sK$ can be identified with $G(F)_+\backslash (G_{\cS}/K\times \sX_+^{\cS})$. The assertion now follows from Prop.\ \ref{prop:uni4} (b) since we have $\cA_R(M, K_f^{\cS}; R)\cong \Maps(G(\bA_f^{\cS})/K_f^{\cS}, \Hom_R(M, R))$. 
\end{proof}

\begin{remarks}
\label{remarks:galaction}
\rm (a) As explained in section \ref{subsection:galcohomology} under the assumption in (a) above the cohomology groups $H_A^{\bu}(X_{\barQ}; M, \sF)$ carry a canonical $\fG:=\fG_L = \Gal(\barQ/L)$ action. Note though that in general it does not commute with the $\Delta$-action induced by the isomorphism \eqref{comphbvb}. 
\medskip

\noi (b) If $d=0$ then $X_S \, = \, X_S^{\an}\, =\, G(F)\backslash \wX_S$ is just a discrete space for every $S\in \sS$. In this case $S\mapsto X_S$ may be viewed as a $0$-dimensional $\sS$-space. It will also be denoted by $X=X(K_f^{\cS})$. It is easy to see that the proof of Prop.\ \ref{prop:drmqch3} can be modified so that it also covers the case $d=n=0$, i.e.\ we have 
\begin{equation*}
\label{comphbvatotdeg}
H_R^0(X^{\an}; M, R)\, \cong\, H^0(G(F), \cA_R(M, K_f^{\cS}; R)).
\end{equation*}
\enddemo
\end{remarks}

We will consider in particular the following type of level $K_f^{\cS}=K_0(\fn)^{\cS}\subseteq G(\bA_f^{\cS})$. Let $\fn\ne (0)$ be an ideal of $\cO_F$ that is relatively prime to $\cS\cup \Ram_D$. Let $\cO_D$ be an Eichler order of level $\fn$ in $D$ (if $D= M_2(F)$ then we choose $\cO_D$ to be the subalgebra $M_2(\fn)\subseteq M_2(\cO_F)$ of matrices that are upper triangular modulo $\fn$).  For a nonarchimedean place $\fq$ of $F$ we put $\cO_{D, \fq} = \cO_D\otimes_{\cO_F} \cO_{\fq}$ and define $K_0(\fn)^{\cS}$ to be the image of 
\begin{equation*}
\label{level}
\wK_0(\fn)^{\cS} \, =\, \prod_{\fq\not\in \cS\cup \cS_{\infty}} \cO_{D, \fq}^*
\end{equation*}
under the projection $\wG(\bA_f^{\cS})\to G(\bA_f^{\cS})$. For $K_f^{\cS} = K_0(\fn)^{\cS}$ we write $X^D_0(\fn)^{\cS}$ for the $\sS$-scheme $X(K_0(\fn)^{\cS})$ and $\sA_R(M, \fn; N)$ for $\sA_R(M, K_f^{\cS}; N)$. The $R$-modules appearing in \eqref{comphbvb} (for $\sF= A_X$) and \eqref{comphbva} are equipped with an action of a spherical prime-to-$(\fn, \cS)$ Hecke algebra
\begin{equation*}
\label{heckealg}
\bT \, = \, \bT^{\cS_0}\, =\, \bZ[G(\bA_f^{\cS_0})// K_f^{\cS_0}]
\end{equation*}
and the isomorphism are $\bT$-equivariant. Here $\cS_0 = \cS \cup \{\fq \mid \fd\fn\}$ and $K_f^{\cS_0} = K_0(\cO_F)^{\cS_0}$ is a maximal compact open subgroup of $G(\bA_f^{\cS_0})$. The multiplication in $\bT$ is given by convolution (the Haar measure on $G(\bA_f^{\cS_0})$ is normalized so that $\Vol(K_f^{\cS_0})=1$). The algebra $\bT$ can be identified with the polynomial ring $\bZ[T_{\fq}\mid\fq \not\in \cS_0]$. Here the variable $T_{\fq}$ corresponds to the characteristic function of the double cosets $K_f^{\cS_0}x_{\fq} K_f^{\cS_0}$ where the idele $x_{\fq} \in G(\bA_f^{\cS_0})$ has the component $=1$ at all places except at $\fq$ where its component lies in $\cO_{\fq}$ and has reduced norm $= \Norm(\fq)$.

 Though the definition of the action of $\bT$ on the cohomology groups appearing in \eqref{comphbvb} and \eqref{comphbva} 
is given by standard arguments, for completeness, we will briefly recall it for the groups $H_R^{\bu}((X^D_0(\fn)^{\cS})^{\an}; M, R)$. For a nonarchimedean place of $F$ with $\fq\not\in S_0$ there exists 
two canonical finite {\'e}tale morphisms of $\sS$-schemes 
\begin{equation*}
\label{heckealg2}
\pr_1, \pr_2: X^D_0(\fn\fq)^{\cS} \, \lra \, X^D_0(\fn)^{\cS}
\end{equation*}
induced by the two canonical maps between the levels $K_0(\fn \fq)^{\cS}$ and $K_0(\fn)^{\cS}$ namely the inclusion and the map induced by conjugation with the idele $x_{\fq}$ respectively. The action of $T_{\fq}$ on $H_R^n((X^D_0(\fn)^{\cS})^{\an}; M, \sF)$ is given by the homomorphism 
\begin{equation*}
\label{comphbv3}
H_R^n((X^D_0(\fn)^{\cS})^{\an}; M, R)\stackrel{\pr_1^*}{\lra}  H_R^n((X^D_0(\fn\fq )^{\cS})^{\an}; M, R)
\stackrel{(\pr_2)_*}{\lra}  H_R^n((X^D_0(\fn)^{\cS})^{\an}; M, R)
\end{equation*}
where the first map is the homomorphism \eqref{extpullback1} (with respect to $\pr_1$)
and the second the homomorphism  \eqref{trace2} (with respect to $\pr_2$).

\subsection{Borel-Serre compactification} 
\label{subsection:borelserre}

In this section we consider the case $D=M_2(F)$ (so that $G={\PGL_2}_{/F}$ and $d=[F:\bQ]$). We fix a level of the form $K_0(\fn)^{\cS}\subseteq G(\bA_f^{\cS})$ where $\fn\ne (0)$ is an ideal of $\cO_F$ that is relatively prime to $\cS$ and let $X^{\an}=X(K_0(\fn)^{\cS})^{\an}$. As in section \ref{section:borelind} we denote by $B$ the standard Borel subgroup of $G$ with Levi decomposition $B=TN$. Let $\sK$ be the set of compact open subgroups of $G_{\cS}$ defined in \ref{df:goodsk} and let 
$\sS=\sS_{G_{\cS}, \sK}$. 

For $S\in \sS$ let $X_S^{\BS}$ denote the Borel-Serre compactification of $X_S^{\an}$. It
can be written as 
\begin{equation*}
\label{borelserre1}
X_S^{\BS}\, =\,  G(F)_+\backslash \left( S \times G(\bA_f^{\cS})/K_0(\fn)^{\cS} \times \fH^{\BS}\right)
\end{equation*}
where $\fH^{\BS}$ denotes the partial Borel-Serre compactification of the symmetric space $G_{\infty, +}/K_{\infty, +} \cong \bH^d$ (see \cite{borelJi}, III.5). The space $X_S^{\BS}$ is a compact real analytic manifold with corners. By (\cite{borelJi}, Prop.\ III.5.14) and Lemma \ref{lemma:uniformsspace} the assignment $\cS \mapsto X_S^{\BS}$ 
defines an exact $\sS$-space as well which will be denoted by $X^{\BS}$. It contains $X^{\an}$ as an open $\sS$-subspace. Indeed, $X_S^{\an}$ is an open dense subspace of $X_S^{\an}$ for every $S\in \sS$ and the inclusion $X_S^{\an}\hra X_S^{\BS}$ is a homotopy equivalence.

The boundary of $X^{\BS}$, i.e.\ the complement of $X^{\an}$ in $X^{\BS}$ will be denoted by $\partial X$. It has the following concrete description. Put $B_{\infty, +} = B(F_{\infty}) \cap G_{\infty, +}$ and $T_{\infty, +} = T(F_{\infty}) \cap G_{\infty, +}$. Let $A$ be the maximal split subtorus of $\Res_{F/\bQ} T$ (so $A\cong (\bG_m)_{\bQ}$) and let $A_{\infty}^0\subseteq T_{\infty, +}$ be the connected component of $1$ of the $\bR$-valued points of $A$. If we identify $T_{\infty, +}$ with the group $(\bR_{>0})^d$ via the isomorphism $T(F_{\infty}) = F_{\infty}^* \cong (\bR^*)^d$ (induced by $\sigma_1, \ldots, \sigma_d$) then $A_{\infty}^0$ corresponds to the subgroup $\bR_{>0}$ embedded diagonally into $(\bR^*)^d$. The quotient 
$\fh:= B_{\infty, +}/A_{\infty}^0$ carries a canonical left $B(F)_+ = B(F)\cap G(F)_+$-action. By (\cite{borelJi}, III.5.8) and Lemma \ref{label:equivariantsspace} we have
\begin{equation*}
\label{borelserre2}
\partial X_S\, =\,  B(F)_+\backslash \left(S \times G(\bA_f^{\cS})/K_0(\fn)^{\cS} \times \fh\right)
\end{equation*}
for every $S\in \sS$. Note that the projection $\fh= B_{\infty, +}/A_{\infty}^0\to \barfh:=T_{\infty, +}/A_{\infty}^0$ is a principal $N(F_{\infty})\cong \bR^d$-bundle. The action of the Hecke algebra $\bT^{\cS_0}=\bZ[T_{\fq}\mid \fq\not\in \cS_0]$ (with $\cS_0=\cS\cup\{\fq\mid \fn\}$) on the cohomology of $X$ extends naturally to an action on the cohomology of $X^{\BS}$ and $\partial X$. 

Put $B_{\cS}=B(F_{\cS})$, $T_{\cS}=T(F_{\cS})$ and $N_{\cS}=N(F_{\cS})$ an let $\iota: B_{\cS}\to G_{\cS}$ denote the inclusion and $\pi: B_{\cS}\to T_{\cS}= B_{\cS}/N_{\cS}$ denote the projection. We let $\sK^B$ and $\sK^T$ denote the set of compact open subgroups of $B_{\cS}$ and $T_{\cS}$ respectively. For $S\in \sS^B$, the group $B(F)_+$ acts properly discontinuously from the left on $S \times G(\bA_f^{\cS})/K_0(\fn)^{\cS} \times \fh$ by part (a) of the following Lemma (applied to $\sY= S\times G(\bA_f^{\cS})/K_0(\fn)^{\cS}$).

\begin{lemma}
\label{lemma:boundaryfib}
Let $\cY$ be a discrete left $B(\bA_f)$-set such that $\Stab_{B(\bA_f)}(y)$ is compact for every $y\in \cY$ and put $\cZ= N(\bA_f)\backslash \cY$. 
\medskip

\noi (a) The group $B(F)_+$ acts properly discontinuously on $\cY\times \fh$.
\medskip

\noi (b) The left $T(\bA_f)$-set $\cZ$ is discrete and $\Stab_{T(\bA_f)}(z)$ is compact for every $z\in \cZ$. 
\medskip

\noi (c) The group $T(F)_+$ acts properly discontinuously on $\cZ\times \barfh$.
\medskip

\noi (d) The canonical map 
\begin{equation}
\label{fiberbund}
B(F)_+\backslash(\cY\times \fh) \lra T(F)_+\backslash(\cZ\times \barfh) 
\end{equation}
is a fibration. The fibers are $d$-dimension (real) tori. 
\end{lemma}

\begin{proof} (a) Recall (see \eqref{semidirect}) that there exists a canonical isomorphism between the semidirect product $\bG_m\ltimes \bG_a$ and the Borel subgroup $B$ of $\PGL_2$. Given a fractional ideal $\fa$ of $F$ we denote by $K_f(\fa)$ the subgroup of $B(\bA_f)$ that corresponds to the subgroup $\prod_{v\nmid \infty} \cO_v^*\ltimes \widehat{\fa}$ of $\bA_f^*\ltimes \bA_f$ (here $\widehat{\fa}:=\fa\otimes_{\cO_F}\prod_{v\nmid \infty} \cO_v$). It is easy to see that any compact open subgroup of $B(\bA_f)$ is contained in one of the groups $K_f(\fa)$. Thus the assertion follows from the fact that the group $B(F)_+\cap K_f(\fa)$ acts properly discontinuously on $\fh$.

(b) If $z= N(\bA_f)\cdot y\in \cZ$, $y\in \cY$ then $\Stab_{T(\bA_f)}(z)$ is the image $\Stab_{B(\bA_f)}(y)$ under the projection $B(\bA_f)\to T(\bA_f)$, so $\Stab_{T(\bA_f)}(z)$ is compact and open.

(c) We identify $T$ with $\bG_m$ via \eqref{splittorus}. Since for any compact open subgroup of $\bA_f^*$ is contained in 
$\prod_{v\nmid \infty} \cO_v^*$ the assertion follows from the fact that the group of totally positive units in $\cO_F^*$ act properly discontinuously on $\barfh$ by Dirichlet's unit theorem. 

(d) Since $N(F)\cong F$ is dense in $N(\bA_f)\cong \bA_f$ we have $N(F)\backslash \cY= \cZ$. Hence the left and therefore also the right vertical map in the diagram 
\[
\begin{CD} 
\pi_0(B(F)_+\backslash(\cY\times \fh))@>\cong >> B(F)_+\backslash \cY\\
@VVV @VVV\\
\pi_0(T(F)_+\backslash(\cZ\times \barfh)) @>\cong >> T(F)_+\backslash \cZ
\end{CD}
\] 
is bijective. Let $y\in \cY$, $z= N(\bA_f)\cdot y\in \cZ$ and put $\Gamma_y= B(F)_+\cap \Stab_{B(\bA_f)}(y)$, $\Gamma_z= T(F)_+\cap \Stab_{T(\bA_f)}(z)$. We have seen that the diagram 
\begin{equation}
\label{fiberbund2}
\begin{CD} 
\Gamma_y\backslash\fh@>>> B(F)_+\backslash(\cY\times \fh)\\
@VVV @VV\eqref{fiberbund}V\\
\Gamma_z\backslash\barfh @>\cong >> T(F)_+\backslash(\cZ\times \barfh)
\end{CD}
\end{equation}
is cartesian and that the horizontal maps are open, closed and injective (they are induced by the obvious maps $\fh \to \fh\times \{y\}$, $\barfh \to \barfh\times \{z\}$).
Note that $\Gamma_y = \Gamma_z\ltimes \Lambda$ where $\Lambda:= N(F)\cap \Stab_{N(\bA_f)}(y)$. Hence 
$\Lambda \cong \bZ^d$ is a lattice in $N(F_{\infty})\cong \bR^d$. It follows that the left vertical map in \eqref{fiberbund2} is a fibration whose fiber can be identified with $N(F_{\infty})/\Lambda\cong \bR^d/\bZ^d$.
\end{proof}

Similar to Prop.\ \ref{prop:hbvar} one shows that by setting
 \begin{equation}
\label{borelserre4}
Y_S \, =\,  B(F)_+\backslash \left(S \times G(\bA_f^{\cS})/K_0(\fn)^{\cS} \times \fh\right)
\end{equation} 
for $S\in \sS^B$ we obtain an exact $\sS^B$-space $Y$. Note that the boundary $\partial X$ is equal to the induced $\sS$-space $\caInd_{B_{\cS}}^{G_{\cS}} \,Y$ (see section \ref{subsection:functinG}). Lemma \ref{lemma:indsshn} thus implies

\begin{lemma}
\label{lemma:borserrebound} 
Let $R$ be a ring. We have
\begin{equation*}
\label{borelserre5}
H^n(\partial X_{\infty}, R)\,=\, \Ind_{B_{\cS}}^{G_{\cS}} H^n(Y_{\infty}, R)
\end{equation*}
for every $n\ge 0$.
\end{lemma}

Under certain assumptions on $R$ and $\cS$ the computation of $H^n(\partial X_{\infty}, R)$ can be further simplified. For that we define the (non-exact) $\sS^T$-space $Z$ by
\begin{equation*}
\label{borelserre6}
Z_{S'} \, =\,  T(F)_+\backslash \left(S' \times N(\bA_f^{\cS})\backslash G(\bA_f^{\cS})/K_0(\fn)^{\cS} \times \barfh\right)
\end{equation*} 
for $S'\in \sS^T$. Note that $Z_{S'}$ is a disjoint union of (real) $d-1$-dimensional tori. 

For $S\in \sS^B$, $S'\in \sS^T$ and a $B_{\cS}$-equivariant map $\xi: S\to S'$ we define
\begin{equation*}
\label{borelserre7}
f_{\xi}: Y_S \lra Z_{S'}, \, \, \, B(F)_+\cdot (s, x) \mapsto T(F)_+\cdot (\xi(s), \pr(x))
\end{equation*} 
where $\pr:G(\bA_f^{\cS})/K_0(\fn)^{\cS} \times \fh\to N(\bA_f^{\cS})\backslash G(\bA_f^{\cS})/K_0(\fn)^{\cS} \times \barfh$ denotes the canonical projection. 
If we view $Y$ and $Z$ as $\sS$-objects in $\spaces_{\bg}$ (see the end of section \ref{subsection:ssspacescheme}) then 
the collection of maps $f_{\xi}$ define a $\pi$-compatible morphism $f:Y\to Z$ in the sense of Def.\ \ref
{df:Smorphphi}. Recall that by Prop.\ \ref{prop:bigsmallcohom} viewing $\partial X$, $Y$ and $Z$ either as $\sS$-objects in $\spaces_{\cl}$ or $\spaces_{\bg}$ does not effect the cohomology groups $H^{\bu}(\partial X_{\infty}, R)$, $H^{\bu}(Y_{\infty}, R)$ and $H^{\bu}(Z_{\infty}, R)$.

As in section \ref{subsection:functinG} for $S\in \sS^B$ we put $f_S = f_{\xi_S}: Y_S\to Z_{\barS}$ where $\barS = N_{\cS}\backslash S$ and where $\xi_S: S\to \barS$ is the projection. Also for $K\in \sK^B$ we put $f_K=f_{B_{\cS}/K}:Y_K\to Z_{\barK}$ so that $f_K$ is the obvious map
\begin{equation*}
\label{borelserre8}
f_K: B(F)_+\backslash (G(\bA_f)/(K\times K_0(\fn)^{\cS}) \times \fh) \lra T(F)_+ \backslash \left( N(\bA_f)\backslash G(\bA_f)/(K \times K_0(\fn)^{\cS}) \times \barfh \right).
\end{equation*}
Note that $f_K$ induces a bijection $\pi_0(Y_K) =\pi_0(Z_{\barK})$ and that $\pi_0(Y_K)$ is finite.

By $C_F^{\fn, \cS}$ we denote idele class group
\begin{equation}
\label{rayclass}
C_F^{\fn, \cS} \,=\, \bA_f^*/(F_+^*\prod_{\fq\not\in \cS} U_{\fq}^{(n_{\fq})})
\end{equation}
where $\fn=\prod_{\fq} \fq^{n_{\fq}}$ denotes the prime factorization of $\fn$. For a nonarchimedian place $\fq\not\in \cS$, $\fq\nmid \fn$ we denote by $[\fq]\in C_F^{\fn, \cS}$ the class of an idele with components $=1$ at all places except at $\fq$ where its component is a local uniformizer. 

Since $T$ is abelian the left $T(\bA_f)$-action on $S'\times N(\bA_f^{\cS})\backslash G(\bA_f^{\cS})/K_0(\fn)^{\cS}$ induces a left $T(\bA_f)$-action on $Z_{S'}$ for every $S'\in \sS^T$. By identifying $T$ with $\bG_m$ via \eqref{splittorus} it is easy to see that the  
$T(\bA_f)$-action on $Z$ induces a $C_F^{\fn, \cS}$-action on the cohomology $H^{\bu}(Z_{\infty}, R)$ for any ring $R$. 

\begin{prop}
\label{prop:borelserrefiber}
Let $A$ be an Artinian local ring with finite residue field of characteristic $p$ and assume that one prime $\fp$ in $\cS$ lies above $p$.
\medskip

\noi (a) We have 
\begin{equation}
\label{borelserre9}
H^n(\partial X_{\infty}, A) \, =\, \Ind_{B_{\cS}}^{G_{\cS}} H^n(Z_{\infty}, A)
\end{equation}
for every $n\ge 0$ (here $\Ind_{B_{\cS}}^{G_{\cS}}: \Mod_A^{\sm}(T_{\cS})\to \Mod_A^{\sm}(G_{\cS})$ denotes parabolic induction). 
\medskip

\noi (b) Under the identification \eqref{borelserre9} the action of the Hecke operator $T_{\fq}$, $\fq\not\in \cS_0$ is induced by the map 
\begin{equation*}
\label{borelserre10}
H^n(Z_{\infty}, A) \lra H^n(Z_{\infty}, A),\,\, x\mapsto [\fq]^{-1}\cdot x + \Norm(\fq)\cdot [\fq] \cdot x.
\end{equation*}
\end{prop}

\begin{proof} (a) We will verify that conditions (i) and (ii) of Prop.\ \ref{prop:basechangeiso} hold for the constant sheaf $\sF=A_Y$ on $Y$, i.e.\ we show that 
\begin{itemize}
\item[(i)] For $\barK_0\in \sK^T$ and $K\in \sK^B$ with $\barK=\barK_0$ the base change morphism $\BC_K: (f_* A_Y)_{\barK_0} \to (f_K)_*A_{Y_K}$ is an isomorphism. 
\item[(ii)] For $K\in \sK^B$ there exists an open subgroup $K'\subseteq K$ such that the base change homomorphism $\BC_{K, K'}^n: \barrho^* R^n (f_K)_*(A)\to R^n (f_{K'})_*(A)$ vanishes 
for every $n\ge 1$. 
\end{itemize}
To prove (i) we fix $\barK_0\in \sK^T$ and put $\cJ = \{K\in \sK^B\mid \pi(K) = \barK_0\}$. Note that $(\cJ, \subseteq)$ is a directed set (i.e.\ for every $K_1, K_2\in \cJ$ there exists $K_3\in \cJ$ with $K_1, K_2 \subseteq K_3$). The collection of sheaves $\{(f_K)_*A_{Y_K}\}_{K\in \cJ}$ together with base change morphism $(f_{K_1})_*(A_{X_{K_1}})\to (f_{K_2})_*(A_{X_{K_2}})$ for $K_2\subseteq K_1$ define a projective system over $(\cJ, \subseteq)$. By unwinding the definitions it is easy to see that 
$(f_* A_Y)_{\barK_0} =\prolim_{K\in \cJ} (f_K)_*A_{Y_K}$. Now for $K\in \cJ$ the map $f_K:Y_K \to Z_{\barK_0}$ 
has connected fibers by Lemma \ref{lemma:boundaryfib} (d), hence we have $(f_K)_* A_{Y_K} = A_{Z_{\barK_0}}$ for every $K\in \sK^B$. In particular for any pair $K_1, K_2\in \cJ$ with $K_2\subseteq K_1$ the base change morphism $(f_{K_1})_*(A_{X_{K_1}})\to (f_{K_2})_*(A_{X_{K_2}})$ is an isomorphism. Hence for fixed $K\in \cJ$ we have 
$(f_* A_Y)_{\barK_0} \cong (f_K)_*A_{Y_K} \cong A_{Z_{\barK_0}}$. Thus we also obtain
\begin{equation}
\label{borelserre11}
f_* A_Y \, = \, A_Z.
\end{equation}
To keep the notation simple we give the proof of (ii) only in the case when $\cS$ consists of the single place $\fp$ lying above $p$ (the following argument can be easily adapted to the case of an arbitrary $\cS$). For $r, s\in \bZ$, $r\ge 1$ let $K(r,s)\in \sK^B$ be the image of $U_{\fp}^{(r)}\ltimes \fp^s \cO_{\fp}$ under the isomorphism 
\eqref{semidirect} and $\barK(r)$ the image of $U_{\fp}^{(r)}$ under the isomorphism \eqref{splittorus}. Since any $K\in \sK^B$ contains a group $K(r,s)$ for $r,s$ sufficiently large it suffices to prove (ii) for $K=K(r,s)$.

For $t\ge s$ and fixed $r$ consider the base change homomorphism 
\begin{equation}
\label{borelserre12}
\BC^n: R^n (f_{K(r, s)})_*(A)\to R^n (f_{K(r, t)})_*(A).
\end{equation}
By Lemma \ref{lemma:boundaryfib} (d) its stalk $(\BC^n)_z: (R^n (f_{K(r, s)})_*(A))_z\to (R^n (f_{K(r, t)})_*(A))_z$ at a point $z\in Z_{\barK(r)}$
can be identified with a homomorphism of the form
\begin{equation}
\label{borelserre13}
\pi^*: H^n(\bR^d/\Lambda_1, A) \lra H^n(\bR^d/\Lambda_2, A)
\end{equation}
where $\Lambda_2\subseteq \Lambda_1$ are cocompact lattices and $\pi: \bR^d/\Lambda_2\to \bR^d/\Lambda_1$ is the projection. For $n=1$ the map \eqref{borelserre13} can be identified with 
\[
\iota^*: \Hom(\Lambda_1, A) \lra \Hom(\Lambda_2, A), \,\,\, \varphi\mapsto \varphi\circ \iota
\]
where $\iota: \Lambda_2\hra \Lambda_1$ is the inclusion. Let $p^N$ be the order of $A$. By choosing $t$ sufficiently large we can arrange that $\Lambda_2$ is contained in $\subseteq p^N \Lambda_1$ and $(\BC^1)_z=0$. In fact since the source and target of \eqref{borelserre12} are locally constant sheaves and since $Z_{\barK(r)}$ has only finitely many connected components, the map $\BC^1$ vanishes for $t$ sufficiently large. Moreover the vanishing of $\BC^1$ also implies the vanishing of $\BC^n$ (since $H^{\bu}(\bR^d/\Lambda, A) = \wedge_A^{\bu} H^1(\bR^d/\Lambda, A)$ for any cocompact lattice $\Lambda \subseteq \bR^d$). 

Now (i) and (ii) together with Lemma \ref{lemma:borserrebound}, Prop.\ \ref{prop:basechangeiso} and \eqref{borelserre11} imply 
\begin{equation*}
\label{borelserre14}
H^n(\partial X_{\infty}, A) \, =\, \Ind_{B_{\cS}}^{G_{\cS}} H^n(Y_{\infty}, A)\, =\, \Ind_{B_{\cS}}^{G_{\cS}} H^n(Z_{\infty}, f_* A)\, =\, \Ind_{B_{\cS}}^{G_{\cS}} H^n(Z_{\infty}, A)
\end{equation*}
(b) is proved by unwinding the definitions. 
\end{proof}

\begin{remark}
\label{remark:borelserrecohom}
\rm Note that the manifold $Z_{S'}$ is a disjoint union of (real) $d-1$-dimensional real tori for every $S'\in \cS^T$. Thus under the assumption of Prop.\ \ref{prop:borelserrefiber} we have 
\begin{equation*}
\label{borelserre14a}
H^n(Z_{\infty}, A) \, =\, 0 \qquad \forall \,\,n\ge d-1.
\end{equation*}
If $\cS\subseteq \cS_p$ consists only of primes lying above the prime number $p$ and if $r$ denotes the $\bZ_p$-rank of the closure of $\cO_F^*\cap \prod_{\fp\in \cS} U_{\fp}^{(1)}$ in $\prod_{\fp\in \cS} U_{\fp}^{(1)}$ then one can show that $H^{d-r}(Z_{\infty}, A) \ne 0$ and 
\begin{equation*}
\label{borelserre14b}
H^n(Z_{\infty}, A) \, =\, 0 \qquad \forall \,\,n\ge d-r.
\end{equation*}
Thus if $\cS=\cS_p$ then the vanishing of $H^n(Z_{\infty}, A)$ for $n\ge 1$ is equivalent to Leopold's conjecture (see \cite{hill}, Cor.\ 5 for the case $\cS=\cS_p$). \enddemo
\end{remark}

\subsection{$\varpi$-adic cohomology and group cohomology}
\label{subsection:hmsvarprincipal}

As in section \ref{section:ordinary} we consider the cohomology groups $H_R^n(X^{\an}; M, \sF)$ and $H_R^n(X_{\barQ}; M, \sF)$ in the case where $M$ is a parabolically induced representation and $\cS$ consists of a single place $\fp$ of $F$ lying above a prime number $p$. We will write $G_{\fp}$, $\bA_f^{\fp}$, $K_0(\fn)^{\fp}$, $X=X^D_0(\fn)^{\fp}$ etc.\ instead of $G_{\cS}$, $\bA_f^{\cS}$, $K_0(\fn)^{\cS}$, $X^D_0(\fn)^{\cS}$ etc.\ As before $\bT= \bZ[T_{\fq}\mid \fq \nmid p\fd\fn]$ denotes the spherical prime-to-$p\fd\fn$ Hecke algebra. We denote by $B_{\fp}$ the standard Borel subgroup of $G_{\fp}$ and by $T_{\fp}$ its maximal torus. We identify $T_{\fp}$ with $F_{\fp}^*$ using the isomorphism \eqref{splittorus} and let $T_{\fp}^+$ (resp.\ $T_{\fp}^0$) denote the submonoid (resp.\ the subgroup) that correspond to $\cO_{\fp}\setminus\{0\}$ (resp.\ to $U_{\fp} = \cO_{\fp}^*$). More generally, as in section \ref{subsection:ladm}, for a closed subgroup $H$ of $T_{\fp}^0$ we put $\barT_{\fp}= T_{\fp}/H$, $\barT_{\fp}^+= T_{\fp}^+/H$ and $\barT_{\fp}^0= T_{\fp}^0/H$. 

Let $\cO$ be a complete discrete valuation ring with maximal ideal $\fP= (\varpi)$, finite residue field $k$ of characteristic $p$ and quotient field $E$ of characteristic $0$. Recall that we considered several types of cohomology groups in section \ref{section:ordinary}. In section \ref{subsection:principal} for a locally admissible $\cO_m[\barT_{\fp}]$-module $W$ and an $\cO_m$-sheaf $\sF$ on $X^{\an}$ and $X$ respectively we have defined
\begin{eqnarray*}
\bH_{\cO_m}^{\bu}(X^{\an}; W, \sF) & = & H_{\cO_m}^n(X^{\an}; \Ind_{B_{\fp}}^{G_{\fp}} W, \sF)\\
\bH_{\cO_m}^{\bu}(X_{\barQ}; W, \sF) & = & H_{\cO_m}^n(X_{\barQ}; \Ind_{B_{\fp}}^{G_{\fp}} W, \sF).
\end{eqnarray*}
In section \ref{subsection:varpiadic} for a $\varpi$-adically admissible $\cO[\barT_{\fp}]$-module $W$ with $W_{\tor}=0$ we have introduced the $\varpi$-adic cohomology $\bH_{\varpi-\ad}^n(X^{\an}; W, \wcdot)$ and $\bH_{\varpi-\ad}^n(X_{\barQ}; W, \wcdot)$. Moreover we introduced cohomology groups $\bH_{\varpi-\ad}^n(X^{\an}; V, E)$ and $\bH_{\varpi-\ad}^n(X_{\barQ}; W, E)$ where $V$ an admissible $E$-Banach space representation of $\barT_{\fp}$.

In this section we establish the following variant of the comparison isomorphism \eqref{comphbva} between $\varpi$-adic cohomology and group cohomolgy. 

\begin{prop}
\label{prop:drmqch4}
Assume that $d\ge 1$ and let $W\in \Mod_{\cO}^{\varpi-\adm}(\barT_{\fp})_{\fl}$. There exists a canonical isomorphism of finitely generated augmented $\cO[\barT_{\fp}]$-modules
\begin{equation}
\label{comphbvc}
\bH_{\varpi-\ad}^n(X^{\an}; W, \cO)\, \cong\, H^n(G(F)_+, \cA_{\cO}(\Ind_{B_{\fp}, \cont}^{G_{\fp}}  W, K_0(\fn)^{\fp}; \cO)).
\end{equation}
for every $n\ge 0$.
\end{prop}

Here $\Ind_{B_{\fp}, \cont}^{G_{\fp}} W:= \{\Phi\in C_{\cont}(G_{\fp}, W)\mid \Phi(tng) = t\Phi(g) \,\,\forall t\in T_{\fp}, n\in N_{\fp}, g\in G_{\fp}\}$ where $C_{\cont}(G_{\fp}, W)$ is the $\cO$-module of maps $G_{\fp} \to W$ that are continuous with respect to the $\varpi$-adic topology on $W$.

\begin{proof} Put $M= \Ind_{B_{\fp},\cont}^{G_{\fp}} W$ so that $M_m = \Ind_{B_{\fp}}^{G_{\fp}} W_m$. By Prop.\ \ref{prop:drmqch3} there exists a canonical isomorphism 
\begin{eqnarray*}
&& \bH_{\varpi-\ad}^n(X^{\an}; W, \cO_m)\, =\, \bH_{\cO_m}^n(X^{\an}; W_m, \cO_m) \, \cong\, H^n(G(F)_+, \cA_{\cO_m}(\Ind_{B_{\fp}}^{G_{\fp}} W_m, K_0(\fn)^{\fp}; \cO_m))\\
&&\, =\,H^n(G(F)_+, \cA_{\cO_m}(M_m, K_0(\fn)^{\fp}; \cO_m))\, =\,H^n(G(F)_+, \cA_{\cO}(M, K_0(\fn)^{\fp}; \cO_m))
\end{eqnarray*}
for every $m \ge 0$. Since the left hand side of \eqref{comphbvc} is the projective 
limit of the inverse system $\{\bH_{\varpi-\ad}^n(X^{\an}; W, \cO_m)\}_m$ it suffices to see that the canonical map
\begin{equation}
\label{comphbv2}
H^n(G(F)_+, \cA_{\cO}(M, K_0(\fn)^{\fp}; \cO)) \lra \prolim_m H^n(G(F)_+, \cA_{\cO}(M, K_0(\fn)^{\fp}; \cO_m))
\end{equation}
is an isomorphism for every $n\ge 0$. Note that
\begin{equation*}
\label{gammacoh}
\cA_{\cO}(M, K_0(\fn)^{\fp}; \cO) \, \cong\, \prolim_m \cA_{\cO}(M, K_0(\fn)^{\fp}; \cO_m)
\end{equation*}
and that the transition maps $\cA_{\cO}(M, K_0(\fn)^{\fp}; \cO_{m+1})\to \cA_{\cO}(M, K_0(\fn)^{\fp}; 
\cO_m)$ are surjective for $m\ge 0$. Indeed, they can be identified with the maps
\begin{equation*}
\Maps(\sX_f^{\fp}, \Hom_{\cO_{m+1}}(M_{m+1}, \cO_{m+1}))\to \Maps(\sX_f^{\fp}, \Hom_{\cO_{m+1}}(M_{m+1}, \cO_m))
\end{equation*}
where $\sX_f^{\fp} = G(\bA_f^{\fp})/K_0(\fn)^{\fp}$, so surjectivity follows from the fact that $M_m$ is free as an $\cO_m$-module for every $m\ge 1$. Hence we have $\prolim^{(1)} \cA_{\cO}(M, K_0(\fn)^{\fp}; \cO_m) = 0$. Note also that by Prop.\ \ref{prop:finiteness1} (b) the inverse system 
\begin{equation}
\label{gammacoh2}
\{ H^n(G(F)_+, \cA_{\cO}(M, K_0(\fn)^{\fp}; \cO_m))\}_m\, \cong\, \{ \bH_{\cO_m}^n(X^{\an}; W_m, \cO_m)\}_n
\end{equation}
consists of finitely generated $\La_{\cO}(U_{\fp}^{(1)})$-modules for every $n\ge 0$. Hence by Lemma \ref{lemma:lim1profinite} this implies that the $\prolim^{(1)}$-term in \eqref{gammacoh2} vanishes.
Now we can apply Lemma \ref{lemma:derfunclim} to deduce that \eqref{comphbv2} is an isomorphism. 
\end{proof}

\subsection{Relation to cuspidal automorphic forms}
\label{subsection:automorph}

We keep the notation and assumption of the end of last section, so $\cS$ consists of a single place $\fp$ of $F$ lying above a prime number $p$, $X$ denotes the $\sS$-scheme $X^D_0(\fn)^{\fp}$ and $X^{\an}$ the associated $\sS$-space. We fix a finite extension of $E/\bQ_p$ contained in $\bC_p$ with valuation ring $\cO$, maximal ideal $(\varpi)$, finite residue field $k$. We fix an isomorphism $\bC_p\cong \bC$ so that we can view $E$ as a subfield of $\bC$. 

Let $\bD$ be the totally ramified incoherent quaternion algebra over $\bA:=\bA_F$ (in the sense \cite{yuanzhangzhang}) with ramification set $\Ram_{\bD} = \Ram_D \cup \cS_{\infty}$. Following (\cite{yuanzhangzhang}, bottom of page 70) for a subfield  $E$ of $\bC$ we define the set $\cA(\bD^*/\bA^*, E)$ of automorphic representations of $\bD^*/\bA^*$ over $E$ as the set isomorphism classes of irreducible representations of $\bD^*/\bA^*$ such that $\pi\otimes_E \bC$ is a sum of automorphic representations of $\bD^*/\bA^*$ of weight $0$ (i.e.\ under the Jacquet-Langlands correspondence they are associated to cuspidal automorphic representations of $\PGL_2(\bA)$ whose archimedian components are the discrete series representation of $\PGL_2(\bR)$ of weight $2$ and whose components at nonarchimedian places $v$ that are ramified in $D$ are square-integrable). For $\pi \in \cA(\bD^*/\bA^*, E)$ the field $E_{\pi}: = \End_{E[\bD^*/\bA*]}(\pi)$ is a finite extension of $E$. We denote the valuation ring of $E_{\pi}$ by $\cO_{\pi}$.

Similar to (\cite{yuanzhangzhang}, Thm.\ 3.4 (3)) one shows that every $\pi \in \cA(\bD^*/\bA^*, E)$ has a decomposition $\pi = \bigotimes_{v}'  \pi_v$ and the $\pi_v$ are irreducible admissible representations of $\bD_v$ defined over $E_{\pi}$. As usual we put $\pi_f= \pi^{\infty} =\bigotimes_{v\nmid \infty}' \pi_v$. Also for an ideal $\fa\ne 0$ of $\cO_F$ we put $\pi_f^{\fa}= \bigotimes_{v\nmid \infty}' \pi_v$. Let $\pi \in \cA(\bD^*/\bA^*, E)$ and assume that the conductor $\ff=\ff(\pi_f^{\fp})$ divides $\fn \fd$, i.e.\ we have $(\pi_f^{\fp})^{K_0(\fn)^{\fp}}\ne 0$. Then we get
\[
\dim_{E_{\pi}} (\pi_f^{\fp\fn\fd})^{K_f^{\cS_0}}\,=\,1
\]
where $\cS_0= \{\fq\mid \fp\fn\fd\}$. The spherical Hecke algebra $\bT_{\cO}$ acts on $\pi_{0,\Omega}^{K_0(\fn)^{\fp\fd}}$ via the Hecke eigenvalue homomorphism 
\begin{equation*}
\label{eigenvalue}
\la_{\pi}:\bT_{\cO}\lra {\cO}_{\pi}.
\end{equation*} 
For a $\bT_E$-module $H$ we denote by $H_{\pi}$ the localization of $H$ with respect to the kernel of $(\la_{\pi})_E: \bT_E\lra E_{\pi}$. Also, for a $\bT_E[\Delta]$-module $H$ and a character $\ep: \Delta\to \{\pm\}$ we denote by $H_{\pi, \ep}$ the localization of $H$ with respect to homomorphism $\bT_E[\Delta] \to E$ induced by the pair $(\la_{\pi}, \ep)$.

The aim of this section is to prove the following 

\begin{prop}
\label{prop:automorphpart}
Let $\chi: T_{\fp}\cong F_{\fp}^* \to \cO^*$ be a quasicharacter and $\ep:\Delta \to \{\pm 1\}$ be a character. Let $\pi \in \cA(\bD^*/\bA^*, E)$ such that the conductor of $\pi_f^{\fp}$ divides $\fn \fd$.
\medskip

\noi (a) We have 
\[
\bH_{\varpi-\ad}^n(X^{\an}, E(\chi), E)\, \cong \, \bH_E^n(X^{\an}, E(\chi), E)
\]
for every $n\ge 0$. 
\medskip

\noi (b) We have 
\[
\bH_{\varpi-\ad}^n(X^{\an}, E(\chi), E)_{\pi, \ep} \, \ne \, 0
\]
if and only if $n\ne d, d+1$ and $\pi_{\fp}$ is a subquotient of $\Ind_{B_{\fp}}^{G_{\fp}} E(\chi)$.
\medskip

\noi (c) If $n\in \{d, d+1\}$ and if $\pi_{\fp}$ is a subquotient\footnote{Recall that the representation $\Ind_{B_{\fp}}^{G_{\fp}} E(\chi)$ is irreducible except if $\chi^2=1$.} of $\Ind_{B_{\fp}}^{G_{\fp}} E(\chi)$ then 
the action of the Hecke algebra $\bT_E$ on $\bH_{\varpi-\ad}^n(X^{\an}, E(\chi), E)_{\pi, \ep}$ factors through 
$\la_{\pi}$ so that $\bH_{\varpi-\ad}^n(X^{\an}, E(\chi), E)_{\pi, \ep}$ can be viewed as $E_{\pi}$-vector space (i.e.\ the Hecke action on $\bH_{\varpi-\ad}^n(X^{\an}, E(\chi), E)_{\pi, \ep}$ is semisimple). The dimension of $\bH_{\varpi-\ad}^n(X^{\an}, E(\chi), E)_{\pi, \ep}$ is independent of $\ep$ and $n\in \{d, d+1\}$. Moreover we have 
$\dim_{E_{\pi}}\bH_{\varpi-\ad}^n(X^{\an}, E(\chi), E)_{\pi, \ep} = 1$
if and only if the conductor of $\pi_f^{\fp\fd}$ equals $\fn$.
\end{prop}

\begin{remarks}
\label{remarks:automorphbcoh}
\rm \noi (a) It is easy that the proof below can be modified so that Prop.\ \ref{prop:automorphpart} (a), (b) holds as well if $d=n=0$, i.e.\ if $D$ is a totally definite quaternion algebra. 
\medskip

\noi (b) If $D$ is a division algebra then the description of $\bH_{\varpi-\ad}^n(X^{\an}, E(\chi), E)$ above for $n\in \{d, d+1\}$ can be strengthen as follows. Let $\cR$ be the set of $\pi\in \cA(\bD^*/\bA^*, E)$ such that the conductor of $\pi_f^{\fp}$ divides $\fn \fd$ and so that $\pi_{\fp}$ is a subquotient of $\Ind_{B_{\fp}}^{G_{\fp}} E_{\pi}(\chi)$. We then have 
\begin{equation*}
\label{shimtower2}
\bH_{\varpi-\ad}^n(X^{\an}, E(\chi), E)_{\ep} \, \cong \, \bigoplus_{\pi\in \cR} \, (\pi_f^{\fp})^{K_0(\fn)^{\fp}}.
\end{equation*}
for $n=d, d+1$ and every character $\ep: \Delta\to \{\pm 1\}$.
\medskip

\noi (c) By (\cite{carayol}, \cite{taylor}) there exists a twodimensional $\varpi$-adic Galois representation $V=V_{\pi}$ associated to $\pi$, i.e.\ there exists a twodimensional $E$-vector space $V$ together with a continuous homomorphism $\rho=\rho_{\pi}: \Gal(\barQ/F)\to \GL(V)$, so that $\rho$ is unramified outside the primes dividing $p\cdot \fn$ and so that 
\begin{equation*}
\label{galrephmfclass}
\Tr(\rho(\Frob_{\fq})) = \la_{\fp}(T_{\fq}), \qquad \det( \rho(\Frob_{\fq})) \, =\, \Norm(\fq).
\end{equation*}
for all $\fq\nmid p\cdot \fn$ (where $\Frob_{\fq} \in \Gal(\barQ/F)$ denotes a Frobenius for the prime $\fq$). Following (\cite{nekovar}, 5.12) we denote by $\Ind_{\Sigma}^{\otimes} \rho$ the (partial) tensor induction of $\rho$ (so $\Ind_{\Sigma}^{\otimes} \rho$ is a $2^d$-dimensional representation of $\fG_L$). If $F^{\gal}$ denotes the normal closure of $F/\bQ$ in $\barQ$ and $\widetilde{\sigma}_1, \ldots, \widetilde{\sigma}_d\in \Gal(\barQ/F^{\gal})$ are lifts of the embeddings $\sigma_1, \ldots, \sigma_d: F\to \barQ$ then the restriction of $\Ind_{\Sigma}^{\otimes} \rho$ to $\Gal(\barQ/F^{\gal})$ is the representation $\bigotimes_{i=1}^d V^{\widetilde{\sigma}_i}$. If $n\ne d, d+1$ and $\pi_{\fp}$ is a subquotient of $\Ind_{B_{\fp}}^{G_{\fp}} E(\chi)$ then one can show that there exists an isomorphism of $E[\fG_L]$-modules
\[
\bH_{\varpi-\ad}^n(X_{\barQ}, E(\chi), E)_{\pi} \, \cong \,(\Ind_{\Sigma}^{\otimes} \rho)^m
\]
where $m=\dim_E(\bH_{\varpi-\ad}^n(X^{\an}, E(\chi), E)_{\pi, \ep})$.
\enddemo
\end{remarks}

We need the following

\begin{lemma}
\label{lemma:cohhmsdlim}
Let $W\in \Mod_{\cO, f}^{\sm}(T_{\fp})$ and assume that $W$ is free as an $\cO$-module. Then the natural map 
\[
\bH_{\cO}^n(X^{\an}; W, \cO)_E \lra \bH_E^n(X^{\an}; W_E, E)
\]
is an isomorphism for every $n\ge 0$.
\end{lemma}

\begin{proof} We first remark that for every $K\in \sK$ the natural map 
\begin{equation}
\label{cohhmsdlim1}
H^n(X_K^{\an}, \cO)_E \lra H^n(X_K^{\an}, E)
\end{equation}
is an isomorphism. This follows from (\cite{iversen}, Ch.\ III, Thm.\ 5.1) if $D$ is a division algebra. If $D$ is the matrix algebra $M_2(F)$, we use the fact that the corresponding statement holds for the Borel-Serre compactification $X_K^{\BS}$ and that the inclusion $\io: X_K^{\an}\to X_K^{\BS}$ is a homotopy equivalence. In fact the first vertical and the horizontal maps in the diagram
\begin{equation*}
\label{cohhmsdlim2}
\begin{CD}
H^n(X_K^{\BS}, \cO)_E @> \io^* >> H^n(X_K^{\an}, \cO)_E \\
@VVV @VV \eqref{cohhmsdlim1} V\\
H^n(X_K^{\BS}, E) @> \io^* >> H^n(X_K^{\an}, E) 
\end{CD}
\end{equation*}
are isomorphism by (\cite{iversen}, Ch.\ III, Thm.\ 5.1) and (\cite{iversen}, Ch.\ IV, Thm.\ 1.1). Hence \eqref{cohhmsdlim1} is an isomorphism as well. 

Next we show that the canonical map 
\begin{equation}
\label{cohhmsdlim3}
\fOrd_{\cO}^{H, n}(X, \cO)_E \lra \fOrd_{\cO}^{H, n}(X, E)
\end{equation}
is an isomorphism for every $n\ge 0$. Here $H$ is an arbitrary open subgroup of $T_{\fp}^0$. Put $\barT_{\fp}^+ =T_{\fp}^+/H$, $\barT_{\fp} =T_{\fp}/H$ and let $m\ge 1$ such that $\delta(U^{(m)})\subseteq H$. Note that the fact that \eqref{cohhmsdlim1} is an isomorphism and that $H^n(X_{K_H(m)}, \cO)$ is finitely generated as an $\cO$-module implies that $H^n(X_{K_H(m)}, E)$ is locally admissible as $\cO[\barT_{\fp}^+]$-module. Hence by Propositions \ref{prop:adTTplus} (d), \ref{prop:rgammaord} (b) and \ref{prop:finiteness2} (a) we have 
\[
\fOrd_{\cO}^{H, n}(X, E)\, \cong\, \Gamma_{\cO}^{\ord}(H^n(X_{K_H(m)}, E))\, \cong\, \Gamma_{\cO}^{\ord}(H^n(X_{K_H(m)}, \cO))_E\, \cong\, \fOrd_{\cO}^{H, n}(X, \cO)_E.
\]
Now choose $H$ with $W^H= W$ and consider the natural morphism of spectral sequences 
\begin{eqnarray}
\label{fordextssmap}&& \left( \, \Ext_{\cO, \barT_{\fp}}^r(W^{\io}, \fOrd_{\cO}^{H, s}(X, \cO))_E \, \Longrightarrow\, \bH_{\cO}^{r+s}(X; W, \cO)_E \, \right) \, \lra \, \\
&& \hspace{3cm} \left( \, \Ext_{\cO, \barT_{\fp}}^r(W^{\io}, \fOrd_{\cO}^{H, s}(X, E)) \, \Longrightarrow\, \bH_{\cO}^{r+s}(X; W, E) \, \right).
\nonumber
\end{eqnarray}
Here the source and target spectral sequence is the spectral sequence \eqref{fordextss} for $\sF=\cO$ (tensored with $E$) and $\sF=E$ respectively. Note that the functors $\Ext_{\cO, \barT_{\fp}}^r(W^{\io}, \wcdot)= \Ext_{\cO[\barT_{\fp}]}^r(W^{\io}, \wcdot)$ for $r\ge 0$ commute with direct limits (since $\cO[\barT_{\fp}]$ is noetherian). Together with the fact that \eqref{cohhmsdlim3} is an isomorphism this implies that \eqref{fordextssmap} is an isomorphism on $E_2$-terms. Hence the first of the canonical maps
\[
\bH_{\cO}^n(X^{\an}; W, \cO)_E \lra \bH_{\cO}^n(X^{\an}; W, E)\lra \bH_E^n(X^{\an}; W_E, E)
\]
is an isomorphism for every $n\ge 0$. That the second is an isomorphism as well follows from Prop.\ \ref{prop:drmqch}.
\end{proof}

We also use the following elementary 

\begin{lemma}
\label{lemma:extgmquasi}
Let $\chi_1, \chi_2: T_{\fp}\cong F_{\fp}^*\to E^*$ be quasicharacters. Then we have 
\begin{equation*}
\Ext_{E, T_{\fp}}^n(E(\chi_1), E(\chi_2))\, =\,  \left\{ \begin{array}{cc} E & \mbox{if $n=0,1$ and $\chi_1=\chi_2$}\\
0 & \mbox{otherwise.}
\end{array}\right.
\end{equation*}
\end{lemma}

\begin{proof}[Proof of Prop.\ \ref{prop:automorphpart}] (a) follows from Prop.\ \ref{prop:compwwhat} and the Lemma \ref{lemma:cohhmsdlim}. 
For (b) and (c) we use the covering spectral sequence (see Prop.\ \ref{prop:hs} (c))
\begin{equation*}
\label{comphbv1}
\Ext_{E, G_{\fp}}^r(\Ind_{B_{\fp}}^{G_{\fp}} E(\chi), H^s(X^{\an}_{\infty}, E)) \Rightarrow  H_E^{r+s}(X^{\an}, \Ind_{B_{\fp}}^{G_{\fp}} E(\chi), E) = \bH_E^{r+s}(X^{\an}, E(\chi), E).
\end{equation*}
By Matsushima's Theorem (see e.g.\ \cite{harder}, \cite{reimann}) we have 
\begin{equation*}
\label{harderreimann}
H^s(X^{\an}_{\infty}, E)_{\pi}  \, =\, \left\{ \begin{array}{cc} (\pi_{\fp}\otimes (\pi_f^{\fp})^{K_0(\fn)^{\fp}}) \otimes_E E[\Delta] & \mbox{if $n=d$,}\\
 0 & \mbox{if $n\ne d$.}
\end{array}\right.
\end{equation*}
Put $V= \pi_{\fp}$, $m= \dim_{E_{\pi}} (\pi_f^{\fp})^{K_0(\fn)^{\fp}}$ and $E'= E_{\pi}$. We obtain 
\begin{eqnarray*}
\label{shimtower3}
&& \bH_E^n(X^{\an}, E(\chi), E)_{\pi, \ep} \, \cong \, \Ext_{E, G_{\fp}}^{n-d}(\Ind_{B_{\fp}}^{G_{\fp}} E(\chi), H^d(X^{\an}_{\infty}, E)_{\pi, \ep})\\
&&  \cong \, \Ext_{E', G_{\fp}}^{n-d}( \Ind_{B_{\fp}}^{G_{\fp}} E'(\chi) , V)^m \,  \cong\, \Ext_{E', T_{\fp}}^{n-d}( E'(\chi^{-1}), \Ord_{E'}(V))^m.
\nonumber\end{eqnarray*}
Therefore by Remark \ref{remarks:emertonord} (e) and Lemma \ref{lemma:extgmquasi} we obtain
\begin{equation*}
\label{shimtower4}
\bH_E^n(X^{\an}, E(\chi), E)_{\pi, \ep} \, \cong \, \left\{ \begin{array}{cc} (E')^m & \mbox{if $n=d,d+1$ and if $\pi_{\fp}$ is a subquotient of $\Ind_{B_{\fp}}^{G_{\fp}} E'(\chi)$,}\\
 0 & \mbox{otherwise.}
\end{array}\right.\nonumber
\end{equation*}
for every character $\ep:\Delta\to \{\pm 1\}$ and $n\ge 0$. This together with (a) implies (b) and (c).  
\end{proof}

\subsection{The case $d=0$ and $d=1$}
\label{subsection:shimuracurve}

In this following we assume that $\cO$ is a complete noetherian local ring of dimension $\le 1$ with finite residue field of characteristic $p$. In this section we mainly consider the case  $d=1$ and $d=0$, i.e.\ $D$ is totally definite in which case $X_K=X_K^{\an}$ is a finite discrete space for each $K\in \sK$. We will study the 
cohomology groups \eqref{ordsheaf2} and \eqref{parindmod} for $n=0$ and $n=1$ with trivial coefficients $\sF=\cO_{X^{\an}}$. We begin with

\begin{lemma}
\label{lemma:ordnull}
Assume that $d\ge 1$.
\medskip

\noi (a) The group $H^0(X_{\infty}^{\an}, \cO)$ is free and of finite rank as an $\cO$-module. Moreover the
$G_{\fp}$-action on $H^0(X_{\infty}^{\an}, \cO)$ is trivial. 
\medskip

\noi (b) We have $\fOrd_{\cO}^{H, 0}(X^{\an}, \cO) =0$ for every closed subgroup $H$ of $T_{\fp}^0$.
\end{lemma}

\begin{proof} (a) By Lemma \ref{lemma:conncomp} we have  
\[
H^0(X_{\infty}^{\an}, \cO)\, =\, \dlim_K \Maps( \pi_0(X_K^{\an}), \cO) \, =\,  \Maps( \pi_0(X_{K_0}^{\an}), \cO)
\]
for $K_0\in \sK$ sufficiently small. Now the assertion follows from the fact that $\pi_0(X_{K_0}^{\an})$ is a finite set.

For (b) note that $\fOrd_{\cO}^{H, 0}(X^{\an}, \cO) = \Ord^H(H^0(X_{\infty}^{\an}, \cO))$ so the assertion follows from (a) and Remark \ref{remarks:emertonord} (d).
\end{proof}

Let $H_1$, $H_2$ be closed subgroups of $T_{\fp}^0$ with $H_2\subseteq H_1$ and put $\Gamma=H_1/H_2$. The functor $\fOrd^{H_1}(X^{\an}, \wcdot)$ is equal to the composite 
\begin{equation*}
\label{fOrddelta}
\begin{CD}
S(X^{\an}, \cO)@> \fOrd^{H_2}(X^{\an}, \wcdot) >> \Mod_{\cO}^{\ladm}(T_{\fp}/H_2) @> H^0(\Gamma, \wcdot) >>  \Mod_{\cO}^{\ladm}(T_{\fp}/H_1)
\end{CD}
\end{equation*}
If $d\ge 1$ then the first functor has an exact left adjoint so there exists a corresponding Grothendieck spectral sequence 
\begin{equation}
\label{fOrddeltass}
E_2^{rs} = H^r(\Gamma, \fOrd_{\cO}^{H_2, s}(X^{\an}, \sF)) \, \Longrightarrow\, \fOrd_{\cO}^{H_1, s}(X^{\an}, \sF).
\end{equation}
Here $H^{\bu}(\Gamma, \wcdot)$ denotes the cohomology of the profinite group $\Gamma$ computed using continuous cochains. Indeed, by Prop.\ \ref{prop:rnladm} the $r$-th derived functor of $\Mod_{\cO}^{\ladm}(T_{\fp}/H_2) \to  \Mod_{\cO}^{\ladm}(T_{\fp}/H_1), W \mapsto W^{\Gamma}$ is isomorphic to the functor $W \mapsto H^r(\Gamma, W)$. 

By Lemma \ref{lemma:ordnull}, Remark \ref{remarks:fordhida} (b), Prop.\ \ref{prop:finiteness2} (a) and by considering the five term exact sequence associated to \eqref{fOrddeltass} and we obtain 

\begin{lemma}
\label{lemma:ordnulleins}
(a) The canonical map 
\begin{equation}
\label{ordnull}
\fOrd_{\cO}^{H_1, n}(X^{\an}, \cO) \lra \fOrd_{\cO}^{H_2, n}(X^{\an}, \cO)^{\Gamma}
\end{equation}
is an isomorphism for $n=0$ or for $n = 1$ if $d\ge 1$.
\medskip

\noi (b) Let $H$ be a closed subgroup of $T_{\fp}^0$. The $\cO[T_{\fp}/H]$-module $\fOrd_{\cO}^{H, n}(X^{\an}, \cO)$ is admissible for $n=0$ and for $n=1$ if $d\ge 1$. 
\end{lemma}

Note that both groups in \eqref{ordnull} vanish if $n=0$ and $d\ge 1$. 

\begin{lemma}
\label{lemma:ordzwei}
Assume that $d=1$ and let $H$ be a closed subgroup of $T_{\fp}^0$. Then 
\[
\fOrd_{\cO}^{H, n}(X^{\an}, \cO)\, =\, 0
\]
for every $n\ge 2$ and every closed subgroup $H$ of $T_{\fp}^0$.  
\end{lemma}

\begin{proof} By Remark \ref{remarks:fordhida} it suffices to consider the case where $n=2$ and $H$ is open in $T_{\fp}^0$. Now if we choose $m\ge 1$ such that $\delta(U_{\fp}^{(m)})\subseteq H$ then we have to show 
\begin{equation*}
\label{fordd1}
\fOrd_{\cO}^{H, 2}(X^{\an}, \cO) \, =\,  H^2(X_{K_H(m)}^{\an}, \cO)^{\ord} \, =\, 0.
\end{equation*}
The assumption $d=1$ implies that either (i) $F=\bQ$ and $D = M_2(\bQ)$ (in which case $X_{K}^{\an}$ is an open modular curve for each $K$) or (ii) $D$ is a division algebra (when each $X_{K}^{\an}$ is a Shimura curve). In the first case we have $H^2(X_{K}^{\an}, \cO)=0$ for 
every $K\in \sK$. It remains to consider the case when $D$ is a division algebra. Then we have 
\[
H^2(X_{K}^{\an}, \cO)\, \cong \, \Maps(\pi_0(X_{K}^{\an}), \cO).
\]
Let $t_0\in T_{\fp}^+$ be the image of a uniformizer of $F_{\fp}$ under the isomorphism \eqref{splittorus} and let $i\ge 0$. Since 
$K_m(H) \cap K_m(H)^{t_0^{-1}} = K_{m+1}(H)$, the action of $t_0$ on $H^i(X_{K_H(m)}^{\an}, \cO)$ is given by 
\begin{equation}
\label{heckemackey}
H^i(X_{K_H(m)}^{\an}, \cO)  \stackrel{\pi_1^*}{\lra}  H^i(X_{K_H(m+1)}^{\an}, \cO)  \stackrel{(\sigma_{t_0})_*}{\lra}  
H^i(X_{K_H(m+1)^{t_0}}^{\an}, \cO) 
 \stackrel{(\pi_2)_*}{\lra}  H^i(X_{K_H(m)}^{\an}, \cO).
\end{equation}
Here $\pi_1$ and $\pi_2$ denote the coverings $X_{K_H(m+1)}^{\an} \to X_{K_H(m)}^{\an}$ and $X_{K_H(m+1)^{t_0}}^{\an}\to X_{K_H(m)}^{\an}$ corresponding to the inclusion $K_H(m+1)\hra K_H(m)$ and $K_H(m+1)^{t_0}\hra K_H(m)$ respectively.

By Lemma \ref{lemma:conncomp} we can identify the group of connected components of $X_{K_H(m)}^{\an}$, $X_{K_H(m+1)}^{\an}$ and $X_{K_H(m+1)^{t_0}}$ with the same finite quotient $\pi_0$ of the group $\bA_f^*/(\bA_f^*)^2$ (since $\barnu(K_H(m)) = \barnu(K_H(m+1))= \barnu(K_H(m+1)^{t_0})$). For $i=2$ the first homomorphism in \eqref{heckemackey} can be identified with the map $\Maps(\pi_0, \cO)\to \Maps(\pi_0, \cO), x\mapsto [K_H(m): K_H(m+1)]\cdot x$ whereas the second and third with the identity map on $\Maps(\pi_0, \cO)$. Since $[K_H(m): K_H(m+1)]= \Norm(\fp)$ is a power of $p$ we conclude that $H^2(X_{K_H(m)}^{\an}, \cO)^{\ord}=0$.
\end{proof}

For the remainder of this section we assume that $d=0$ or $d=1$ and also that $\cO$ is a complete discrete valuation ring with maximal ideal $\fm= (\varpi)$, finite residue field $k$ of characteristic $p$ and quotient field $E$ of characteristic $0$. As before for an $\cO$-module $N$ and $m\ge 1$ we put $N_m= N\otimes_{\cO} \cO_m$ where $\cO_m = \cO/\fm^m$. 

\begin{prop}
\label{prop:ordvarpiadm}
Assume that $d=0$ or $d=1$. Let $H$ be a closed subgroup of $T_{\fp}^0$ and put $\barT_{\fp} = T_{\fp}/H$.
\medskip

\noi (a) If $d=1$ then there are short exact sequences 
\begin{equation}
\label{fordkummer}
\begin{CD}
0 @>>> \fOrd_{\cO}^{H, n}(X^{\an}, \cO) @> \varpi^m \cdot >>  \fOrd_{\cO}^{H, n}(X^{\an}, \cO) @>>> \fOrd_{\cO_m}^{H, n}(X^{\an}, \cO_m) @>>> 0
\end{CD}
\end{equation}
for every $m\ge 1$ and $n\ge 0$. If $d=0$ then \eqref{fordkummer} is exact for every $m\ge 1$ and $n=0$.
\medskip

\noi (b) The $\cO[\barT_{\fp}]$-module 
\begin{equation}
\label{fordpiadic}
\fOrd_{\varpi-\ad}^{H, d}(X^{\an}, \cO): =\, \prolimm \fOrd_{\cO_m}^{H, d}(X^{\an}, \cO_m) 
\end{equation}
is $\varpi$-adically admissible. It is isomorphic to the $\varpi$-adic completion of $\fOrd_{\cO}^{H, d}(X^{\an}, \cO)$ and we have $\fOrd_{\varpi-\ad}^{H, d}(X^{\an}, \cO)_{\tor} =0$. Moreover if $d=1$ then $\fOrd_{\varpi-\ad}^n(X^{\an}, \cO)=0$ for all $n\ne 1$.
\medskip

\noi (c) The $\cO[\barT_{\fp}]$-module $\fOrd_{\cO}^{H, d}(X^{\an}, \cO)$ is admissible. It is free as an $\cO$-module.
\medskip

\noi (d) If $H$ is open in $T_{\fp}$ and if $n\ge 1$ such that $\delta(U_{\fp}^{(n)})\subseteq H$ then we have
\begin{equation*}
\label{fordpiadic2}
\fOrd_{\varpi-\ad}^{H, d}(X^{\an}, \cO) \,= \, \fOrd_{\cO}^{H, d}(X^{\an}, \cO) \, =\, H^d(X_{K_H(n)}^{\an}, \cO)^{\ord}.
\end{equation*}
\noi (e) Assume that $d=1$. For $W\in \Mod_{\cO}^{\varpi-\adm}(\barT_{\fp})$ we have 
\begin{equation*}
\label{ordancomplete1}
\bH_{\varpi-\ad}^n(X^{\an}; W, \cO) \, = \, \prolimm \Ext_{\cO_m, \barT_{\fp}}^{n-1}(W_m^{\io}, \fOrd^1(X^{\an}, \cO_m))
\nonumber 
\end{equation*}
for every $n\ge 0$. In particular if $V\in \Ban_E^{\adm}(\barT_{\fp})$ then we get 
\begin{equation*}
\label{ordancomplete2}
\bH_{\varpi-\ad}^1(X_{\barQ}; V, E)\, =\,  \Hom_{\Ban_E(\barT_{\fp})}(V^{\io},\fOrd_{\varpi-\ad}^1(X^{\an}, \cO)_E).
\end{equation*}
\end{prop}

\begin{proof} (a) First assume $d=1$. The short exact sequence $0\lra \cO\stackrel{\varpi^m \cdot}{\lra}\cO\lra \cO_m\lra 0$ induces a long exact sequence 
\[
\ldots \to \fOrd_{\cO}^{H, n}(X^{\an}, \cO) \to  \fOrd_{\cO}^{H, n}(X^{\an}, \cO) \to \fOrd_{\cO}^{H, n}(X^{\an}, \cO_m) \to \fOrd_{\cO}^{H, n+1}(X^{\an}, \cO) \to \ldots 
\]
Moreover if $d=1$ then by Cor.\ \ref{corollary:ordnxcoeff} we have $\fOrd_{\cO}^{H, n}(X^{\an}, \cO_m)=\fOrd_{\cO_m}^{H, n}(X^{\an}, \cO_m)$ for every $m\ge 1$ and $n\ge 0$. Thus the assertion follows from Lemmas \ref{lemma:ordnull}, \ref{lemma:ordzwei} and Remark \ref{remarks:fordhida} (b) 
since the groups $\fOrd_{\cO}^{H, n}(X^{\an}, \cO)$ and $\fOrd_{\cO_m}^{H, n}(X^{\an}, \cO_m)$ vanish except for $n=1$. For $d=0$ the assertion can be proved easily using formula \eqref{fordlim2} (for $n=0$ of Remark \ref{remarks:fordhida} (a).

(b) follows immediately from (a) and Prop.\ \ref{prop:finiteness1}, (c) follows from Lemmas \ref{lemma:ordnulleins} and \ref{lemma:torsionfree} and (d) follows from (b) and Prop.\ \ref{prop:finiteness2} (a).

For (e) note that by (b) and Lemma \ref{prop:ford} (b) we have
\begin{eqnarray*}
\label{ordancomplete3}
&& \bH_{\varpi-\ad}^n(X^{\an}; W, \cO) \, = \, \prolimm \bH_{\varpi-\ad}^n(X^{\an}; W, \cO_m) \, =\, \prolimm \bH_{\cO_m}^n(X^{\an}; W_m, \cO_m)\\
&& =\, \prolimm \Ext_{\cO_m, T_{\fp}}^{n-1}(W_m^{\io}, \fOrd^1(X^{\an}, \cO_m))
\nonumber 
\end{eqnarray*}
for every $n\ge 0$. For $n=1$ and $V\in \Ban_E(\barT_{\fp})$ we choose $W\in \Mod_{\cO}^{\varpi-\adm}(\barT_{\fp})$ with $V=W_E$. By Lemma \ref{lemma:locbanach} we have 
\begin{eqnarray*}
\label{ordancomplete4}
&& \bH_{\varpi-\ad}^1(X^{\an}; V, E)\, = \, \bH_{\varpi-\ad}^1(X^{\an}; W, \cO)_E \, =\, \left(\prolimm \Hom_{\cO[\barT_{\fp}]}(W^{\io}, \fOrd^1(X^{\an}, \cO_m))\right)_E \\
&&  = \, \Hom_{\cO[T_{\fp}]}(W^{\io},\fOrd_{\varpi-\ad}^1(X_{\barQ}, \cO))_E\, =\, \Hom_{\Ban_E(\barT_{\fp})}(V^{\io},\fOrd_{\varpi-\ad}^1(X^{\an}, \cO)_E).
\end{eqnarray*}
\end{proof}

As in section \ref{subsection:admbanach} we can consider the dual \eqref{pontrajagin3} of $\fOrd_{\varpi-\ad}^{H, n}(X^{\an}, \cO)$ and the dual \eqref{pontrajagin4} of the admissible $E$-Banach space representation $\fOrd_{\varpi-\ad}^{H, n}(X^{\an}, \cO)_E$ of $\barT_{\fp}$.  

\begin{prop}
\label{prop:freetotdeg}
Assume that $d\le 1$ and that $D$ is a division algebra. Suppose also that $H$ is a closed subgroup of $T_{\fp}^0$ such that $\barT_{\fp}^0 = T_{\fp}^0/H$ is a pro-$p$-group. If we put $\La= \La_{\cO}(\barT_{\fp}^0)$ then we have
\medskip

\noi (a) The dual $\cD(\fOrd_{\varpi-\ad}^{H, n}(X^{\an}, \cO)_E)$ is a finitely generated projective $\La_E$-module. 
\medskip

\noi (b) If moreover $K_H(1)\in \sK$ then $\sD(\fOrd_{\varpi-\ad}^{H, d}(X^{\an}, \cO))$ is free of finite rank as a $\La$-module.
\end{prop}

The proof requires some preparation (see also \cite{emerton1} for a similar line of arguments). 

\begin{lemma}
\label{lemma:cohomfree}
Let $\Gamma$ be a finite abelian $p$-group and let $A$ be $\cO[\Gamma]$-module such that

\noi (i) $\whH^0(\Gamma, A_1)=0$.

\noi (ii) $A$ is free as an $\cO$-module.

\noi Then $A$ is a free $\cO[\Gamma]$-module.
\end{lemma}

Here $\whH^{\bu}(\Gamma, A)$ denotes Tate cohomology.

\begin{proof} This result is a variant of (\cite{atiyah-wall}, Thm.\ 8). Firstly, since $\cO[\Gamma]$ is a local ring it suffices to see that $A$ is a projective $\cO[\Gamma]$-module. By (\cite{atiyah-wall}, Thm.\ 6), assumption (i) implies that $A_1$ is a free $k[\Gamma]$-module. By (\cite{atiyah-wall}, Thm. 7) and condition (ii) we deduce that $A$ is cohomological trivial, i.e.\ $\whH^q(\Gamma', A)=0$ for every $q\in \bZ$ and every subgroup $\Gamma'$ of $\Gamma$. We can now argue as in the proof of (\cite{atiyah-wall}, Theorem 8): We choose an exact sequence 
\begin{equation}
\label{resolution}
0\lra Q\lra F\lra A\lra 0
\end{equation} 
where $F$ is a free $\cO[\Gamma]$-module. By condition (ii) the sequence
\[
0\lra \Hom_{\cO}(A, Q) \lra  \Hom_{\cO}(A, F) \lra  \Hom_{\cO}(A, A) \lra 0
\]
is exact. The fact that $Q$ is free as an $\cO$-module and that $A_1$ is a free $k[\Gamma]$-module implies that 
\[
\Hom_{\cO}(Q,A)\otimes_\cO k\, \cong\, \Hom_k(Q_1, A_1)
\]
is an induced $k[\Gamma]$-module, hence the $\Gamma$-module $\Hom_{\cO}(A, Q)$ is cohomologically trivial (again by \cite{atiyah-wall}, Thm.\ 7). 
Therefore $H^1(\Gamma, \Hom_{\cO}(A, Q))=0$, so it follows that $\Hom_{\Gamma}(F, A) \to \Hom_{\Gamma}(A,A)$ is surjective. Consequently, the sequence \eqref{resolution} splits and $A$ is a direct summand of $F$. 
\end{proof}

\begin{lemma}
\label{lemma:cofree}
Let $H$ be a closed subgroup of $T_{\fp}^0$ such that $\barT_{\fp}^0 = T_{\fp}^0/H$ is a pro-$p$-group. Let $W$ be an admissible $\cO[\barT_{\fp}]$-module and put $\hatW = \prolimm W_m$. Assume that
\medskip

\noi (i) $W$ is free as an $\cO$-module.

\noi (ii) $W_1$ is an admissible $k[\barT_{\fp}]$-module.

\noi (iii) $W^U$ is free as an $\cO[T_{\fp}^0/U]$-module for every open subgroup $U$ of $T_{\fp}^0$ with $U\supseteq H$.

\noi Then we have 
\medskip

\noi (a) $\hatW$ is a $\varpi$-adically admissible $\cO[\barT_{\fp}]$-module and the inclusion $W\hra \hatW$ induces an isomorphism $\sD(\hatW)\to \sD(W) : =\Hom_{\cO}(W, \cO)$.
\medskip

\noi (b) $\sD(\hatW)$ is free of finite rank as a $\La$-module. 
\end{lemma}

\begin{proof} (a) The first statement is obvious and for the second note that $\hatW_m = W_m$ for every $m\ge 1$ hence
\begin{equation*}
\label{fordcomplete3}
\sD(\hatW)\,  =\,  \prolimm \Hom_{\cO}(\hatW, \cO_m)  \, = \, \prolimm \Hom_{\cO}(W , \cO_m) \, =\, \Hom_{\cO}(W, \cO)\, = : \sD(W).
\end{equation*}
For (b) we choose a decreasing sequence $U_1 \supseteq U_2 \supseteq \ldots \supseteq U_n  \supseteq \ldots$ of open subgroups of $T_{\fp}^0$ containing $H$ with $\bigcap_{n \ge 1} U_n= H$.
For $n\ge 1$ we put $\La_n = \cO[T_{\fp}^0/U_n]$, $W^n = W^{U_n}$ and $r_n=\rank_{\La_n} W^n$. Then $\Hom_{\cO}(W^n, \cO) \cong \Hom_{\La_n}(W^n, \La_n)$ is also a free $\La_n$-module of rank $r_n$ for every $n\ge 1$. Since $W^{n+1}/W^n$ is a torsionfree $\cO$-module the transition maps of the inverse system of $\La$-modules 
\[
\ldots\lra \Hom_{\cO}(W^n, \cO) \lra \ldots \lra \Hom_{\cO}(W^2, \cO) \lra \Hom_{\cO}(W^1, \cO) 
\]
are surjective. Moreover its projective limit $\sD(W)$ is a finitely generated $\La$-module. Hence $\sD(\hatW)=\sD(W)$ is a free $\La$-module of finite rank. Indeed, the fact that the canonical map $\sD(W) \otimes_{\La} \La_n \to \Hom_{\cO}(W^n, \cO)$ is surjective implies that the increasing sequence of ranks $\left\{r_n\right\}_{n\ge 1}$ is bounded. So for sufficiently large $n$ the sequence $r_n$ is constant. Coherent choices of bases of the $\La_n$-modules $\Hom_{\cO}(W^n, \cO)$ then yield a $\La$-basis of $\sD(W)$. 
\end{proof}

\begin{lemma}
\label{lemma:keyordedge2}
Let $U_2 \subseteq U_1$ be open subgroups of $T_{\fp}^0$ and put $\Gamma= U_1/U_2$. Then under the assumptions (i) -- (iii) of Prop.\ \ref{prop:freetotdeg} 
the $\cO[\Gamma]$-module $\fOrd_{\cO}^{U_2, d}(X^{\an}, \cO)$ is free of finite rank. 
\end{lemma}

\begin{proof} By Prop.\ \ref{prop:ordvarpiadm} (c) and Lemma \ref{lemma:cohomfree} it suffices to verify 
condition (i) of \ref{lemma:cohomfree} for the $k[\Gamma]$-module $\fOrd_k^{U_2, d}(X^{\an}, k)$, i.e.\ it suffices to verify that the norm map 
\begin{equation*}
\label{normord} 
N_{\Gamma}: \fOrd_k^{U_2, d}(X^{\an}, k) \lra \fOrd_k^{U_2, d}(X^{\an}, k)^{\Gamma}
\end{equation*}
is surjective. By Lemma \ref{lemma:ordnulleins} (a) we have $\fOrd_k^{U_2, d}(X^{\an}, k)^{\Gamma} = \fOrd_k^{U_1, d}(X^{\an}, k)$. If we choose $m\ge 1$ large enough such that $\delta(U_{\fp}^{(m)})\subseteq U_2$ then the group $\fOrd_k^{U_i, d}(X^{\an}, k)$ for $i=1,2$ can be identified with 
$H^d(X_{K_{U_i}(m)}^{\an}, k)^{\ord}$. Hence it suffices to see that the map
\begin{equation}
\label{trace} 
(\pi)_* : H^d(X_{K_{U_2}(m)}^{\an}, k)^{\ord}\lra H^d(X_{K_{U_1}(m)}^{\an}, k)^{\ord}
\end{equation}
is surjective where $\pi: X_{K_{U_2}(m)}\to X_{K_{U_1}(m)}^{\an}$ denotes the projection (note that the latter is a Galois covering with Galois group $K_{U_1}(m)/K_{U_2}(m)\cong U_1/U_2$). By using 
Poincar{\'e} duality for $X_{K_{U_1}(m)}^{\an}$ and $X_{K_{U_1}(m)}^{\an}$ this map can be identified with the dual of the map
\[
\fOrd_k^{U_1, d}(X^{\an}, k) = H^d(X_{K_{U_1}(m)}^{\an}, k)^{\ord} \stackrel{\pi^*}{\lra} H^d(X_{K_{U_2}(m)}^{\an}, k)^{\ord}=\fOrd_k^{U_2, d}(X^{\an}, k)
\]
(here we use assumption (i) of Prop.\ \ref{prop:freetotdeg}). 
This map is a monomorphism by Lemma \ref{lemma:ordnulleins} (a). Hence \eqref{trace} is surjective and $\whH^0(\Gamma, \fOrd_k^{U_2, d}(X^{\an}, k))=0$. 
\end{proof}

\begin{proof}[Proof of Prop.\ \ref{prop:freetotdeg}] Firstly, we prove (b). Put 
\begin{equation*}
\label{fordcomplete2}
W :=\, \fOrd_{\cO}^{H, d}(X^{\an}, \cO),\qquad \hatW := \prolimm W_m \, =\, \fOrd_{\varpi-\ad}^{H, d}(X^{\an}, \cO).
\end{equation*}
By Lemma \ref{lemma:cofree} (b) it suffices to verify the conditions (i) -- (iii) for $W$ and $\hatW$. These follow from Prop.\ \ref{prop:ordvarpiadm} (b), (c) and Lemma \ref{lemma:keyordedge2}. 

If $K_H(1)\in \sK$ then (a) follows from (b) and Lemma \ref{lemma:locbanach}. In the general case we can choose a sufficiently small ideal $\fn'\subseteq \fn$ of $\cO_F$ so that by replacing $X=X^D_0(\fn)^{\fp}$ with $X'=X^D_0(\fn')^{\fp}$ we can apply (a), i.e.\ if we put $\hatW':=\fOrd_{\varpi-\ad}^{H, d}(X^{\an}, \cO)$ then $\cD(\hatW'_E)= \sD(\hatW')_E$ is a free $\La_E$-module of finite rank. Since $\hatW_E$ is a direct summand of $\hatW'_E$ we conclude that  $\cD(\hatW_E)$ is a finitely generated projective $\La_E$-module.
\end{proof}

\subsection{Cohomology of $\sS$-Shimura curves}
\label{subsection:shimuracurve2}

As in the last section $\cO$ denotes a complete discrete valuation ring with maximal ideal $\fP= (\varpi)$, finite residue field $k$ of characteristic $p$ and quotient field $E$ of characteristic $0$.
For a nonarchimedean place $\fq$ of $F$ we denote by $\fG_{\fq} \cong \Gal(\overline{F_{\fq}}/F_{\fq}) \subseteq \fG= \Gal(\barQ/F)$ the decomposition group of a prime of $\barQ$ above $\fq$. By identifying $T_{\fp}$ with $F_{\fp}^*$ via \eqref{splittorus} we can view the local reciprocity map as a homomorphism $\rec: T_{\fp} \to \fG_{\fp}^{\ab}$. Let $\chi_{\cycl}: \fG \to \bZ_p^*=\Aut(\mu_{p^n}(\barQ))$ denote the cyclotomic character. For an $\cO[\fG]$-module $A$ the Tate twist $A(1)$ is given by $A(1)  = A\otimes_{\bZ_p} \bZ_p(1)\cong A(\chi_{\cycl})$ where $\bZ_p(1)= \prolimn \mu_{p^n}(\barQ)$. Let $\alpha: T_{\fp}\to \bZ_p^*$ be the continuous character \eqref{grosschar}, i.e.\ $\alpha$ is given by 
\begin{equation*}
\label{grosschar2}
\alpha\left( \begin{pmatrix} x_1 & 0\\0 & x_2\end{pmatrix}\!\! \!\!\!\!\mod Z\right)\, =\, \Norm_{F_{\fp}/\bQ_p}(x_1/x_2) \Norm(\fp)^{-v_F(x_1) + v_F(x_2)}
\end{equation*}
Recall that $\chi_{\cycl}\circ \rec = \alpha^{-1}$.

We now assume $d=1$. Then $X_K$ is either a Shimura curve defined over $F$ (if $D$ is a division algebra) or $X_K$ is an open modular curve (if $F=\bQ$ and $D=M_2(\bQ)$).
As before for a closed subgroup $H$ of $T_{\fp}^0$ we consider the $\cO[T_{\fp}/H\times \fG]$-module 
\begin{equation}
\label{ordetcomplete}
\fOrd_{\varpi-\ad}^{H, 1}(X_{\barQ}, \cO)(1) := \, \prolimm \fOrd_{\cO_m}^{H, 1}(X_{\barQ}, \cO_m)(1).
\end{equation}
As an $\cO[T_{\fp}/H]$-module it is isomorphic to the module $\fOrd_{\varpi-\ad}^{H, 1}(X^{\an}, \cO)$ defined in \eqref{fordpiadic} so it is in particular a torsionfree $\varpi$-adically admissible $\cO[T_{\fp}/H]$-module. Note that if $H$ is open in $T_{\fp}$ and if $n\ge 1$ is large enough such that $\delta(U_{\fp}^{(n)})\subseteq H$ we have by Prop.\ \ref{prop:ordvarpiadm} (d)
\begin{equation*}
\label{ordetcomplete2}
\fOrd_{\varpi-\ad}^{H, 1}(X_{\barQ}, \cO)(1) \, =\, H^1((X_{K_H(n)})_{\barQ}, \cO)(1)^{\ord}.
\end{equation*}
hence it is free of finite rank as an $\cO$-module.

We also define the $\cO[T_{\fp}\times \fG]$-module
\begin{equation*}
\label{ordetcomplete3}
\fOrd^1(X_{\barQ}, \cO)(1) :=  \dlim_{U \open,\, U \le T_{\fp}^0} \fOrd_{\varpi-\ad}^{U, 1}(X_{\barQ}, \cO)(1)   = \dlim_n H^1((X_{K_1(n)})_{\barQ}, \cO)(1)^{\ord}.
\end{equation*}
It is isomorphic to $\fOrd_{\cO}^1(X^{\an}, \cO)$ as an $\cO[T_{\fp}]$-module so it is admissible as $\cO[T_{\fp}]$-module and free as an $\cO$-module. Moreover by Prop.\ \ref{prop:ordvarpiadm} (b) the $\varpi$-adic completion of $\fOrd^1(X_{\barQ}, \cO)(1)$ is the module $\fOrd_{\varpi-\ad}^1(X_{\barQ}, \cO)(1)$. 

For the remainder of this section we consider the case of Shimura curves, i.e.\ we assume that $D$ is a division algebra (the case of modular curves will be studied in the next section). For the action of the decomposition group $\fG_{\fp}$ on \eqref{ordetcomplete} we have 

\begin{prop}
\label{prop:ordetlocgal}
The $\cO[T_{\fp}\times \fG_{\fp}]\otimes \bT$-module $\fOrd_{\varpi-\ad}^1(X_{\barQ}, \cO)(1)$ contains a submodule \break $\fOrd_{\varpi-\ad}^1(X_{\barQ}, \cO)(1)^0$ with the following properties:

\noi (i) The $\cO$-modules $\fOrd_{\varpi-\ad}^1(X_{\barQ}, \cO)(1)^0$ and 
\[
\fOrd_{\varpi-\ad}^1(X_{\barQ}, \cO)(1)^{\et}:=\fOrd_{\varpi-\ad}^1(X_{\barQ}, \cO)(1)/\fOrd_{\varpi-\ad}^1(X_{\barQ}, \cO)(1)^0
\] 
are torsionfree. They are $\varpi$-adically admissible as $\cO[T_{\fp}]$-modules.

\noi (ii) The action of $\fG_{\fp}$ on $\fOrd_{\varpi-\ad}^1(X_{\barQ}, \cO)(1)^0$ and $\fOrd_{\varpi-\ad}^1(X_{\barQ}, \cO)(1)^{\et}$ factors through $\fG_{\fp}^{\ab}$ and we have
\[
\rec(t^{-1}) \cdot x = \alpha(t)\, t\cdot x\qquad \mbox{and} \qquad \rec(t) \cdot y = t \cdot y
\]
for every $t\in T_{\fp}$ and $x\in \fOrd_{\varpi-\ad}^1(X_{\barQ}, \cO)(1)^0$ and $y\in \fOrd_{\varpi-\ad}^1(X_{\barQ}, \cO)(1)^{\et}$. 
\end{prop}

The proof is based on the following Lemma.

\begin{lemma}
\label{lemma:ordetlocgal2}
For $n\ge 1$ the $\cO[T_{\fp}\times \fG_{\fp}]\otimes \bT$-module $H^1((X_{K_1(n)})_{\barQ}, \cO)(1)^{\ord}$ contains a submodule $\cF_n= \cF$ with the following properties:

\noi (i) The action of $\fG_{\fp}$ on $\cF_n$ factors through $\fG_{\fp}^{\ab}$ and we have $\rec(t^{-1}) \cdot x = \alpha(t)\, t\cdot x$ for every $t\in T_{\fp}$ and $x\in \cF$.

\noi (ii) $\cF$ is maximal isotropic with respect to the cup product pairing restricted to ordinary parts
\begin{equation}
\label{cupord}
H^1((X_{K_1(n)})_{\barQ}, \cO)(1)^{\ord}\times H^1((X_{K_1(n)})_{\barQ}, \cO)(1)^{\ord}\lra H^2((X_{K_1(n)})_{\barQ}, \cO)(2) \cong \cO(1).
\end{equation}
\noi (iii) The action of $\fG_{\fp}$ on $H^1((X_{K_1(n)})_{\barQ}, \cO)^{\ord}/\cF$ factors through $\fG_{\fp}^{\ab}$ and we have 
$\rec(t) \cdot y = t \cdot y$ for every $y\in H^1((X_{K_1(n)})_{\barQ}, \cO)^{\ord}/\cF$. Moreover as an $\cO$-module $H^1((X_{K_1(n)})_{\barQ}, \cO)^{\ord}/\cF$ is free. 
\end{lemma}

\begin{proof} Since $(t\cdot x) \cup y = x\cup (t\cdot y)$ for all $x,y\in H^1((X_{K_1(n)})_{\barQ}, \cO)(1)$ and $t\in T_{\fp}^+$ we see that \eqref{cupord} is a perfect pairing. Also because $H^1((X_{K_1(n)})_{\barQ}, \cO)(1)=H^1((X_{K_1(n)})_{\barQ}, \bZ_p)(1)\otimes_{\bZ_p} \cO$ it suffices to consider the case $\cO =\bZ_p$.

 We will show below using (\cite{wiles1}, Thm.\ 2.2) that $H^1((X_{K_1(n)})_{\barQ}, \bZ_p)(1)^{\ord}$ contains a submodule $\cF'$ satisfying (i) such that
\begin{equation}
\label{cupisotrop}
  \rank_{\bZ_p} \cF' \,=\, 1/2 \cdot \rank_{\bZ_p} H^1((X_{K_1(n)})_{\barQ}, \bZ_p)(1)^{\ord}.
\end{equation}

This implies the assertion. Indeed, if we define 
\[
    \cF\, =\, \left\{ x\in \left(H^1((X_{K_1(n)})_{\barQ}, \bZ_p)(1)^{\ord}\right)^{[\fG_{\fp}, \fG_{\fp}]}\mid \rec(t^{-1}) \cdot x = \alpha(t) \,t\cdot x\,\, \forall \, t\in T_{\fp}\right\}
  \]
($[\fG_{\fp}, \fG_{\fp}]$ denotes the commutator subgroup of $\fG_{\fp}$) then $\cF' \subseteq \cF$. The subspace $\cF$ of $H^1((X_{K_1(n)})_{\barQ},\break \bZ_p)(1)^{\ord}$ is isotropic because of
\[
  \alpha(t)  x\cup y =\rec(t^{-1}) \cdot ( x\cup y) = (\rec(t^{-1}) \cdot x)\cup (\rec(t^{-1})\cdot y) =\alpha(t)^2(t\cdot x \cup t\cdot y) = x\cup y
\]
for $t\in \delta(U^{(m)})\subseteq T_{\fp}$ and $x, y\in \cF$. Since we can choose $t\in \delta(U^{(m)})$ such that $\alpha(t)\ne 1$ we get $x\cup y =0$ for all $x, y\in \cF$. Also note that the fact that $H^1((X_{K_1(n)})_{\barQ}, \bZ_p)(1)^{\ord}$ is a torsionfree $\bZ_p$-module implies that $H^1((X_{K_1(n)})_{\barQ}, \bZ_p)^{\ord}/\cF$ is torsionfree as well. Thus assuming \eqref{cupisotrop} we see that $\cF$ is a maximal isotropic subspace of $H^1((X_{K_1(n)})_{\barQ}, \bZ_p)(1)^{\ord}$ hence
\[
  H^1((X_{K_1(n)})_{\barQ}, \bZ_p)^{\ord}/\cF\, \cong\, \Hom_{\bZ_p}(\cF, \bZ_p)(1)
  \]
 is an isomorphism of $\bZ_p[T_{\fp}\times \fG_{\fp}]$-modules. This yields (iii).
  
To construct $\cF'$ we need to review the decomposition of $H^1((X_{K_1(n)})_{\barQ}, \bQ_p)(1)$ in terms of automorphic representations.
Let $\bD$ be the totally ramified incoherent quaternion algebra over $\bA:=\bA_F$ (in the sense \cite{yuanzhangzhang}) with ramification set $\Ram_{\bD} = \Ram_D \cup \{\infty_1\}$. Following (\cite{yuanzhangzhang}, bottom of page 70) for a subfield  $E$ of $\bC$ we define the set $\cA(\bD^*/\bA^*, E)$ of automorphic representations of $\bD^*/\bA^*$ over $E$ as the set isomorphism classes of irreducible representations of $\bD^*/\bA^*$ such that $\pi\otimes_E \bC$ is a sum of automorphic representations of $\bD^*/\bA^*$ of weight $0$.

By  (\cite{yuanzhangzhang}, p.\ 78) the Jacobian $J_{K_1(n)}$ of $X_{K_1(n)}$ admits a decomposition in the category of abelian varieties over $F$ up to isogeny 
\begin{equation}
\label{jacshimura}
J_{K_1(n)}\, \simeq \, \bigoplus_{\pi\in \cA(\bD^*/\bA^*, \bQ)} \, \pi^{K_1(n) \times K_0(\fn)^{\fp}} \otimes_{\End(\pi)} A_{\pi}.
\end{equation}
Here the abelian variety $A_{\pi}$ corresponds to $\pi$ under the bijection of (\cite{yuanzhangzhang}, Thm.\ 3.8). It is defined over $F$ and it is of strict $\GL_2$-type in the sense of (\cite{yuanzhangzhang}, 3.2.1). In particular we have $\End(\pi)=\End(A_{\pi})=\bQ_{\pi}$ is a totally real number field of degree $[\bQ_{\pi}:\bQ]=\dim(A_{\pi})$ generated by the spherical Hecke eigenvalues of $\pi$. From \eqref{jacshimura} we deduce that 
\begin{equation}
\label{jacshimura2}
H^1((X_{K_1(n)})_{\barQ}, \bQ_p)(1) \, = \, T_p(J_{K_1(n)})\otimes_{\bZ_p} \bQ_p \, = \, \bigoplus_{\pi\in \cA(\bD^*/\bA^*, \bQ)} \, \pi^{K_1(n) \times K_0(\fn)^{\fp}} \otimes_{\bQ_{\pi}} V(\pi)
\end{equation}
where $V(\pi) := T_p(A_{\pi}) \otimes_{\bZ_p} \bQ_p$. 

We fix an isomorphism $\bC\cong \bC_p$ so that we can view $\bQ_p$ as a subfield of $\bC$. For fixed $\pi\in \cA(\bD^*/\bA^*, \bQ)$ we decompose $\bQ_{\pi}\otimes_{\bQ} \bQ_p$ as a product of fields $\bQ_{\pi}\otimes_{\bQ} \bQ_p = \prod_{i} E_i$.
We obtain corresponding decompositions 
\[
\pi\otimes_{\bQ} \bQ_p = \prod_{i} \pi_i, \qquad V(\pi) = \prod_{i} V(\pi_i)
\]
with $\End_{\bQ_p}(\pi_i) = E_i = \End_{\bQ_p[\fG_{\fp}]}(V(\pi_i))$. Thus the decomposition \eqref{jacshimura2} can be refined to 
\begin{equation}
\label{jacshimura3}
H^1((X_{K_1(n)})_{\barQ}, \bQ_p)(1) \,  = \, \bigoplus_{\pi\in \cA(\bD^*/\bA^*, \bQ_p)} \, \pi^{K_1(n) \times K_0(\fn)^{\fp}} \otimes_{E_{\pi}} V(\pi)
\end{equation}
where $E_{\pi} = \End_{\bQ_p}(\pi) = \End_{\bQ_p[\fG_{\fp}]}(V(\pi))$ is generated (over $\bQ_p$) by spherical Hecke eigenvalues. 

For $\pi \in \cA(\bD^*/\bA^*, \bQ_p)$ with $\pi^{K_1(n) \times K_0(\fn)^{\fp}}\ne 0$ the $E_{\pi}[\fG]$-module $V(\pi)$ is twodimensional as an $E_{\pi}$-vector space. Furthermore the $\fG$-action on $V(\pi)$ is unramified outside the set primes dividing $p\cdot \fd \cdot \fn$
and we have for every prime $\fq\nmid p\fd\fn$ (Eichler-Shimura relations)
\begin{equation}
\label{galrephmf}
\Tr(\Frob_{\fq}; V(\pi) \to V(\pi)) = \la_{\fp}(T_{\fq}), \qquad \det(\Frob_{\fq}: V(\pi)\to V(\pi)) \, =\, \Norm(\fq)
\end{equation}
where $\Frob_{\fq} \in \Gal(F_{\fq}/F)$ is a Frobenius and $\la_{\pi}: \bT^{p\fd \fn} \to E_{\pi}$ is the Hecke eigenvalue homomorphism (i.e.\ it describes the action of the spherical prime-to-$p\fd\fn$ Hecke algebra on the set of spherical vectors in $\pi$). If $\cO_{\pi}$ denotes the valuation ring of $E_{\pi}$ then there exists a $\fG$-stable $\cO_{\pi}$-lattice $T(\pi)$ in $V(\pi)$.

Similar to (\cite{yuanzhangzhang}, Thm.\ 3.4 (3)) one shows that every $\pi \in \cA(\bD^*/\bA^*, \bQ_p)$ has a decomposition $\pi = \bigotimes_{v}'  \pi_v$ and the $\pi_v$ are irreducible admissible representations of $\bD_v$ defined over $E_{\pi}$. Thus for $\pi \in \cA(\bD^*/\bA^*, \bQ_p)$ we have 
\begin{equation*}
\label{autdecomp2}
\pi^{K_1(n) \times K_0(\fn)^{\fp}}  \,  = \, \pi_{\fp}^{K_1(n)} \otimes_{E_{\pi}} (\pi^{\fp})^{K_0(\fn)^{\fp}}.
\end{equation*}
We denote by $\cR_1$ the set of $\pi\in \cA(\bD^*/\bA^*, \bQ_p)$ contributing to \eqref{jacshimura3}, i.e.\ those $\pi$ such that $\pi_{\fp}^{K_1(n)}\ne 0 \ne (\pi^{\fp})^{K_0(\fn)^{\fp}}$. By using the fact that $H^1((X_{K_1(n)})_{\barQ}, \bZ_p)(1)$ is a $\bZ_p$-lattice in $H^1((X_{K_1(n)})_{\barQ}, \bQ_p)(1)$ it is easy to see that the $E_{\pi}$-module $\pi_{\fp}^{K_1(n)}$ (resp.\ $(\pi^{\fp})^{K_0(\fn)^{\fp}}$) admits a $T_{\fp}^+$-stable (resp.\ a $\bT$-stable) $\cO_{\pi}$-lattice $M_{\pi, \fp}$ (resp.\ $M_{\pi}^{\fp}$) for every $\pi \in \cR_1$. Furthermore we may assume that the image of $\bigoplus_{\pi\in \cR_1} \, M_{\pi, \fp} \otimes_{\cO_{\pi}} M_{\pi}^{\fp} \otimes_{\cO_{\pi}} T(\pi)$ under the isomorphism \eqref{jacshimura3} is contained in $H^1((X_{K_1(n)})_{\barQ}, \bZ_p)(1)$ so that \eqref{jacshimura3} is induced by an injective map 
\begin{equation*}
\label{jacshimura4}
\bigoplus_{\pi\in \cR_1} \, M_{\pi, \fp} \otimes_{\cO_{\pi}} M_{\pi}^{\fp} \otimes_{\cO_{\pi}} T(\pi) \lra
H^1((X_{K_1(n)})_{\barQ}, \bZ_p)(1) 
\end{equation*}
with finite cokernel. Passing to ordinary parts yields a monomorphism with finite cokernel
\begin{equation*}
\label{jacshimura5}
\bigoplus_{\pi\in \cR_2} \, M_{\pi, \fp}^{\ord} \otimes_{\cO_{\pi}} M_{\pi}^{\fp} \otimes_{\cO_{\pi}} T(\pi)\lra H^1((X_{K_1(n)})_{\barQ}, \bZ_p)(1)^{\ord}
\end{equation*}
where $\cR_2$ denotes the subset of $\pi\in \cR_1$ such that $M_{\pi, \fp}^{\ord}\ne 0$. 

We fix $\pi\in \cR_2$. By Prop.\ \ref{prop:classifyordjacq} there exists a quasicharacter $\chi: F_{\fp}^* \to \cO^*$ such that $M_{\pi, \fp}^{\ord} \cong \cO_{\pi}(\chi)$ and $\pi_{\fp} = \Ind_{B_{\fp}}^{G_{\fp}} \chi^{-1}$ or $\pi_{\fp} = \sigma(\chi)$ (depending on $\chi^2\ne 1$ or $\chi^2 =1$ respectively). The $\cO_{\pi}$-module $T(\pi)$ contains a $\fG_{\fp}$-stable $\cO_{\pi}$-submodule $\cF_{\pi}'$ of rank $1$ so that we have
\begin{equation}
\label{ordpigal}
\rec(t^{-1}) \cdot x = (\alpha\chi)(t)\, x \qquad \forall\, x\in \cF_{\pi}', \, t\in T_{\fp}. 
\end{equation}
(note that since $\rank_{\cO_{\pi}} \cF_{\pi}' =1$ the $\fG_{\fp}$-action on $\cF_{\pi}'$ factors through $\fG_{\fp}^{\ab}$). This implies that the submodule $M_{\pi, \fp}^{\ord} \otimes_{\cO_{\pi}} M_{\pi}^{\fp} \otimes_{\cO_{\pi}} \cF_{\pi} \cong \cO_{\pi}(\chi) \otimes_{\cO_{\pi}} M_{\pi}^{\fp} \otimes_{\cO_{\pi}} \cF_{\pi}$ of $H^1((X_{K_1(n)})_{\barQ}, \bZ_p)(1)^{\ord}$ is contained in $\cF$.

As explained on page 130 of \cite{hida2} the existence of the submodule $\cF_{\pi}'$ can deduced from (\cite{wiles1}, Thm.\ 2.2; see also \cite{wiles2}, Lemma 2.1.5). Indeed, we can choose a character $\xi: \bA/F^*\to \cO^*$ of finite order with $\xi_v=1$ for all archimedean places $v$ and $\xi_{\fp}|_{U_{\fp}}= \chi|_{U_{\fp}}$. If we view $\pi$ as a representation of $\wG(\bA)$ via the projection $\wG(\bA)\to G(\bA)$ and $\xi$ as a character $\wG(\bA)\to \cO^*$ via the reduced norm then for the Galois representation $V(\pi\otimes \xi)$ associated to $\pi\otimes \xi$ we have $V(\pi\otimes \xi^{-1})=V(\pi)\otimes_{E_{\pi}} E_{\pi}(\xi'^{-1})$ where $\xi':\fG^{\ab}\to \cO^*$ is the character corresponding to $\xi$ via global class field theory. We can now apply (\cite{wiles1}, Thm.\ 2.2 and \cite{wiles2}, Lemma 2.1.5) to the automorphic representation $\pi\otimes \xi^{-1}$ to deduce the existence of $\cF_{\pi}'$ such that \eqref{ordpigal} holds.\footnote{Strictly speaking we cannot apply (\cite{wiles1}, Thm.\ 2.2) directly as stated since we have not assumed that $[F:\bQ]$ is odd; however the proof in loc.\ cit.\ covers the case of the Galois representation associated to the automorphic representation $\pi\otimes \xi^{-1}$ at hand since the latter is a direct summand of the Tate module of the Jacobian of a Shimura curve asscociated to the quaternion algebra $D$.}

Thus if we set
\[
\cF': =\, \bigoplus_{\pi\in \cR_2} \, M_{\pi, \fp}^{\ord} \otimes_{\cO_{\pi}} M_{\pi}^{\fp} \otimes_{\cO_{\pi}} \cF_{\pi}'
\]
then $\cF'$ is contained in $\cF$. Using the fact that $\rank_{\cO_{\pi}} \cF_{\pi}' = 1/2\cdot \rank_{\cO_{\pi}} T(\pi)$ for $\pi\in \cR_2$ we see that \eqref{cupisotrop} holds as well. 
\end{proof}

\begin{proof}[Proof of Prop.\ \ref{prop:ordetlocgal}] Note that the restriction map 
\[
H^1((X_{K_1(n)})_{\barQ}, \cO)(1)^{\ord}\to H^1((X_{K_1(n+1)})_{\barQ}, \cO)(1)^{\ord}
\] 
is injective and that we have $\cF_n = H^1((X_{K_1(n)})_{\barQ}, \cO)(1)^{\ord}\cap \cF_{n+1}$ for every $n\ge 1$. Thus in the short exact sequence 
\[
0 \lra \dlim_n \cF_n \lra \dlim_n H^1((X_{K_1(n)})_{\barQ}, \cO)(1)^{\ord} \lra \dlim_n H^1((X_{K_1(n)})_{\barQ}, \cO)(1)^{\ord}/\cF_n\lra 0
\]
the second term is $=\fOrd^1(X_{\barQ}, \cO)(1)$ and the third is torsionfree as an $\cO$-module. It follows that sequence remains exact after passing to $\varpi$-adic completions. We define $\fOrd_{\varpi-\ad}^1(X_{\barQ}, \cO)(1)^0$ as the $\varpi$-adic completion of $\dlim_n \cF_n$. It is an $\cO[T_{\fp}\times \fG_{\fp}]$-submodule of $\fOrd_{\varpi-\ad}^1(X_{\barQ}, \cO)(1)$. The assertion now follows from Lemma \ref{lemma:ordetlocgal2}.
\end{proof}

Let $W$ be a $\varpi$-adically admissible $\cO[T]$-module that is free and of finite rank as an $\cO$-module. 
Recall that the cohomology group 
\begin{equation*}
\label{laextgald1}
\bH_{\varpi-\ad}^n(X_{\barQ}; W, \cO) \, = \, \prolimm \bH_{\varpi-ad}^n(X_{\barQ}; W, \cO_m)
\end{equation*}
is an augmented $\cO[T_{\fp}]$-module that is is finitely generated as an $\cO$-module. It carries additionally a continuous $\fG$-action commuting with the $\La_{\cO}(T_{\fp})$- and Hecke action. Similar to Prop.\ \ref{prop:ordvarpiadm} (e) we have
\begin{equation}
\label{ordetcomplete5}
\bH_{\varpi-\ad}^n(X_{\barQ}; W, \cO) \, =\, \prolimm \Ext_{\cO_m, T_{\fp}}^{n-1}(W_m^{\io}, \fOrd^1(X_{\barQ}, \cO_m)) \nonumber 
\end{equation}
for every $n\ge 0$. Also if $V\in \Ban_E^{\adm}(T_{\fp})$ with $\dim_E(V)<\infty$ then we have 
\begin{eqnarray*}
\label{laextgald2}
\bH_{\varpi-\ad}^1(X_{\barQ}; V, E) & = & \Hom_{\Ban_E^{\adm}(T_{\fp})}(V^{\io},\fOrd_{\varpi-\ad}^1(X_{\barQ}, \cO)_E)\\
& = & \Hom_{E[T_{\fp}]}(V^{\io},\fOrd_{\varpi-\ad}^1(X_{\barQ}, \cO)_E).\nonumber
\end{eqnarray*}
Now assume that $V=A(\chi)$ where $A$ is an Artinian local $E$-algebra with maximal ideal $\fM$ such that $A/\fM=E$ and where $\chi: F_{\fp}^*\cong T_{\fp} \to A^*$ is a bounded continuous homomorphism (so that $A(\chi) \in \Ban_E^{\adm}(T_{\fp})$; see Remark \ref{remarks:lamod} (b)). 

If $V'$ is any admissible $E$-Banach space representation of $T_{\fp}$ equipped with a continuous $E$-linear $\fG_{\fp}$-action commuting with the $T_{\fp}$-action then we equip 
$\Hom_{E[T_{\fp}]}(A(\chi), V')$ with the following $A[\fG_{\fp}]$-module structure 
\begin{equation*}
\label{laextgald3}
(a \cdot \Phi)(b) \, = \, \Phi(ab), \qquad (\sigma \cdot \Phi)(b) \, = \, \sigma \cdot \Phi(b)\nonumber
\end{equation*}
for all $\Phi\in \Hom_{E[T_{\fp}]}(A(\chi), V')$, $a\in A$, $b\in A(\chi)$ and $\sigma\in \fG_{\fp}$. 

We define 
\begin{eqnarray*} 
\label{laextgald4}
&& \bH_{\varpi-\ad}^1(X_{\barQ}; A(\chi), E)(1)^0 := \Hom_{E[T_{\fp}]}(A(\chi^{-1}),  \fOrd_{\varpi-\ad}^1(X_{\barQ}, \cO)(1)^0_E)\\
&& \bH_{\varpi-\ad}^1(X_{\barQ}; A(\chi), E)(1)^{\et}:= \Hom_{E[T_{\fp}]}(A(\chi^{-1}), \fOrd_{\varpi-\ad}^1(X_{\barQ}, \cO)^{\et}_E) \nonumber
\end{eqnarray*}
so that we obtain an exact sequence of $A[\fG_{\fp}]\otimes \bT$-modules
\begin{eqnarray} 
\label{laextgald5}
&&0\lra \bH_{\varpi-\ad}^1(X_{\barQ}; A(\chi), E)(1)^0 \lra \bH_{\varpi-\ad}^1(X_{\barQ}; A(\chi), E)(1) \\
&& \hspace{2cm} \lra \bH_{\varpi-\ad}^1(X_{\barQ}; A(\chi), E)^{\et} \stackrel{\delta}{\lra} \Ext_{E[T_{\fp}]}^1(A(\chi^{-1}), \fOrd_{\varpi-\ad}^1(X_{\barQ}, \cO)(1)_E^0). \nonumber
\end{eqnarray}
Note that all three terms are finitely generated $A$-modules by Remark \ref{remark:locdualpply}. Let $\barchi: F_{\fp}^*\cong T_{\fp} \to (A/\fM)^*=E^*$ denote the character $\barchi(x) := \chi(x) \!\!\mod \fM$.

\begin{theorem}
\label{theorem:reciprocity}
(a) The action of $\fG_{\fp}$ on the first and third term in \eqref{laextgald5} factors through $\fG_{\fp}^{\ab}$ and we have 
\begin{equation*}
\label{galdecfpact}
\rec(t) \cdot x \, =\, \alpha^{-1}(t) \cdot \chi(t) \cdot x\qquad \mbox{and} \qquad \rec(t)\cdot y\, =\, \chi^{-1}(t) \cdot y
\end{equation*}
for every $t\in T_{\fp}$, $x\in \bH_{\varpi-\ad}^1(X_{\barQ}; A(\chi), E)(1)^0$ and $y\in \bH_{\varpi-\ad}^1(X_{\barQ}; A(\chi), E)(1)^{\et}$. Moreover if $\barchi^2\ne \alpha$ holds then the sequence 
\begin{equation*} 
\label{laextgald5a}
0\to \bH_{\varpi-\ad}^1(X_{\barQ}; A(\chi), E)(1)^0 \to \bH_{\varpi-\ad}^1(X_{\barQ}; A(\chi), E)(1) \to \bH_{\varpi-\ad}^1(X_{\barQ}; A(\chi), E)^{\et}\to 0
\end{equation*}
is exact. 
\medskip

\noi (b) The $\fG$-action on $\bH_{\varpi-\ad}^n(X_{\barQ}; A(\chi), E)(1)$ is unramified outside the set of primes dividing $p\cdot \fn$ and the Eichler-Shimura relations
\begin{equation*}
\label{eichlershimura}
\Frob_{\fq}^2 + T_{\fq}  \Frob_{\fq} + \Norm(\fq) \, =\, 0
\end{equation*}
hold on $\bH_{\varpi-\ad}^n(X_{\barQ}; A(\chi), E)(1)$ for all $\fq\nmid p \fd \fn$ and $n\ge 0$.
\end{theorem}

\begin{proof} (a) follows immediately from Prop.\ \ref{prop:ordetlocgal} except for the last claim. For that we have to show that if the connecting homomorphism $\delta$ in \eqref{laextgald5} is $\ne 0$ then $\barchi^2= \alpha$. Let $y$ be an element $\ne 0$ in the image of $\delta$. By Prop.\ \ref{prop:ordetlocgal}, for $t\in T_{\fp}$ the element $\rec(t)$ of $\fG_{\fp}$ acts on the source of $\delta$ by multiplication with $\alpha^{-1}(t) \cdot \chi(t)$ and on the target by multiplication with $\chi^{-1}(t)$. It follows $\alpha^{-1}(t) \cdot \chi(t) \cdot y = \chi^{-1}(t) \cdot y$, hence $(\alpha^{-1} \chi^2)(t) -1 \in \Ann_{A}(y)\subseteq \fM$ for every $t\in T_{\fp}$, whence $\barchi^2= \alpha$.

For (b) note that by \eqref{jacshimura3} and \eqref{galrephmf} (or more directly by work of Ihara \cite{ihara}) the assertion holds for the action of $\fG$ and $\bT$ on $H^1((X_{K_1(n)})_{\barQ}, \bQ_p)(1)$, hence also on $H^1((X_{K_1(n)})_{\barQ}, \cO)(1)^{\ord}  \subseteq H^1((X_{K_1(n)})_{\barQ}, \bQ_p)(1)_E$ and therefore also for the action of $\fG$ and $\bT$ on $\fOrd^1(X_{\barQ}, \cO)(1)$ and on its $\varpi^m$-torsion subgroup $\fOrd_{\cO_m}^1(X_{\barQ}, \cO_m)(1)$. Using \eqref{ordetcomplete5} we deduce that the assertion holds for the action of $\fG$ and $\bT$ on $\bH_{\varpi-\ad}^n(X_{\barQ}; V, E)(1)$ for every $n\ge 0$ and any admissible $E$-Banach space representation $V$ of $T_{\fp}$. 
\end{proof}

\subsection{Cohomology of $\sS$-modular curves}
\label{section:modcurves}

As before $\cO$ denotes a complete discrete valuation ring with maximal ideal $\fP= (\varpi)$, finite residue field $k$ of characteristic $p$ and quotient field $E$ of characteristic $0$.
We now assume that $F=\bQ$ and that $D= M_2(\bQ)$. The ideals $\fp$, $\fn$ will be written as $\fp=p\bZ$, $N\bZ=\fn$ for $p$ a prime and an integer $N\ge 1$ with $p\nmid N$. We write $T_p$, $B_p$, $G_p$, $\fG_p$ etc.\ for $T_{\fp}$, $B_{\fp}$, $G_{\fp}$, $\fG_{\fp}$ etc.\ As in \ref{subsection:borelserre} we let $X^{\BS}$ be the Borel-Serre compactification of $X^{\an}$, we let $\partial X$ be its boundary and let $Y$ be the $\sS^B$-space defined in \eqref{borelserre4}. 

\begin{lemma}
\label{lemma:boundmodcurve}
There exists a canonical long exact sequence 
\begin{equation}
\label{borelserrord}
\ldots \lra \fOrd_{\cO, c}^n(X^{\an}, \cO) \lra  \fOrd_{\cO}^n(X^{\an}, \cO) \lra \fOrd_{\cO}^n(\partial X, \cO) \lra \fOrd_{\cO, c}^{n+1}(X^{\an}, \cO)\lra\ldots
\end{equation}
\end{lemma}

\begin{proof} By Lemma \ref{lemma:csres} we have an exact sequence \eqref{borelserrord} with the groups  $\fOrd_{\cO}^{\bu}(X^{\an}, \cO)$ replaced by $\fOrd_{\cO}^{\bu}(X^{\BS}, \cO)$. So it suffices to see that $\fOrd_{\cO}^n(X^{\an}, \cO)\cong \fOrd_{\cO}^n(X^{\BS}, \cO)$. For that, note that $X_{K_1(m)}^{\an}\hra X_{K_1(m)}^{\BS}$ is a homotopy equivalence so we have 
$\fOrd_{\cO}^n(X^{\an}, \cO)\cong  \dlim_m H^n(X_{K_1(m)}, \cO)^{\ord} \break \cong \dlim_m H^n(X_{K_1(m)}^{\BS}, \cO)^{\ord} =  \fOrd_{\cO}^n(X^{\BS}, \cO)$.
\end{proof}

We consider the following $\cO[T_p]$-submodule of $\fOrd^1(X^{\an}, \cO)$ 
\begin{equation*}
\label{ordoneint}
\fOrd_{\cO, !}^1(X^{\an}, \cO) : =\, \image(\fOrd_{\cO, c}^n(X^{\an}, \cO) \lra  \fOrd_{\cO}^n(X^{\an}, \cO))
\end{equation*}
Note that we have 
\begin{equation*}
\label{ordoneint2}
\fOrd_{\cO, !}^1(X^{\an}, \cO) \, =\, \dlim_m H_{!}^1(X_{K_1(m)}, \cO)^{\ord}
\end{equation*}
where $H_{!}^n(X_{K_1(m)}, \cO):=\image(H_c^1(X_{K_1(m)}, \cO)\to H^1(X_{K_1(m)}, \cO))$.

\begin{lemma}
\label{lemma:boundmodcurve2}
For every $n\ge 0$ we have
\medskip

\noi (a) $\fOrd_{\cO, c}^n(X^{\an}, \cO)_{\tor}=0$ (i.e.\ $\fOrd_{\cO, c}^n(X^{\an}, \cO)$ is torsionfree as an $\cO$-module). 
\medskip

\noi (b) $\fOrd_{\cO}^n(\partial X, \cO)_{\tor}=0$. 
\medskip

\noi(c) The sequences
\begin{equation*}
\label{ordcsmod}
\begin{CD}
0 @>>> \fOrd_{\cO, c}^n(X^{\an}, \cO) @> \varpi^m \cdot >>  \fOrd_{\cO, c}^n(X^{\an}, \cO) @>>> \fOrd_{\cO, c}^n(X^{\an}, \cO_m) @>>> 0
\end{CD}
\end{equation*}
and 
\begin{equation*}
\label{ordrandmod}
\begin{CD}
0 @>>> \fOrd_{\cO}^n(\partial X, \cO) @> \varpi^m \cdot >>  \fOrd_{\cO}^n(\partial X, \cO) @>>> \fOrd_{\cO}^n(\partial X, \cO_m) @>>> 0
\end{CD}
\end{equation*}
are exact for every $m\ge 1$. 
\end{lemma}

\begin{proof} (a) By definition of $\fOrd_{\cO, c}^n(X^{\an}, \cO)$ it suffices to see that the $\cO$-module $H_c^n(X_{K_1(m)}^{\an}, \cO)$ is free and of finite rank. Poincar{\'e} duality for the non-compact Riemann surface $X_{K_1(m)}^{\an}$ yields
\[
H_c^n(X_{K_1(m)}^{\an}, \cO) \, \cong \, \Hom_{\cO}(H^{2-n}(X_{K_1(m)}^{\an}, \cO), \cO) \, \cong \, \Hom_{\bZ}(H^{2-n}(X_{K_1(m)}^{\an}, \bZ), \cO).
\]
Since $H^n(X_{K_1(m)}^{\an}, \bZ)$ is a free abelian group of finite rank (and $=0$ if $n\ge 2$) we see that 
$H_c^n(X_{K_1(m)}^{\an}, \cO)$ is a free $\cO$-module of finite rank (and $=0$ if $n\ne 1,2$). 

(b) Let $\sS^B$ (resp.\ $\sS^T$) denote the set of discrete left $B_p$-sets $S$ (resp.\ $T_p$-sets $S'$) such that $\Stab_{B_p}(s)$ (resp.\ $\Stab_{T_p}(s')$ is compact for every $s\in S$ (resp.\ $s'\in S'$). Let $f:Y\to Z$ denote the morphism between the $\sS^B$-space $Y$ and $\sS^T$-space $Z$ compatible with the canonical projection $\pi: B_p\to B_p/N_p=T_p$ introduced in section \ref{subsection:borelserre}. For $S\in \sS^B$ having finitely many $B_p$-orbits the map $f_S: Y_S\to Z_{\barS}$ (with $\barS:=N_p\backslash S\in \sS^T$) is a fibration whose fibers are onedimensional real tori and whose base $Z_{\barS}$ is a discrete space consisting of finitely many points. It follows that $H^n(Y_S, \cO)$ is a free $\cO$-module of rank $\# Z_{\barS}$ for $n=0,1$ and $H^n(Y_S, \cO)=0$ for $n\ge 2$. 

Since $\partial X_{K_1(m)} = Y_S$ with $S= G_p/K_1(m)$ (viewed as an element of $\sS^B$) it follows that 
$H^n(\partial X_{K_1(m)}, \cO)$ is a free $\cO$-module of finite rank for every $m\ge 1$ (and $=0$ if $n\ge 2$). Now the assertion follows from Remark \ref{remarks:fordhida} (a).

(c) follows immediately from (a) and (b).
\end{proof}

Consider the long exact sequence of inverse systems 
\begin{eqnarray}
\label{borelserrord2}
&& \ldots \lra \{\fOrd_{\cO, c}^n(X^{\an}, \cO_m)\}_{m\ge 1} \lra  \{\fOrd_{\cO}^n(X^{\an}, \cO_m)\}_{m\ge 1} \lra \{\fOrd_{\cO}^n(\partial X, \cO_m)\}_{m\ge1}\\
&& \hspace{3cm} \lra \{\fOrd_{\cO, c}^{n+1}(X^{\an}, \cO_m)\}_{m\ge 1}\lra\ldots\nonumber
\end{eqnarray}
Prop.\ \ref{prop:ordvarpiadm} (a) and Lemma \ref{lemma:boundaryfib} (c) imply that the transition maps 
of each system are surjective, so the sequence remains exact after passing to the limit.
Note that every $\cO_m[T_p]$-module occurring in \eqref{borelserrord2} is admissible. Indeed by Prop.\ 
\ref{prop:finiteness1} the $\cO_m[T_p]$-modules $\fOrd_{\cO}^n(X^{\an}, \cO_m)$ and $\fOrd_{\cO}^n(\partial X, \cO_m)$ are. Since $\Mod_{\cO_m}^{\adm}(T_p)$ is a Serre subcategory of $\Mod_{\cO_m}^{\sm}(T_p)$
(see \cite{emerton3}, Prop.\ 2.2.13) we conclude that the modules $\fOrd_{\cO, c}^n(X^{\an}, \cO_m)$ are admissible as well.

The sequence of the projective limits will be denoted by
\begin{eqnarray}
\label{ordpiadicbs}
&& \ldots \lra \fOrd_{\varpi-\ad, c}^n(X^{\an}, \cO) \lra \fOrd_{\varpi-\ad}^n(X^{\an}, \cO)  \lra \fOrd_{\varpi-\ad}^n(\partial X, \cO)\\
&& \hspace{3cm} \lra \fOrd_{\varpi-\ad, c}^{n+1}(X^{\an}, \cO)\lra \ldots\nonumber
\end{eqnarray}
We define 
\begin{equation*}
\label{ordoneint2a}
\fOrd_{\varpi-\ad, !}^1(X^{\an}, \cO) :=\, \prolimm \fOrd_{\cO, !}^1(X^{\an}, \cO_m)\, =\, \image(\fOrd_{\varpi-\ad, c}^1(X^{\an}, \cO) \to \fOrd_{\varpi-\ad}^1(X^{\an}, \cO)).
\end{equation*}
As a consequence of Prop.\ \ref{prop:ordvarpiadm} (a) and Lemma \ref{lemma:boundaryfib} (c) we get

\begin{lemma}
\label{lemma:borelserrepiadic}
(a) The sequence \eqref{ordpiadicbs} is an exact sequence of torsionfree $\varpi$-adically admissible $\cO[T_p]$-modules.
\medskip

\noi (b) The $\cO[T_p]$-module $\fOrd_{\varpi-\ad, !}^1(X^{\an}, \cO)$ is $\varpi$-adically admissible and torsionfree. It is the $\varpi$-adic completion of the admissible $\cO[T_p]$-module $\fOrd_{\cO, !}^1(X^{\an}, \cO)$.
\end{lemma}

Next we will give an explicit description of the $\varpi$-adically admissible $\cO[T_p]$-module $\fOrd_{\varpi-\ad}^n(\partial X, \break \cO)$ for $n=0, 1$ with its action of the Hecke algebra $\bT=\bZ[T_{\ell}\mid \ell\nmid pN]$. Let $C(\bQ_p^* \times (\bZ/N\bZ)^*, \cO)$ denote the $\cO$-module of continuous maps $f: \bQ_p^* \times (\bZ/N\bZ)^*\to \cO$ and let
\begin{equation*}
\label{ordrandpiadic}
C^0\, =\, \{ f\in C(\bQ_p^* \times (\bZ/N\bZ)^*, \cO)\mid  f(px,\barp y) = f(x,y)\,\forall (x,y)\in \bQ_p^* \times (\bZ/N\bZ)^*\}
\end{equation*}
(where $\barp:=p\mod N\in (\bZ/N\bZ)^*$). We define a $T_p$- and a $\bT$- action on $C^0$ by
\begin{eqnarray*}
\label{ordrandpiadic2}
&& (\delta(a) \cdot f)(x,y) \, = \, f(a^{-1} x, y)\\
&& (T_{\ell} f)(x,y) \, = \, f(\ell^{-1} x, y) + \ell f(\ell x, y)\nonumber
\end{eqnarray*}
where $f\in C^0$, $(x,y)\in \bQ_p^* \times (\bZ/N\bZ)^*$, $a\in \bQ_p^*$ and $\ell\nmid pN$ is a prime number ($\delta: \bQ_p^*\to T_p$ denotes the isomorphism \eqref{splittorus}). It is easy to see that $C^0$ is a $\varpi$-adically admissible $\cO[T_p]$-module equipped with an equivariant Hecke action. We define $C^1$ by $C^1=C^0$ with the same Hecke action but with $T_p$-action given by
\begin{equation*}
\label{ordrandpiadic3}
(t \cdot f)(x,y) \, = \, \alpha(t) f(\delta^{-1}(t) x, y)
\end{equation*}
for all $t\in T_p$ (i.e.\ $(\delta(a) \cdot f)(x,y) = \frac{a}{p^{v_p(a)}}  f(a x, y)$ for all $a\in \bQ_p^*$).

\begin{prop}
\label{prop:ordrandpiadic}
There exists a canonical Hecke equivariant isomorphism of $\varpi$-adically admissible $\cO[T_p]$-modules
\begin{equation}
\label{ordrandpiadic4}
\fOrd_{\varpi-\ad}^n(\partial X, \cO) \, \cong \, C^n
\end{equation}
for $n=0, 1$. Moreover $\fOrd_{\varpi-\ad}^n(\partial X, \cO)=0$ for $n\ge 2$.
\end{prop}

\begin{proof} Using Prop.\ \ref{prop:borelserrefiber}, Remark \ref{remark:borelserrecohom}, Cor.\ \ref{corollary:ordind} and the spectral sequence \eqref{ordinftyss} we obtain
\begin{equation}
\label{ordrandpiadic5}
\fOrd_{\cO_m}^n(\partial X, \cO_m) \, \cong \, \Ord^n( \Ind_{B_p}^{G_p} H^0(Z_{\infty}, \cO_m)) \, = \, \left\{\begin{array}{cc} H^0(Z_{\infty}, \cO_m)^{\io} & \mbox{if $n=0$,}\\
H^0(Z_{\infty}, \cO_m)(\alpha) &  \mbox{if $n=1$,}\\
0 & \mbox{if $n>1$,}
\end{array}\right.
\end{equation}
for every $m\ge 1$. By identifying the torus $T$ with $\bG_m$ via \eqref{splittorus} it is easy to see that $H^0(Z_{\infty}, \cO_m)$ can be identified with the
\[
\dlim_n \Maps(  \bA_f^*/(\bQ_+^*(U_p^{(n)} \times \prod_{\ell\ne p} U_{\ell}^{(n_{\ell})})), \cO_m)\, \cong \, \dlim_n \Maps( (\bQ_p^*/U_p^{(n)}\times (\bZ/N\bZ)^*)/\langle p\rangle, \cO_m)
\]
Thus we see that $\prolimm H^0(Z_{\infty}, \cO_m)$ is isomorphic to $C^0$ viewed solely as an $\cO$-module. The assertion follows by keeping track of the $T_p$ and Hecke action under the isomorphism \eqref{ordrandpiadic5} and by using Prop.\ \ref{prop:borelserrefiber} (b).
\end{proof}

\begin{remark}
\label{remark:randeisenstein}
\rm Using Prop.\ \ref{prop:ordrandpiadic} we can give an explicit description of the $\bT_E$-module 
\begin{equation*}
\label{randeisenstein1}
\Hom_{E[T_p]}(E(0), \fOrd_{\varpi-\ad}^n(\partial X, \cO)_E)
\end{equation*}
namely it is an $E$-vector space whose dimension is the order of the group $(\bZ/N\bZ)^*/\langle \barp\rangle$ and where $\bT$ acts via the homomorphism $\eis: \bT\to \bZ, \, T_{\ell} \mapsto \ell +1$. Indeed, by Prop.\ \ref{prop:ordrandpiadic} we have 
\[
\Hom_{E[T_p]}(E(0), \fOrd_{\varpi-\ad}^n(\partial X, \cO)_E)\, \cong \, \Hom_{E[T_p]}(E(0), C^1_E)\,\cong \, \{f\in C^1_E\mid t\cdot f = f\}
\]
and a simple computation shows that the map 
\begin{eqnarray*}
\label{randeisenstein2}
&& \{f\in C^1_E\mid t\cdot f = f\} \lra \{\wf\in \Maps((\bZ/N\bZ)^*, E)\mid \wf(\barp y) =\wf(y)\, \forall \,y\in (\bZ/N\bZ)^*\}\\
&&  f \mapsto \left( \wf: (\bZ/N\bZ)^*\to E, \,\, y\mapsto \wf(y) := f(1, y)\right)\nonumber
\end{eqnarray*}
is an isomorphism and that the $\bT_E$-action on the source corresponds to the $\bT_E$-action via the homomorphism $\eis$ on the target. \enddemo
\end{remark}

We are going to study the $\cO[T_p\times \fG]$-modules $\fOrd^1(X_{\barQ}, \cO)(1)$, $\fOrd_{\varpi-\ad}^1(X_{\barQ}, \cO)(1)$ and \break $H_{\varpi-\ad}^1(X_{\barQ}; W, \cO)(1)$. For that we also have to consider the $\cO[T_p\times \fG]$-module (see Remark \ref{remark:fpordcsuppet})
\begin{eqnarray*}
\label{ordetonec}
&& \fOrd_c^1(X_{\barQ}, \cO)(1) \, =\,  \dlim_n H_c^1((X_{K_1(n)})_{\barQ}, \cO)(1)^{\ord} \, =\, \dlim_n\left( \prolimm H_c^1((X_{K_1(n)})_{\barQ}, \cO_m)(1)^{\ord}\right)\\
&& \fOrd_{\varpi-\ad, c}^1(X_{\barQ}, \cO)(1) \, =\,  \prolimm \fOrd_c^1(X_{\barQ}, \cO_m)(1) \, =\, \prolimm\left(\dlim_n H_c^1((X_{K_1(n)})_{\barQ}, \cO_m)(1)^{\ord}\right)
\end{eqnarray*}
As $\cO[T_p]$- and Hecke modules we have $\fOrd_c^1(X_{\barQ}, \cO)(1)$ and $\fOrd_{\varpi-\ad, c}^1(X_{\barQ}, \cO)(1)$ can be identified with $\fOrd_c^1(X^{\an}, \cO)$ and $\fOrd_{\varpi-\ad, c}^1(X^{\an}, \cO)$ respectively.
There are obvious $T_p$-, Hecke and Galois equivariant homomorphism 
\begin{equation*}
\label{ordetonec2}
\fOrd_{c}^1(X_{\barQ}, \cO)(1) \lra \fOrd^1(X_{\barQ}, \cO)(1), \quad \fOrd_{\varpi-\ad, c}^1(X_{\barQ}, \cO)(1) \lra \fOrd_{\varpi-\ad}^1(X_{\barQ}, \cO)(1)\end{equation*}
and we denote their images by 
\begin{equation*}
\label{ordetoneint}
\fOrd_{!}^1(X_{\barQ}, \cO)(1) \qquad \mbox{and} \qquad \fOrd_{\varpi-\ad, !}^1(X_{\barQ}, \cO)(1) 
\end{equation*}
respectively. As $\cO[T_p]$- and Hecke modules we have again
\begin{equation*}
\label{ordetoneint2}
\fOrd_{!}^1(X_{\barQ}, \cO)(1)=\fOrd_{!}^1(X^{\an}, \cO), \qquad \fOrd_{\varpi-\ad, !}^1(X_{\barQ}, \cO)(1) =\fOrd_{\varpi-\ad, !}^1(X^{\an}, \cO).
\end{equation*}

\begin{remark}
\label{remark:intcohom}
\rm Note that 
\[
H_{!}^1((X_{K_1(n)})_{\barQ}, \cO)(1)= \image(H_c^1((X_{K_1(n)})_{\barQ}, \cO)(1)\to H^1((X_{K_1(n)})_{\barQ}, \cO)(1)
\]
can be identified with the group $H^1((X_{K_1(n)}^*)_{\barQ}, \cO)(1)$, were for any $K\in \sK$ we let $X_K^*$ denote the completion of the open (adelic) modular curve $X_K$ (i.e.\ $X_K^*$ is the Baily-Borel compactification of $X_K$). So we have 
\[
\fOrd_{!}^1(X_{\barQ}, \cO)(1)\, \cong \, \dlim_n H^1((X_{K_1(n)}^*)_{\barQ}, \cO)(1)^{\ord}.
\] 
\end{remark} 

We will see that an analogue of Prop.\ \ref{prop:ordetlocgal} holds for the "cuspidal part" $\fOrd_{\varpi-\ad, !}^1(X_{\barQ}, \cO)(1)$ of $\fOrd_{\varpi-\ad}^1(X_{\barQ}, \cO)(1)$.
For that let $\bD$ be the totally ramified incoherent quaternion algebra over $\bA=\bA_{\bQ}$ with ramification set $\Ram_{\bD} = \{\infty\}$. As in \cite{yuanzhangzhang} and the proof of Lemma \ref{lemma:ordetlocgal2} the set of automorphic representations of $\bD^*/\bA^*$ over a subfield $E$ of $\bC\cong \bC_p$ is denoted by $\cA(\bD^*/\bA^*, E)$. Note that $\cA(\bD^*/\bA^*, \bC)$ can be identified with the set of cuspidal automorphic representations $\pi= \bigotimes'_v \pi_v$ of $G(\bA)$ such that $\pi_{\infty}$ is the discrete series representation of weight $2$.

\begin{prop}
\label{prop:ordetlocgal3}
There exists an $\cO[T_p\times \fG_p]\otimes \bT$-submodule $\fOrd_{\varpi-\ad, !}^1(X_{\barQ}, \cO)(1)^0$ of \break $\fOrd_{\varpi-\ad, !}^1(X_{\barQ}, \cO)(1)$ with the the following properties:
\medskip

\noi (i) As $\cO$- (resp.\ as $\cO[T_{\fp}]$-) modules $\fOrd_{\varpi-\ad, !}^1(X_{\barQ}, \cO)(1)^0$ and 
\[
\fOrd_{\varpi-\ad, !}^1(X_{\barQ}, \cO)(1)^{\et}:=\,\fOrd_{\varpi-\ad, !}^1(X_{\barQ}, \cO)(1)/\fOrd_{\varpi-\ad, !}^1(X_{\barQ}, \cO)(1)^0
\]
are torsionfree (resp.\ $\varpi$-adically admissible). 

\noi (ii) The action of $\fG_p$ on $\fOrd_{\varpi-\ad, !}^1(X_{\barQ}, \cO)(1)^0$ and $\fOrd_{\varpi-\ad, !}^1(X_{\barQ}, \cO)(1)^{\et}$ factors through $\fG_{\fp}^{\ab}$ and we have
\[
\rec(t^{-1}) \cdot x = \alpha(t)\, t\cdot x\qquad \mbox{and} \qquad \rec(t) \cdot y = t \cdot y
\]
for every $t\in T_p$ and $x\in \fOrd_{\varpi-\ad, !}^1(X_{\barQ}, \cO)(1)^0$ and $y\in \fOrd_{\varpi-\ad, !}^1(X_{\barQ}, \cO)(1)^{\et}$. 
\end{prop} 

\begin{proof} This can be proved in essentially the same way as Prop.\ \ref{prop:ordetlocgal} and Lemma \ref{lemma:ordetlocgal2} so we omit the details. The proof requires adjustments only at two places. Firstly, by work of Eichler and Shimura and because of Remark \ref{remark:intcohom} we now have a decomposition only for the submodule $H_{!}^1((X_{K_1(n)})_{\barQ}, \bQ_p)(1)$ (rather than $H^1((X_{K_1(n)})_{\barQ}, \bQ_p)(1)$) that is similar to \eqref{jacshimura3}
\begin{equation*}
\label{jacmodcurve}
H_{!}^1((X_{K_1(n)})_{\barQ}, \bQ_p)(1) \,  = \,\bigoplus_{\pi\in \cA(\bD^*/\bA^*, \bQ_p)} \, \pi^{K_1(n) \times K_0(\fn)^{\fp}} \otimes_{E_{\pi}} V(\pi). 
\end{equation*}
Secondly, for the existence of the submodule $\cF_{\pi}'$ we use (\cite{mazur-wiles}, Prop.\ on p.\ 243; see also \cite{tilouine}, Cor.\ 4.2) instead of (\cite{wiles1}, Thm.\ 2.2) and again Hida's trick.
\end{proof} 

As in the previous section we can use Propositions \ref{prop:ordrandpiadic} and \ref {prop:ordetlocgal3} in order to define a filtration on the subgroup 
\begin{equation}
\label{h1modcurve2}
\bH_{\varpi-\ad, !}^1(X_{\barQ}; A(\chi), E) := \,\Hom_{E[T_{\fp}]}(A(\chi^{-1}),\fOrd_{\varpi-\ad, !}^1(X_{\barQ}, \cO)_E). 
\end{equation}
of the cohomology group $\bH_{\varpi-\ad}^1(X_{\barQ}; A(\chi), E)= \Hom_{E[T_{\fp}]}(A(\chi^{-1}),\fOrd_{\varpi-\ad}^1(X_{\barQ}, \cO)_E)$. Here as before $A$ is an Artinian local $E$-algebra with maximal ideal $\fM$, residue field $A/\fM=E$ and $\chi: \bQ_p^*\cong T_p \to A^*$ is a bounded continuous homomorphism so that $A(\chi) \in \Ban_E^{\adm}(T_{\fp})$.

Note that as a $\bT_A$-module \eqref{h1modcurve2} can be identified with $\Hom_{E[T_{\fp}]}(A(\chi^{-1}),\fOrd_{\varpi-\ad}^1(X^{\an}, \cO)_E)$ so  \eqref{ordpiadicbs} together with  Prop.\ \ref{ordrandpiadic4} implies that there exists an exact sequence of $\bT_A$-modules 
\begin{equation*}
\label{h1modcurve3}
0\lra \bH_{\varpi-\ad, !}^1(X_{\barQ}; A(\chi), E) \lra \bH_{\varpi-\ad}^1(X_{\barQ}; A(\chi), E) \lra \Hom_{E[T_{\fp}]}(A(\chi^{-1}),C^1_E).
\end{equation*}
Moreover there is an exact sequence of $A[\fG_{\fp}]\otimes \bT$-modules
\begin{eqnarray} 
\label{h1modcurve4}
&&0\lra \bH_{\varpi-\ad, !}^1(X_{\barQ}; A(\chi), E)(1)^0 \lra \bH_{\varpi-\ad, !}^1(X_{\barQ}; A(\chi), E)(1) \\
&& \hspace{4cm} \lra \bH_{\varpi-\ad, !}^1(X_{\barQ}; A(\chi), E)^{\et} \stackrel{\delta}{\lra} \Ext_{E[T_{\fp}]}^1(A(\chi^{-1}),\cF_E) \nonumber
\end{eqnarray}
where as before we put
\begin{eqnarray*} 
\label{h1modcurve5}
&& \bH_{\varpi-\ad, !}^1(X_{\barQ}; A(\chi), E)(1)^0 :=\, \Hom_{E[T_{\fp}]}(A(\chi^{-1}), \cF_E) \\ 
&& \bH_{\varpi-\ad, !}^1(X_{\barQ}; A(\chi), E)(1)^{\et}:=\, \Hom_{E[T_{\fp}]}(A(\chi^{-1}), \fOrd_{\varpi-\ad, !}^1(X_{\barQ}, \cO), _E/\cF_E). \nonumber
\end{eqnarray*}
Similarly to Thm.\ \ref{theorem:reciprocity} we have

\begin{theorem}
\label{theorem:reciprocity2}
(a) The action of $\fG_p$ on the first and third term in \eqref{h1modcurve4} factors through $\fG_p^{\ab}$ and we have 
$\rec(t) \cdot x = \alpha^{-1}(t) \cdot \chi(t) \cdot x$ and $\rec(t)\cdot y =\chi^{-1}(t) \cdot y$ for every $t\in T_p$, $x\in \bH_{\varpi-\ad, !}^1(X_{\barQ}; A(\chi), E)(1)^0$ and $y\in \bH_{\varpi-\ad, !}^1(X_{\barQ}; A(\chi), E)(1)^{\et}$. Moreover if $\barchi^2\ne \alpha$ then the sequence 
\begin{equation*} 
\label{h1modcurve4a}
0\to \bH_{\varpi-\ad, !}^1(X_{\barQ}; A(\chi), E)(1)^0 \to \bH_{\varpi-\ad, !}^1(X_{\barQ}; A(\chi), E)(1) \to \bH_{\varpi-\ad, !}^1(X_{\barQ}; A(\chi), E)^{\et}\to 0
\end{equation*}
is exact. 

\noi (b) The $\fG$-action on $\bH_{\varpi-\ad, !}^n(X_{\barQ}; A(\chi), E)(1)$ is unramified outside the set of primes dividing $pN$ and the Eichler-Shimura relations $\Frob_{\ell}^2 + T_{\ell}  \Frob_{\ell} + \ell = 0$ hold on $\bH_{\varpi-\ad, !}^n(X_{\barQ}; A(\chi), E)(1)$ for all $\ell\nmid p N$ and $n\ge 0$.
\end{theorem}

\addappendix




In the appendix we assemble a few results of homological algebra that are seemingly not be covered by the literature.

\paragraph{Right derived functors and inverse limits} As in \cite{jannsen} for an abelian category $\sA$ we denote by $\sA^{\bN}$ the category of inverse systems indexed by the set $\bN$ of natural numbers, so objects in $\sA^{\bN}$ are inverse systems 
\[
(A_n, d_n):  \ldots \lra A_{n+1} \stackrel{d_{n+1}}{\lra} A_n \lra \ldots \lra A_1
\] 
and morphisms are morphisms of inverse systems. 

\begin{lemma}
\label{lemma:derfunclim}
Let $F: \sA\to \sB$ be a left exact functor between abelian categories that has an exact left adjoint. Assume that $\sA$ and $\sB$ satisfy (AB3*) and (AB4*) and that they have enough injectives. Let $(A_n, d_n)$ be an object of $\sA^{\bN}$ that satisfy 
the following conditions 

\noi (i) $\prolimn^{(1)} (A_n, d_n) = 0$.

\noi (ii) $\prolimn^{(1)}(R^q F(A_n), R^q F(d_n)) =0$ for every $q\ge 0$.

\noi Then the canonical homomorphism 
\[
R^q F( \prolimn A_n) \, \lra \, \prolimn R^q F(A_n)
\]
is an isomorphism for every $q\ge 0$.
\end{lemma}

Note that (AB3*) implies that projective limits exists in $\sA$ and $\sB$.

\begin{proof} The functor $F$ induces a functor $F^{\bN}: \sA^{\bN}\to \sB^{\bN}$ that has again an exact left adjoint. By (\cite{jannsen}, 1.1 and 1.2) the category $\sA^{\bN}$ has enough injectives as well and we have $R^q F^{\bN} = (R^q F)^{\bN}$ for every $q\ge 0$.
We denote its composite with $\prolimn : \sB^{\bN}\to \sB$ by $\prolimn F: \sA^{\bN} \to \sB$. Since $F$ has a left adjoint, $\prolimn F$ can also be decomposed as $F\circ \prolimn: \sA^{\bN}\to \sA\to \sB$. So there are two Grothendieck spectral sequences 
\begin{eqnarray*}
&& ^{(1)} E_2^{rs} = R^s F \circ \prolimn^{(r)}(A_n, d_n)\, \Longrightarrow\, (R^{r+s} \prolimn F)(A_n, d_n).\\
&& ^{(2)} E_2^{rs} = \prolimn^{(r)} \circ (R^s F)^{\bN}(A_n, d_n)\, \Longrightarrow\, (R^{r+s} \prolimn F)(A_n, d_n).
\end{eqnarray*}
Conditions (i), (ii) and (AB4*) together with (\cite{weibel}, Corollary 3.5.4) imply that both spectral sequences degenerate, i.e.\ we have $^{(1)} E_2^{rs}=0$ if $s\ge 1$ and $^{(2)} E_2^{rs}=0$ if $r\ge 1$. Thus we get  
\[
R^q F(\prolimn A_n) \, \cong \, (R^q \prolimn F)(A_n, d_n)\, \cong \, \prolimn R^q F(A_n)
\]
for every $q\ge 0$.
\end{proof}

\begin{lemma}
\label{lemma:derfunclim2}
Let $F: \sA\to \sB$, $\sA$ and $\sB$ be as in Lemma \ref{lemma:derfunclim} and let $(A_n, d_n)$ be an object of $\sA^{\bN}$. Assume that $d_n: A_n \to A_{n-1}$ is surjective and that $\ker(d_n)$ is $F$-acyclic for every $n\ge 2$. We also assume that $A_1$ is $F$-acyclic. Then $\prolimn A_n$ is $F$-acyclic as well. 
\end{lemma}

\begin{proof} Note that assumptions imply that $A_n$ is $F$-acyclic for every $n\ge 1$. Also by (\cite{weibel}, Prop.\ 3.5.7) conditions (i) and (ii) of Lemma \ref{lemma:derfunclim} hold. Hence we have 
$R^q F( \prolimn A_n) \cong \prolimn R^q F(A_n)=0$ for every $q\ge 1$. 
\end{proof}

\paragraph{$R$-linear abelian categories} Let $R$ be a ring and let $\sA$ be an $R$-linear category, i.e.\ $\sA$ is an additive category together with a ringhomomorphism $\mu: R\to \End(\id_{\sA})$. The later induces an $R$-module structure on $\Hom_{\sA}(A,B)$ for all $A,B\in \sA$, so that the composition law $\Hom_{\sA}(B,C)\times \Hom_{\sA}(A,B)\to \Hom_{\sA}(A,C), (\beta,\alpha)\mapsto \beta\circ \alpha$ is $R$-bilinear. A functor $F: \sA\to \sB$ will be called $R$-linear, if $F: \Hom_{\sA}(A_1,A_2)\to \Hom_{\sA}(F(A_1), F(A_2)), \alpha\mapsto F(\alpha)$ is a homomorphism of $R$-modules for all $A_1, A_2\in \sA$. 

We now assume that $R$ is noetherian and that $\sA$ is an $R$-linear abelian category. For $M\in \Mod_{R,f}$ and $A\in \sA$ we define the objects $A\otimes_R M$, $\Hom_R(M, A)$ of $\sA$ as follows. Let $R^m \stackrel{f}{\lra} R^n \to M\to 0$ be a free resolution of $M$. The map $f$ is given by left multiplication with an $n\times m$-matrix $X$ with entries in $R$. It thus defines morphism $X\cdot : A^m \to A^n$ (resp.\ $X^t: A^n \to A^m$) in $\sA$ and we define $A\otimes_R M\in \sA$ (resp.\ $\Hom_R(M, A)$) as its cokernel (resp.\ kernel). It is easy to see that these objects are independent of the choice of the free resolution of $M$ (up to canonical isomorphism).
In particular for an ideal $\fa$ of $R$ and an object $A$ of $\sA$ we can consider the following objects in $\sA$: 
\begin{equation*}
\label{fatortens}
A[\fa]=\Hom_R(R/\fa, A), \qquad \fa A= \image(A\otimes_R \fa\to A), \qquad  A\otimes_R R/\fa.
\end{equation*} 
Let $F: \sA\to \sB$ be an $R$-linear additive covariant functor between $R$-linear abelian categories. For $M\in \Mod_{R,f}$ and $A\in \sA$ there exists canonical morphism $F(A)\otimes_R M \to F(A\otimes_R M)$ and $F(\Hom_R(M, A))\to \Hom_R(M, F(A))$. In particular morphisms for an ideal $\fa$ of $R$ we obtain morphisms
\begin{equation}
\label{addtenshom}
F(A)\otimes_R R/\fa \lra F(A\otimes_R R/\fa),\qquad F(A[\fa])\lra F(A)[\fa]
\end{equation} 
These are isomorphism if $F$ is exact. Similarly, if $F$ is contravariant there are canonical morphisms $F(A\otimes_R M)\to \Hom(M, F(A))$ and $F(A)\otimes_R M \to F(\Hom(M,A))$. In particular there are morphisms 
\begin{equation}
\label{addtenshom2}
F(A\otimes_R R/\fa)\lra F(A)[\fa], \qquad F(A)\otimes_R R/\fa \lra F(A[\fa])
\end{equation} 
for every ideal $\fa$ of $R$ that are isomorphisms if $F$ is exact.

\paragraph{Localization of $\delta$-functors} 
Let $S$ be a given multiplicative subset $R$. We denote by $S^{-1} \sA$ the $S^{-1}R$-linear additive category with the same objects as $\sA$ and with morphisms
\[
\Hom_{S^{-1} \sA}(A, B) \, =\, S^{-1} \Hom_{\sA}(A, B)\, =\,  \Hom_{\sA}(A, B)\otimes_R S^{-1}R.
\]
If $R$ is an integral domain with field of fractions $E=S^{-1}R$ (with $S=\setminus \{0\}$) then we will also write $\sA_E$ instead of $S^{-1} \sA$. 

This construction is a special case of the localization of a category with respect to multiplicative system of morphism (see e.g.\ \cite{weibel}, \S 10.3). Let
\begin{equation}
\label{slocal}
\iota: \sA \lra S^{-1} \sA
\end{equation} 
be the canonical functor. We sometimes denote the image of $A\in \sA$ under \eqref{slocal} by
$S^{-1}A$ or $A\otimes_R S^{-1} R$. If $R$ is noetherian and $M$ is a finitely generated $R$-module then we have
\begin{equation}
\label{slocal2}
S^{-1} A \otimes_{S^{-1} R} S^{-1} M\, =\, S^{-1}(A\otimes_R M),\qquad \Hom_{S^{-1}R}(S^{-1} M, S^{-1} A) \, =\, S^{-1}\Hom_R(M, A).
\end{equation} 
For $s\in S$ and $A\in \sA$ we put $sA:=\image(s\cdot 1_A)$.

\begin{lemma}
\label{lemma:locab}
For $A\in \sA$, the following conditions are equivalent

\noi (i) $\iota(A) = 0$.

\noi (ii) The covariant functor $S^{-1}\Hom_{\sA}(A, \wcdot) : \sA\to \Mod_{S^{-1}R}$ is trivial, i.e.\ $S^{-1}\Hom_{\sA}(A, B)=0$ for all $B\in \sA$.

\noi (iii) The contravariant functor $S^{-1}\Hom_{\sA}(\wcdot, A) : \sA\to \Mod_{S^{-1}R}$ is trivial.

\noi (iv) There exists $s\in S$ such that $sA = 0$. 
\end{lemma}

\begin{proof} The equivalence of (i) and (ii) and (i) and (iii) is an immediate consequence of the Yoneda Lemma. If (iii) holds then we have in particular $S^{-1}\Hom_{\sA}(A, A)=0$. Hence there exists $s\in S$ with $s1_A = s \cdot 0 =0 \in \Hom_{\sA}(A, A)$ and we get (iv). Conversely, (iv) implies $0 = \iota(1_A)=1_{\iota(A)}$ in $\End_{S^{-1}\sA}(\iota(A))$ hence (i) holds.
\end{proof}

\begin{lemma}
\label{lemma:locab2}
\noi (a) Let $\sB$ be an $S^{-1}R$-linear category and let $F: \sA\to \sB$ be an $R$-linear functor. Then $F$ factors uniquely through $\iota$.
\medskip

\noi (b) If $\sA$ is abelian then $S^{-1} \sA$ is abelian as well and $\iota: \sA\to S^{-1} \sA$ is exact.
\medskip

\noi (c) If $\sA$ is abelian category then any short exact sequence in $S^{-1} \sA$ is isomorphic to the image of a short exact sequence in $\sA$ under the functor $\iota: \sA\to S^{-1} \sA$. 
\end{lemma}

\begin{proof} (a) is obvious. For (b) consider the full subcategory $\sB = \ker(\iota) = \{A\in \sA\mid \iota(A)=0\}$ of $\sA$. The characterization of objects in $\sB$ given in Lemma \ref{lemma:locab} above implies easily that $\sB$ is a Serre subcategory, i.e.\ $\sB$ is closed under the formation of subobjects, quotients and extensions.  
Let $\Sigma$ be the family of morphisms $f$ in $\sA$ such that $\ker(f)$ and $\coker(f)$ lie in $\sB$. 
For $s\in S$ and $A\in \sA$ we have $s\cdot 1_A\in \Sigma$. Hence the canonical functor $\sA \to \Sigma^{-1} \sA$ factors uniquely in the form $\sA\stackrel{\iota}{\lra} S^{-1} \sA \to \Sigma^{-1} \sA$. Since $\Sigma^{-1} \sA=\sA/\sB$ is abelian (see e.g.\ \cite{weibel}, 10.3.2) it suffices to prove that the second functor is an equivalence or, equivalently, that $\iota(f)$ is an isomorphism for every morphism $f: A_1\to A_2$ in $\Sigma$. For that it suffices to consider the two cases (M) $f$ is a monomorphism and (E) $f$ is an epimorphism. 

For case (M) let $s\in S$ such that $s\cdot \coker(f) =0$, hence $sA_2\subseteq \image(f)$. Let $g: sA_2\to A_1$ be the composite $sA_2 \stackrel{\incl}{\lra} \image(f)\stackrel{f^{-1}}{\lra} A_2$. Then $s^{-1}\cdot \iota(g)$ is the inverse of $\iota(f)$ in $S^{-1} \sA$. In case (E) there exists $s\in S$ with $s\cdot \ker(f) =0$, hence $\ker(f) \subseteq \ker(s\cdot :A_1\to sA_1)$. Hence $s\cdot :A_1\to sA_1$ factors over $f$, i.e.\ there exists a morphism $h: A_2\to sA_1$ such that $h\circ f= s\cdot: A_1\to sA_1$. Then $s^{-1} \iota(h)$ is the inverse of $\iota(f)$.  

(c) A sequence of morphisms $A\stackrel{f}{\lra} B \stackrel{g}{\lra} C$ in $\sA$ will be called $S$-exact if $g\circ f=0$ and if the cohomology $\ker(g)/\image(g)$ lies in $\sB$. The latter condition is equivalent to require that $A\stackrel{f}{\lra} B \stackrel{g}{\lra} C$ is exact when viewed as a sequence in $S^{-1}\sA$.

We prove the assertion in two steps. Firstly, we show that any short exact sequence 
\begin{equation}
\label{sexact}
\begin{CD}
0 @>>> A @> f >> B @> g >> C @>>> 0\end{CD}
\end{equation}
in $S^{-1}\sA$ is isomorphic to the image under $\iota$ of a short $S$-exact sequence in $\sA$. Indeed, there exists $s_1, s_2\in S$ and morphisms $f': A\to B$, $g': B\to C$ such that $f=s_1^{-1} f'$ and $g= s_2^{-1} g'$. Since $\iota(g'\circ f') = (s_1 s_2) \cdot (g\circ f)=0$ we have $\image(g'\circ f')\in \sB$, i.e.\ there exists $s\in S$ with $s\cdot (g'\circ f')=0$. Hence by replacing $f'$ with $s\cdot f$ and $s_1$ with $s s_1$ we may assume that $g'\circ f' =0$. Consider the commutative diagram
\[
\begin{CD} 
0 @>>> A @> f >> B @> g >> C @>>> 0\\
@.@VV s_1 s_2 1_A V @VV s_2 1_B V @VV 1_C V @.\\
0 @>>> A @> f' >> B @> g' >> C @>>> 0
\end{CD}
\]
Viewed as morphisms in $S^{-1}\sA$ the vertical arrows are isomorphisms. Hence the lower sequence is $S$-exact and isomorphic to \eqref{sexact}.

In the second step we assume that \eqref{sexact} is a short $S$-exact sequence in $\sA$. Since $g\circ f = 0$ there exists $\ep: \coker(f) \to C$ with $\ep\circ q= g$ where $q:B \to \coker(f)$ is the projection. The exact sequence 
\[
A= \ker(g\circ f) \stackrel{f}{\lra} \ker(g)\lra \coker(f) \stackrel{\ep}{\lra} C =\coker(g\circ f) \lra \coker(g)\lra 0
\]
togerther with the fact that $\ker(g)/\image(f)$ and $\coker(g)$ lie in $\sB$ imply that $\iota(\ep)$ is an isomorphism. 
Hence the vertical arrows in the diagram 
\[
\begin{CD} 
0 @>>> A @> f >> B @> g >> C @>>> 0\\
@.@VV f V @VV 1_B V @VV \iota(\ep)^{-1} V @.\\
0 @>>> \image(f) @> \incl >> B @> q >> \coker(f) @>>> 0
\end{CD}
\]
are isomorphisms $S^{-1}\sA$. The lower row is exact in $\sA$. Hence \eqref{sexact} is isomorphic to a short exact sequence in $\sA$.
\end{proof}

Lemma \ref{lemma:locab} implies in particular that an $R$-linear functor $F: \sA\to \sB$ between $R$-linear categories extends uniquely to an $S^{-1}R$-linear functor $F: S^{-1} \sA\to S^{-1} \sB$.

\begin{lemma}
\label{lemma:locdelta}
Let $\sA$ (resp.\ $\sB$) be an $R$- (resp.\ $S^{-1}R$-) linear abelian category and let
\[
H^n: \sA\lra \sB,\qquad n\ge 0
\]
be an $R$-linear $\delta$-functor. Let 
\[
0 \to A_1 \stackrel{f_1}{\lra} A_2 \stackrel{f_2}{\lra} A_3 \to 0, \qquad 0 \to B_1 \stackrel{g_1}{\lra} B_2 \stackrel{g_2}{\lra} B_3 \to 0
\]
be two short exact sequences in $\sA$ and let $\alpha_i\in \Hom_{S^{-1}\sA}(\iota(A_i), \iota(B_i))$, $i=1,2,3$ so that 
\begin{equation}
\label{deltaloc1}
\begin{CD} 
0 @>>> \iota(A_1) @> \iota(f_1) >> \iota(A_2) @> \iota(f_2) >>\iota(A_3) @>>> 0\\
@.@VV \alpha_1 V @VV \alpha_2 V @VV \alpha_3 V @.\\
0 @>>> \iota(B_1) @> \iota(g_1) >> \iota(B_2) @> \iota(g_2) >>\iota(B_3) @>>> 0
\end{CD}
\end{equation}
commutes. Then the diagram 
\begin{equation}
\label{deltaloc2}
\begin{CD} 
H^n(A_3) @> \delta^n >> H^{n+1}(A_1) \\
@VV (\alpha_3)_* V@VV (\alpha_1)_* V\\
H^n(B_3) @> \delta^n >> H^{n+1}(B_1)
\end{CD}
\end{equation}
commutes for every $n\ge 0$.  
\end{lemma}

\begin{proof} There exists a common "denominator" $s\in S$ of $\alpha_1$, $\alpha_2$ and $\alpha_3$, i.e.\ there exists morphisms $\beta_i\in \Hom_{\sA}(A_i, B_i)$ such that $\alpha_i = s^{-1} \iota(\beta_i)$. By multiplying the vertical morphism  in \eqref{deltaloc1} with $s$ we see that $\iota(\beta_{i+1} \circ f_i) = \iota(g_i\circ \beta_i)$ for $i=1,2$. Thus the image of $\beta_{i+1} \circ f_i- g_i\circ \beta_i: A_i \to B_{i+1}$ lies in the subcategory $\ker(\iota)$, i.e.\ there exists $s'\in S$ with $s'\cdot (\beta_{i+1} \circ f_i- g_i\circ \beta_i) = 0$. Therefore if we replace the morphisms $\alpha_i$ with $s's\alpha_i$ then \eqref{deltaloc1} is the image under $\iota$ of a commutative diagram with exact rows in $\sA$. Hence \eqref{deltaloc2} commutes after multiplying the vertical arrows with $s s'$. However since $s s'1_B$ is an isomorphism for every $B\in \sB$ the original diagram commutes as well. 
\end{proof}

\begin{prop}
\label{prop:locdelta2}
Let $\sA$ be an $R$-linear abelian category, let $\sB$ be an $S^{-1}R$-linear abelian category and let
\[
H^n: \sA\lra \sB,\qquad n\ge 0
\]
be an $R$-linear $\delta$-functor.
Then $(H^n)_{n\ge 0}$ extends uniquely to an $S^{-1}R$-linear $\delta$-functor
\[
\wH^n: S^{-1}\sA \lra \sB,\qquad n\ge 0.
\]
\end{prop}

\begin{proof} By Lemma \ref{lemma:locab2} (a) for $n\ge 0$ the $R$-linear functors $H^n$ extends uniquely to $S^{-1}R$-linear functors $\wH^n: S^{-1} \sA\to \sB$.  

Let $0 \to A_1' \stackrel{f_1'}{\lra} A_2' \stackrel{f_2'}{\lra} A_3'' \to 0$ be a short exact sequence in $S^{-1}\sA$. We have to show that there exists canonical connecting homomorphisms 
\begin{equation}
\label{deltaloc}
\widetilde{\delta}^n: \wH^n(A_3') \lra \wH^{n+1}(A_1') \qquad \forall \,\, n\ge 0
\end{equation}
so that  
\[
\ldots \stackrel{(f_2')_*}{\lra}  \wH^{n-1}(A_3') \stackrel{\widetilde{\delta}}{\lra} \wH^n(A_1')\stackrel{(f_1')_*}{\lra} \wH^n(A_2')\stackrel{(f_2')_*}{\lra} \wH^n(A_3')\stackrel{\widetilde{\delta}}{\lra}\ldots
\]
is exact. By Lemma \ref{lemma:locab2} there exists a commutative diagram in $S^{-1}\sA$
\begin{equation}
\label{deltaloc3}
\begin{CD} 
0 @>>> \iota(A_1) @> \iota(f_1) >> \iota(A_2) @> \iota(f_2) >>\iota(A_3) @>>> 0\\
@.@VV \alpha_1 V @VV \alpha_2 V @VV \alpha_3 V @.\\
0 @>>> A_1' @> f_1' >> A_2' @> f_2' >> A_3' @>>> 0
\end{CD}
\end{equation}
with isomorphisms $\alpha_1$, $\alpha_2$ and $\alpha_3$. We define \eqref{deltaloc} so that the diagram 
\[
\begin{CD} 
\ldots@>>> H^{n-1}(A_3) @> \delta >> H^n(A_1) @> (f_1')_*>> H^n(A_2) @> (f_2')_* >> H^n(A_3) @>>> \ldots\\
@.@V\cong V (\alpha_3)_* V @V \cong V (\alpha_1)_* V @V\cong V (\alpha_2)_*  V @V\cong V (\alpha_3)_*  V @.\\
\ldots@>>> \wH^{n-1}(A_3') @> \widetilde{\delta} >> \wH^n(A_1') @> (f_1')_*>> \wH^n(A_2') @> (f_2')_* >> \wH^n(A_3') @>>> \ldots 
\end{CD}
\]
i.e.\ we put $\widetilde{\delta}= (\alpha_1)_*\circ \delta \circ (\alpha_3^{-1})_*$. By Lemma \ref{lemma:locdelta} this definition does not depend on the choice of the diagram \eqref{deltaloc3}.
\end{proof}

\begin{corollary}
\label{corollary:deltaloc}
Any $R$-linear $\delta$-functor $H=(H^n)_{n\ge 0}: \sA\to \sB$ between $R$-linear abelian categories extends uniquely to a $S^{-1}R$-linear $\delta$-functor $S^{-1}H= (S^{-1}H^n)_{n\ge 0}: S^{-1} \sA\to S^{-1} \sB$.
\end{corollary}

In the companion paper \cite{spiess2} we use the following simple 

\begin{prop}
\label{prop:deltafp}
Let $R$ be a noetherian ring and let $\fq\subseteq R$ be a primary ideal whose radical $\fm= \sqrt{\fq}$ is a regular maximal ideal of $R$ of height 1. Let $H=(H^n)_{n\ge 0}: \sA\to \sB$ be an $R$-linear $\delta$-functor between $R$-linear abelian categories and let $A\in \sA$. 
\medskip

\noi (a) Assume that $(H^n)_{n\ge 0}$ is covariant and that $A[\fm] = 0$. Then there exists natural short exact sequences
\begin{equation*}
\label{deltaloc4}
\begin{CD}
0@>>> H^n(A)\otimes_R R/\fq @> \eqref{addtenshom} >>  H^n(A\otimes_R R/\fq) @>>> H^{n+1}(A)[\fq]@>>> 0
\end{CD}
\end{equation*}
for every $n\ge 0$. 

\noi (b) Assume that $(H^n)_{n\ge 0}$ is contravariant and that $A\otimes_R R/\fm = 0$. Then there exists short exact sequences
\begin{equation*}
\label{deltaloc4b}
\begin{CD}
0@>>> H^n(A)\otimes_R R/\fq @> \eqref{addtenshom2} >>  H^n(A[\fq]) @>>> H^{n+1}(A)[\fq]@>>> 0
\end{CD}
\end{equation*}
for every $n\ge 0$. 
\end{prop}

\begin{proof} We prove only (a) since the proof of part (b) is very similar. Using \eqref{slocal2} and Cor.\ \ref{corollary:deltaloc} it is easy to see that the assertion is "local in $\fm$", i.e.\ we may pass to the localization of $R$ and of the $\delta$-functor $H=(H^n)_{n\ge 0}$ with respect to $S=R\setminus\fm$ without changing kernel and cokernel of \eqref{addtenshom} (and the object $H^{n+1}(A)[\fq]$). Thus we may assume that $R$ is a regular local ring of dimension $1$ with maximal ideal $\fm$ so that $\fq= \fm^s$ for some $s\ge 1$. Let $\varpi$ be a generator of $\fm$. The assumption $A[\fm] = 0$ implies that the sequence
\[
\begin{CD}
0 @>>> A @> \varpi^s >> A @>>> A\otimes_R R/\fq @>>> 0
\end{CD}
\]
is exact in $\sA$. The long exact sequence
\[
\begin{CD}
\ldots \to H^n(A) @> \varpi^s \cdot >>  H^n(A) @>>> H^n(A\otimes_R R/\fq) @>>> H^{n+1}(A) @> \varpi^s  >>  H^{n+1}(A)\to \ldots 
\end{CD}
\]
therefore yields the short exact sequences 
\[
\begin{CD}
0 @>>> H^n(A)\otimes_R R/\fq @>>> H^n(A\otimes_R R/\fq) @>>> H^{n+1}(A)[\fq]  @>>> 0
\end{CD}
\]
for every $n\ge 0$.
\end{proof}

 \textsc{Fakult{\"a}t f{\"u}r Mathematik, Universit{\"a}t Bielefeld, Germany} 

 \textit{E-mail} \text{mspiess@math.uni-bielefeld.de}

\end{document}